\documentstyle[12pt,epic,eepic]{article}
\title{{\bf Chiral Structure of
Modular Invariants for Subfactors}}

\author{{\sc Jens B\"ockenhauer} and {\sc David E. Evans}\\
School of Mathematics\\
University of Wales, Cardiff\\
PO Box 926, Senghennydd Road\\
Cardiff CF2 4YH, Wales, U.K.\\
e-mail: {\tt BockenhauerJM@cf.ac.uk, EvansDE@cf.ac.uk}  \\
\vphantom{X}\\
{\sc Yasuyuki Kawahigashi}\\
Department of Mathematical Sciences\\
University of Tokyo, Komaba, Tokyo, 153-8914, JAPAN\\
e-mail: {\tt yasuyuki@ms.u-tokyo.ac.jp}}
\date{July 23, 1999}

\begin{document}
\maketitle

\input amssym.def

\newsymbol\rtimes 226F

\def\Ad            {{\rm{Ad}}}
\def\Aut           {{\rm{Aut}}}
\def\bbC           {\Bbb{C}}
\def\bbM           {\Bbb{M}}
\def\bbN           {\Bbb{N}}
\def\bbNo          {\Bbb{N}_0}
\def\bbR           {\Bbb{R}}
\def\bbT           {\Bbb{T}}
\def\bbZ           {\Bbb{Z}}
\def\be            {\begin{equation}}
\def\bearl         {\begin{array}{l}}
\def\bearll        {\begin{array}{ll}}
\def\bearlll       {\begin{array}{lll}}
\def\bearrl        {\begin{array}{rl}}
\def\bea           {\begin{eqnarray}}
\def\beaa          {\begin{eqnarray*}}
\def\bfe           {{\bf1}}
\def\can           {\gamma}
\def\canr          {\theta}
\def\cA            {{\cal{A}}}
\def\cC            {{\cal{C}}}
\def\cCA           {\frak{C}}
\def\cD            {{\cal{D}}}
\def\cE            {{\cal{E}}}
\def\cF            {{\cal{F}}}
\def\cG            {{\cal{G}}}
\def\cH            {{\cal{H}}}
\def\cJ            {{\cal{J}}}
\def\cK            {{\cal{K}}}
\def\cL            {{\cal{L}}}
\def\cM            {{\cal{M}}}
\def\cN            {{\cal{N}}}
\def\cO            {{\cal{O}}}
\def\cP            {{\cal{P}}}
\def\cR            {{\cal{R}}}
\def\cS            {{\cal{S}}}
\def\cT            {{\cal{T}}}
\def\cV            {{\cal{V}}}
\def\cW            {{\cal{W}}}
\def\cX            {{\cal{X}}}
\def\cY            {{\cal{Y}}}
\def\cZ            {{\cal{Z}}}
\newcommand\co[1]  {\bar{{#1}}}
\def\diag          {{\rm{diag}}}
\def\dim           {{\rm{dim}}}
\newcommand\del[2] {\delta_{{#1},{#2}}}
\def\E             {{\rm{e}}}
\def\ee            {\end{equation}}
\def\eear          {\end{array}}
\def\eea           {\end{eqnarray}}
\def\eeaa          {\end{eqnarray*}}
\def\End           {{\rm{End}}}
\newcommand\eps[2] {\varepsilon({#1},{#2})}
\newcommand\epsm[2]{\varepsilon^-({#1},{#2})}
\newcommand\epsp[2]{\varepsilon^+({#1},{#2})}
\newcommand\epspm[2]{\varepsilon^\pm({#1},{#2})}
\newcommand\epsmp[2]{\varepsilon^\mp({#1},{#2})}
\newcommand\Eps[2] {{\cal E}({#1},{#2})}
\newcommand\Epsm[2]{{\cal E}^-({#1},{#2})}
\newcommand\Epsp[2]{{\cal E}^+({#1},{#2})}
\newcommand\Epspm[2]{{\cal E}^\pm({#1},{#2})}
\newcommand\Epsmp[2]{{\cal E}^\mp({#1},{#2})}
\newcommand\Epsr[2]{{\cal E}_{{\rm{r}}}({#1},{#2})}
\newcommand\erf[1] {Eq.\ (\ref{#1})}
\def\Exp           {{\rm{Exp}}}
\def\ext           {{\rm{ext}}}
\def\Gtwo          {{\rm{G}}_2}
\def\Hom           {{\rm{Hom}}}
\def\I             {{\rm{i}}}
\def\id            {{\rm{id}}}
\def\iotab         {{\co\iota}}
\def\Jz            {\cal{J}_z}
\def\lan           {\langle}
\def\lab           {{\co \lambda}}
\newcommand\laend[2]{\lambda_{{#1},{#2}}}
\def\LG            {{\it{LG}}}
\def\LH            {{\it{LH}}}
\def\LIG           {{\it{L}}_I{\it{G}}}
\def\LIcG          {{\it{L}}_{\Ic}{\it{G}}}
\def\LIH           {{\it{L}}_I{\it{H}}}
\def\LIcSUn        {{\it{L}}_{I'}{\it{SU}}(n)}
\def\LISUn         {{\it{L}}_I{\it{SU}}(n)}
\def\LISUk         {{\it{L}}_I{\it{SU}}(k)}
\def\LISUnk        {{\it{L}}_I{\it{SU}}(nk)}
\def\LIcSUz        {{\it{L}}_{I'}{\it{SU}}(2)}
\def\LISUz         {{\it{L}}_I{\it{SU}}(2)}
\newcommand\ls[1]  {[\lambda_{{#1}}]}
\def\LSE           {L^2(S^1)}
\def\LSn           {L^2(S^1;\Bbb{C}^n)}
\def\LSUd          {{\it{LSU}}(3)}
\def\LSUk          {{\it{LSU}}(k)}
\def\LSUn          {{\it{LSU}}(n)}
\def\LSUnk         {{\it{LSU}}(nk)}
\def\LSUz          {{\it{LSU}}(2)}
\def\Mat           {{\rm{Mat}}}
\def\Mor           {{\rm{Mor}}}
\def\mub           {{\co{\mu}}}
\def\mult          {{\rm{mult}}}
\def\MXN           {{}_M {\cal X}_N}
\def\MXM           {{}_M {\cal X}_M}
\def\MXMo          {{}_M^{} {\cal X}_M^0}
\def\MXMp          {{}_M^{} {\cal X}_M^+}
\def\MXMm          {{}_M^{} {\cal X}_M^-}
\def\MXMpm         {{}_M^{} {\cal X}_M^\pm}
\newcommand\N[3]   {N_{{#1},{#2}}^{{#3}}}
\def\Nres          {\tilde{N}}
\def\NXN           {{}_N {\cal X}_N}
\def\NXM           {{}_N {\cal X}_M}
\def\nub           {\overline{\nu}}
\def\oto           {=0,1,2,\ldots}
\def\pio           {\pi_0}
\def\PSLZ          {{\it{PSL}}(2;\bbZ)}
\def\PSU           {{\it{PSU}}(1,1)}
\def\reso          {|_{\cH_0}}
\def\rmA           {{\rm{A}}}
\def\rmD           {{\rm{D}}}
\def\rmE           {{\rm{E}}}
\def\Sect          {{\rm{Sect}}}
\def\SLC           {{\it{SL}}(2;\bbC)}
\def\SLnC          {{\it{SL}}(n;\bbC)}
\def\SLZ           {{\it{SL}}(2;\bbZ)}
\def\Ssys          {{\Sigma(\sys)}}
\def\SOf           {{\it{SO}}(5)}
\def\son           {\frak{so}(N)}
\def\SON           {{\it{SO}}(N)}
\def\SUd           {{\it{SU}}(3)}
\def\SUk           {{\it{SU}}(k)}
\def\SUm           {{\it{SU}}(m)}
\def\SUn           {{\it{SU}}(n)}
\def\SUnk          {{\it{SU}}(nk)}
\def\SUz           {{\it{SU}}(2)}
\def\SUzk          {{\it{SU}}(2k)}
\def\SUf           {{\it{SU}}(4)}
\def\sud           {\frak{su}(3)}
\def\sun           {\frak{su}(n)}
\def\suz           {\frak{su}(2)}
\def\sudh          {\widehat{\frak{su}}(3)}
\def\sunh          {\widehat{\frak{su}}(n)}
\def\suzh          {\widehat{\frak{su}}(2)}
\def\sys           {{\Delta}}
\def\tr            {{\rm{tr}}}
\def\Tr            {{\rm{Tr}}}
\def\Un            {\it{U}(n)}
\newcommand\V[3]   {V_{{#1};{#2}}^{#3}}
\def\Vir           {\frak{Vir}}


\def\qed{{\unskip\nobreak\hfil\penalty50
\hskip2em\hbox{}\nobreak\hfil  $\Box$      
\parfillskip=0pt \finalhyphendemerits=0\par}\medskip}
\def\proof{\trivlist \item[\hskip \labelsep{\bf Proof.\ }]}
\def\endproof{\null\hfill\qed\endtrivlist}

\def\equi{\sim}
\def\isom{\cong}
\def\ti{\tilde}
\def\lan{\langle}
\def\ran{\rangle}
\def\a{\alpha}
\def\de{\delta}
\def\e{\varepsilon}
\def\ga{\gamma}
\def\Ga{\Gamma}
\def\la{\lambda}
\def\La{\Lambda}
\def\th{\theta}
\def\om{\omega}
\def\Om{\Omega}
  
\newcommand\dta{\begin{picture}(12,10)
\thicklines
\path(2,4)(6,8)(10,4)(2,4)(6,0)(10,4)
\end{picture}}

\newcommand\dtap{\begin{picture}(12,10)\thicklines
\path(6,8)(10,4)(6,0)(6,8)(2,4)(6,0)\end{picture}}


\def\thinlines{\allinethickness{0.3pt}}
\def\thicklines{\allinethickness{1.0pt}}
\def\Thicklines{\allinethickness{2.0pt}}


\newtheorem{theorem}{Theorem}[section]
\newtheorem{lemma}[theorem]{Lemma}
\newtheorem{conjecture}[theorem]{Conjecture}
\newtheorem{corollary}[theorem]{Corollary}
\newtheorem{definition}[theorem]{Definition}
\newtheorem{assumption}[theorem]{Assumption}
\newtheorem{proposition}[theorem]{Proposition}
\newtheorem{remark}[theorem]{Remark}
\newtheorem{example}[theorem]{Example}

\begin{abstract}
In this paper we further analyze modular invariants
for subfactors, in particular the structure of the
chiral induced systems of $M$-$M$ morphisms.
The relative braiding between the chiral
systems restricts to a proper braiding on their
``ambichiral'' intersection, and we show that
the ambichiral braiding is non-degenerate if the
original braiding of the $N$-$N$ morphisms is.
Moreover, in this case the dimensions of
the irreducible representations of the
chiral fusion rule algebras are given by the
chiral branching coefficients which describe
the ambichiral contribution in the irreducible
decomposition of $\alpha$-induced sectors.
We show that modular invariants come along naturally
with several non-negative integer valued matrix
representations of the original $N$-$N$ Verlinde
fusion rule algebra, and we completely determine
their decomposition into its characters.
Finally the theory is illustrated by various examples,
including the treatment of all $\SUz_k$ modular
invariants, some $\SUd$ conformal inclusions
and the chiral conformal Ising model.
\end{abstract}

\newpage

\tableofcontents

\section{Introduction}

An important step towards complete classification of
rational conformal field theory would be an exhaustive
list of all modular invariant partition functions of
WZW models based on simple Lie groups $G$. In such models
one deals with a chiral algebra which is given by a
semi-direct sum of the affine Lie algebra of $G$ and the
associated Virasoro algebra arising from the Sugawara
construction. Fixing the level $k=1,2,...$, which specifies
the multiplier of the central extension, the chiral algebra
possesses a certain finite spectrum of representations
acting on (pre-) Hilbert spaces $\cH_\la$, labelled by
``admissible weights'' $\la$.
The characters
\[ \chi_\la(\tau;z_1,z_2,\ldots,z_\ell;u)
=\E^{2\pi\I ku} \, \tr_{\cH_\la}(\E^{2\pi\I\tau(L_0-c/24)}
\E^{2\pi\I(z_1 H_1+z_2 H_2+\ldots+z_\ell H_\ell)}) \,, \]
with $\rm{Im}(\tau)>0$, $L_0$ being the conformal Hamiltonian,
$c$ the central charge and $H_r$,
$r=1,2,...,\ell={{\rm rank}}(G)$, Cartan subalgebra
generators, transform unitarily under the action of
the (double cover of the) modular group, defined
by re-substituting the arguments as
\[
(\tau;\vec{z};u) \;\; \longmapsto \;\;
g(\tau;\vec{z};u) =
\left( \frac{a\tau+b}{c\tau+d} ; \frac{\vec{z}}{c\tau+d} ;
u-\frac{c(z_1^2+z_2^2+\ldots+z_\ell^2)}{2(c\tau+d)} \right)
\]
for $g=\left({a\atop c}{b\atop d}\right)\in\SLZ$,
see e.g.\ \cite{Kc}.
A modular invariant partition function is then
a sesqui-linear expression
$Z= \sum_{\la,\mu} Z_{\la,\mu}\chi_\la \chi_\mu^*$
which is is invariant under the $\SLZ$ action,
$Z(g(\tau;\vec{z};u))=Z(\tau;\vec{z};u)$, and subject to
\begin{equation}
\label{massmatZ}
Z_{\la,\mu} = 0,1,2,\ldots \,, \qquad\qquad Z_{0,0}=1 \,.
\end{equation}
Here the label ``0'' refers to the ``vacuum'' representation,
and the condition $Z_{0,0}=1$ reflects the physical concept of
uniqueness of the vacuum state. For the canonical generators
$\cS=\left({0\atop 1}{-1\atop 0}\right)$ and
$\cT=\left({1\atop 0}{1\atop 1}\right)$ of $\SLZ$ we
obtain the unitary Kac-Peterson matrices $S=[S_{\la,\mu}]$
and $T=[T_{\la,\mu}]$ transforming the characters,
where $T$ is diagonal and $S$ is symmetric
as well as $S_{\la,0}\ge S_{0,0}>0$.
Then the classification of modular invariants can be
rephrased like this:
Find all the matrices $Z$ subject to the conditions in
\erf{massmatZ} and commuting with $S$ and $T$.
This problem turns out to be a rather difficult one;
a complete list is known for all
simple Lie groups at low levels, however, a list covering
all levels is known to be complete only for Lie groups
$\SUz$ and $\SUd$.

Let us consider the $\SUz$ case. For $\SUz$ at level $k$,
the admissible weights are just spins $\la=0,1,2,...,k$. 
The Kac-Peterson matrices are given explicitly as
\[ S_{\la,\mu}=\sqrt{\frac 2{k+2}} \sin \left(
\frac{\pi(\la+1)(\mu+1)}{k+2} \right) , \qquad
T_{\la,\mu}=\del\la\mu \exp \left(
\frac{\pi\I(\la+1)^2}{2k+4} - \frac {\pi\I}4 \right) , \]
with $\la,\mu=0,1,...,k$.
A list of $\SUz$ modular invariants was given in
\cite{CIZ1} and proven to be complete in \cite{CIZ2,Kt},
the celebrated A-D-E classification of $\SUz$ modular
invariants. The A-D-E pattern arises as follows.
The eigenvalues of the (adjacency matrices of the) A-D-E
Dynkin diagrams are all of the form
$2\cos(m\pi/h)$ with $h=3,4,...$ being the (dual) Coxeter
number and $m$ running over a subset of $\{1,2,...,h-1\}$,
the Coxeter exponents of the diagram.
The bijection between the modular invariants $Z$
in the list of \cite{CIZ1} and Dynkin diagrams is
then such that the diagonal
entries $Z_{\la,\la}$ are given exactly by the
multiplicity of the eigenvalue
$2\cos(\pi(\la+1)/k+2)$ of one of the A-D-E Dynkin
diagrams with Coxeter number $h=k+2$. In particular,
the trivial modular invariants, $Z_{\la,\mu}=\del\la\mu$,
correspond to the diagrams $\rmA_{k+1}$. Note that the
adjacency matrix of the $\rmA_{k+1}$ diagram is given by the
level $k$ fusion matrix $N_1$ of the spin $\la=1$ representation.
Here $N_\la=[N_{\la,\mu}^\nu]$, and the
(non-negative integer valued) fusion rules $N_{\la,\mu}^\nu$
are generically (e.g.\ for all $\SUn$)
given by the Verlinde formula
\begin{equation}
\label{verlinde}
N_{\la,\mu}^\nu = \sum_\rho
\frac{S_{\rho,\la}}{S_{\rho,0}}
S_{\rho,\mu} S_{\rho,\nu}^* \,.
\end{equation}
As we have
$N_\la N_\mu = \sum_\nu N_{\la,\mu}^\nu N_\nu$,
we may interpret the $\rmA_{k+1}$ matrix is the
spin one representation matrix in the regular
representation of the fusion rules.
The meaning of the D and E diagrams, however, remained
obscure, and this has been regarded as a
``mystery of the A-D-E classification'' \cite{Itz}.
In fact, the adjacency matrices of the D-E diagrams
turned out to be only the spin
$\la=1$ matrices $G_1$ of a whole family
of non-negative integer valued matrices $G_\la$ providing
a representation of the original $\SUz_k$ fusion rules:
$G_\la G_\mu = \sum_\nu N_{\la,\mu}^\nu G_\nu$.
By the Verlinde formula, \erf{verlinde}, the representations
of the commutative fusion rule algebra are given by
the characters $\chi_\rho(\la)=S_{\rho,\la}/S_{\rho,0}$,
and therefore the multiplicities of the Coxeter
exponents just reflect the multiplicity of the
character $\chi_\rho$ in the representation given
by the matrices $G_\la$. Di Francesco, Petkova and
Zuber similarly observed \cite{DZ1,DZ2,PZ}
that there are non-negative integer valued matrix
representations (nimreps, for short)
of the $\SUn_k$ fusion rules which decompose into the
characters matching the diagonal part of some non-trivial
$\SUn_k$ modular invariants (mainly $\SUd$).
Graphs are then obtained by reading the matrices
$G_\la$ as adjacency matrices, with $\la$ now chosen
among the fundamental weights of $\SUn$ generalizing
appropriately the spin 1 weight for $\SUz$.
The classification of $\SUd$ modular invariants \cite{G2}
shows a similar pattern as the $\SUz$ case, called
$\cA$-$\cD$-$\cE$, $\cA$ referring to the diagonal invariants,
$\cD$ to ``simple current invariants'' and $\cE$ to
exceptionals. Again, it is the nimreps
associated to the $\cD$ and $\cE$ invariants
which call for an explanation whereas the $\cA$ invariants
just correspond to the original fusion rules: $G_\la=N_\la$.
Why are there graphs and, even more, nimreps of the Verlinde
fusion rules associated to modular invariants?
This question has not been answered for a long time.
Nahm found a relation between the diagonal part of
$\SUz$ modular invariants and  Lie group exponents using
quaternionic coset spaces \cite{N}, however, his construction
does not explain the appearance of nimreps of
fusion rules and seems impossible to be extended to other
Lie groups e.g.\ $\SUd$.

A first step in associating systematically nimreps
of the Verlinde fusion rules was done by F.\ Xu \cite{X1}
using nets of subfactors \cite{LR} arising from conformal
inclusions of $\SUn$ theories.
However, only a small number of modular invariants
comes from conformal inclusions, e.g.\ the $\rmD_4$,
$\rmE_6$ and $\rmE_8$ invariants for $\SUz$.
Developing systematically the $\alpha$-induction machinery
\cite{BE1,BE3} for nets of subfactors, a notion originally
introduced by Longo and Rehren \cite{LR}, such nimreps
where shown in \cite{BE2,BE3} to arise similarly from all
(local) simple current extensions \cite{SY} of $\SUn$
theories, thus covering in particular the $\rmD_{{\rm{even}}}$
series for $\SUz$. Yet, type \nolinebreak II invariants
(cf.\ $\rmD_{{\rm{odd}}}$ and $\rmE_7$ for $\SUz$)
were not treated in \cite{BE2,BE3}.

In \cite{BEK1} we have constructed modular invariants from
braided subfactors in a very general approach which unifies
and develops further the ideas of $\alpha$-induction
\cite{LR,X1,BE1,BE2,BE3} and Ocneanu's double triangle
algebras \cite{O7}. We started with a von Neumann factor $N$
endowed with a system $\NXN$ of braided endomorphisms
(``$N$-$N$ morphisms'').
Such a braiding defines ``statistics'' matrices $S$ and $T$
\cite{R0,FG} which, as shown by Rehren \cite{R0}, provide a
unitary representation of $\SLZ$ if it is non-degenerate.
The statistics matrices are analogous to the Kac-Peterson
matrices: $T$ is diagonal, $S$ is symmetric and
$S_{\la,0}\ge S_{0,0}>0$. (The label ``$0$'' now
refers to the identity morphism $\id\in\NXN$ which
corresponds to the vacuum in applications.)
Moreover, the endomorphism
fusion rules are diagonalized by
the statistical S-matrix, i.e.\ obey the Verlinde
formula \erf{verlinde} in the non-degenerate case.
We then studied embeddings $N\subset M$ in larger
factors $M$ which are in a certain sense compatible
with the braided system of endomorphisms; namely,
such a subfactor $N\subset M$ is essentially
given by specifying a canonical endomorphism
within the system $\NXN$. Then one can apply
$\alpha$-induction which associates to an
$N$-$N$ morphism $\la$ two $M$-$M$ morphisms, $\a^+_\la$
and $\a^-_\la$. Motivated by the analysis in \cite{BE3},
we {\em defined} a matrix $Z$ with entries
\begin{equation}
\label{Zdef}
Z_{\la,\mu} = \langle \a^+_\la,\a^-_\mu \rangle \,,
\qquad \la,\mu\in\NXN \,,
\end{equation}
where the brackets denote the dimension of the
intertwiner space $\Hom(\a^+_\la,\a^-_\mu)$.
Then $Z$ automatically fulfills the conditions of
\erf{massmatZ} and we showed that it commutes with
$S$ and $T$ \cite[Thm.\ 5.7]{BEK1}.
The inclusion $N\subset M$ associates
to $\NXN$ a system $\MXM$ of $M$-$M$ morphisms as well
as ``intermediate'' systems $\NXM$ and $\MXN$ where the
latter are related by conjugation. In turn, one obtains
a (graded) fusion rule algebra from the sector products.
Decomposing the induced morphisms $\a^\pm_\la$ into
irreducibles yields ``chiral'' subsystems of $M$-$M$ morphisms,
and it was shown that the whole system $\MXM$ is generated
by the chiral systems whenever the original braiding
is non-degenerate \cite[Thm.\ 5.10]{BEK1}.
We showed that each non-zero entry
$Z_{\la,\mu}$ labels one of the irreducible representations
of the $M$-$M$ fusion rules and its dimensions is
exactly given by $Z_{\la,\mu}$ \cite[Thm.\ 6.8]{BEK1}.
Moreover, we showed that the irreducible decomposition of the
representation obtained by multiplying $M$-$M$ morphisms on
$M$-$N$ morphisms corresponds exactly to the diagonal part
of the modular invariant \cite[Thm.\ 6.12]{BEK1}.

In this paper we take the analysis further and
investigate the chiral induced systems.
The matrix entry of \erf{Zdef} can be written as
\[ Z_{\la,\mu} = \sum\nolimits_\tau
b_{\tau,\la}^+ b_{\tau,\mu}^- \,,\]
where the sum runs over morphisms $\tau$ in the
``ambichiral'' intersection of the chiral systems, and
$b_{\tau,\la}^\pm=\langle\tau,\a^\pm_\la\rangle$
are the chiral branching coefficients. Analogous
to the second interpretation of the entries of $Z$,
we show that the chiral branching coefficients are
at the same time the dimensions of the irreducible
representations of the chiral fusion rules. We can evaluate the
induced morphisms $\a^\pm_\la$ in all these representations
of the chiral or full $M$-$M$ fusion rule algebra.
The representation which decomposes according to the
diagonal part of the modular invariant is the one
obtained by multiplying $M$-$M$ morphisms
on the $M$-$N$ system. By evaluating $\a^+_\la$
(here $\a^-_\la$ yields the same) we obtain a
family of matrices $G_\la$. Since $\alpha$-induction
preserves the fusion rules, this provides a matrix
representation of the original $N$-$N$ (Verlinde)
fusion rule algebra,
$G_\la G_\mu = \sum_\nu N_{\la,\mu}^\nu G_\nu$,
which therefore must decompose into the characters
given in terms of the S-matrix.
Moreover, as the $G_\la$'s are just fusion
matrices (i.e.\ each entry is the dimension of a
finite-dimensional intertwiner space), we have in fact
obtained nimreps here. We are able to compute the
eigenvalues of the matrices and thus we determine
the multiplicities of the characters,
proving that $\chi_\la$ appears in it exactly
with multiplicity $Z_{\la,\la}$.

The structure of the induced $M$-$M$ system
is quite different from the original braided $N$-$N$
system. In general, neither the full system $\MXM$ nor
the chiral induced subsystems are braided, they can
even have non-commutative fusion. In fact, our
results show that the entire $M$-$M$ fusion
algebra (respectively a chiral fusion algebra) is
non-commutative if and only if an entry of $Z$
(respectively a chiral branching coefficient)
is strictly larger than one.
However, as constructed in \cite{BE3},
there is a relative braiding between the chiral
induced systems which restricts to a proper braiding
on the ambichiral system. We show that the ambichiral
braiding is non-degenerate provided that the original
braiding on $\NXN$ is.

Contact with conformal field theory, in particular
with $\SUn$ WZW models, is made through Wassermann's
loop group construction \cite{W2}. The factor $N$ can
be viewed as $\pi_0(\LISUn)''$, a local loop group in
the level $k$ vacuum representation. Wassermann's bimodules
corresponding to the positive energy representations
yield the system of $N$-$N$ morphisms, labelled by the the
$\SUn$ level $k$ admissible weights and obeying the
$\SUn_k$ fusion rules by \cite{W2}. The statistics matrices
$S$ and $T$ are then forced to coincide with the $\SUn_k$
Kac-Peterson matrices, so that
$Z_{\la,\mu} = \langle \a^+_\la,\a^-_\mu \rangle$
produced from subfactors $N\subset M$ will in fact give
modular invariants of the $\SUn_k$ WZW models.

Can any modular invariant of, say, $\SUn$ models, be
realized from some subfactor? We tend to believe that this
is true. A systematic construction of canonical endomorphisms
is available for all modular invariants arising form
conformal inclusions \cite{X1,BE2,BE3} or by simple currents
\cite{BE2,BE3}; the canonical endomorphism for modular
invariants from non-local simple currents (with fractional
conformal dimensions) can be obtained in the same way as
in the local case \cite{BE2,BE3} since the ``chiral locality
condition'' is no longer required to hold in our
general framework. Maybe not too surprising for experts
in modular invariants, it is the few
--- in Gannon's language --- $\cE_7$ type invariants for
which we do not (yet?) have a systematic construction.
Nevertheless we can realize the complete list of
$\SUz$ modular invariants, including $\rmE_7$. We
can determine the structure of the induced systems
completely and we can draw the simultaneous fusion
graphs of the left and right chiral generators.
For $\rmD_{{\rm{even}}}$, $\rmE_6$ and $\rmE_8$ this
was already presented in \cite{BE3}, and here we present
the remaining cases $\rmD_{{\rm{odd}}}$ and $\rmE_7$.
As in \cite{BE3} we obtain Ocneanu's pictures for
the ``quantum symmetries of Coxeter graphs''
\cite{O7}, a coincidence which reflects the
identification of $\a$-induced sectors
with chiral generators in the double triangle
algebra \cite[Thm.\ 5.3]{BEK1}.

This paper is organized as follows. In Sect.\ \ref{prelim}
we recall some basic facts and notations from \cite{BEK1}
and introduce more intertwining braiding fusion symmetry.
In Sect.\ \ref{chiralanal} we introduce basic notions
and we start to analyze the structure of the chiral
induced system. As a by-product, we show in our setting
that $Z_{\la,0}=\del\la0$ implies that $Z$ is a
permutation matrix corresponding to a fusion rule automorphism,
even if the braiding is degenerate.
Sect.\ \ref{Zpm} contains the main
results. We assume non-degeneracy of the braiding on
$\NXN$ and show that then the ambichiral braiding
is non-degenerate.
We decompose the chiral parts of the center of the double
triangle algebra into simple matrix blocks, corresponding
to a ``diagonalization'' of the chiral fusion rule algebras.
We evaluate the chiral generators in the simple matrix
blocks, corresponding to the evaluation of the induced
morphisms in the irreducible representations of the
chiral fusion rule algebras.
Sects.\ \ref{ADESUz} and \ref{morex} are devoted to
examples. In Sect.\ \ref{ADESUz} we realize the
remaining $\SUz$ invariants $\rmD_{{\rm{odd}}}$ and $\rmE_7$,
and we give an overview over all A-D-E cases.
We also discuss the nimreps of the Verlinde fusion
rules and the problems in finding an underlying fusion
rule structure for the type \nolinebreak II invariants,
a problem  which was noticed by
Di Francesco and Zuber \cite{DZ1,DZ2},
based on an observation of Pasquier \cite{Pa} who noticed
that for Dynkin diagrams $\rmA$, $\rmD_{{\rm{even}}}$,
$\rmE_6$ and $\rmE_8$ there exist positive fusion rules,
but not for $\rmD_{{\rm{odd}}}$ and $\rmE_7$.
In Sect.\ \ref{morex} we present more
examples arising from conformal inclusions of $\SUd$.
We also discuss in detail the chiral conformal Ising model
as an example of a non-trivial canonical endomorphism
producing the trivial modular invariant.
Finally we discuss degenerate examples.

\section{Preliminaries}
\label{prelim}

Let $A$, $B$ be infinite factors. We denote by $\Mor(A,B)$
the set of unital $\ast$-ho\-mo\-mor\-phisms from
$A$ to $B$. The statistical dimension of $\rho\in\Mor(A,B)$
is defined as $d_\rho=[B:\rho(A)]^{1/2}$ where $[B:\rho(A)]$
is the minimal index \cite{J,Ko}. A morphism
$\rho\in\Mor(A,B)$ is called irreducible if
$\rho(A)\subset B$ is irreducible, i.e.\ 
$\rho(A)'\cap B=\bbC\bfe_B$. Two morphisms
$\rho,\rho'\in\Mor(A,B)$ are called equivalent
if there is a unitary $u\in B$ such that
$\rho'=\Ad(u)\circ\rho$. The unitary equivalence
class $[\rho]$ of a morphism $\rho\in\Mor(A,B)$
is called a $B$-$A$ sector. For sectors we
have a notion of sums, products and conjugates
(cf.\ \cite[Sect.\ 2]{BEK1} and the references
therein for more details).
For $\rho,\tau\in\Mor(A,B)$ we denote
$\Hom(\rho,\tau)=\{t\in B:t\rho(a)=\tau(a)t,\,a\in A\}$
and $\langle\rho,\tau\rangle=\dim\,\Hom(\rho,\tau)$.
Let $N$ be a type \nolinebreak III factor equipped with
a system $\sys\subset\Mor(N,N)$ of endomorphisms
in the sense of \cite[Def.\ 2.1]{BEK1}. This
means essentially that the morphisms in $\sys$
are irreducible and have finite statistical dimension
and, as sectors, they are different and form a
closed fusion rule algebra. Then
$\Ssys\subset\Mor(N,N)$ denotes the set of morphisms
which decompose as sectors into finite sums of
elements in $\sys$. We assume the system $\sys$
to be braided in the sense of \cite[Def.\ 2.2]{BEK1}
and we extend the braiding to $\Ssys$
(see \cite[Subsect.\ 2.2]{BEK1}). We then
consider a subfactor $N\subset M$, i.e.\ $N$
embedded into another type \nolinebreak III factor $M$,
of that kind that the dual canonical endomorphism
sector $[\canr]$ decomposes in a finite sum of sectors
of morphisms in $\sys$, i.e.\ $\canr\in\Ssys$.
Here $\canr=\co\iota\iota$ with
$\iota:N\hookrightarrow M$ being the injection map
and $\co\iota\in\Mor(M,N)$ being a conjugate morphism.
Note that this forces the statistical dimension of
$\canr$ and thus the index of $N\subset M$ to be
finite, $d_\canr=[M:N]<\infty$.
Then we can define $\a$-induction \cite{BE1} along
the lines of \cite{BEK1} just by using the extension
formula of Longo and Rehren \cite{LR}, i.e.\ by putting
\[ \a_\la^\pm = \co\iota^{\,-1} \circ \Ad (\epspm \la\canr)
\circ \la \circ \co\iota \]
for $\la\in\Ssys$, using braiding operators
$\epspm \la\canr \in\Hom(\la\canr,\canr\la)$.
Then $\a^+_\la$ and $\a^-_\la$ are morphisms in $\Mor(M,M)$
satisfying in particular $\a^\pm_\la\iota=\iota\la$.

\begin{figure}[htb]
\begin{center}
\unitlength 0.6mm
\begin{picture}(160,40)
\Thicklines
\put(20,40){\line(0,-1){5.858}}
\put(30,20){\line(-1,1){7.071}}
\put(30,20){\line(1,-1){7.071}}
\put(30,34.142){\arc{20}{2.356}{3.142}}
\put(40,5.858){\vector(0,-1){5.858}}
\put(30,5.858){\arc{20}{5.498}{0}}
\put(40,40){\line(0,-1){5.858}}
\put(32,22){\line(1,1){5.071}}
\put(28,18){\line(-1,-1){5.071}}
\put(30,34.142){\arc{20}{0}{0.785}}
\put(20,5.858){\vector(0,-1){5.858}}
\put(30,5.858){\arc{20}{3.142}{3.927}}
\put(13,5){\makebox(0,0){$\beta_-$}}
\put(47,5){\makebox(0,0){$\beta_+$}}
\put(120,40){\line(0,-1){5.858}}
\put(128,22){\line(-1,1){5.071}}
\put(132,18){\line(1,-1){5.071}}
\put(130,34.142){\arc{20}{2.356}{3.142}}
\put(140,5.858){\vector(0,-1){5.858}}
\put(130,5.858){\arc{20}{5.498}{0}}
\put(140,40){\line(0,-1){5.858}}
\put(130,20){\line(1,1){7.071}}
\put(130,20){\line(-1,-1){7.071}}
\put(130,34.142){\arc{20}{0}{0.785}}
\put(120,5.858){\vector(0,-1){5.858}}
\put(130,5.858){\arc{20}{3.142}{3.927}}
\put(113,5){\makebox(0,0){$\beta_+$}}
\put(147,5){\makebox(0,0){$\beta_-$}}
\end{picture}
\end{center}
\caption{Wire diagrams for $\Epsr {\beta_+}{\beta_-}$
and $\Epsr {\beta_+}{\beta_-} ^*$}
\label{relcross}
\end{figure}
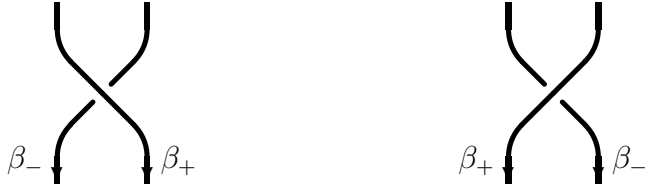

\begin{figure}[htb]
\begin{center}
\unitlength 0.6mm
\begin{picture}(95,40)
\Thicklines
\put(0,20){\arc{50}{5.356}{0.927}}
\put(40,20){\arc{50}{2.214}{2.418}}
\put(40,20){\arc{50}{2.578}{3.705}}
\put(40,20){\arc{50}{3.865}{4.069}}
\put(8,5){\makebox(0,0){$\beta_+$}}
\put(25,17){\vector(0,-1){0}}
\put(15,17){\vector(0,-1){0}}
\put(32,5){\makebox(0,0){$\beta_-$}}
\put(50,20){\makebox(0,0){$=$}}
\put(70,0){\line(0,1){40}}
\put(80,0){\line(0,1){40}}
\put(70,17){\vector(0,-1){0}}
\put(80,17){\vector(0,-1){0}}
\put(63,5){\makebox(0,0){$\beta_+$}}
\put(87,5){\makebox(0,0){$\beta_-$}}
\end{picture}
\end{center}
\caption{Unitarity of relative braiding operators}
\label{Reid2}
\end{figure}
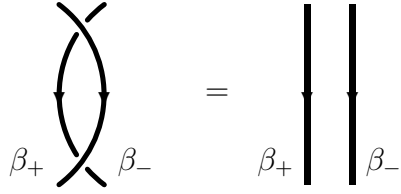

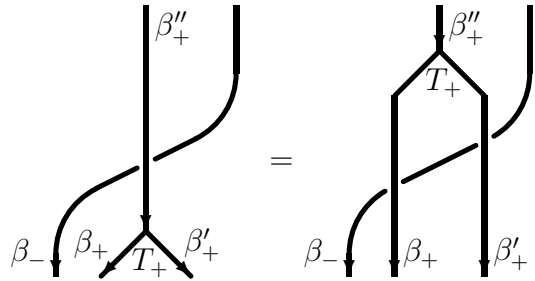
\begin{figure}[htb]
\begin{center}
\unitlength 0.6mm
\begin{picture}(117,60)
\Thicklines
\put(26.180,5){\arc{32.361}{3.142}{4.249}}
\put(10,5){\vector(0,-1){5}}
\put(28,24){\line(-2,-1){9.1}}   
\put(33.820,45){\arc{32.361}{0}{1.107}}
\put(50,60){\line(0,-1){15}}
\put(32,26){\line(2,1){9.1}}
\put(30,60){\line(0,-1){35}}
\put(30,25){\vector(0,-1){15}}
\put(30,10){\vector(1,-1){10}}
\put(30,10){\vector(-1,-1){10}}
\put(31,3){\makebox(0,0){$T_+$}}
\put(36,55){\makebox(0,0){$\beta_+''$}}
\put(18,7){\makebox(0,0){$\beta_+$}}
\put(42.5,7){\makebox(0,0){$\beta_+'$}}
\put(4,5){\makebox(0,0){$\beta_-$}}
\put(60,25){\makebox(0,0){$=$}}
\put(91.180,5){\arc{32.361}{3.142}{4.179}}    
\put(75,5){\vector(0,-1){5}}
\put(95,25){\line(-2,-1){8.1}}   
\put(98.820,45){\arc{32.361}{0}{1.037}} 
\put(115,60){\line(0,-1){15}}
\put(95,25){\line(2,1){8.1}}
\put(95,60){\vector(0,-1){10}}
\put(95,50){\line(1,-1){10}}
\put(95,50){\line(-1,-1){10}}
\put(85,40){\line(0,-1){20}}
\put(105,40){\line(0,-1){10}}
\put(85,10){\line(0,1){10}}
\put(105,10){\line(0,1){20}}
\put(85,10){\vector(0,-1){10}}
\put(105,10){\vector(0,-1){10}}
\put(96,43){\makebox(0,0){$T_+$}}
\put(101,55){\makebox(0,0){$\beta_+''$}}
\put(91,5){\makebox(0,0){$\beta_+$}}
\put(111,5){\makebox(0,0){$\beta_+'$}}
\put(69,5){\makebox(0,0){$\beta_-$}}
\end{picture}
\end{center}
\caption{The first braiding fusion relation
for the relative braiding}
\label{wireErbfe}
\end{figure}

In \cite[Subect.\ 3.3]{BE3}, a relative braiding
between representative endomorphisms of subsectors of
$[\a_\la^+]$ and $[\a_\mu^-]$ was introduced.
Namely, if $\beta_+,\beta_-\in\Mor(M,M)$ are such that
$[\beta_+]$ and $[\beta_-]$ are subsectors of
$[\a_\la^+]$ and $[\a_\mu^-]$ for some $\la,\mu\in\Ssys$,
respectively, then
\[ \Epsr {\beta_+}{\beta_-} = S^* \a_\mu^-(T^*)
\epsp \la\mu \a_\la^+(S)T
\in \Hom (\beta_+\beta_-,\beta_-\beta_+) \]
is unitary where $T\in\Hom(\beta_+,\a_\la^+)$ and
$S\in\Hom_M(\beta_-,\a_\mu^-)$ are isometries.
It was shown that $\Epsr {\beta_+}{\beta_-}$ does not
depend on $\lambda,\mu\in\Ssys$ and not on the isometries
$S,T$, in the sense that, if
there are isometries $X\in\Hom(\beta_+,\a_\nu^+)$ and
$Y\in\Hom(\beta_-,\a_\rho^-)$ with some $\nu,\rho\in\Ssys$,
then
$\Epsr{\beta_+}{\beta_-} = Y^* \a_\rho^-(X^*)
\epsp\nu\rho \a_\nu^+(Y)X$.
Moreover, it was shown\footnote{The proof of
\cite[Prop.\ 3.12]{BE3} is actually formulated in
the context of nets of subfactors. However, the
proof is exactly the same in the setting of braided
subfactors and it does not depend on the chiral
locality condition.} in \cite[Prop.\ 3.12]{BE3}
that the system of unitaries $\Epsr{\beta_+}{\beta_-}$
provides a relative braiding between representative
endomorphisms of subsectors
of $[\a_\la^+]$ and $[\a_\mu^-]$ in the sense that, if
$\beta_+,\beta_-,\beta_+',\beta_-'\in\Mor(M,M)$ are such that
$[\beta_+],[\beta_-],[\beta_+'],[\beta_-']$ are subsectors of
$[\a_\la^+],[\a_\mu^-],[\a_\nu^+],[\a_\rho^-]$, respectively,
$\lambda,\mu,\nu,\rho\in\Ssys$, then we have
``initial conditions''
$\Epsr{\id}{\beta_-}=\Epsr{\beta_+}{\id}=\bfe$,
``composition rules''
\begin{equation}
\label{Ercomp}
\bearl
\Epsr {\beta_+\beta_+'}{\beta_-} = \Epsr {\beta_+}{\beta_-}
\, \beta_+ (\Epsr {\beta_+'}{\beta_-}) \,, \\[.4em]
\Epsr {\beta_+}{\beta_-\beta_-'} = \beta_-
(\Epsr {\beta_+}{\beta_-'}) \, \Epsr {\beta_+}{\beta_-} \,,
\eear
\end{equation}
and whenever $Q_+\in \Hom(\beta_+,\beta_+')$
and $Q_-\in \Hom(\beta_-,\beta_-')$ then we have
``naturality''
\begin{equation}
\label{Ernat}
\beta_-(Q_+) \, \Epsr {\beta_+}{\beta_-}
= \Epsr {\beta_+'}{\beta_-} \, Q_+ \,, \qquad
Q_- \, \Epsr {\beta_+}{\beta_-}
= \Epsr {\beta_+}{\beta_-'} \, \beta_+ (Q_-) \,.
\end{equation}
Now let also $\beta_\pm''\in\Mor(M,M)$ and
$T_\pm\in\Hom(\beta_\pm'',\beta_\pm \beta_\pm')$ be an
intertwiner. From Eqs.\ (\ref{Ercomp}) and (\ref{Ernat})
we obtain the following braiding fusion relations:
\begin{equation}
\label{Erbfe}
\bearrl
\beta_-(T_+) \, \Epsr {\beta_+''}{\beta_-}
&=\,\, \Epsr {\beta_+}{\beta_-}
\beta_+(\Epsr {\beta_+'}{\beta_-}) \, T_+ \\[.4em]
T_- \, \Epsr {\beta_+}{\beta_-''}
&=\,\, \beta_-(\Epsr {\beta_+}{\beta_-'})
\Epsr {\beta_+}{\beta_-} \, \beta_+(T_-) \\[.4em]
\beta_-(T_+)^* \, \Epsr {\beta_+}{\beta_-}
\beta_+(\Epsr {\beta_+'}{\beta_-})
&=\,\, \Epsr {\beta_+''}{\beta_-} \, T_+^* \\[.4em]
T_-^* \, \beta_-(\Epsr {\beta_+}{\beta_-'})
\Epsr {\beta_+}{\beta_-}
&=\,\, \Epsr {\beta_+}{\beta_-''} \, \beta_+(T_-)^*
\,. \eear
\end{equation}
We can include the relative braiding operators in
the ``graphical intertwiner calculus'' along the
lines of \cite{BEK1} where isometric intertwiners
(with certain prefactors) realizing ``fusion channels''
and unitary braiding operators are diagrammatically
represented by trivalent vertices and crossings,
respectively. Again, the symmetry relations
fulfilled by the relative braiding operators
determine topological moves of the corresponding
wire diagrams. For $\Epsr {\beta_+}{\beta_-}$
and $\Epsr {\beta_+}{\beta_-} ^*$ we draw wire
diagrams as in Fig.\ \ref{relcross}.
%
%
The unitarity of $\Epsr {\beta_+}{\beta_-}$ then gives
a (vertical) ``Reidemeister move of type \nolinebreak II''
as displayed in Fig.\ \ref{Reid2}.
%
Next, the braiding fusion equations of \erf{Erbfe}
correspond to ``crossings moved over trivalent
vertices''. We display here only the move obtained
from the first relation in Fig.\ \ref{wireErbfe}.
%

If $a\in\Mor(M,N)$ is such that $[a]$ is a subsector
of $\mu\co\iota$ for some $\mu$ in $\Ssys$ then
$a\iota\in\Ssys$. Hence the braiding operators
$\epspm\la{a\iota}$ are well defined for $\la\in\Ssys$.
We showed in \cite[Prop.\ 3.1]{BEK1} that
$\epspm\la{a\iota}\in\Hom(\la a, a\a^\pm_\la)$.
If $\co b\in\Mor(N,M)$ is such that
$[\co b]$ is a subsector of $[\iota\co\nu]$
for some $\co\nu\in\Ssys$ and
$T\in\Hom(\co b,\iota\co\nu)$ is an isometry,
then we showed that
$\Epspm \la{\co b}=T^*\epspm\la{\co\nu}\a^\pm_\la(T)$
(and $\Epspm{\co b}\la=\Epsmp\la{\co b} ^*$) are
independent of the particular choice of $\co\nu$
and $T$ and are unitaries in
$\Hom(\a^\pm_\la \co b,\co b \la)$
(respectively $\Hom(\co b \la,\a^\mp_\la \co b)$).
These operators obey certain symmetry relations
\cite[Prop.\ 3.3]{BEK1} which we called
``intertwining braiding fusion relations''
(IBFE's), and they can nicely be represented
graphically by ``mixed crossings'' which
involve ``thick wires'' representing $N$-$M$
morphisms (see \cite[Fig.\ 30]{BEK1}).
We will now complete the picture by relating
their braiding symmetry to the relative braiding
by means of additional IBFE's.

\begin{lemma}
Let $\la,\mu,\nu\in\Ssys$, $\beta_\pm\in\Mor(M,M)$
$a,b\in\Mor(M,N)$
such that $[\beta_\pm]$, $[a]$, $[b]$
are subsectors of
$[\a^\pm_\la]$, $[\mu\co\iota]$
and $[\nu\co\iota]$ respectively.
Let also $\co a,\co b\in\Mor(N,M)$
be conjugates of $a,b$, respectively.
Then we have
\begin{equation}
\label{comp1}
\Epsp {\beta_+ \co b}\rho =
\Epsr {\beta_+}{\a^-_\rho}\beta_+(\Epsp {\co b}\rho) \,,
\quad \Epsm {\beta_- \co b}\rho =
\Epsr {\a^+_\rho}{\beta_-} ^* \beta_-(\Epsm {\co b}\rho) \,,
\end{equation}
and
\begin{equation}
\label{comp2}
\epsp \rho{b\beta_-\iota} =
b(\Epsr{\a^+_\rho}{\beta_-})\epsp \rho{b\iota} \,,
\qquad \epsm \rho{b\beta_+\iota} =
b(\Epsr{\beta_+}{\a^-_\rho})^*\epsm \rho{b\iota} \,,
\end{equation}
for all $\rho\in\Ssys$.
\end{lemma}

\begin{proof}
Let $S_\pm\in\Hom(\beta_\pm,\a^\pm_\la)$ and
$T\in\Hom(\co b,\iota\co\nu)$ be isometries.
Then we have
\[ \Epsr {\beta_+}{\a^-_\rho} = \a^-_\rho(S_+)^*
\epsp \la\rho S_+ \,, \qquad
\Epsr {\a^+_\rho}{\beta_-} ^* = \a^+_\rho(S_-)^*
\epsm \la\rho S_- \,, \]
and
\[ \Epspm {\co b}\rho = \a^\mp_\rho(T)^*
\epspm {\co\nu}\rho T \,. \]
Since
$\a^\pm_\la(T)S_\pm \in
\Hom(\beta_\pm\co b,\a^\pm_\la \iota\co\nu)
\equiv\Hom(\beta_\pm\co b, \iota\la\co\nu)$
is an isometry we can compute
\[ \bearll
\Epspm {\beta_\pm\co b}\rho &= \a^\mp_\rho
(S_\pm^* \a^\pm_\la(T)^*) \epspm {\la\co\nu}\rho
\a^\pm_\la(T)S_\pm \\[.4em]
&= \a^\mp_\rho (S_\pm)^* \a^\mp_\rho\a^\pm_\la(T)^* 
\epspm \la\rho \la(\epspm {\co\nu}\rho )
\a^\pm_\la(T)S_\pm \\[.4em]
&= \a^\mp_\rho (S_\pm)^* \epspm \la\rho
\a^\pm_\la \a^\mp_\rho(T)^* \la(\epspm {\co\nu}\rho)
\a^\pm_\la(T)S_\pm \\[.4em]
&= \a^\mp_\rho (S_\pm)^* \epspm \la\rho S_\pm \,
\beta_\pm (\a^\mp_\rho(T)^* \epspm {\co\nu}\rho T) \,,
\eear \]
which gives the desired \erf{comp1}. Here we have used
that
$\epspm\la\rho\in\Hom(\a^\pm_\la \a^\mp_\rho,
\a^\mp_\rho \a^\pm_\la)$ by \cite[Lemma 3.24]{BE1}.
Since
$b(S_\pm)^*\in\Hom(b\a^\pm_\la\iota,b\beta_\pm\iota)
\equiv\Hom(b\iota\la,b\beta_\pm\iota)$ we can compute
by virtue of naturality (cf.\ \cite[Eq.\ (8)]{BEK1})
\[ \bearll
\epspm \rho{b\beta_\pm\iota} &= 
\epspm \rho{b\beta_\pm\iota} \rho b(S_\pm^* S_\pm)
= b(S_\pm)^* \epspm \rho{b\iota\la}
\rho b(S_\pm) \\[.4em]
&= b(S_\pm)^* b(\epspm \rho\la) \epspm \rho{b\iota}
\rho b(S_\pm) \\[.4em]
&= b(S_\pm)^* b(\epspm \rho\la) b \a^\pm_\rho(S_\pm) \,
\epspm \rho{b\iota} \,,
\eear \]
which gives the desired \erf{comp2}.
\end{proof}

From the naturality equations for the braiding
operators \cite[Lemma 3.2]{BEK1} and
\cite[Lemma 3.25]{BE1} we then obtain the following

\begin{corollary}
For $X_\pm\in\Hom(\co a,\beta_\pm\co b)$
and $x_\pm\in\Hom(a,\beta_\pm b)$ we have
IBFE's 
\begin{equation}
\label{newibfe1}
\bearll
\a^-_\rho(X_+) \Epsp {\co a}\rho &= \,\,
\Epsr {\beta_+}{\a^-_\rho} \beta_+
(\Epsp {\co b}\rho) X_+ \,,\\[.4em]
\a^+_\rho(X_-) \Epsm {\co a}\rho &= \,\,
\Epsr {\a^+_\rho}{\beta_-} ^* \beta_-
(\Epsm {\co b}\rho) X_- \,, \eear
\end{equation}
and
\begin{equation}
\label{newibfe2}
\bearll
x_+ \epsm \rho{a\iota} &= \,\,
b(\Epsr {\beta_+}{\a^-_\rho})^*
\epsm \rho{b\iota} \rho(x_+) \,,\\[.4em]
x_- \epsp \rho{a\iota} &= \,\,
b(\Epsr {\a^+_\rho}{\beta_-})
\epsp \rho{b\iota} \rho(x_-) \,.
\eear
\end{equation}
\end{corollary}
These IBFE's can again be visualized in diagrams. We display
the first relation of \erf{newibfe1} in Fig.\ \ref{IBFE1}
%
\begin{figure}[htb]
\begin{center}
\unitlength 0.6mm
\begin{picture}(128,60)
\Thicklines
\put(26.180,5){\arc{32.361}{3.142}{4.249}}
\put(10,5){\vector(0,-1){5}}
\put(28,24){\line(-2,-1){9.1}}
\put(30,10){\vector(-1,-1){10}}
\thinlines
\put(33.820,45){\arc{32.361}{0}{1.107}}
\put(50,60){\line(0,-1){15}}
\put(32,26){\line(2,1){9.1}}
\thicklines
\put(30,10){\vector(0,1){50}}
\put(30,10){\vector(1,-1){10}}
\put(30.5,3){\makebox(0,0){\footnotesize{$X_+$}}}
\put(35,55){\makebox(0,0){$a$}}
\put(55,55){\makebox(0,0){$\rho$}}
\put(18,6){\makebox(0,0){$\beta_+$}}
\put(42,5){\makebox(0,0){$b$}}
\put(3,5){\makebox(0,0){$\a^-_\rho$}}
\put(65,25){\makebox(0,0){$=$}}
\Thicklines
\put(96.180,5){\arc{32.361}{3.142}{4.179}}    
\put(80,5){\vector(0,-1){5}}
\put(90,40){\vector(0,-1){40}}
\put(100,50){\line(-1,-1){10}}
\thinlines
\put(103.820,45){\arc{32.361}{0}{1.057}} 
\put(120,60){\line(0,-1){15}}
\thicklines
\put(100,50){\vector(0,1){10}}
\put(110,40){\vector(-1,1){10}}
\put(110,40){\line(0,-1){40}}
\Thicklines
\put(100,25){\line(2,1){8}}
\put(100,25){\line(-2,-1){8}}
\put(100.5,43){\makebox(0,0){\footnotesize{$X_+$}}}
\put(105,55){\makebox(0,0){$a$}}
\put(96,5){\makebox(0,0){$\beta_+$}}
\put(115,5){\makebox(0,0){$b$}}
\put(73,5){\makebox(0,0){$\a^-_\rho$}}
\put(125,55){\makebox(0,0){$\rho$}}
\end{picture}
\end{center}
\caption{The first intertwining braiding fusion relation of \protect\erf{newibfe1}}
\label{IBFE1}
\end{figure}
and the second relation of \erf{newibfe2}
in Fig.\ \ref{IBFE2}.
%
\begin{figure}[htb]
\begin{center}
\unitlength 0.6mm
\begin{picture}(132,60)
\thinlines
\put(26.180,45){\arc{32.361}{2.034}{3.142}}
\put(10,60){\line(0,-1){15}}
\put(30,25){\line(-2,1){11.1}}   
\Thicklines
\put(33.820,5){\arc{32.361}{5.176}{6.283}}
\put(50,5){\vector(0,-1){5}}
\put(30,25){\line(2,-1){11.1}}
\thicklines
\put(30,27){\line(0,1){33}}
\put(30,23){\vector(0,-1){13}}
\put(30,10){\vector(-1,-1){10}}
\Thicklines
\put(30,10){\vector(1,-1){10}}
\put(5,55){\makebox(0,0){$\rho$}}
\put(30.5,3){\makebox(0,0){$x_-$}}
\put(35,55){\makebox(0,0){$a$}}
\put(18,5){\makebox(0,0){$b$}}
\put(42,7){\makebox(0,0){$\beta_-$}}
\put(57,5){\makebox(0,0){$\a^+_\rho$}}
\put(65,25){\makebox(0,0){$=$}}
\thinlines
\put(96.180,45){\arc{32.361}{2.034}{3.142}}
\put(80,60){\line(0,-1){15}}
\put(90,30){\line(-2,1){1.1}}
\Thicklines
\put(103.820,5){\arc{32.361}{5.176}{6.283}}
\put(120,5){\vector(0,-1){5}}
\put(110,20){\line(2,-1){1.1}}
\thicklines
\put(100,60){\vector(0,-1){10}}
\put(90,40){\line(1,1){10}}
\put(90,40){\line(0,-1){8}}
\put(90,28){\vector(0,-1){28}}
\Thicklines
\put(110,40){\line(-1,1){10}}
\put(110,40){\line(0,-1){18}}
\put(110,10){\line(0,1){8}}
\put(110,10){\vector(0,-1){10}}
\put(100,25){\line(-2,1){10}}
\put(100,25){\line(2,-1){10}}
\put(100.5,43){\makebox(0,0){$x_-$}}
\put(105,55){\makebox(0,0){$a$}}
\put(85,5){\makebox(0,0){$b$}}
\put(105,5){\makebox(0,0){$\beta_-$}}
\put(127,5){\makebox(0,0){$\a^+_\rho$}}
\put(75,55){\makebox(0,0){$\rho$}}
\end{picture}
\end{center}
\caption{The second intertwining braiding fusion relation of
\protect\erf{newibfe2}}
\label{IBFE2}
\end{figure}

Next we recall our definition of Ocneanu's
double triangle algebra. For the above considerations
we did not need finiteness of the system $\sys$.
For the definition of the double triangle algebra
we do need such a finiteness assumption but
it does not rely on the braiding. Therefore
we start again and work for the rest of this paper
with the following

\begin{assumption}
\label{assbasic}{\rm
Let $N\subset M$ be a type \nolinebreak III subfactor
of finite index. We assume that we have a finite system
of endomorphisms $\NXN\subset\Mor(N,N)$ in the sense of
\cite[Def.\ 2.1]{BEK1} such that
$\canr=\co\iota\iota\in\Sigma(\NXN)$ for the
injection map $\iota:N\hookrightarrow M$ and a
conjugate $\co\iota\in\Mor(M,N)$.
We choose sets of morphisms $\NXM\subset\Mor(M,N)$,
$\MXN\subset\Mor(N,M)$ and $\MXM\subset\Mor(M,M)$
consisting of representative endomorphisms of irreducible
subsectors of sectors of the form $[\la\co\iota]$,
$[\iota\la]$ and $[\iota\la\co\iota]$, $\la\in\NXN$,
respectively. We choose $\id\in\Mor(M,M)$
representing the trivial sector in $\MXM$.
}\end{assumption}
Then the the double triangle algebra $\dta$ is
given as a linear space by
\[\dta=\bigoplus_{a,b,c,d\in\NXM}\Hom (a\co b,c\co d)\]
and is equipped with two different multiplications;
the horizontal product $*_h$ and the vertical product
$*_v$ (cf.\ \cite[Sect.\ 4]{BEK1}). The center
$\cZ_h$ of $(\dta,*_h)$ is closed under the vertical
product. In fact, the algebra $(\cZ_h,*_v)$ is
isomorphic to the fusion rule algebra associated
to the system $\MXM$ (cf.\ \cite[Thm.\ 4.4]{BEK1}).
This fact provides a useful tool since in examples
the system $\NXN$ is typically the known part of
the theory whereas the dual system $\MXM$ is the
unknown part. To determine the structure of the
fusion rule algebra of $\MXM$, i.e.\
of $(\cZ_h,*_v)$, completely is often a rather
difficult problem. However, a braiding on $\NXN$
forces a lot of symmetry structure within the entire
set $\cX$ which can in turn be enough to
determine the whole $M$-$M$ fusion table completely.
For the rest of this paper we therefore now impose
the following

\begin{assumption}
\label{assbraid}
{\rm In addition to Assumption \ref{assbasic} we now assume that
the system $\NXN$ is braided in the sense of \cite[Def.\ 2.2]{BEK1}.
}\end{assumption}
In particular we then have the notion of $\a$-induction.
The relation $\a_\la^\pm\iota=\iota\la$ implies that
for any $\la\in\NXN$ each irreducible subsector of
$[\a^\pm_\la]$ is of the form $[\beta]$ for some
$\beta\in\MXM$. By $\MXMpm\subset\MXM$ we denote the
subsets corresponding to subsectors of $[\a^\pm_\la]$
when $\la$ varies in $\NXN$. By virtue of the
homomorphism property of $\a$-induction, the
sets $\MXMpm$ must in fact be systems of endomorphism
themselves. We call $\MXMp$ and $\MXMm$ the
{\sl chiral systems}. Clearly, another system is
obtained by taking the intersection $\MXMo=\MXMp\cap\MXMm$
which we call the {\sl ambichiral system}.
In this paper, we will make special use of the relative
braiding between $\MXMp$ and $\MXMm$. Note that the
relative braiding restricts to a proper braiding on
$\MXMo$. The relative braiding symmetry also gives rise
to new useful graphical identities.
Let $\beta_\pm,\beta_\pm'\in\MXMpm$
and $V_\pm\in\Hom(\beta_\pm,\a^\pm_\la)$.
From naturality \erf{Ernat} we obtain
\begin{equation}
\bearl
\beta_-'(V_+) \, \Epsr {\beta_+}{\beta_-'}
=  \Epsr {\a^+_\la}{\beta_-'} \, V_+ \,,\\[.4em]
V_- \,  \Epsr {\beta_+'}{\beta_-}
=  \Epsr {\beta_+'}{\a^-_\la} \, \beta_+'(V_-) \,.
\eear
\end{equation}
Graphically these equations are displayed in
Figs.\ \ref{natrelui}
%
\begin{figure}[htb]
\begin{center}
\unitlength 0.6mm
\begin{picture}(130,50)
\Thicklines
\put(26.180,10){\arc{32.361}{3.142}{4.249}}
\put(10,10){\vector(0,-1){10}}
\put(28,29){\line(-2,-1){9.1}}
\put(33.820,50){\arc{32.361}{0}{1.107}}
\put(32,31){\line(2,1){9.1}}
\put(30,50){\vector(0,-1){30}}
\put(30,10){\vector(0,-1){10}}
\thinlines
\put(25,10){\dashbox{2}(10,10){$V_+$}}
\put(25,45){\makebox(0,0){$\beta_+$}}
\put(37,5){\makebox(0,0){$\a^+_\la$}}
\put(5,5){\makebox(0,0){$\beta_-'$}}
\put(65,25){\makebox(0,0){$=$}}
\Thicklines
\put(96.180,0){\arc{32.361}{3.142}{4.249}}
\put(80,0){\vector(0,-1){0}}
\put(98,19){\line(-2,-1){9.1}}
\put(103.820,40){\arc{32.361}{0}{1.107}}
\put(120,50){\line(0,-1){10}}
\put(102,21){\line(2,1){9.1}}
\put(100,50){\vector(0,-1){10}}
\put(100,30){\vector(0,-1){30}}
\thinlines
\put(95,30){\dashbox{2}(10,10){$V_+$}}
\put(95,45){\makebox(0,0){$\beta_+$}}
\put(107,5){\makebox(0,0){$\a^+_\la$}}
\put(75,5){\makebox(0,0){$\beta_-'$}}
\end{picture}
\end{center}
\caption{Naturality move for relative braiding}
\label{natrelui}
\end{figure}
and \ref{natreloi}.
%
\begin{figure}[htb]
\begin{center}
\unitlength 0.6mm
\begin{picture}(130,50)
\Thicklines
\put(26.180,50){\arc{32.361}{2.034}{3.142}}
\put(30,30){\line(-2,1){11.1}}   
\put(33.820,10){\arc{32.361}{5.176}{6.283}}
\put(50,10){\vector(0,-1){10}}
\put(30,30){\line(2,-1){11.1}}
\put(30,50){\line(0,-1){18}}
\put(30,28){\vector(0,-1){8}}
\put(30,10){\vector(0,-1){10}}
\thinlines
\put(25,10){\dashbox{2}(10,10){$V_-$}}
\put(37,45){\makebox(0,0){$\beta_-$}}
\put(23,5){\makebox(0,0){$\a^-_\la$}}
\put(56,4){\makebox(0,0){$\beta_+'$}}
\put(65,25){\makebox(0,0){$=$}}
\Thicklines
\put(80,50){\line(0,-1){10}}
\put(96.180,40){\arc{32.361}{2.034}{3.142}}
\put(100,20){\line(-2,1){11.1}}   
\put(103.820,0){\arc{32.361}{5.176}{6.283}}
\put(120,0){\vector(0,-1){0}}
\put(100,20){\line(2,-1){11.1}}
\put(100,50){\vector(0,-1){10}}
\put(100,30){\line(0,-1){8}}
\put(100,18){\vector(0,-1){18}}
\thinlines
\put(95,30){\dashbox{2}(10,10){$V_-$}}
\put(106,45){\makebox(0,0){$\beta_-$}}
\put(93,5){\makebox(0,0){$\a^-_\la$}}
\put(126,4){\makebox(0,0){$\beta_+'$}}
\end{picture}
\end{center}
\caption{Naturality move for relative braiding}
\label{natreloi}
\end{figure}

Recall that we defined \cite[Def.\ 5.5]{BEK1}
a matrix $Z$ by setting
\[   Z_{\la,\mu}=\langle\a^+_\la,\a^-_\mu\rangle \,,
\qquad \la,\mu\in\NXN \,, \]
and we showed in \cite[Thm.\ 5.7]{BEK1} that it
commutes with Rehren's monodromy matrix $Y$ and
statistics T-matrix which have matrix elements
\begin{equation}
\label{YT}
Y_{\la,\mu}=\sum_{\nu}\frac{\om_\la\om_\mu}{\om_\nu}
N_{\la,\mu}^\nu d_\nu \,,\qquad
T_{\la,\mu}=\del\la\mu \E^{-\I\pi c/12} \om_\la \,,
\qquad \la,\mu\in\NXN \,,
\end{equation}
where $c=4\arg(\sum_\nu \om_\nu d_\nu^2)/\pi$.
As $Z$ has by definition
non-negative integer entries and satisfies
$Z_{0,0}=1$ (the label ``0'' stands as usual for
the identity morphism $\id\in\NXN$), it therefore
constitutes a modular invariant in the sense of
conformal field theory whenever the braiding is
non-degenerate because matrices $S=w^{-1/2}Y$ and $T$
obey the modular Verlinde algebra in that case
\cite{R0} (see also \cite{FG,FRS2} or our review
in \cite[Subsect.\ 2.2]{BEK1}).

\section{Chiral analysis}
\label{chiralanal}

In this section we begin to analyze the structure of the
chiral systems $\MXMpm$. So far the analysis will be carried
out without an assumption of non-degeneracy of the braiding,
and in fact several structures appear independently of it.

\subsection{Chiral horizontal projectors and chiral global indices}

Let $w_\pm=\sum_{\beta\in\MXMpm} d_\beta^2$. We call
$w_+$ and $w_-$ the {\sl chiral global indices}.
In the double triangle algebra, we define
$P^\pm=\sum_{\beta\in\MXMpm} e_\beta$.
We call (slightly different from Ocneanu's definition)
$P^+$ and $P^-$ chiral horizontal projectors.

\begin{proposition}
\label{sumcg}
In the $M$-$M$ fusion rule algebra we have
\begin{equation}
\label{firstsumcg}
\sum_{\la\in\NXN} d_\la [\a_\la^\pm]
= \frac w{w_\pm} \sum_{\beta\in\MXMpm} d_\beta [\beta]
\end{equation}
and consequently
$\sum_{\la\in\NXN} p_\la^\pm = w w_\pm^{-1} P^\pm$ in the double
triangle algebra. Moreover, the chiral global indices coincide
and are given by
\begin{equation}
w_+ = w_- = \frac w{\sum_{\la\in\NXN} d_\la Z_{\la,0}} \,.
\end{equation}
\end{proposition}

\begin{proof}
Put $\Gamma_{\la;\beta}^{\beta',\pm}=\langle\beta\a_\la^\pm,
\beta'\rangle$ for $\la\in\NXN$ and $\beta,\beta'\in\MXMpm$.
This defines square matrices $\Gamma_\la^\pm$ and we have
$\Gamma_\la^\pm=\sum_{\beta\in\MXMpm} \langle\beta,
\a_\la^\pm\rangle N_\beta$ where the $N_\beta$'s are the
fusion matrices of $\beta$ within $\MXMpm$.
With these, the matrices $\Gamma_\la^\pm$ therefore
share the simultaneous eigenvector $\vec{d}^\pm$,
defined by entries $d_\beta$, $\beta\in\MXMpm$,
with respective eigenvalues $d_\la$.
Note that the sum matrix $Q^\pm=\sum_\la \Gamma_\la^\pm$
is irreducible since each $[\beta]$ with $\beta\in\MXMpm$
is a subsector of some $[\a_\la^\pm]$ by definition.
Now define another vector $\vec{v}^\pm$ with entries
$v^\pm_\beta=\sum_\la d_\la \langle \beta,\a^\pm_\la \rangle$,
$\beta\in\MXMpm$. Note that all entries are positive.
We now compute
\[ \bearll
(\Gamma_\la^\pm \vec{v}^\pm)_\beta
&= \sum_{\beta'\in\MXMpm}
\sum_{\nu\in\NXN} \langle \beta \a^\pm_\la, \beta' \rangle
d_\nu \langle \beta', \a^\pm_\nu \rangle
=  \sum_{\nu\in\NXN} d_\nu \langle \beta \a^\pm_\la,
\a^\pm_\nu \rangle \\[.4em]
&= \sum_{\mu,\nu\in\NXN} N_{\nu,\co\la}^\mu d_\nu
\langle \beta, \a^\pm_\mu \rangle
= \sum_{\mu\in\NXN} d_\la d_\mu \langle \beta,
\a^\pm_\mu \rangle
= d_\la v^\pm_\beta \,,
\eear \]
i.e.\ $\Gamma_\la^\pm \vec{v}^\pm = d_\la \vec{v}^\pm$.
Hence $\vec{v}^\pm$ is another eigenvector of $Q$
with the same eigenvalue $\sum_\la d_\la$. By uniqueness
of the Perron-Frobenius eigenvector it follows
$v^\pm_\beta = \zeta_\pm d_\beta$ for all $\beta\in\MXMpm$
with some number $\zeta_\pm\in\bbC$. We can determine
this number in two different ways. We first find that now
$\sum_\la d_\la [\a^\pm_\la]
= \sum_{\beta\in\MXMpm} v^\pm_\beta [\beta]
= \sum_{\beta\in\MXMpm} \zeta_\pm d_\beta [\beta]$.
By computing the dimension we obtain
$w=\zeta_\pm w_\pm$, establishing \erf{firstsumcg}.
On the other hand we can compare the norms:
For $\vec{d}^\pm$ we have $\|\vec{d}^\pm\|=w_\pm$.
For $\vec{v}^\pm$ we compute
\[ \bearll
\|\vec{v}^\pm\|^2 &= \sum_{\beta\in\MXMpm} \sum_{\mu,\nu\in\NXN}
d_\mu d_\nu \langle \beta,\a^\pm_\mu \rangle
\langle \beta, \a^\pm_\nu \rangle
= \sum_{\mu,\nu\in\NXN} d_\mu d_\nu \langle \a^\pm_\mu,
\a^\pm_\nu \rangle \\[.4em]
&= \sum_{\la,\mu,\nu\in\NXN}
d_\mu d_\nu N_{\mu,\co\nu}^\la
\langle \a^\pm_\la,\id \rangle
= \sum_{\la,\mu\in\NXN} d_\mu^2 d_\la
\langle \a^\pm_\la,\id \rangle \,,
\eear \]
hence $\|\vec{v}^+\|^2=w\sum_{\la} d_\la Z_{\la,0}$
whereas $\|\vec{v}^-\|^2=w\sum_{\la} d_\la Z_{0,\la}$.
But note that
\[ \sum_{\la\in\NXN} d_\la Z_{\la,0}=(YZ)_{0,0}=(ZY)_{0,0}
=\sum_{\la\in\NXN} Z_{0,\la} d_\la \,. \]
We have found
\[ \frac{w^2}{w_\pm^2} = \zeta_\pm^2 =
\frac{\|\vec{v}^\pm\|^2}{\|\vec{d}^\pm\|^2} =
\frac w{w_\pm} \sum_{\la\in\NXN} d_\la Z_{\la,0}\,, \]
and this proves the proposition.
\end{proof}

Note that the $\Hom(a\co a,b\co b)$ part of the equality
$\sum_\la p_\la^+ = w w_+^{-1} P^+$ gives us the graphical
identity of Fig. \ref{scg=cp}.
%
\thinlines
\begin{figure}[htb]
\begin{center}
\unitlength 0.6mm
\begin{picture}(143,40)
\thicklines
\put(5,17){\makebox(0,0){$\displaystyle\sum_{\la}$}}
\put(25,5){\line(0,1){8}}
\put(25,17){\line(0,1){18}}
\put(45,5){\line(0,1){8}}
\put(45,17){\line(0,1){18}}
\Thicklines
\put(25,15){\line(1,0){20}}
\put(37,15){\vector(1,0){0}}
\thinlines
\put(25,20){\arc{10}{1.571}{4.712}}
\put(45,20){\arc{10}{4.712}{1.571}}
\put(25,25){\arc{5}{1.571}{4.712}}
\put(45,25){\arc{5}{4.712}{1.571}}
\put(20,35){\makebox(0,0){$a$}}
\put(50,35){\makebox(0,0){$a$}}
\put(20,5){\makebox(0,0){$b$}}
\put(50,5){\makebox(0,0){$b$}}
\put(35,21){\makebox(0,0){$\a_\la^+$}}
\thicklines
\put(85,17){\makebox(0,0){$=\;\displaystyle \frac w{w_+}
\sum_{\beta\in\MXMp}$}}
\put(115,5){\line(0,1){30}}
\put(135,5){\line(0,1){30}}
\Thicklines
\put(115,20){\line(1,0){20}}
\put(127,20){\vector(1,0){0}}
\thinlines
\put(135,20){\arc{5}{1.571}{4.712}}
\put(115,20){\arc{5}{4.712}{1.571}}
\put(110,35){\makebox(0,0){$a$}}
\put(140,35){\makebox(0,0){$a$}}
\put(110,5){\makebox(0,0){$b$}}
\put(140,5){\makebox(0,0){$b$}}
\put(125,25){\makebox(0,0){$\beta$}}
\end{picture}
\end{center}
\caption{Chiral generators sum up to chiral horizontal projectors}
\label{scg=cp}
\end{figure}
(And we obtain a similar identity for ``$-$''.)

We next claim the following

\begin{proposition}
\label{Zperm}
The following conditions are equivalent:
\begin{enumerate}
\item We have $Z_{0,\la}=\del \la0$.
\item We have $Z_{\la,0}=\del \la0$.
\item $Z$ is a permutation matrix,
$Z_{\la,\mu}=\del\la{\pi(\mu)}$ where
$\pi$ is a permutation of $\NXN$ satisfying
$\pi(0)=0$ and defining a fusion rule
automorphism of the $N$-$N$ fusion rule
algebra.
\end{enumerate}
\end{proposition}

\begin{proof}
The implication $1. \Rightarrow 2.$ follows
again from
$\sum_\la d_\la Z_{\la,0} = \sum_\la Z_{0,\la} d_\la$
arising from $[Y,Z]=0$. We next show the implication
$2. \Rightarrow 3.$: From Proposition \ref{sumcg}
we here obtain $w_+=w$. Because we have in general
$\MXMpm\subset\MXM$ we here find $\MXMpm=\MXM$.
Consequently,
$\sum_{\la\in\NXN} d_\la [\a^\pm_\la]=\sum_{\beta\in\MXM}
d_\beta[\beta]$ in the $M$-$M$ fusion rule algebra.
Assume for contradiction that some $[\a^\pm_\la]$ is
reducible. Then $d_\beta<d_\la$ if $[\beta]$ is a
subsector. But $[\beta]$ appears on the left hand
side with a coefficient larger or equal $d_\la$
whereas with coefficient $d_\beta$ on the right hand
side which cannot be true. Hence all $[\a_\la^\pm]$'s
are irreducible and as $w=w_+$, they must also be distinct.
Therefore $\#\MXMpm=\#\NXN$, and consequently
$Z$ must be a permutation matrix:
$Z_{\la,\mu}=\del\la{\pi(\mu)}$ with $\pi(0)=0$
as $Z_{0,0}=0$. Moreover, by virtue of the
homomorphism property of $\a$-induction we
have two isomorphisms
$\vartheta_\pm:[\la]\mapsto[\a^\pm_\la]$
from the $N$-$N$ into the $M$-$M$ fusion rule
algebra and consequently
$\vartheta_+^{-1}\circ\vartheta_-([\mu])=[\pi(\mu)]$
defines an automorphism of the $N$-$N$ fusion rules.
Finally, the implication $3. \Rightarrow 1.$ is trivial.
\end{proof}

Note that the statement of Proposition \ref{Zperm} is
well known for modular invariants in conformal field
theory \cite{GS,G1} (see also \cite{MS}).
However, it is remarkable that our statement does not
rely on the non-degeneracy of the braiding, i.e.\
it holds even if there is no representation of the
modular group around. An analogous statement has also
been derived recently for the coupling matrix
arising from the embedding of left and right chiral
observables into a ``canonical tensor product
subfactor'', not relying on modularity either
\cite{R4}. Yet our result turns up by considering
chiral observables only.

\subsection{Chiral branching coefficients}

We will now introduce the chiral branching coefficients
which play an important (twofold) role for the chiral
systems, analogous to the role of the entries of the
matrix $Z$ for the entire system.

\begin{lemma}
We have
\begin{equation}
\langle \beta, \a^\pm_\la \rangle = \frac w{d_\la d_\beta}
\varphi_h (p_\la^\pm *_h e_\beta)
\end{equation}
for any $\la\in\NXN$ and any $\beta\in\MXM$.
\end{lemma}

\begin{proof}
By \cite[Thm.\ 5.3]{BEK1} we have
\[ \frac 1{d_\la} p^\pm_\la = \sum_{\beta\in\MXM}
\frac 1{d_\beta} \langle\beta,\a^\pm_\la\rangle e_\beta \,,\]
hence
\[ \frac 1{d_\la d_\beta} p^\pm_\la *_h e_\beta =
\frac 1{d_\beta^2} \langle\beta,\a^\pm_\la\rangle e_\beta \,.\]
Application of $\varphi_h$ now yields the claim
since $\varphi_h(e_\beta)=d_\beta^2/w$
by \cite[Lemma 4.7]{BEK1}
\end{proof}

Hence the number $\langle\beta,\a^+_\la\rangle$
(and similarly $\langle\beta,\a^-_\la\rangle$) can be displayed
graphically as in Fig.\ \ref{bgraph} (cf.\ the argument to get
the picture for $Z_{\la,\mu}$ in \cite[Thm.\ 5.6]{BEK1}).
%
\begin{figure}[htb]
\begin{center}
\unitlength 0.6mm
\begin{picture}(110,50)
\thinlines
\put(32,23){\makebox(0,0){$\langle\beta,a^+_\la\rangle \;= \;
\displaystyle\sum_{b,c} \frac{d_b d_c}{w d_\la d_\beta}$}}
\put(80,5){\line(0,1){10}}
\put(80,45){\line(0,-1){10}}
\put(80,45){\arc{5}{0}{3.142}}
\put(100,35){\arc{5}{0}{3.142}}
\put(80,5){\arc{5}{3.142}{0}}
\put(100,15){\arc{5}{3.142}{0}}
\Thicklines
\put(80,15){\line(0,1){20}}
\put(100,15){\line(0,1){20}}
\put(80,23){\vector(0,-1){0}}
\put(100,27){\vector(0,1){0}}
\thicklines
\put(77,5){\line(1,0){26}}
\put(77,45){\line(1,0){26}}
\put(77,15){\line(1,0){1}}
\put(77,35){\line(1,0){1}}
\put(82,15){\line(1,0){21}}
\put(82,35){\line(1,0){21}}
\put(77,40){\arc{10}{1.571}{4.712}}
\put(103,40){\arc{10}{4.712}{1.571}}
\put(77,10){\arc{10}{1.571}{4.712}}
\put(103,10){\arc{10}{4.712}{1.571}}
\put(90,2){\makebox(0,0){$c$}}
\put(90,48){\makebox(0,0){$c$}}
\put(90,12){\makebox(0,0){$b$}}
\put(90,38){\makebox(0,0){$b$}}
\put(74,25.4){\makebox(0,0){$\a_\la^+$}}
\put(107,24.6){\makebox(0,0){$\beta$}}
\end{picture}
\end{center}
\caption{Graphical representation of $\langle \beta, \a^+_\la \rangle$}
\label{bgraph}
\end{figure}
For $\tau\in\MXMo$ we call the numbers
$b^\pm_{\tau,\la}=\langle\tau,\a^\pm_\la\rangle$
{\sl chiral branching coefficients}. Note that from
$Z_{\la,\mu}=\langle\a^+_\la,\a^-_\mu\rangle$
we obtain the formula
\begin{equation}
\label{Z=bb}
Z_{\la,\mu} = \sum_{\tau\in\MXMo}b^+_{\tau,\la}b^-_{\tau,\mu} \,.
\end{equation}
Introducing rectangular matrices $b^\pm$ with
entries $(b^\pm)_{\tau,\la}=b^\pm_{\tau,\la}$ we
can thus write $Z={}^{{\rm{t}}}\!b^+ b^-$. The
name ``chiral branching coefficients'' is motivated
from the case where chiral locality condition holds.
The canonical sector restriction \cite{LR} of
some morphism $\beta\in\Mor(M,M)$ is given by
$\sigma_\beta=\co\iota\beta\iota\in\Mor(N,N)$
and was named ``$\sigma$-restriction'' in
\cite{BE1}. Now suppose $\beta\in\MXM$. Then
$\sigma_\beta\in\Sigma(\NXN)$. We put
$b_{\tau,\la}=\langle\la,\sigma_\beta\rangle$
for $\la\in\NXN$. The following proposition is
just the version of $\alpha\sigma$-reciprocity
\cite[Thm.\ 3.21]{BE1} in our setting of
braided subfactors.

\begin{proposition}
\label{asreci}
Whenever the chiral locality condition
$\epsp\canr\canr\can(v)=\can(v)$ holds then
we have $b_{\tau,\la}^+=b_{\tau,\la}^-=b_{\tau,\la}$
for all $\tau\in\MXMo$, $\la\in\NXN$.
\end{proposition}

\begin{proof}
Using chiral locality, it was proven in
\cite[Prop.\ 3.3]{BE3} that
$\langle\a^\pm_\la,\beta\rangle=
\langle \la,\sigma_\beta\rangle$
whenever $[\beta]$ is a subsector of some
$[\a^\pm_\mu]$. Hence
\[ \langle\a^+_\la,\tau\rangle=
\langle \la,\sigma_\tau\rangle=
\langle\a^-_\la,\tau\rangle \]
for $\tau\in\MXMo$.
\end{proof}

Note that, with chiral locality, the modular invariant
matrix is written as
$Z_{\la,\mu}=\sum_{\tau\in\MXMo}b_{\tau,\la}b_{\tau,\mu}$,
and this is exactly the expression which characterizes
``block-diagonal'' or ``type \nolinebreak I'' invariants.
In fact, in the net of subfactor setting, the numbers
$\langle\la,\sigma_\beta\rangle$ describe the decomposition
of restricted representations $\pi_0\circ\beta$ as
established in \cite{LR}. For conformal inclusions
or simple current extensions treated in \cite{BE2},
the $b_{\tau,\la}$'s are exactly the branching coefficients
because the ambichiral system corresponds to the DHR
morphisms of the extended theory by the results of
\cite{BE3,BEK1}. Without chiral locality we only
have $b^\pm_{\tau,\la}\le b_{\tau,\la}$ similar
to the inequality
$\langle\a^\pm_\la,\a^\pm_\mu\rangle
\le\langle\canr\la,\mu\rangle$ which replaces
the ``main formula'' of \cite[Thm.\ 3.9]{BE1}.

\subsection{Chiral vertical algebras}

We define for each $\tau\in\MXMo$ a vector space
\[ \cA_\tau = \bigoplus_{a,b\in\NXM}
\Hom (a\tau \co a, b\tau\co b) \,, \]
and we endow it, similar to the double triangle algebra,
with a vertical product $\star_v$ defined graphically
in Fig.\ \ref{Atvprod}.
%
\begin{figure}[htb]
\begin{center}
\unitlength 0.6mm
\begin{picture}(210,50)
\thinlines
\put(5,20){\dashbox{2}(40,10){$X$}}
\thicklines
\put(10,30){\line(0,1){10}}
\put(40,30){\line(0,1){10}}
\put(10,10){\line(0,1){10}}
\put(40,10){\line(0,1){10}}
\Thicklines
\put(25,40){\vector(0,-1){10}}
\put(25,20){\vector(0,-1){10}}
\put(5,37){\makebox(0,0){$a$}}
\put(30,37){\makebox(0,0){$\tau$}}
\put(45,37){\makebox(0,0){$a$}}
\put(5,13){\makebox(0,0){$b$}}
\put(30,13){\makebox(0,0){$\tau$}}
\put(45,13){\makebox(0,0){$b$}}
\put(60,25){\makebox(0,0){$\star_v$}}
\thinlines
\put(75,20){\dashbox{2}(40,10){$Y$}}
\thicklines
\put(80,30){\line(0,1){10}}
\put(110,30){\line(0,1){10}}
\put(80,10){\line(0,1){10}}
\put(110,10){\line(0,1){10}}
\Thicklines
\put(95,40){\vector(0,-1){10}}
\put(95,20){\vector(0,-1){10}}
\put(75,37){\makebox(0,0){$c$}}
\put(100,37){\makebox(0,0){$\tau$}}
\put(115,37){\makebox(0,0){$c$}}
\put(75,13){\makebox(0,0){$d$}}
\put(100,13){\makebox(0,0){$\tau$}}
\put(115,13){\makebox(0,0){$d$}}
\put(141,25){\makebox(0,0){$= \, \del ad \,d_a$}}
\thinlines
\put(165,30){\dashbox{2}(40,10){$Y$}}
\put(165,10){\dashbox{2}(40,10){$X$}}
\thicklines
\put(170,40){\line(0,1){10}}
\put(200,40){\line(0,1){10}}
\put(170,20){\line(0,1){10}}
\put(200,20){\line(0,1){10}}
\put(170,0){\line(0,1){10}}
\put(200,0){\line(0,1){10}}
\Thicklines
\put(185,50){\vector(0,-1){10}}
\put(185,30){\vector(0,-1){10}}
\put(185,10){\vector(0,-1){10}}
\put(165,47){\makebox(0,0){$c$}}
\put(190,47){\makebox(0,0){$\tau$}}
\put(205,47){\makebox(0,0){$c$}}
\put(165,25){\makebox(0,0){$a$}}
\put(190,25){\makebox(0,0){$\tau$}}
\put(205,25){\makebox(0,0){$a$}}
\put(165,3){\makebox(0,0){$b$}}
\put(190,3){\makebox(0,0){$\tau$}}
\put(205,3){\makebox(0,0){$b$}}
\end{picture}
\end{center}
\caption{Vertical product for $\cA_\tau$}
\label{Atvprod}
\end{figure}
Then it is not hard to see that a complete set of matrix units
is given by elements $f_{\la;b,d,k,l}^{a,c,i,j}$ as defined in
Fig.\ \ref{matunAtfig}.
%
\begin{figure}[htb]
\begin{center}
\unitlength 0.6mm
\begin{picture}(117,75)
\thicklines
\put(32,37.5){\makebox(0,0){$f_{\la;b,d,i,k}^{a,c,j,l}\;
=\;\displaystyle\frac 1{d_a d_b} \sqrt{\frac{d_\la}{d_\tau}}$}}
\put(80,0){\line(0,1){5}}
\put(110,0){\line(0,1){20}}
\put(80,5){\line(1,1){10}}
\put(90,15){\line(0,1){5}}
\put(90,20){\line(1,1){10}}
\put(110,20){\line(-1,1){10}}
\put(80,75){\line(0,-1){5}}
\put(110,75){\line(0,-1){20}}
\put(80,70){\line(1,-1){10}}
\put(90,60){\line(0,-1){5}}
\put(90,55){\line(1,-1){10}}
\put(110,55){\line(-1,-1){10}}
\Thicklines
\put(100,5){\vector(0,-1){5}}
\put(100,75){\line(0,-1){5}}
\put(100,5){\line(-1,1){10}}
\put(100,70){\vector(-1,-1){10}}
\thinlines
\put(100,30){\line(0,1){15}}
\put(100,35.5){\vector(0,-1){0}}
\put(76,71){\makebox(0,0){$a$}}
\put(114,71){\makebox(0,0){$a$}}
\put(100,65){\makebox(0,0){$\tau$}}
\put(90,50){\makebox(0,0){$c$}}
\put(105,37.5){\makebox(0,0){$\la$}}
\put(76,4){\makebox(0,0){$b$}}
\put(114,4){\makebox(0,0){$b$}}
\put(100,10){\makebox(0,0){$\tau$}}
\put(90,25){\makebox(0,0){$d$}}
\put(90,5){\makebox(0,0){\footnotesize{$T_{b,\tau}^{d;i}$}}}
\put(100,20){\makebox(0,0){\footnotesize{$t_{d,\co b}^{\la;k}$}}}
\put(90,71){\makebox(0,0){\footnotesize{$(T_{b,\tau}^{d;j})^*$}}}
\put(100,56){\makebox(0,0){\footnotesize{$(t_{d,\co b}^{\la;l})^*$}}}
\end{picture}
\end{center}
\caption{Matrix units for $\cA_\tau$}
\label{matunAtfig}
\end{figure}
They obey
\[ f_{\la;b,d,i,k}^{a,c,j,l} \star_v f_{\la';b',d',i',k'}^{a',c',j',l'}
= \del \la{\la'} \del a{b'} \del c{d'} \del j{i'} \del l{k'} \,
f_{\la;b,d,i,k}^{a',c',j',l'} \,. \]
We define a functional $\psi^\tau_v$ as in Fig.\ \ref{traceAt}.
%
\begin{figure}[htb]
\begin{center}
\unitlength 0.6mm
\begin{picture}(197,40)
\thinlines
\put(7,20){\makebox(0,0){$\psi^\tau_v \; :$}}
\put(25,15){\dashbox{2}(40,10){$X$}}
\thicklines
\put(30,25){\line(0,1){10}}
\put(60,25){\line(0,1){10}}
\put(30,5){\line(0,1){10}}
\put(60,5){\line(0,1){10}}
\Thicklines
\put(45,35){\vector(0,-1){10}}
\put(45,15){\vector(0,-1){10}}
\put(25,32){\makebox(0,0){$a$}}
\put(50,32){\makebox(0,0){$\tau$}}
\put(65,32){\makebox(0,0){$a$}}
\put(25,8){\makebox(0,0){$b$}}
\put(50,8){\makebox(0,0){$\tau$}}
\put(65,8){\makebox(0,0){$b$}}
\thinlines
\put(100,20){\makebox(0,0){$\longmapsto \; \del ab \, d_a$}}
\put(125,15){\dashbox{2}(40,10){$X$}}
\thicklines
\put(165,25){\arc{10}{3.142}{0}}
\put(170,15){\line(0,1){10}}
\put(165,15){\arc{10}{0}{3.142}}
\put(145,25){\arc{30}{3.142}{4.712}}
\put(145,40){\line(1,0){30}}
\put(175,25){\arc{30}{4.712}{0}}
\put(190,15){\line(0,1){10}}
\put(175,15){\arc{30}{0}{1.571}}
\put(145,0){\line(1,0){30}}
\put(145,15){\arc{30}{1.571}{3.142}}
\Thicklines
\put(155,25){\arc{20}{3.142}{4.712}}
\put(155,35){\line(1,0){15}}
\put(170,25){\arc{20}{4.712}{0}}
\put(180,15){\line(0,1){10}}
\put(180,20){\vector(0,1){0}}
\put(170,15){\arc{20}{0}{1.571}}
\put(155,5){\line(1,0){15}}
\put(155,15){\arc{20}{1.571}{3.142}}
\put(174,20){\makebox(0,0){$a$}}
\put(184,20){\makebox(0,0){$\tau$}}
\put(194,20){\makebox(0,0){$a$}}
\end{picture}
\end{center}
\caption{Trace for $\cA_\tau$}
\label{traceAt}
\end{figure}
It fulfills
$\psi^\tau_v(f_{\la;b,d,i,k}^{a,c,j,l})
= \del ab \del cd \del ij \del kl d_\la$,
and therefore it is a faithful (un-normalized)
trace on $\cA_\tau$.
We next define vector spaces $\cH_{\tau,\la}$ by
\[ \cH_{\tau,\la} = \bigoplus_{a\in\NXM} \Hom (\la,a\tau\co a) \,,\]
and special vectors
$\omega^{\tau,\la,+}_{b,c,t,X}\in\cH_{\tau,\la}$ and
$\omega^{\tau,\la,-}_{b,c,t,X}\in\cH_{\co\tau,\co\la}$
as given in Fig.\ \ref{ombctsX}.
%
\begin{figure}[htb]
\begin{center}
\unitlength 0.6mm
\begin{picture}(219,51)
\thicklines
\put(22,24){\makebox(0,0){$\omega^{\tau,\la,+}_{b,c,t,X}\;
=\;\displaystyle\sum_a$}}
\put(52,0){\line(0,1){26}}
\put(62,26){\arc{20}{3.142}{4.712}}
\put(62,36){\line(1,0){6}}
\put(72,36){\line(1,0){6}}
\put(78,26){\arc{20}{4.712}{0}}
\put(88,0){\line(0,1){26}}
\put(62,16){\line(0,1){10}}
\put(67,26){\arc{10}{3.142}{4.712}}
\put(67,31){\line(1,0){1}}
\put(72,31){\line(1,0){1}}
\put(73,26){\arc{10}{4.712}{0}}
\put(78,16){\line(0,1){10}}
\put(73,16){\arc{10}{0}{1.571}}
\put(67,11){\line(1,0){6}}
\put(67,16){\arc{10}{1.571}{3.142}}
\Thicklines
\put(70,0){\line(0,1){11}}
\put(70,3.5){\vector(0,-1){0}}
\put(70,31){\line(0,1){5}}
\thinlines
\put(62,21){\line(1,0){3}}
\put(65,26){\arc{10}{0}{1.571}}
\put(70,26){\line(0,1){5}}
\put(70,36){\line(0,1){15}}
\put(70,41.5){\vector(0,-1){0}}
\put(48,4){\makebox(0,0){$a$}}
\put(62,10){\makebox(0,0){$c$}}
\put(82,21){\makebox(0,0){$b$}}
\put(75,4){\makebox(0,0){$\tau$}}
\put(75,43.5){\makebox(0,0){$\la$}}
\put(57,21){\makebox(0,0){$t$}}
\put(70,15){\makebox(0,0){$X^*$}}
\thicklines
\put(152,24){\makebox(0,0){$\omega^{\tau,\la,-}_{b,c,t,X}\;
=\;\displaystyle\sum_a$}}
\put(182,0){\line(0,1){26}}
\put(192,26){\arc{20}{3.142}{4.712}}
\put(192,36){\line(1,0){16}}
\put(208,26){\arc{20}{4.712}{0}}
\put(218,0){\line(0,1){26}}
\put(192,16){\line(0,1){10}}
\put(197,26){\arc{10}{3.142}{4.712}}
\put(197,31){\line(1,0){6}}
\put(203,26){\arc{10}{4.712}{0}}
\put(208,16){\line(0,1){10}}
\put(203,16){\arc{10}{0}{1.571}}
\put(197,11){\line(1,0){6}}
\put(197,16){\arc{10}{1.571}{3.142}}
\Thicklines
\put(200,0){\line(0,1){11}}
\put(200,7.5){\vector(0,1){0}}
\put(200,33){\line(0,1){1}}
\thinlines
\put(208,21){\line(-1,0){3}}
\put(205,26){\arc{10}{1.571}{3.142}}
\put(200,26){\line(0,1){3}}
\put(200,38){\line(0,1){13}}
\put(200,44.5){\vector(0,1){0}}
\put(178,4){\makebox(0,0){$a$}}
\put(208,10){\makebox(0,0){$c$}}
\put(188,21){\makebox(0,0){$b$}}
\put(195,4){\makebox(0,0){$\tau$}}
\put(205,43.5){\makebox(0,0){$\la$}}
\put(213,21){\makebox(0,0){$t^*$}}
\put(200,15){\makebox(0,0){$X$}}
\end{picture}
\end{center}
\caption{The vectors 
$\omega^{\tau,\la,+}_{b,c,t,X}\in\cH_{\tau,\la}$ and
$\omega^{\tau,\la,-}_{b,c,t,X}\in\cH_{\co\tau,\co\la}$}
\label{ombctsX}
\end{figure}
Note that such vectors may be linearly dependent.
Let $H_{\tau,\la}^+\subset\cH_{\tau,\la}$ respectively
$H_{\tau,\la}^-\subset\cH_{\co\tau,\co\la}$ be the
subspaces spanned by vectors  $\omega^{\tau,\la,+}_{b,c,t,X}$
respectively $\omega^{\tau,\la,-}_{b,c,t,X}$
where $b,c\in\NXM$ and $t\in\Hom(\la,b\co c)$ and
$X\in\Hom(\tau,\co c b)$ are isometries. Now take
such vectors $\omega^{\tau,\la,\pm}_{b,c,t,X}$ and
$\omega^{\tau,\la,\pm}_{b',c',t',X'}$. We define
an element 
$|\omega^{\tau,\la,+}_{b',c',t',X'}\rangle
\langle\omega^{\tau,\la,+}_{b,c,t,X}|\in\cA_\tau$
by the diagram in Fig.\ \ref{|om><om|}.
%
\begin{figure}[htb]
\begin{center}
\unitlength 0.6mm
\begin{picture}(129,87)
\thicklines
\put(38,41){\makebox(0,0){$|\omega^{\tau,\la,+}_{b',c',t',X'}
\rangle\langle\omega^{\tau,\la,+}_{b,c,t,X}|\;
=\;\displaystyle\sum_{a,d}$}}
\put(92,0){\line(0,1){26}}
\put(92,87){\line(0,-1){26}}
\put(102,26){\arc{20}{3.142}{4.712}}
\put(102,61){\arc{20}{1.571}{3.142}}
\put(102,36){\line(1,0){6}}
\put(112,36){\line(1,0){6}}
\put(102,51){\line(1,0){6}}
\put(112,51){\line(1,0){6}}
\put(118,26){\arc{20}{4.712}{0}}
\put(118,61){\arc{20}{0}{1.571}}
\put(128,0){\line(0,1){26}}
\put(128,87){\line(0,-1){26}}
\put(102,16){\line(0,1){10}}
\put(102,71){\line(0,-1){10}}
\put(107,26){\arc{10}{3.142}{4.712}}
\put(107,61){\arc{10}{1.571}{3.142}}
\put(107,31){\line(1,0){1}}
\put(112,31){\line(1,0){1}}
\put(107,56){\line(1,0){1}}
\put(112,56){\line(1,0){1}}
\put(113,26){\arc{10}{4.712}{0}}
\put(113,61){\arc{10}{0}{1.571}}
\put(118,16){\line(0,1){10}}
\put(118,71){\line(0,-1){10}}
\put(113,16){\arc{10}{0}{1.571}}
\put(113,71){\arc{10}{4.712}{0}}
\put(107,11){\line(1,0){6}}
\put(107,76){\line(1,0){6}}
\put(107,16){\arc{10}{1.571}{3.142}}
\put(107,71){\arc{10}{3.142}{4.712}}
\Thicklines
\put(110,0){\line(0,1){11}}
\put(110,3.5){\vector(0,-1){0}}
\put(110,31){\line(0,1){5}}
\put(110,87){\line(0,-1){11}}
\put(110,80.5){\vector(0,-1){0}}
\put(110,56){\line(0,-1){5}}
\thinlines
\put(102,21){\line(1,0){3}}
\put(105,26){\arc{10}{0}{1.571}}
\put(110,26){\line(0,1){5}}
\put(110,36){\line(0,1){15}}
\put(102,66){\line(1,0){3}}
\put(105,61){\arc{10}{4.712}{0}}
\put(110,61){\line(0,-1){5}}
\put(110,41.5){\vector(0,-1){0}}
\put(88,4){\makebox(0,0){$d$}}
\put(102,9){\makebox(0,0){$c'$}}
\put(122,21){\makebox(0,0){$b'$}}
\put(115,4){\makebox(0,0){$\tau$}}
\put(115,43.5){\makebox(0,0){$\la$}}
\put(97,21){\makebox(0,0){$t'$}}
\put(110.5,15){\makebox(0,0){\footnotesize{$(X')^*$}}}
\put(88,83){\makebox(0,0){$a$}}
\put(102,77){\makebox(0,0){$c$}}
\put(122,66){\makebox(0,0){$b$}}
\put(115,83){\makebox(0,0){$\tau$}}
\put(97,66){\makebox(0,0){$t^*$}}
\put(110,72){\makebox(0,0){\footnotesize{$X$}}}
\end{picture}
\end{center}
\caption{The elements $|\omega^{\tau,\la,+}_{b',c',t',X'}\rangle
\langle\omega^{\tau,\la,+}_{b,c,t,X}|\in\cA_\tau$}
\label{|om><om|}
\end{figure}
Analogously we define
$|\omega^{\tau,\la,-}_{b',c',t',X'}\rangle
\langle\omega^{\tau,\la,-}_{b,c,t,X}|\in\cA_{\co\tau}$.
Choosing orthonormal bases of isometries
$t_{b,\co c}^{\la;i}\in\Hom(\la,b\co c)$ and
$X_{\co c,b}^{\tau,j}\in\Hom(\tau,\co c b)$
we sometimes abbreviate
$\omega^{\tau,\la,\pm}_{b,c,i,j}
=\omega^{\tau,\la,\pm}_{b,c,t_{b,\co c}^{\la;i},X_{\co c,b}^{\tau,j}}$
and we also use the notation
$\omega^{\tau,\la,\pm}_\xi=\omega^{\tau,\la,\pm}_{b,c,i,j}$
with some multi-index $\xi=(b,c,i,j)$.
For vectors $\varphi^{\tau,\la,\pm}_\ell\in H_{\tau,\la}^\pm$
with expansions
$\varphi^{\tau,\la,\pm}_\ell
=\sum_\xi c^\xi_{\ell,\pm} \omega^{\tau,\la,\pm}_\xi$,
$\ell=1,2$, we define elements
$|\varphi^{\tau,\la,+}_1\rangle\langle\varphi^{\tau,\la,+}_2|
\in\cA_\tau$
and
$|\varphi^{\tau,\la,-}_1\rangle\langle\varphi^{\tau,\la,-}_2|
\in\cA_{\co\tau}$
by
\begin{equation}
|\varphi^{\tau,\la,\pm}_1\rangle\langle\varphi^{\tau,\la,\pm}_2|
=\sum_{\xi,\xi'} c^\xi_{1,\pm}(c^{\xi'}_{2,\pm})^*
|\omega^{\tau,\la,\pm}_\xi\rangle
\langle\omega^{\tau,\la,\pm}_{\xi'}| \,,
\label{|phi><phi|}
\end{equation}
and scalars
$\langle\varphi^{\tau,\la,\pm}_2,\varphi^{\tau,\la,\pm}_1\rangle
\in\bbC$
by
\begin{equation}
\bearl
\langle\varphi^{\tau,\la,+}_2,\varphi^{\tau,\la,+}_1\rangle
= \displaystyle\frac 1{d_\la}
\psi^\tau_v(|\varphi^{\tau,\la,+}_1\rangle
\langle\varphi^{\tau,\la,+}_2|) \,, \\[.8em]
\langle\varphi^{\tau,\la,-}_2,\varphi^{\tau,\la,-}_1\rangle
= \displaystyle\frac 1{d_\la}
\psi^{\co\tau}_v(|\varphi^{\tau,\la,-}_1\rangle
\langle\varphi^{\tau,\la,-}_2|) \,.
\eear
\label{<phi,phi>}
\end{equation}

\begin{lemma}
\label{positive}
\erf{|phi><phi|} extends to positive definite sesqui-linear maps
$H_{\tau,\la}^+ \times H_{\tau,\la}^+ \rightarrow \cA_\tau$ and
$H_{\tau,\la}^- \times H_{\tau,\la}^- \rightarrow \cA_{\co\tau}$.
Consequently, \erf{<phi,phi>} defines scalar products turning
$H_{\tau,\la}^\pm$ into Hilbert spaces.
\end{lemma}

\begin{proof}
(Similar to the proof of \cite[Lemma 6.1]{BEK1}.)
We only show it for ``$+$''; the proof for ``$-$'' is
analogous. As in particular
$\varphi^{\tau,\la,+}_\ell\in\cH_{\tau,\la}$,
we can write $\varphi^{\tau,\la,+}_\ell=
\bigoplus_a (\varphi^{\tau,\la,+}_\ell)_a$ with
$(\varphi^{\tau,\la,+}_\ell)_a\in\Hom(\la,a\tau\co a)$
according to the direct sum structure of
$\cH_{\tau,\la}$, $\ell=1,2$.
Assume $\varphi^{\tau,\la,+}_1=0$.
Then clearly $(\varphi^{\tau,\la,+}_1)_a=0$ for all $a$.
Now the $\Hom(a\tau\co a,b \tau\co b)$ part of
$|\varphi^{\tau,\la,+}_1\rangle
\langle\varphi^{\tau,\la,+}_2|\in\cA_\tau$
is given by 
$(\varphi^{\tau,\la,+}_1)_b (\varphi^{\tau,\la,+}_2)_a^*$,
hence
$|\varphi^{\tau,\la,+}_1\rangle
\langle\varphi^{\tau,\la,+}_2|=0$.
A similar argument applies to $\varphi^{\tau,\la,+}_2$,
and hence the element
$|\varphi^{\tau,\la,+}_1\rangle
\langle\varphi^{\tau,\la,+}_2|\in\cA_\tau$
is independent of the linear expansions of the
$\varphi^{\tau,\la,+}_\ell$'s.
Therefore \erf{|phi><phi|} defines a sesqui-linear map
$H_{\tau,\la}^+\times H_{\tau,\la}^+\rightarrow\cA_\tau$.
Now assume 
$|\varphi^{\tau,\la,+}_1\rangle
\langle\varphi^{\tau,\la,+}_1|=0$.
Then in particular
$(\varphi^{\tau,\la,+}_1)_a (\varphi^{\tau,\la,+}_1)_a^*=0$
for all $a\in\NXM$,
and hence $\varphi^{\tau,\la,+}_1=0$,
proving strict positivity. That the
sesqui-linear form $\langle\cdot,\cdot\rangle$ on $H_{\tau,\la}^+$
is non-degenerate follows now from positive definiteness of
$\psi_v^\tau$.
\end{proof}

Note that the scalar product
$\langle\omega^{\tau,\la,+}_{b,c,t,X},
\omega^{\tau,\la,+}_{b',c',t',X'}\rangle$
is given graphically as in Fig.\ \ref{<om,om>}.
%
\begin{figure}[htb]
\begin{center}
\unitlength 0.6mm
\begin{picture}(141,70)
\thicklines
\put(41,33){\makebox(0,0){$\langle\omega^{\tau,\la,+}_{b,c,t,X},
\omega^{\tau,\la,+}_{b',c',t',X'}\rangle
\;=\;\displaystyle\frac w{d_\la}$}}
\put(102,15){\line(0,1){10}}
\put(102,55){\line(0,-1){10}}
\put(107,25){\arc{10}{3.142}{4.712}}
\put(107,45){\arc{10}{1.571}{3.142}}
\put(107,30){\line(1,0){1}}
\put(112,30){\line(1,0){1}}
\put(107,40){\line(1,0){1}}
\put(112,40){\line(1,0){1}}
\put(113,25){\arc{10}{4.712}{0}}
\put(113,45){\arc{10}{0}{1.571}}
\put(118,15){\line(0,1){10}}
\put(118,55){\line(0,-1){10}}
\put(113,15){\arc{10}{0}{1.571}}
\put(113,55){\arc{10}{4.712}{0}}
\put(107,10){\line(1,0){6}}
\put(107,60){\line(1,0){6}}
\put(107,15){\arc{10}{1.571}{3.142}}
\put(107,55){\arc{10}{3.142}{4.712}}
\Thicklines
\put(130,10){\line(0,1){50}}
\put(130,37){\vector(0,1){0}}
\put(110,30){\line(0,1){10}}
\put(110,33){\vector(0,-1){0}}
\put(120,10){\arc{20}{0}{3.142}}
\put(120,60){\arc{20}{3.142}{0}}
\thinlines
\put(102,20){\line(1,0){3}}
\put(105,25){\arc{10}{0}{1.571}}
\put(110,25){\line(0,1){5}}
\put(102,50){\line(1,0){3}}
\put(105,45){\arc{10}{4.712}{0}}
\put(110,45){\line(0,-1){5}}
\put(102,8){\makebox(0,0){$c'$}}
\put(122,20){\makebox(0,0){$b'$}}
\put(135,35){\makebox(0,0){$\tau$}}
\put(103,35){\makebox(0,0){$\a^+_\la$}}
\put(97,20){\makebox(0,0){$t'$}}
\put(110.5,14){\makebox(0,0){\footnotesize{$(X')^*$}}}
\put(102,61){\makebox(0,0){$c$}}
\put(122,50){\makebox(0,0){$b$}}
\put(97,50){\makebox(0,0){$t^*$}}
\put(110,56){\makebox(0,0){\footnotesize{$X$}}}
\end{picture}
\end{center}
\caption{The scalar product $\langle\omega^{\tau,\la,+}_{b,c,t,X},
\omega^{\tau,\la,+}_{b',c',t',X'}\rangle$}
\label{<om,om>}
\end{figure}
Here we pulled out a closed wire $a$ so that the summation
over $a$ produced together with the prefactor $d_a$ just
the global index $w$.
We define subspaces $A^+_{\tau,\la}\subset\cA_\tau$ respectively
$A^-_{\tau,\la}\subset\cA_{\co\tau}$ given
as the linear span of elements
$|\omega^{\tau,\la,+}_{b',c',t',X'}\rangle
\langle\omega^{\tau,\la,+}_{b,c,t,X}|$ respectively
$|\omega^{\tau,\la,-}_{b',c',t',X'}\rangle
\langle\omega^{\tau,\la,-}_{b,c,t,X}|$.

\begin{lemma}
\label{lemchirkey}
We have the identity of Fig.\ \ref{chirkey} for
intertwiners in $\Hom(\la',\la)$. An analogous
identity can be established using vectors
$\omega^{\tau,\la,-}_{b,c,t,X}\in H_{\tau,\la}^-$.
\end{lemma}
%
\begin{figure}[htb]
\begin{center}
\unitlength 0.6mm
\begin{picture}(180,80)
\thicklines
\put(11,37){\makebox(0,0){$\displaystyle\sum_a \; d_a$}}
\put(32,20){\line(0,1){40}}
\put(42,60){\arc{20}{3.142}{4.712}}
\put(42,20){\arc{20}{1.571}{3.142}}
\put(42,10){\line(1,0){6}}
\put(52,10){\line(1,0){6}}
\put(42,70){\line(1,0){6}}
\put(52,70){\line(1,0){6}}
\put(58,60){\arc{20}{4.712}{0}}
\put(58,20){\arc{20}{0}{1.571}}
\put(68,20){\line(0,1){40}}
\put(42,20){\line(0,1){10}}
\put(42,60){\line(0,-1){10}}
\put(47,30){\arc{10}{3.142}{4.712}}
\put(47,50){\arc{10}{1.571}{3.142}}
\put(47,15){\line(1,0){1}}
\put(52,15){\line(1,0){1}}
\put(47,65){\line(1,0){1}}
\put(52,65){\line(1,0){1}}
\put(53,30){\arc{10}{4.712}{0}}
\put(53,50){\arc{10}{0}{1.571}}
\put(58,20){\line(0,1){10}}
\put(58,60){\line(0,-1){10}}
\put(53,20){\arc{10}{0}{1.571}}
\put(53,60){\arc{10}{4.712}{0}}
\put(47,35){\line(1,0){6}}
\put(47,45){\line(1,0){6}}
\put(47,20){\arc{10}{1.571}{3.142}}
\put(47,60){\arc{10}{3.142}{4.712}}
\Thicklines
\put(50,35){\line(0,1){10}}
\put(50,38){\vector(0,-1){0}}
\put(50,10){\line(0,1){5}}
\put(50,70){\line(0,-1){5}}
\thinlines
\put(42,55){\line(1,0){3}}
\put(45,60){\arc{10}{0}{1.571}}
\put(50,60){\line(0,1){5}}
\put(50,70){\line(0,1){10}}
\put(50,0){\line(0,1){10}}
\put(42,25){\line(1,0){3}}
\put(45,20){\arc{10}{4.712}{0}}
\put(50,20){\line(0,-1){5}}
\put(50,73){\vector(0,-1){0}}
\put(50,3){\vector(0,-1){0}}
\put(28,20){\makebox(0,0){$a$}}
\put(42,44){\makebox(0,0){$c'$}}
\put(62,55){\makebox(0,0){$b'$}}
\put(55,40){\makebox(0,0){$\tau$}}
\put(55,4){\makebox(0,0){$\la$}}
\put(55,76){\makebox(0,0){$\la'$}}
\put(37,55){\makebox(0,0){$t'$}}
\put(50.5,49){\makebox(0,0){\footnotesize{$(X')^*$}}}
\put(42,35){\makebox(0,0){$c$}}
\put(62,25){\makebox(0,0){$b$}}
\put(37,25){\makebox(0,0){$t^*$}}
\put(50,31){\makebox(0,0){\footnotesize{$X$}}}
\put(125,40){\makebox(0,0){$=\; \del \la{\la'}\,
\langle\omega^{\tau,\la,+}_{b,c,t,X},
\omega^{\tau,\la,+}_{b',c',t',X'}\rangle$}}
\thinlines
\put(170,0){\line(0,1){80}}
\put(170,38){\vector(0,-1){0}}
\put(175,40){\makebox(0,0){$\la$}}
\end{picture}
\end{center}
\caption{An identity in $\Hom(\la',\la)$}
\label{chirkey}
\end{figure}
\begin{proof}
(Similar to the proof of \cite[Lemma 6.2]{BEK1}.)
It is clear that we obtain a scalar which is
zero unless $\la=\la'$. To compute the scalar, we put
$\la=\la'$ and then we can close the wire $\la$ on the left
hand side, what has to be compensated by a factor
$1/d_\la$. We can now open the wire $a$ on the right and
close it on the left, and this way we can pull out
the wire $a$, yielding a closed loop. Hence the summation
over $a$ gives the global index, and the resulting
picture is regularly isotopic to Fig.\ \ref{<om,om>}.
\end{proof}

\begin{corollary}
\label{matunAt}
The subspaces $A_{\tau,\la}^+\subset\cA_\tau$ and
$A_{\tau,\la}^-\subset\cA_{\co\tau}$ are in fact
subalgebras. Moreover, in $\cA_\tau$ respectively
$\cA_{\co\tau}$ we have multiplication rules
\begin{equation}
|\varphi^{\tau,\la,\pm}_1\rangle
\langle\varphi^{\tau,\la,\pm}_2| \star_v
|\varphi^{\tau,\mu,\pm}_3\rangle
\langle\varphi^{\tau,\mu,\pm}_4| = \del \la\mu
\langle\varphi^{\tau,\la,\pm}_2,
\varphi^{\tau,\la,\pm}_3\rangle \,
|\varphi^{\tau,\la,\pm}_1\rangle
\langle\varphi^{\tau,\la,\pm}_4| \,,
\end{equation}
$\varphi^{\tau,\la,\pm}_\ell\in H_{\tau,\la}^\pm$, $\ell=1,2,3,4$.
Consequently, we have subalgebras
$A_{\tau}^+\subset\cA_\tau$ and $A_{\tau}^-\subset\cA_{\co\tau}$
given as the direct sums $A_\tau^\pm=\bigoplus_\la A_{\tau,\la}^\pm$.
We can choose orthonormal bases
$\{ u^{\tau,\la,\pm}_i \}_{i=1}^{\dim H_{\tau,\la}^\pm}$
of $H_{\tau,\la}^\pm$ to obtain systems of matrix units
$\{ |u^{\tau,\la,\pm}_i \rangle
\langle u^{\tau,\la,\pm}_j| \}_{\la,i,j}$
in $A_\tau^\pm$.
\end{corollary}
We call the algebras $A_\tau^\pm$, $\tau\in\MXMo$,
{\sl chiral vertical algebras}. Next we define elements
$I_\tau^+\in\cA_\tau$ and $I_\tau^-\in\cA_{\co\tau}$
by the diagrams in Fig.\ \ref{chirmunit} and we call them
{\sl chiral multiplicative units}
(for reasons given below).
%
\begin{figure}[htb]
\begin{center}
\unitlength 0.6mm
\begin{picture}(236,30)
\thicklines
\put(31,13){\makebox(0,0){$I_\tau^+\;=\;\displaystyle\frac 1{w_+}
\sum_{\beta\in\MXMm} \sum_{a,b}$}}
\put(70,0){\line(0,1){30}}
\put(100,0){\line(0,1){30}}
\thinlines
\put(70,15){\arc{5}{4.712}{1.571}}
\put(100,15){\arc{5}{1.571}{4.712}}
\Thicklines
\put(85,30){\vector(0,-1){30}}
\put(70,15){\line(1,0){13}}
\put(87,15){\vector(1,0){13}}
\put(93,20){\makebox(0,0){$\beta$}}
\put(89,4){\makebox(0,0){$\tau$}}
\put(66,26){\makebox(0,0){$a$}}
\put(104,26){\makebox(0,0){$a$}}
\put(66,4){\makebox(0,0){$b$}}
\put(104,4){\makebox(0,0){$b$}}
\thicklines
\put(161,13){\makebox(0,0){$I_\tau^-\;=\;\displaystyle\frac 1{w_+}
\sum_{\beta\in\MXMp} \sum_{a,b}$}}
\put(200,0){\line(0,1){30}}
\put(230,0){\line(0,1){30}}
\thinlines
\put(200,15){\arc{5}{4.712}{1.571}}
\put(230,15){\arc{5}{1.571}{4.712}}
\Thicklines
\put(200,15){\vector(1,0){30}}
\put(215,17){\vector(0,1){13}}
\put(215,13){\line(0,-1){13}}
\put(223,20){\makebox(0,0){$\beta$}}
\put(211,26){\makebox(0,0){$\tau$}}
\put(196,26){\makebox(0,0){$a$}}
\put(234,26){\makebox(0,0){$a$}}
\put(196,4){\makebox(0,0){$b$}}
\put(234,4){\makebox(0,0){$b$}}
\end{picture}
\end{center}
\caption{Chiral multiplicative units $I_\tau^\pm$}
\label{chirmunit}
\end{figure}
We then claim the following

\begin{lemma}
\label{chirality}
We have
\begin{equation}
I_\tau^\pm = \frac 1{w^2\sqrt{d_\tau}} \sum_{\la,\xi} \sqrt{d_\la}
|\omega^{\tau,\la,\pm}_\xi\rangle\langle\omega^{\tau,\la,\pm}_\xi|\,.
\end{equation}
\end{lemma}

\begin{proof}
We compute the sum
\[ \sum_{\la,b,c,i,j}\frac{\sqrt{d_\la}}{w^2 \sqrt{d_\tau}}
|\omega^{\tau,\la,+}_{b,c,i,j}
\rangle\langle\omega^{\tau,\la,+}_{b,c,i,j}| \]
graphically. The proof for ``$-$'' is analogous.
This sum is given by the left hand side of
Fig.\ \ref{chirsum1}.
%
\begin{figure}[htb]
\begin{center}
\unitlength 0.6mm
\begin{picture}(224,87)
\thicklines
\put(23,41){\makebox(0,0){$\displaystyle\sum_{a,b,c,d,\la,j}
\frac{\sqrt{d_bd_c}}{w^2\sqrt{d_\tau}}$}}
\put(62,0){\line(0,1){26}}
\put(62,87){\line(0,-1){26}}
\put(72,26){\arc{20}{3.142}{4.712}}
\put(72,61){\arc{20}{1.571}{3.142}}
\put(72,36){\line(1,0){6}}
\put(82,36){\line(1,0){6}}
\put(72,51){\line(1,0){6}}
\put(82,51){\line(1,0){6}}
\put(88,26){\arc{20}{4.712}{0}}
\put(88,61){\arc{20}{0}{1.571}}
\put(98,0){\line(0,1){26}}
\put(98,87){\line(0,-1){26}}
\put(72,16){\line(0,1){10}}
\put(72,71){\line(0,-1){10}}
\put(77,26){\arc{10}{3.142}{4.712}}
\put(77,61){\arc{10}{1.571}{3.142}}
\put(77,31){\line(1,0){1}}
\put(82,31){\line(1,0){1}}
\put(77,56){\line(1,0){1}}
\put(82,56){\line(1,0){1}}
\put(83,26){\arc{10}{4.712}{0}}
\put(83,61){\arc{10}{0}{1.571}}
\put(88,16){\line(0,1){10}}
\put(88,71){\line(0,-1){10}}
\put(83,16){\arc{10}{0}{1.571}}
\put(83,71){\arc{10}{4.712}{0}}
\put(77,11){\line(1,0){6}}
\put(77,76){\line(1,0){6}}
\put(77,16){\arc{10}{1.571}{3.142}}
\put(77,71){\arc{10}{3.142}{4.712}}
\Thicklines
\put(80,0){\line(0,1){11}}
\put(80,3.5){\vector(0,-1){0}}
\put(80,31){\line(0,1){5}}
\put(80,87){\line(0,-1){11}}
\put(80,80.5){\vector(0,-1){0}}
\put(80,56){\line(0,-1){5}}
\thinlines
\put(72,21){\arc{5}{4.712}{1.571}}
\put(72,66){\arc{5}{4.712}{1.571}}
\put(72,21){\line(1,0){3}}
\put(75,26){\arc{10}{0}{1.571}}
\put(80,26){\line(0,1){5}}
\put(80,36){\line(0,1){15}}
\put(72,66){\line(1,0){3}}
\put(75,61){\arc{10}{4.712}{0}}
\put(80,61){\line(0,-1){5}}
\put(80,41.5){\vector(0,-1){0}}
\put(58,4){\makebox(0,0){$d$}}
\put(72,9){\makebox(0,0){$c$}}
\put(92,21){\makebox(0,0){$b$}}
\put(85,4){\makebox(0,0){$\tau$}}
\put(85,43.5){\makebox(0,0){$\la$}}
\put(80.5,15){\makebox(0,0){\tiny{$(X_{\co c,b}^{\tau;j})^*$}}}
\put(58,83){\makebox(0,0){$a$}}
\put(72,77){\makebox(0,0){$c$}}
\put(92,66){\makebox(0,0){$b$}}
\put(85,83){\makebox(0,0){$\tau$}}
\put(80,72){\makebox(0,0){\tiny{$X_{\co c,b}^{\tau;j}$}}}
\thicklines
\put(143,41){\makebox(0,0){$=\;\displaystyle
\sum_{a,b,c,d,\atop\la,\nu,\rho,j}
\frac{\sqrt{d_bd_c}}{w^2\sqrt{d_\tau}}$}}
\put(182,0){\line(0,1){36}}
\put(182,87){\line(0,-1){36}}
\put(187,36){\arc{10}{3.142}{0}}
\put(187,51){\arc{10}{0}{3.142}}
\put(213,36){\arc{10}{3.142}{0}}
\put(213,51){\arc{10}{0}{3.142}}
\put(218,0){\line(0,1){36}}
\put(218,87){\line(0,-1){36}}
\put(192,16){\line(0,1){20}}
\put(192,71){\line(0,-1){20}}
\put(208,16){\line(0,1){20}}
\put(208,71){\line(0,-1){20}}
\put(203,16){\arc{10}{0}{1.571}}
\put(203,71){\arc{10}{4.712}{0}}
\put(197,11){\line(1,0){6}}
\put(197,76){\line(1,0){6}}
\put(197,16){\arc{10}{1.571}{3.142}}
\put(197,71){\arc{10}{3.142}{4.712}}
\Thicklines
\put(200,0){\line(0,1){11}}
\put(200,3.5){\vector(0,-1){0}}
\put(200,87){\line(0,-1){11}}
\put(200,80.5){\vector(0,-1){0}}
\thinlines
\put(192,21){\arc{5}{4.712}{1.571}}
\put(192,66){\arc{5}{4.712}{1.571}}
\put(192,31){\arc{5}{4.712}{1.571}}
\put(192,56){\arc{5}{4.712}{1.571}}
\put(208,31){\arc{5}{1.571}{4.712}}
\put(208,56){\arc{5}{1.571}{4.712}}
\put(192,31){\line(1,0){6}}
\put(192,56){\line(1,0){6}}
\put(208,31){\line(-1,0){6}}
\put(208,56){\line(-1,0){6}}
\put(198.5,31){\vector(1,0){0}}
\put(198.5,56){\vector(1,0){0}}
\put(192,21){\line(1,0){3}}
\put(195,26){\arc{10}{0}{1.571}}
\put(200,26){\line(0,1){5}}
\put(200,31){\line(0,1){25}}
\put(192,66){\line(1,0){3}}
\put(195,61){\arc{10}{4.712}{0}}
\put(200,61){\line(0,-1){5}}
\put(200,41.5){\vector(0,-1){0}}
\put(178,4){\makebox(0,0){$d$}}
\put(222,4){\makebox(0,0){$d$}}
\put(192,9){\makebox(0,0){$c$}}
\put(188,26){\makebox(0,0){$b$}}
\put(212,21){\makebox(0,0){$b$}}
\put(205,4){\makebox(0,0){$\tau$}}
\put(205,43.5){\makebox(0,0){$\la$}}
\put(200.5,15){\makebox(0,0){\tiny{$(X_{\co c,b}^{\tau;j})^*$}}}
\put(178,83){\makebox(0,0){$a$}}
\put(222,83){\makebox(0,0){$a$}}
\put(192,77){\makebox(0,0){$c$}}
\put(188,61){\makebox(0,0){$b$}}
\put(212,66){\makebox(0,0){$b$}}
\put(196,35){\makebox(0,0){$\rho$}}
\put(196,51){\makebox(0,0){$\nu$}}
\put(205,83){\makebox(0,0){$\tau$}}
\put(200,72){\makebox(0,0){\tiny{$X_{\co c,b}^{\tau;j}$}}}
\end{picture}
\end{center}
\caption{The sum $\sum_{\la,\xi} w^{-2}\protect\sqrt{d_\la/d_\tau}
|\omega^{\tau,\la,\pm}_\xi\rangle\langle\omega^{\tau,\la,\pm}_\xi|$}
\label{chirsum1}
\end{figure}
Using the expansion of the identity (cf.\ \cite[Lemma 4.3]{BEK1})
for the parallel wires $a,b$ on the top and $d,b$
on the bottom we obtain the
right hand side of Fig.\ \ref{chirsum1}. Using such an expansion
now the other way round for the summation over $\la$ we arrive
at the left hand side of Fig.\ \ref{chirsum2}.
%
\begin{figure}[htb]
\begin{center}
\unitlength 0.6mm
\begin{picture}(219,87)
\thicklines
\put(19,41){\makebox(0,0){$\displaystyle
\sum_{a,b,c,d,\atop\nu,\rho,j}
\frac{\sqrt{d_bd_c}}{w^2\sqrt{d_\tau}}$}}
\put(47,0){\line(0,1){36}}
\put(47,87){\line(0,-1){36}}
\put(52,36){\arc{10}{3.142}{0}}
\put(52,51){\arc{10}{0}{3.142}}
\put(57,26){\line(0,1){10}}
\put(57,61){\line(0,-1){10}}
\put(61,61){\arc{8}{3.142}{0}}
\put(61,26){\arc{8}{0}{3.142}}
\put(65,26){\line(0,1){35}}
\put(75,16){\line(0,1){55}}
\put(96,36){\arc{10}{3.142}{0}}
\put(96,51){\arc{10}{0}{3.142}}
\put(101,0){\line(0,1){36}}
\put(101,87){\line(0,-1){36}}
\put(91,16){\line(0,1){20}}
\put(91,71){\line(0,-1){20}}
\put(86,16){\arc{10}{0}{1.571}}
\put(86,71){\arc{10}{4.712}{0}}
\put(80,11){\line(1,0){6}}
\put(80,76){\line(1,0){6}}
\put(80,16){\arc{10}{1.571}{3.142}}
\put(80,71){\arc{10}{3.142}{4.712}}
\Thicklines
\put(83,0){\line(0,1){11}}
\put(83,3.5){\vector(0,-1){0}}
\put(83,87){\line(0,-1){11}}
\put(83,80.5){\vector(0,-1){0}}
\put(67,31){\line(1,0){6}}
\put(72,31){\vector(1,0){0}}
\put(67,56){\line(1,0){6}}
\put(72,56){\vector(1,0){0}}
\thinlines
\put(57,31){\arc{5}{4.712}{1.571}}
\put(57,56){\arc{5}{4.712}{1.571}}
\put(91,31){\arc{5}{1.571}{4.712}}
\put(91,56){\arc{5}{1.571}{4.712}}
\put(57,31){\line(1,0){6}}
\put(57,56){\line(1,0){6}}
\put(91,31){\line(-1,0){14}}
\put(91,56){\line(-1,0){14}}
\put(43,4){\makebox(0,0){$d$}}
\put(105,4){\makebox(0,0){$d$}}
\put(61,43.5){\makebox(0,0){$b$}}
\put(95,21){\makebox(0,0){$b$}}
\put(88,4){\makebox(0,0){$\tau$}}
\put(83.5,15){\makebox(0,0){\tiny{$(X_{\co c,b}^{\tau;j})^*$}}}
\put(43,83){\makebox(0,0){$a$}}
\put(105,83){\makebox(0,0){$a$}}
\put(79,43.5){\makebox(0,0){$c$}}
\put(95,66){\makebox(0,0){$b$}}
\put(70,25){\makebox(0,0){$\a^-_\rho$}}
\put(70,62){\makebox(0,0){$\a^-_\nu$}}
\put(88,83){\makebox(0,0){$\tau$}}
\put(83,72){\makebox(0,0){\tiny{$X_{\co c,b}^{\tau;j}$}}}
\thicklines
\put(138,41){\makebox(0,0){$=\;\displaystyle
\sum_{a,b,c,d,\la,j,\atop \beta\in\MXMm}
\frac{\sqrt{d_bd_c}}{w^2\sqrt{d_\tau}}$}}
\put(172,0){\line(0,1){87}}
\put(187,41){\line(0,1){5}}
\put(192,46){\arc{10}{3.142}{4.712}}
\put(192,41){\arc{10}{1.571}{3.142}}
\put(192,36){\line(1,0){1}}
\put(197,36){\line(1,0){6}}
\put(192,51){\line(1,0){1}}
\put(197,51){\line(1,0){6}}
\put(203,26){\arc{20}{4.712}{0}}
\put(203,61){\arc{20}{0}{1.571}}
\put(213,0){\line(0,1){26}}
\put(213,87){\line(0,-1){26}}
\put(187,16){\line(0,1){10}}
\put(187,71){\line(0,-1){10}}
\put(192,26){\arc{10}{3.142}{4.712}}
\put(192,61){\arc{10}{1.571}{3.142}}
\put(192,31){\line(1,0){1}}
\put(197,31){\line(1,0){1}}
\put(192,56){\line(1,0){1}}
\put(197,56){\line(1,0){1}}
\put(198,26){\arc{10}{4.712}{0}}
\put(198,61){\arc{10}{0}{1.571}}
\put(203,16){\line(0,1){10}}
\put(203,71){\line(0,-1){10}}
\put(198,16){\arc{10}{0}{1.571}}
\put(198,71){\arc{10}{4.712}{0}}
\put(192,11){\line(1,0){6}}
\put(192,76){\line(1,0){6}}
\put(192,16){\arc{10}{1.571}{3.142}}
\put(192,71){\arc{10}{3.142}{4.712}}
\Thicklines
\put(172,43.5){\line(1,0){15}}
\put(181.5,43.5){\vector(1,0){0}}
\put(195,0){\line(0,1){11}}
\put(195,3.5){\vector(0,-1){0}}
\put(195,31){\line(0,1){5}}
\put(195,87){\line(0,-1){11}}
\put(195,80.5){\vector(0,-1){0}}
\put(195,56){\line(0,-1){5}}
\thinlines
\put(187,21){\arc{5}{4.712}{1.571}}
\put(187,66){\arc{5}{4.712}{1.571}}
\put(187,21){\line(1,0){3}}
\put(190,26){\arc{10}{0}{1.571}}
\put(195,26){\line(0,1){5}}
\put(195,36){\line(0,1){15}}
\put(187,66){\line(1,0){3}}
\put(190,61){\arc{10}{4.712}{0}}
\put(195,61){\line(0,-1){5}}
\put(195,41.5){\vector(0,-1){0}}
\put(168,4){\makebox(0,0){$d$}}
\put(217,4){\makebox(0,0){$d$}}
\put(187,9){\makebox(0,0){$c$}}
\put(207,21){\makebox(0,0){$b$}}
\put(200,4){\makebox(0,0){$\tau$}}
\put(200,43.5){\makebox(0,0){$\la$}}
\put(195.5,15){\makebox(0,0){\tiny{$(X_{\co c,b}^{\tau;j})^*$}}}
\put(168,83){\makebox(0,0){$a$}}
\put(217,83){\makebox(0,0){$a$}}
\put(187,77){\makebox(0,0){$c$}}
\put(207,66){\makebox(0,0){$b$}}
\put(200,83){\makebox(0,0){$\tau$}}
\put(179.5,38.5){\makebox(0,0){$\beta$}}
\put(195,72){\makebox(0,0){\tiny{$X_{\co c,b}^{\tau;j}$}}}
\end{picture}
\end{center}
\caption{The sum $\sum_{\la,\xi} w^{-2}\protect\sqrt{d_\la/d_\tau}
|\omega^{\tau,\la,\pm}_\xi\rangle\langle\omega^{\tau,\la,\pm}_\xi|$}
\label{chirsum2}
\end{figure}
The crucial point is now the observation that left and right part
of this wire diagram are only connected by wires $\a^-_\nu$ and
$\a^-_\rho$. Let us start again with the original picture, namely the
left hand side of Fig.\ \ref{chirsum1}, and make an expansion
for the open ending wires $a$ and $d$ on the left side with a
summation over wires $\beta\in\MXM$. Then it follows that only the
wires with $\beta\in\MXMm$ contribute because
$\Hom(\beta,\a^-_\rho\a^-_\nu)$ is always zero unless $\beta\in\MXMm$.
This establishes equality with the right hand side of 
Fig.\ \ref{chirsum2}. The wire $\beta$ can now be pulled in and
application of the naturality move of Fig.\ \ref{natrelui}
for the relative braiding yields the left hand side of
Fig.\ \ref{chirsum3}.
%
\begin{figure}[htb]
\begin{center}
\unitlength 0.6mm
\begin{picture}(225,87)
\thicklines
\put(23,41){\makebox(0,0){$\displaystyle
\sum_{a,b,c,d,\la,j,\atop \beta\in\MXMm}
\frac{\sqrt{d_bd_c}}{w^2\sqrt{d_\tau}}$}}
\put(62,0){\line(0,1){87}}
\put(98,0){\line(0,1){87}}
\put(72,26){\line(0,1){10}}
\put(72,71){\line(0,-1){10}}
\put(77,36){\arc{10}{3.142}{4.712}}
\put(77,61){\arc{10}{1.571}{3.142}}
\put(77,41){\line(1,0){1}}
\put(82,41){\line(1,0){1}}
\put(77,56){\line(1,0){1}}
\put(82,56){\line(1,0){1}}
\put(83,36){\arc{10}{4.712}{0}}
\put(83,61){\arc{10}{0}{1.571}}
\put(88,26){\line(0,1){10}}
\put(88,71){\line(0,-1){10}}
\put(83,26){\arc{10}{0}{1.571}}
\put(83,71){\arc{10}{4.712}{0}}
\put(77,21){\line(1,0){6}}
\put(77,76){\line(1,0){6}}
\put(77,26){\arc{10}{1.571}{3.142}}
\put(77,71){\arc{10}{3.142}{4.712}}
\Thicklines
\put(62,11){\line(1,0){16}}
\put(82,11){\line(1,0){16}}
\put(73,11){\vector(1,0){0}}
\put(80,0){\line(0,1){21}}
\put(80,3.5){\vector(0,-1){0}}
\put(80,41){\line(0,1){15}}
\put(80,46.5){\vector(0,-1){0}}
\put(80,87){\line(0,-1){11}}
\put(80,80.5){\vector(0,-1){0}}
\thinlines
\put(62,11){\arc{5}{4.712}{1.571}}
\put(98,11){\arc{5}{1.571}{4.712}}
\put(72,31){\arc{5}{4.712}{1.571}}
\put(72,66){\arc{5}{4.712}{1.571}}
\put(72,31){\line(1,0){3}}
\put(75,36){\arc{10}{0}{1.571}}
\put(80,36){\line(0,1){5}}
\put(72,66){\line(1,0){3}}
\put(75,61){\arc{10}{4.712}{0}}
\put(80,61){\line(0,-1){5}}
\put(58,4){\makebox(0,0){$d$}}
\put(102,4){\makebox(0,0){$d$}}
\put(72,19){\makebox(0,0){$c$}}
\put(92,31){\makebox(0,0){$b$}}
\put(85,4){\makebox(0,0){$\tau$}}
\put(87,48.5){\makebox(0,0){$\a^+_\la$}}
\put(80.5,25){\makebox(0,0){\tiny{$(X_{\co c,b}^{\tau;j})^*$}}}
\put(58,83){\makebox(0,0){$a$}}
\put(102,83){\makebox(0,0){$a$}}
\put(72,77){\makebox(0,0){$c$}}
\put(92,66){\makebox(0,0){$b$}}
\put(85,83){\makebox(0,0){$\tau$}}
\put(70,6){\makebox(0,0){$\beta$}}
\put(80,72){\makebox(0,0){\tiny{$X_{\co c,b}^{\tau;j}$}}}
\thicklines
\put(143,41){\makebox(0,0){$=\;\displaystyle
\sum_{a,b,c,d,j,\atop \beta\in\MXMm}
\frac{\sqrt{d_bd_c}}{ww_+\sqrt{d_\tau}}$}}
\put(182,0){\line(0,1){87}}
\put(218,0){\line(0,1){87}}
\put(192,26){\line(0,1){10}}
\put(192,71){\line(0,-1){10}}
\put(197,36){\arc{10}{3.142}{4.712}}
\put(197,61){\arc{10}{1.571}{3.142}}
\put(197,41){\line(1,0){6}}
\put(197,56){\line(1,0){6}}
\put(203,36){\arc{10}{4.712}{0}}
\put(203,61){\arc{10}{0}{1.571}}
\put(208,26){\line(0,1){10}}
\put(208,71){\line(0,-1){10}}
\put(203,26){\arc{10}{0}{1.571}}
\put(203,71){\arc{10}{4.712}{0}}
\put(197,21){\line(1,0){6}}
\put(197,76){\line(1,0){6}}
\put(197,26){\arc{10}{1.571}{3.142}}
\put(197,71){\arc{10}{3.142}{4.712}}
\Thicklines
\put(182,11){\line(1,0){16}}
\put(202,11){\line(1,0){16}}
\put(193,11){\vector(1,0){0}}
\put(200,0){\line(0,1){21}}
\put(200,3.5){\vector(0,-1){0}}
\put(200,41){\line(0,1){15}}
\put(200,46.5){\vector(0,-1){0}}
\put(200,87){\line(0,-1){11}}
\put(200,80.5){\vector(0,-1){0}}
\thinlines
\put(182,11){\arc{5}{4.712}{1.571}}
\put(218,11){\arc{5}{1.571}{4.712}}
\put(200,41){\arc{5}{3.142}{0}}
\put(200,56){\arc{5}{0}{3.142}}
\put(178,4){\makebox(0,0){$d$}}
\put(222,4){\makebox(0,0){$d$}}
\put(188,31){\makebox(0,0){$c$}}
\put(212,31){\makebox(0,0){$b$}}
\put(205,4){\makebox(0,0){$\tau$}}
\put(205,48.5){\makebox(0,0){$\tau$}}
\put(200.5,26){\makebox(0,0){\tiny{$(X_{\co c,b}^{\tau;j})^*$}}}
\put(178,83){\makebox(0,0){$a$}}
\put(222,83){\makebox(0,0){$a$}}
\put(188,66){\makebox(0,0){$c$}}
\put(212,66){\makebox(0,0){$b$}}
\put(205,83){\makebox(0,0){$\tau$}}
\put(190,6){\makebox(0,0){$\beta$}}
\put(200,71){\makebox(0,0){\tiny{$X_{\co c,b}^{\tau;j}$}}}
\end{picture}
\end{center}
\caption{The sum $\sum_{\la,\xi} w^{-2}\protect\sqrt{d_\la/d_\tau}
|\omega^{\tau,\la,\pm}_\xi\rangle\langle\omega^{\tau,\la,\pm}_\xi|$}
\label{chirsum3}
\end{figure}
Then, using the graphical identity of Fig.\ \ref{scg=cp}
gives us the right hand side of Fig.\ \ref{chirsum3},
as only the wire $\tau$ survives in the sum of the chiral
horizontal projector. The two ``bulbs'' yield just a scalar
factor $\sqrt{d_bd_c/d_\tau}$, but due to the summation
over the fusion channels $j$ it appears with multiplicity
$N_{\co c,b}^\tau$. Hence the total prefactor is
calculated as
\[ \sum_{b,c} \frac{d_b d_c N_{\co c,b}^\tau}{ww_+ d_\tau}
= \sum_b \frac{d_b^2}{ww_+} = \frac 1{w_+} \,, \]
and this is the prefactor of $I_\tau^+$.
\end{proof}

Now let us the consider the case $\tau=\id$: Note
that $\cA_{\id}$ is a subspace of the double triangle
algebra $\dta$ containing the horizontal center
$\cZ_h$. Then the sum
$\frac 1{w^2} \sum_{\la,\xi} \sqrt{d_\la}
|\omega^{0,\la,\pm}_\xi\rangle\langle\omega^{0,\la,\pm}_\xi|$
gives graphically exactly the picture ($+$) for
$\sum_\la q_{\la,0}$ respectively ($-$) for
$\sum_\la q_{0,\la}$, where $q_{\la,\mu}\in\cZ_h$
are the vertical projectors of \cite[Def.\ 6.7]{BEK1}.
Hence we obtain the following

\begin{corollary}
In the double triangle algebra
we have $w_+\sum_\la q_{\la,0}=P^-$ and
$w_+\sum_\la q_{0,\la}=P^+$.
\end{corollary}

Next we establish some kind of trivial action of $\MXMpm$ on
$H_{\tau,\la}^\mp$.

\begin{lemma}
\label{chirtriv}
For $\beta\in\MXMm$ we have the identity of Fig.\ \ref{ebomdom}.
An analogous identity holds when we choose $\beta\in\MXMp$ acting
on $\omega^{\tau,\la,-}_{b,c,t,X}$.
\end{lemma}
%
\begin{figure}[htb]
\begin{center}
\unitlength 0.6mm
\begin{picture}(121,60)
\thicklines
\put(10,27){\makebox(0,0){$\displaystyle\sum_{a,d} \, d_a$}}
\put(32,0){\line(0,1){35}}
\put(42,35){\arc{20}{3.142}{4.712}}
\put(42,45){\line(1,0){6}}
\put(52,45){\line(1,0){6}}
\put(58,35){\arc{20}{4.712}{0}}
\put(68,0){\line(0,1){35}}
\put(42,25){\line(0,1){10}}
\put(47,35){\arc{10}{3.142}{4.712}}
\put(47,40){\line(1,0){1}}
\put(52,40){\line(1,0){1}}
\put(53,35){\arc{10}{4.712}{0}}
\put(58,25){\line(0,1){10}}
\put(53,25){\arc{10}{0}{1.571}}
\put(47,20){\line(1,0){6}}
\put(47,25){\arc{10}{1.571}{3.142}}
\Thicklines
\put(50,20){\vector(0,-1){20}}
\put(50,40){\line(0,1){5}}
\put(32,10){\line(1,0){16}}
\put(62,10){\vector(1,0){0}}
\put(52,10){\line(1,0){16}}
\thinlines
\put(42,30){\line(1,0){3}}
\put(45,35){\arc{10}{0}{1.571}}
\put(50,35){\line(0,1){5}}
\put(50,45){\line(0,1){15}}
\put(32,10){\arc{5}{4.712}{1.571}}
\put(68,10){\arc{5}{1.571}{4.712}}
\put(50,50.5){\vector(0,-1){0}}
\put(28,4){\makebox(0,0){$d$}}
\put(72,4){\makebox(0,0){$d$}}
\put(28,18){\makebox(0,0){$a$}}
\put(42,19){\makebox(0,0){$c$}}
\put(62,30){\makebox(0,0){$b$}}
\put(55,4){\makebox(0,0){$\tau$}}
\put(59,14){\makebox(0,0){$\beta$}}
\put(55,52.5){\makebox(0,0){$\la$}}
\put(37,30){\makebox(0,0){$t$}}
\put(50,24){\makebox(0,0){$X^*$}}
\put(102,30){\makebox(0,0){$=\;d_\beta^2 \,
\omega^{\tau,\la,+}_{b,c,t,X}$}}
\end{picture}
\end{center}
\caption{The trivial action of $\MXMm$ on $H_{\tau,\la}^+$}
\label{ebomdom}
\end{figure}

\begin{proof}
Starting with Fig.\ \ref{ebomdom} we can slide around the, say, left
trivalent vertex of the wire $\beta$ to obtain the left hand side of
Fig.\ \ref{pebomdom}.
%
\begin{figure}[htb]
\begin{center}
\unitlength 0.6mm
\begin{picture}(175,65)
\thicklines
\put(10,27){\makebox(0,0){$\displaystyle\sum_{a,d} \, d_a$}}
\put(27,0){\line(0,1){50}}
\put(37,50){\arc{20}{3.142}{4.712}}
\put(37,60){\line(1,0){11}}
\put(52,60){\line(1,0){6}}
\put(58,50){\arc{20}{4.712}{0}}
\put(68,0){\line(0,1){50}}
\put(42,25){\line(0,1){10}}
\put(47,35){\arc{10}{3.142}{4.712}}
\put(47,40){\line(1,0){1}}
\put(52,40){\line(1,0){1}}
\put(53,35){\arc{10}{4.712}{0}}
\put(58,25){\line(0,1){10}}
\put(53,25){\arc{10}{0}{1.571}}
\put(47,20){\line(1,0){6}}
\put(47,25){\arc{10}{1.571}{3.142}}
\Thicklines
\put(50,20){\vector(0,-1){20}}
\put(50,50.5){\vector(0,-1){0}}
\put(50,40){\line(0,1){20}}
\put(42,10){\line(1,0){6}}
\put(42,45){\line(1,0){6}}
\put(32,20){\line(0,1){15}}
\put(42,35){\arc{20}{3.142}{4.712}}
\put(42,20){\arc{20}{1.571}{3.142}}
\put(32,25.5){\vector(0,-1){0}}
\put(52,10){\line(1,0){16}}
\put(52,45){\line(1,0){16}}
\thinlines
\put(42,30){\line(1,0){3}}
\put(45,35){\arc{10}{0}{1.571}}
\put(50,35){\line(0,1){5}}
\put(50,65){\line(0,-1){5}}
\put(68,45){\arc{5}{1.571}{4.712}}
\put(68,10){\arc{5}{1.571}{4.712}}
\put(23,4){\makebox(0,0){$d$}}
\put(72,4){\makebox(0,0){$d$}}
\put(72,27.5){\makebox(0,0){$a$}}
\put(42,19){\makebox(0,0){$c$}}
\put(62,30){\makebox(0,0){$b$}}
\put(55,4){\makebox(0,0){$\tau$}}
\put(34,7){\makebox(0,0){$\beta$}}
\put(43,52.5){\makebox(0,0){$\a^+_\la$}}
\put(37,30){\makebox(0,0){$t$}}
\put(50,24){\makebox(0,0){$X^*$}}
\thicklines
\put(105,27){\makebox(0,0){$=\;\displaystyle\sum_{a,d} \, d_\beta$}}
\put(127,0){\line(0,1){50}}
\put(137,50){\arc{20}{3.142}{4.712}}
\put(137,60){\line(1,0){6}}
\put(147,60){\line(1,0){11}}
\put(158,50){\arc{20}{4.712}{0}}
\put(168,0){\line(0,1){5}}
\put(163,5){\arc{10}{4.712}{0}}
\put(160.5,10){\line(1,0){2.5}}
\put(160.5,25){\line(1,0){2.5}}
\put(163,30){\arc{10}{0}{1.571}}
\put(168,30){\line(0,1){22}}
\put(160.5,10){\line(0,1){15}}
\put(137,35){\line(0,1){10}}
\put(142,45){\arc{10}{3.142}{4.712}}
\put(142,50){\line(1,0){1}}
\put(147,50){\line(1,0){1}}
\put(148,45){\arc{10}{4.712}{0}}
\put(153,35){\line(0,1){10}}
\put(148,35){\arc{10}{0}{1.571}}
\put(142,30){\line(1,0){6}}
\put(142,35){\arc{10}{1.571}{3.142}}
\Thicklines
\put(145,30){\vector(0,-1){30}}
\put(145,53){\vector(0,-1){0}}
\put(145,50){\line(0,1){10}}
\put(158,10){\line(1,0){2.5}}
\put(158,15){\arc{10}{1.571}{3.142}}
\put(153,15){\line(0,1){5}}
\put(153,15.5){\vector(0,-1){0}}
\put(158,20){\arc{10}{3.142}{4.712}}
\put(158,25){\line(1,0){2.5}}
\thinlines
\put(137,40){\line(1,0){3}}
\put(140,45){\arc{10}{0}{1.571}}
\put(145,45){\line(0,1){5}}
\put(145,65){\line(0,-1){5}}
\put(160.5,25){\arc{5}{0}{3.142}}
\put(160.5,10){\arc{5}{3.142}{0}}
\put(123,4){\makebox(0,0){$d$}}
\put(172,4){\makebox(0,0){$d$}}
\put(165,17.5){\makebox(0,0){$a$}}
\put(137,29){\makebox(0,0){$c$}}
\put(157,40){\makebox(0,0){$b$}}
\put(140,4){\makebox(0,0){$\tau$}}
\put(153,5){\makebox(0,0){$\beta$}}
\put(138,55){\makebox(0,0){$\a^+_\la$}}
\put(132,40){\makebox(0,0){$t$}}
\put(145,34){\makebox(0,0){$X^*$}}
\end{picture}
\end{center}
\caption{Proof of the trivial action of $\MXMm$ on $H_{\tau,\la}^+$}
\label{pebomdom}
\end{figure}
Using now the relative braiding naturality move of Fig.\ \ref{natrelui}
and turning around the small arcs, giving a factor $d_\beta/d_a$,
yields the right hand side of Fig.\ \ref{pebomdom}. We now see
that the summation over the wire $a$ is just an expansion of the
identity which can be replaced by parallel wires $\beta$ and $d$
(cf.\ \cite[Lemma 4.3]{BEK1}). Hence we obtain a closed
loop $\beta$ which is just another factor $d_\beta$, and we are
left with the original diagram for $\omega^{\tau,\la,+}_{b,c,t,X}$,
together with a prefactor $d_\beta^2$.
\end{proof}

We now obtain immediately the following corollary which finally
justifies the name ``chiral multiplicative units'' for elements
$I_\tau^\pm\in A_\tau^\pm$.

\begin{corollary}
\label{unitpm}
In $\cA_\tau$ we have
\begin{equation}
I_\tau^\pm \star_v |\omega^{\tau,\la,\pm}_{b',c',t',X'}\rangle
\langle\omega^{\tau,\la,\pm}_{b,c,t,X}| =
|\omega^{\tau,\la,\pm}_{b',c',t',X'}\rangle
\langle\omega^{\tau,\la,\pm}_{b,c,t,X}| \star_v I_\tau^\pm =
|\omega^{\tau,\la,\pm}_{b',c',t',X'}\rangle
\langle\omega^{\tau,\la,\pm}_{b,c,t,X}| \,.
\end{equation}
\end{corollary}

Then we define elements $I_{\tau,\la}^\pm\in A_{\tau,\la}^\pm$ by
\begin{equation}
I_{\tau,\la}^\pm = \frac 1{w^2}\sqrt{\frac{d_\la}{d_\tau}}
\sum_\xi |\omega^{\tau,\la,\pm}_\xi\rangle
\langle\omega^{\tau,\la,\pm}_\xi|\,,
\end{equation}
so that $I_\tau^\pm=\sum_\la I_{\tau,\la}^\pm$.
We now claim

\begin{lemma}
\label{Iexpan}
We have the expansion in matrix units
\begin{equation}
I_{\tau,\la}^\pm = \sum_{i=1}^{\dim H_{\tau,\la}^\pm}
|u_i^{\tau,\la,\pm}\rangle \langle u_i^{\tau,\la,\pm}|\,.
\end{equation}
\end{lemma}

\begin{proof}
Using Lemma \ref{chirality}, Corollary \ref{matunAt}
and Corollary \ref{unitpm} we compute
\[ \bearl
\displaystyle\frac 1{w^2} \sqrt{\frac{d_\la}{d_\tau}} \sum_\xi
\langle u_i^{\tau,\la,\pm}, \om_\xi^{\tau,\la,\pm}\rangle
\langle \om_\xi^{\tau,\la,\pm}, u_j^{\tau,\la,\pm}\rangle
|u_i^{\tau,\la,\pm}\rangle \langle u_j^{\tau,\la,\pm}| = \\[.4em]
\qquad\qquad\qquad=
|u_i^{\tau,\la,\pm}\rangle \langle u_i^{\tau,\la,\pm}| \star_v
I_\tau^\pm \star_v |u_j^{\tau,\la,\pm}\rangle \langle u_j^{\tau,\la,\pm}|
= \del ij |u_i^{\tau,\la,\pm}\rangle \langle u_i^{\tau,\la,\pm}| \,.
\eear \]
On the other hand we obtain by expanding the vectors
$\om_\xi^{\tau,\la,\pm}$ in basis vectors $u_i^{\tau,\la,\pm}$
\[ I_{\tau,\la}^\pm = \frac 1{w^2} \sqrt{\frac{d_\la}{d_\tau}} \sum_{\xi,i,j}
\langle u_i^{\tau,\la,\pm}, \om_\xi^{\tau,\la,\pm}\rangle
\langle \om_\xi^{\tau,\la,\pm}, u_j^{\tau,\la,\pm}\rangle
|u_i^{\tau,\la,\pm}\rangle \langle u_j^{\tau,\la,\pm}| \,, \]
hence
$I_{\tau,\la}^\pm =  \sum_{i,j} \del ij
|u_i^{\tau,\la,\pm}\rangle \langle u_i^{\tau,\la,\pm}|$.
\end{proof}

\begin{lemma}
\label{chirdim}
The dimensions of the Hilbert spaces $H_{\tau,\la}^\pm$ 
are given by the chiral branching coefficients:
$\dim H_{\tau,\la}^\pm=b_{\tau,\la}^\pm$.
\end{lemma}

\begin{proof}
We show $\dim H_{\tau,\la}^+=b_{\tau,\la}^+$; the ``$-$'' case
is analogous.
The dimensions $\dim H_{\tau,\la}^+$ are counted as
\[ \dim H_{\tau,\la}^+= \sum_{i=1}^{\dim H_{\tau,\la}^+}
\langle u_i^{\tau,\la,+} , u_i^{\tau,\la,+} \rangle =
\frac 1{d_\la} \psi_v^\tau (I^+_{\tau,\la}) =
\frac 1{w^2} \sqrt{\frac{d_\la}{d_\tau}} \sum_\xi
\langle \om_\xi^{\tau,\la,+}, \om_\xi^{\tau,\la,+} \rangle \,.\]
Using now the graphical representation of the
scalar product in Fig.\ \ref{<om,om>}, then we obtain
with the normalization convention for the small semicircular
wires exactly the wire diagram for $b^+_{\tau,\la}$,
cf.\ Fig.\ \ref{bgraph}.
\end{proof}

\subsection{Chiral representations}

Recall that the horizontal center $\cZ_h$ of the
double triangle algebra $\dta$ is spanned by
the elements $e_\beta$ with $\beta\in\MXM$.
Denote $\cZ_h^\pm=P^\pm *_h \cZ_h$. Since the $e_\beta$'s
are projections with respect to the horizontal product,
$\cZ_h^\pm\subset\cZ_h$ are the subspaces spanned by
elements $e_{\beta_\pm}$ with $\beta_\pm\in\MXMpm$.
As $(\cZ_h,*_v)$ is isomorphic to the $M$-$M$ fusion rule
algebra (cf.\ \cite[Thm.\ 4.4]{BEK1}) and since
$\MXMpm\subset\MXM$ are subsystems,
$\cZ_h^\pm\subset\cZ_h$ are in fact vertical subalgebras.
We are now going to construct representations of these
chiral vertical algebras $(\cZ_h^\pm,*_v)$.

\begin{lemma}
\label{chirrep}
For $\beta_\pm\in\MXMpm$ let
$\pi_{\tau,\la}^+(e_{\beta_+})
\omega^{\tau,\la,+}_{b,c,t,X}\in\cH_{\tau,\la}$ and
$\pi_{\tau,\la}^-(e_{\beta_-})
\omega^{\tau,\la,-}_{b,c,t,X}\in\cH_{\co\tau,\co\la}$,
respectively, denote the vectors defined graphically by
the left respectively right hand side of Fig.\ \ref{pitlb}.
Then in fact
$\pi_{\tau,\la}^\pm(e_{\beta_\pm})
\omega^{\tau,\la,\pm}_{b,c,t,X}\in H_{\tau,\la}^\pm$.
\end{lemma}
%
\begin{figure}[htb]
\begin{center}
\unitlength 0.6mm
\begin{picture}(190,60)
\thicklines
\put(10,27){\makebox(0,0){$\displaystyle\sum_{a,d} \; d_a$}}
\put(32,0){\line(0,1){35}}
\put(42,35){\arc{20}{3.142}{4.712}}
\put(42,45){\line(1,0){6}}
\put(52,45){\line(1,0){6}}
\put(58,35){\arc{20}{4.712}{0}}
\put(68,0){\line(0,1){35}}
\put(42,25){\line(0,1){10}}
\put(47,35){\arc{10}{3.142}{4.712}}
\put(47,40){\line(1,0){1}}
\put(52,40){\line(1,0){1}}
\put(53,35){\arc{10}{4.712}{0}}
\put(58,25){\line(0,1){10}}
\put(53,25){\arc{10}{0}{1.571}}
\put(47,20){\line(1,0){6}}
\put(47,25){\arc{10}{1.571}{3.142}}
\Thicklines
\put(50,8){\vector(0,-1){8}}
\put(50,40){\line(0,1){5}}
\put(32,10){\line(1,0){36}}
\put(61,10){\vector(1,0){0}}
\put(50,12){\line(0,1){8}}
\thinlines
\put(42,30){\line(1,0){3}}
\put(45,35){\arc{10}{0}{1.571}}
\put(50,35){\line(0,1){5}}
\put(50,45){\line(0,1){15}}
\put(32,10){\arc{5}{4.712}{1.571}}
\put(68,10){\arc{5}{1.571}{4.712}}
\put(50,50.5){\vector(0,-1){0}}
\put(28,4){\makebox(0,0){$d$}}
\put(72,4){\makebox(0,0){$d$}}
\put(28,18){\makebox(0,0){$a$}}
\put(42,19){\makebox(0,0){$c$}}
\put(62,30){\makebox(0,0){$b$}}
\put(55,4){\makebox(0,0){$\tau$}}
\put(59,14){\makebox(0,0){$\beta_+$}}
\put(55,52.5){\makebox(0,0){$\la$}}
\put(37,30){\makebox(0,0){$t$}}
\put(50,24){\makebox(0,0){$X^*$}}
\thicklines
\put(130,27){\makebox(0,0){$\displaystyle\sum_{a,d} \; d_a$}}
\put(152,0){\line(0,1){35}}
\put(162,35){\arc{20}{3.142}{4.712}}
\put(162,45){\line(1,0){16}}
\put(178,35){\arc{20}{4.712}{0}}
\put(188,0){\line(0,1){35}}
\put(162,25){\line(0,1){10}}
\put(167,35){\arc{10}{3.142}{4.712}}
\put(167,40){\line(1,0){6}}
\put(173,35){\arc{10}{4.712}{0}}
\put(178,25){\line(0,1){10}}
\put(173,25){\arc{10}{0}{1.571}}
\put(167,20){\line(1,0){6}}
\put(167,25){\arc{10}{1.571}{3.142}}
\Thicklines
\put(170,0){\vector(0,1){20}}
\put(170,42){\line(0,1){1}}
\put(152,10){\line(1,0){16}}
\put(163,10){\vector(1,0){0}}
\put(172,10){\line(1,0){16}}
\thinlines
\put(178,30){\line(-1,0){3}}
\put(175,35){\arc{10}{1.571}{3.142}}
\put(170,35){\line(0,1){3}}
\put(170,47){\line(0,1){13}}
\put(152,10){\arc{5}{4.712}{1.571}}
\put(188,10){\arc{5}{1.571}{4.712}}
\put(170,54.5){\vector(0,1){0}}
\put(148,4){\makebox(0,0){$d$}}
\put(192,4){\makebox(0,0){$d$}}
\put(148,18){\makebox(0,0){$a$}}
\put(178,19){\makebox(0,0){$c$}}
\put(158,30){\makebox(0,0){$b$}}
\put(165,4){\makebox(0,0){$\tau$}}
\put(161,14){\makebox(0,0){$\beta_-$}}
\put(175,52.5){\makebox(0,0){$\la$}}
\put(183,30){\makebox(0,0){$t^*$}}
\put(170,24){\makebox(0,0){$X$}}
\end{picture}
\end{center}
\caption{The vectors
$\pi_{\tau,\la}^+(e_{\beta_+})
\omega^{\tau,\la,+}_{b,c,t,X}\in\cH_{\tau,\la}$ and
$\pi_{\tau,\la}^-(e_{\beta_-})
\omega^{\tau,\la,-}_{b,c,t,X}\in\cH_{\co\tau,\co\la}$}
\label{pitlb}
\end{figure}

\begin{proof}
We prove
$\pi_{\tau,\la}^+(e_{\beta_+})
\omega^{\tau,\la,+}_{b,c,t,X}\in H_{\tau,\la}^+$ for
$\beta_+\in\MXMp$. The proof of
$\pi_{\tau,\la}^-(e_{\beta_-})
\omega^{\tau,\la,-}_{b,c,t,X}\in H_{\tau,\la}^-$ for
$\beta_-\in\MXMm$ is analogous.
First we can turn around the small arcs at the trivalent
vertices of the wire $\beta_+$ which gives us a factor
$d_{\beta_+}/d_a$. The we use the expansion of the identity
(cf.\ \cite[Lemma 4.3]{BEK1}) for the parallel wires
$a$ and $b$. This we way we obtain the left hand side of
Fig.\ \ref{6jproof1}.
%
\begin{figure}[htb]
\begin{center}
\unitlength 0.6mm
\begin{picture}(237,70)
\thicklines
\put(12,32){\makebox(0,0){$\displaystyle\sum_{a,d,\nu} \; d_{\beta_+}$}}
\put(32,0){\line(0,1){35}}
\put(37,35){\arc{10}{3.142}{4.712}}
\put(37,40){\line(1,0){3}}
\put(40,40){\line(0,1){5}}
\put(45,45){\arc{10}{3.142}{4.712}}
\put(45,50){\line(1,0){10}}
\put(55,45){\arc{10}{4.712}{0}}
\put(60,25){\line(0,1){20}}
\put(65,25){\arc{10}{1.571}{3.142}}
\put(65,20){\line(1,0){10}}
\put(75,25){\arc{10}{0}{1.571}}
\put(80,25){\line(0,1){20}}
\put(85,45){\arc{10}{3.142}{4.712}}
\put(85,50){\line(1,0){10}}
\put(95,45){\arc{10}{4.712}{0}}
\put(100,40){\line(0,1){5}}
\put(100,40){\line(1,0){3}}
\put(103,35){\arc{10}{4.712}{0}}
\put(108,0){\line(0,1){35}}
\Thicklines
\put(40,40){\line(1,0){3}}
\put(43,35){\arc{10}{4.712}{0}}
\put(48,15){\line(0,1){20}}
\put(53,15){\arc{10}{1.571}{3.142}}
\put(53,10){\line(1,0){34}}
\put(87,15){\arc{10}{0}{1.571}}
\put(92,15){\line(0,1){20}}
\put(97,35){\arc{10}{3.142}{4.712}}
\put(97,40){\line(1,0){3}}
\put(83,10){\vector(1,0){0}}
\put(70,20){\line(0,-1){8}}
\put(70,8){\vector(0,-1){8}}
\thinlines
\dottedline{4}(29,35)(111,35)
\put(60,30){\line(1,0){5}}
\put(65,35){\arc{10}{0}{1.571}}
\put(70,35){\line(0,1){35}}
\put(70,50){\vector(0,-1){0}}
\put(50,50){\line(0,1){5}}
\put(55,55){\arc{10}{3.142}{4.712}}
\put(55,60){\line(1,0){13}}
\put(64,60){\vector(1,0){0}}
\put(85,60){\line(-1,0){13}}
\put(85,55){\arc{10}{4.712}{0}}
\put(90,50){\line(0,1){5}}
\put(40,40){\arc{5}{3.142}{0}}
\put(50,50){\arc{5}{3.142}{0}}
\put(90,50){\arc{5}{3.142}{0}}
\put(100,40){\arc{5}{3.142}{0}}
\put(28,4){\makebox(0,0){$d$}}
\put(112,4){\makebox(0,0){$d$}}
\put(35,45){\makebox(0,0){$a$}}
\put(105,45){\makebox(0,0){$a$}}
\put(56,42){\makebox(0,0){$b$}}
\put(84,42){\makebox(0,0){$b$}}
\put(60,19){\makebox(0,0){$c$}}
\put(62,65){\makebox(0,0){$\nu$}}
\put(75,50){\makebox(0,0){$\la$}}
\put(65,4){\makebox(0,0){$\tau$}}
\put(82,4){\makebox(0,0){$\beta_+$}}
\put(56,30){\makebox(0,0){$t$}}
\put(70,24){\makebox(0,0){$X^*$}}
\thicklines
\put(140,32){\makebox(0,0){$=\;\displaystyle\sum_{a,d,\nu}
\; d_{\beta_+}$}}
\put(163,0){\line(0,1){45}}
\put(168,45){\arc{10}{3.142}{4.712}}
\put(168,50){\line(1,0){6}}
\put(174,45){\arc{10}{4.712}{0}}
\put(179,40){\line(0,1){5}}
\put(179,40){\line(1,0){3}}
\put(182,35){\arc{10}{4.712}{0}}
\put(187,25){\line(0,1){10}}
\put(192,25){\arc{10}{1.571}{3.142}}
\put(192,20){\line(1,0){10}}
\put(202,25){\arc{10}{0}{1.571}}
\put(207,25){\line(0,1){10}}
\put(212,35){\arc{10}{3.142}{4.712}}
\put(212,40){\line(1,0){3}}
\put(215,40){\line(0,1){5}}
\put(220,45){\arc{10}{3.142}{4.712}}
\put(220,50){\line(1,0){6}}
\put(226,45){\arc{10}{4.712}{0}}
\put(231,0){\line(0,1){45}}
\Thicklines
\put(176,40){\line(1,0){3}}
\put(176,35){\arc{10}{3.142}{4.712}}
\put(171,15){\line(0,1){20}}
\put(176,15){\arc{10}{1.571}{3.142}}
\put(176,10){\line(1,0){42}}
\put(218,15){\arc{10}{0}{1.571}}
\put(223,15){\line(0,1){20}}
\put(218,35){\arc{10}{4.712}{0}}
\put(215,40){\line(1,0){3}}
\put(210,10){\vector(1,0){0}}
\put(197,20){\line(0,-1){8}}
\put(197,8){\vector(0,-1){8}}
\thinlines
\dottedline{4}(160,35)(234,35)
\put(187,30){\line(1,0){5}}
\put(192,35){\arc{10}{0}{1.571}}
\put(197,35){\line(0,1){35}}
\put(197,50){\vector(0,-1){0}}
\put(171,50){\line(0,1){5}}
\put(176,55){\arc{10}{3.142}{4.712}}
\put(176,60){\line(1,0){19}}
\put(191,60){\vector(1,0){0}}
\put(218,60){\line(-1,0){19}}
\put(218,55){\arc{10}{4.712}{0}}
\put(223,50){\line(0,1){5}}
\put(179,40){\arc{5}{3.142}{0}}
\put(171,50){\arc{5}{3.142}{0}}
\put(223,50){\arc{5}{3.142}{0}}
\put(215,40){\arc{5}{3.142}{0}}
\put(159,4){\makebox(0,0){$d$}}
\put(235,4){\makebox(0,0){$d$}}
\put(184,45){\makebox(0,0){$a$}}
\put(210,45){\makebox(0,0){$a$}}
\put(190,38){\makebox(0,0){$b$}}
\put(204,38){\makebox(0,0){$b$}}
\put(187,19){\makebox(0,0){$c$}}
\put(189,65){\makebox(0,0){$\nu$}}
\put(202,50){\makebox(0,0){$\la$}}
\put(192,4){\makebox(0,0){$\tau$}}
\put(209,4){\makebox(0,0){$\beta_+$}}
\put(183,30){\makebox(0,0){$t$}}
\put(197,24){\makebox(0,0){$X^*$}}
\end{picture}
\end{center}
\caption{The vector
$\pi_{\tau,\la}^+(e_{\beta_+})
\omega^{\tau,\la,+}_{b,c,t,X}\in\cH_{\tau,\la}$}
\label{6jproof1}
\end{figure}
Now let us look at the part of the picture above the
dotted line. In a suitable Frobenius annulus, this part
can be read for fixed $\nu$ and $d$ as
$\sum_i \la (t_i) \epsm\nu\la t_i^*$, and the sum
runs over a full orthonormal basis of isometries
$t_i\in\Hom(\nu,b\co{\beta_+}\co d)$. Next we look
at the part above the dotted line on the right
hand side of Fig.\ \ref{6jproof1}. In the same
Frobenius annulus, this can be similarly read as
$\sum_j \la (s_j) \epsm\nu\la s_j^*$ where the sum runs
over another orthonormal basis of isometries
$s_j\in\Hom(\nu,b\co{\beta_+}\co d)$. Since such bases
are related by a unitary matrix (``unitarity of
$6j$-symbols''), we conclude that both diagrams
represent the same vector in $\cH_{\tau,\la}$.
Now turning around the small arcs at the trivalent
vertices of the wire $\beta_+$ and using the expansion
the identity in the reverse way leads us to the left
hand side of Fig.\ \ref{6jproof2}.
%
\begin{figure}[htb]
\begin{center}
\unitlength 0.6mm
\begin{picture}(221,60)
\thicklines
\put(10,32){\makebox(0,0){$\displaystyle\sum_{a,d} \; d_a$}}
\put(30,0){\line(0,1){40}}
\put(40,40){\arc{20}{3.142}{4.712}}
\put(40,50){\line(1,0){18}}
\put(62,50){\line(1,0){18}}
\put(80,40){\arc{20}{4.712}{0}}
\put(90,0){\line(0,1){40}}
\put(50,25){\line(0,1){15}}
\put(55,40){\arc{10}{3.142}{4.712}}
\put(55,45){\line(1,0){3}}
\put(62,45){\line(1,0){3}}
\put(65,40){\arc{10}{4.712}{0}}
\put(70,25){\line(0,1){15}}
\put(65,25){\arc{10}{0}{1.571}}
\put(55,20){\line(1,0){10}}
\put(55,25){\arc{10}{1.571}{3.142}}
\Thicklines
\put(70,35){\line(1,0){5}}
\put(75,30){\arc{10}{4.712}{0}}
\put(80,15){\line(0,1){15}}
\put(45,15){\arc{10}{1.571}{3.142}}
\put(45,10){\line(1,0){30}}
\put(75,15){\arc{10}{0}{1.571}}
\put(40,15){\line(0,1){15}}
\put(45,30){\arc{10}{3.142}{4.712}}
\put(45,35){\line(1,0){5}}
\put(73,10){\vector(1,0){0}}
\put(60,20){\line(0,-1){8}}
\put(60,8){\vector(0,-1){8}}
\put(60,45){\line(0,1){5}}
\thinlines
\dottedline{4}(35,40)(85,40)(85,5)(35,5)(35,40)
\put(50,28){\line(1,0){5}}
\put(55,33){\arc{10}{0}{1.571}}
\put(60,33){\line(0,1){12}}
\put(60,50){\line(0,1){10}}
\put(60,53){\vector(0,-1){0}}
\put(50,35){\arc{5}{1.571}{4.712}}
\put(70,35){\arc{5}{4.712}{1.571}}
\put(26,4){\makebox(0,0){$d$}}
\put(47,45){\makebox(0,0){$a$}}
\put(54,33){\makebox(0,0){$b$}}
\put(70,18){\makebox(0,0){$b$}}
\put(50,19){\makebox(0,0){$c$}}
\put(65,55){\makebox(0,0){$\la$}}
\put(55,4){\makebox(0,0){$\tau$}}
\put(72,4){\makebox(0,0){$\beta_+$}}
\put(46,28){\makebox(0,0){$t$}}
\put(60,24){\makebox(0,0){$X^*$}}
\thicklines
\put(139,21){\makebox(0,0){$\longleftrightarrow\;
\displaystyle\sum_{c',i,j} \; {{\rm{coeff}}}_{(c',i,j)}$}}
\put(205,20){\line(0,1){20}}
\put(200,20){\arc{10}{0}{1.571}}
\put(190,15){\line(1,0){10}}
\put(190,20){\arc{10}{1.571}{3.142}}
\put(185,20){\line(0,1){20}}
\Thicklines
\put(195,15){\vector(0,-1){10}}
\thinlines
\dottedline{4}(170,40)(220,40)(220,5)(170,5)(170,40)
\put(185,28){\line(1,0){5}}
\put(190,33){\arc{10}{0}{1.571}}
\put(195,33){\line(0,1){7}}
\put(195,34){\vector(0,-1){0}}
\put(181,36){\makebox(0,0){$a$}}
\put(209,36){\makebox(0,0){$a$}}
\put(185,13){\makebox(0,0){$c'$}}
\put(199,35){\makebox(0,0){$\la$}}
\put(200,9){\makebox(0,0){$\tau$}}
\put(177,28){\makebox(0,0){$t_{a,\co{c'}}^{\la;i}$}}
\put(195.5,21){\makebox(0,0){\footnotesize{$(X_{\co{c'},a}^{\tau;j})^*$}}}
\end{picture}
\end{center}
\caption{The vector
$\pi_{\tau,\la}^+(e_{\beta_+})
\omega^{\tau,\la,+}_{b,c,t,X}\in\cH_{\tau,\la}$}
\label{6jproof2}
\end{figure}
Then we look at the part of the picture inside the dotted box.
In a suitable Frobenius annulus, this can be read as an
intertwiner in $\Hom(\co a \la,\tau \co a)$. Since any
element in this space can be expanded in the basis basis
given in the dotted box on the right hand side of
Fig.\ \ref{6jproof2}, we conclude that 
$\pi_{\tau,\la}^+(e_{\beta_+})\omega^{\tau,\la,+}_{b,c,t,X}$
is in fact a linear combination of
$\omega^{\tau,\la,+}_\xi$'s, hence it is in $H_{\tau,\la}^+$.
\end{proof}

Since it is just intertwiner multiplication in each
$\Hom(\la,a\tau\co a)$ block, the prescription
$\omega^{\tau,\la,+}_{b,c,t,X} \mapsto
\pi^+_{\tau,\la}(e_{\beta_+}) \omega^{\tau,\la,+}_{b,c,t,X}$
clearly defines a linear map
$\pi^+_{\tau,\la}(e_{\beta_+}):H_{\tau,\la}^+\mapsto\cH_{\tau,\la}$
for each $\beta_+\in\MXMp$.
From Lemma \ref{chirrep} we now learn that
$\pi^+_{\tau,\la}(e_{\beta_+})$ is in fact
a linear operator on $H_{\tau,\la}^+$. Similarly
$\pi^-_{\tau,\la}(e_{\beta_-})$ is a linear operator
on $H_{\tau,\la}^-$ for each $\beta_-\in\MXMm$.
We therefore obtain linear maps
$\pi^\pm_{\tau,\la}:\cZ_h^\pm\rightarrow B(H_{\tau,\la}^\pm)$
by linear extension of
$e_{\beta_\pm}\mapsto\pi^\pm_{\tau,\la}(e_{\beta_\pm})$,
$\beta_\pm\in\MXMpm$.

\begin{lemma}
\label{lemchirrep}
The maps
$\pi^\pm_{\tau,\la}:\cZ_h^\pm\rightarrow B(H_{\tau,\la}^\pm)$
are in fact linear representations.
\end{lemma}

\begin{proof}
We prove the representation property of
$\pi^+_{\tau,\la}$; the proof for $\pi^-_{\tau,\la}$
is analogous.
For $\beta_+,\beta_+'\in\MXMp$, the vector
$\pi_{\tau,\la}^+(e_{\beta_+}) (\pi_{\tau,\la}^+(e_{\beta_+'})
\omega^{\tau,\la,+}_{b,c,t,X})$
is given graphically by the left hand side of
Fig.\ \ref{proofcr}.
%
\begin{figure}[htb]
\begin{center}
\unitlength 0.6mm
\begin{picture}(226,65)
\thicklines
\put(15,30){\makebox(0,0){$\displaystyle\sum_{a,a',d}
\; d_a d_{a'}$}}
\put(42,0){\line(0,1){45}}
\put(52,45){\arc{20}{3.142}{4.712}}
\put(52,55){\line(1,0){6}}
\put(62,55){\line(1,0){6}}
\put(68,45){\arc{20}{4.712}{0}}
\put(78,0){\line(0,1){45}}
\put(52,35){\line(0,1){10}}
\put(57,45){\arc{10}{3.142}{4.712}}
\put(57,50){\line(1,0){1}}
\put(62,50){\line(1,0){1}}
\put(63,45){\arc{10}{4.712}{0}}
\put(68,35){\line(0,1){10}}
\put(63,35){\arc{10}{0}{1.571}}
\put(57,30){\line(1,0){6}}
\put(57,35){\arc{10}{1.571}{3.142}}
\Thicklines
\put(60,8){\vector(0,-1){8}}
\put(60,18){\line(0,-1){6}}
\put(60,50){\line(0,1){5}}
\put(42,20){\line(1,0){36}}
\put(71,20){\vector(1,0){0}}
\put(42,10){\line(1,0){36}}
\put(71,10){\vector(1,0){0}}
\put(60,22){\line(0,1){8}}
\thinlines
\put(52,40){\line(1,0){3}}
\put(55,45){\arc{10}{0}{1.571}}
\put(60,45){\line(0,1){5}}
\put(60,55){\line(0,1){10}}
\put(42,10){\arc{5}{4.712}{1.571}}
\put(78,10){\arc{5}{1.571}{4.712}}
\put(42,20){\arc{5}{4.712}{1.571}}
\put(78,20){\arc{5}{1.571}{4.712}}
\put(60,58){\vector(0,-1){0}}
\put(38,4){\makebox(0,0){$d$}}
\put(82,4){\makebox(0,0){$d$}}
\put(38,15){\makebox(0,0){$a'$}}
\put(82,15){\makebox(0,0){$a'$}}
\put(38,28){\makebox(0,0){$a$}}
\put(52,29){\makebox(0,0){$c$}}
\put(72,40){\makebox(0,0){$b$}}
\put(55,4){\makebox(0,0){$\tau$}}
\put(69,14){\makebox(0,0){$\beta_+$}}
\put(69,24){\makebox(0,0){$\beta_+'$}}
\put(65,60){\makebox(0,0){$\la$}}
\put(47,40){\makebox(0,0){$t$}}
\put(60,34){\makebox(0,0){$X^*$}}
\thicklines
\put(125,30){\makebox(0,0){$=\;\displaystyle\sum_{a,a',d,\beta_+''}
\; d_a d_{a'}$}}
\put(165,0){\arc{10}{3.142}{4.712}}
\put(165,10){\arc{10}{0}{1.571}}
\put(170,10){\line(0,1){21}}
\put(165,31){\arc{10}{4.712}{0}}
\put(165,41){\arc{10}{1.571}{3.142}}
\put(160,41){\line(0,1){4}}
\put(170,45){\arc{20}{3.142}{4.712}}
\put(170,55){\line(1,0){18}}
\put(192,55){\line(1,0){18}}
\put(210,45){\arc{20}{4.712}{0}}
\put(220,41){\line(0,1){4}}
\put(215,0){\arc{10}{4.712}{0}}
\put(215,10){\arc{10}{1.571}{3.142}}
\put(210,10){\line(0,1){21}}
\put(215,31){\arc{10}{3.142}{4.712}}
\put(215,41){\arc{10}{0}{1.571}}
\put(182,35){\line(0,1){10}}
\put(187,45){\arc{10}{3.142}{4.712}}
\put(187,50){\line(1,0){1}}
\put(192,50){\line(1,0){1}}
\put(193,45){\arc{10}{4.712}{0}}
\put(198,35){\line(0,1){10}}
\put(193,35){\arc{10}{0}{1.571}}
\put(187,30){\line(1,0){6}}
\put(187,35){\arc{10}{1.571}{3.142}}
\Thicklines
\put(190,18.5){\vector(0,-1){18.5}}
\put(190,30){\line(0,-1){7.5}}
\put(190,50){\line(0,1){5}}
\put(180,20.5){\line(1,0){20}}
\put(197,20.5){\vector(1,0){0}}
\put(170,13){\line(1,0){5}}
\put(180,20.5){\vector(0,1){0}}
\put(205,13){\line(1,0){5}}
\put(210,13){\vector(1,0){0}}
\put(170,28){\line(1,0){5}}
\put(180,20.5){\vector(0,-1){0}}
\put(205,28){\line(1,0){5}}
\put(210,28){\vector(1,0){0}}
\put(175,23){\arc{10}{4.712}{0}}
\put(175,18){\arc{10}{0}{1.571}}
\put(205,23){\arc{10}{3.142}{4.712}}
\put(205,18){\arc{10}{1.571}{3.142}}
\put(180,18){\line(0,1){5}}
\put(200,18){\line(0,1){5}}
\thinlines
\put(182,40){\line(1,0){3}}
\put(185,45){\arc{10}{0}{1.571}}
\put(190,45){\line(0,1){5}}
\put(190,55){\line(0,1){10}}
\put(170,13){\arc{5}{4.712}{1.571}}
\put(210,13){\arc{5}{1.571}{4.712}}
\put(170,28){\arc{5}{4.712}{1.571}}
\put(210,28){\arc{5}{1.571}{4.712}}
\put(180,20.5){\arc{5}{4.712}{1.571}}
\put(200,20.5){\arc{5}{1.571}{4.712}}
\put(190,58){\vector(0,-1){0}}
\put(156,4){\makebox(0,0){$d$}}
\put(224,4){\makebox(0,0){$d$}}
\put(165,20.5){\makebox(0,0){$a'$}}
\put(215,20.5){\makebox(0,0){$a'$}}
\put(170,59){\makebox(0,0){$a$}}
\put(182,29){\makebox(0,0){$c$}}
\put(202,40){\makebox(0,0){$b$}}
\put(185,4){\makebox(0,0){$\tau$}}
\put(175,8){\makebox(0,0){$\beta_+$}}
\put(175,23){\makebox(0,0){$\beta_+'$}}
\put(205,8){\makebox(0,0){$\beta_+$}}
\put(205,23){\makebox(0,0){$\beta_+'$}}
\put(196,15){\makebox(0,0){$\beta_+''$}}
\put(195,60){\makebox(0,0){$\la$}}
\put(177,40){\makebox(0,0){$t$}}
\put(190,34){\makebox(0,0){$X^*$}}
\end{picture}
\end{center}
\caption{The vector
$\pi_{\tau,\la}^+(e_{\beta_+}) \pi_{\tau,\la}^+(e_{\beta_+'})
\omega^{\tau,\la,+}_{b,c,t,X}\in H_{\tau,\la}^+$}
\label{proofcr}
\end{figure}
Next we use the expansion of the identity
(cf.\ \cite[Lemma 4.3]{BEK1}) for the parallel
wires $\beta_+$ and $\beta_+'$ on, say, the
left hand side of the crossings with the wire $\tau$.
Note that only $\beta_+''\in\MXMp$ can contribute
because $\Hom(\beta_+\beta_+',\beta_+'')=0$ otherwise.
Application of the braiding fusion relation
for the relative braiding, Fig.\ \ref{wireErbfe},
yields the right hand side of Fig.\ \ref{proofcr}.
Using expansions of the identity also for the
parallel pieces of the wires $a$ and $d$ on the left and
on the right, we obtain a picture where the bottom
part coincides with the wire diagram in
\cite[Fig.\ 42]{BEK1}, up to the crossing with the wire
$\tau$. In fact we can use the same argument (``unitarity
of $6j$-symbols'') as in the proof of
\cite[Thm.\ 4.4]{BEK1} to obtain the desired result
\[ \pi_{\tau,\la}^+(e_{\beta_+}) (\pi_{\tau,\la}^+(e_{\beta_+'})
\omega^{\tau,\la,+}_{b,c,t,X} ) = \sum_{\beta_+''\in\MXMp}
\frac{d_{\beta_+} d_{\beta_+'}}{d_{\beta_+''}}
N_{\beta_+,\beta_+'}^{\beta_+''} \,\pi_{\tau,\la}^+(e_{\beta_+'})
\omega^{\tau,\la,+}_{b,c,t,X} \,. \]
As the prefactors coincide with those in the decomposition
of the vertical product
$e_{\beta_+} *_v e_{\beta_+'}$ into $e_{\beta_+''}$'s,
the claim is proven.
\end{proof}

\section{Chiral structure of the center $\cZ_h$}
\label{Zpm}

In this section we will analyze the chiral systems
$\MXMpm$ in the non-degenerate case, i.e.\ from now on
we impose the following

\begin{assumption}
\label{assnondeg}{\rm
In addition to Assumption \ref{assbraid}, we now assume 
that the braiding on $\NXN$ is non-degenerate in the
sense of \cite[Def.\ 2.3]{BEK1}.
}\end{assumption}

\subsection{Non-degeneracy of the ambichiral braiding}

We define $w_0=\sum_{\beta\in\MXMo}d_\beta^2$ and call it the
{\sl ambichiral global index}.

\begin{theorem}
\label{ambinondeg}
The braiding on the ambichiral system $\MXMo$ arising from the
relative braiding of the chiral systems is non-degenerate.
Moreover, the ambichiral global index is given by $w_0=w_+^2/w$.
\end{theorem}

\begin{proof}
From Lemma \cite[Thm.\ 6.8]{BEK1} we obtain
$\sum_{\la,\mu} q_{\la,\mu} *_h e_\tau =
\del \tau 0 e_0$. The left hand side is displayed
graphically by the left hand side of Fig.\ \ref{ndab1}.
%
\thinlines
\begin{figure}[htb]
\begin{center}
\unitlength 0.6mm
\begin{picture}(214,50)
\thicklines
\put(18,25){\makebox(0,0){$\displaystyle\sum_{a,b,c,d,\la,\mu}
\frac{d_b d_c}{w^2}$}}
\put(42,0){\line(0,1){10}}
\put(80.5,22.5){\line(0,1){5}}
\put(52,10){\arc{20}{3.142}{4.712}}
\put(78,27.5){\arc{5}{4.712}{0}}
\put(52,10){\arc{10}{1.571}{4.712}}
\put(78,10){\arc{10}{4.712}{1.571}}
\put(52,40){\arc{20}{1.571}{3.142}}
\put(78,22.5){\arc{5}{0}{1.571}}
\put(52,40){\arc{10}{1.571}{4.712}}
\put(78,40){\arc{10}{4.712}{1.571}}
\put(52,5){\line(1,0){26}}
\put(52,15){\line(1,0){1}}
\put(57,15){\line(1,0){21}}
\put(52,20){\line(1,0){1}}
\put(57,20){\line(1,0){21}}
\put(42,50){\line(0,-1){10}}
\put(96,0){\line(0,1){50}}
\put(52,30){\line(1,0){1}}
\put(57,30){\line(1,0){21}}
\put(52,35){\line(1,0){1}}
\put(57,35){\line(1,0){21}}
\put(52,45){\line(1,0){26}}
\thinlines
\put(55,5){\line(0,1){10}}
\put(55,20){\line(0,1){10}}
\put(75,5){\line(0,1){8}}
\put(75,22){\line(0,1){6}}
\put(55,23){\vector(0,-1){0}}
\put(75,27){\vector(0,1){0}}
\put(55,45){\line(0,-1){10}}
\put(75,45){\line(0,-1){8}}
\put(55,5){\arc{5}{3.142}{0}}
\put(75,5){\arc{5}{3.142}{0}}
\put(55,45){\arc{5}{0}{3.142}}
\put(75,45){\arc{5}{0}{3.142}}
\put(80.5,25){\arc{5}{4.712}{1.571}}
\put(96,25){\arc{5}{1.571}{4.712}}
\Thicklines
\put(80.5,25){\line(1,0){15.5}}
\put(90.25,25){\vector(1,0){0}}
\put(75,17){\line(0,1){1}}
\put(55,15){\line(0,1){5}}
\put(75,33){\line(0,-1){1}}
\put(55,35){\line(0,-1){5}}
\put(38,4){\makebox(0,0){$d$}}
\put(100,4){\makebox(0,0){$d$}}
\put(65,2){\makebox(0,0){$c$}}
\put(65,11){\makebox(0,0){$b$}}
\put(38,46){\makebox(0,0){$a$}}
\put(100,46){\makebox(0,0){$a$}}
\put(65,48){\makebox(0,0){$c$}}
\put(65,39){\makebox(0,0){$b$}}
\put(50,25){\makebox(0,0){$\la$}}
\put(70,24){\makebox(0,0){$\mu$}}
\put(88.25,20){\makebox(0,0){$\tau$}}
\thicklines
\put(133,25){\makebox(0,0){$=\;\displaystyle\sum_{a,b,c,d,\la,\mu}
\frac{d_b d_c}{w^2}$}}
\put(162,0){\line(0,1){50}}
\put(172,10){\arc{10}{1.571}{4.712}}
\put(198,10){\arc{10}{4.712}{1.571}}
\put(172,40){\arc{10}{1.571}{4.712}}
\put(198,40){\arc{10}{4.712}{1.571}}
\put(172,5){\line(1,0){26}}
\put(162,50){\line(0,-1){10}}
\put(208,0){\line(0,1){50}}
\put(172,15){\line(1,0){1}}
\put(177,15){\line(1,0){21}}
\put(172,35){\line(1,0){1}}
\put(177,35){\line(1,0){21}}
\put(172,45){\line(1,0){26}}
\thinlines
\put(175,5){\line(0,1){10}}
\put(195,5){\line(0,1){8}}
\put(175,45){\line(0,-1){10}}
\put(195,45){\line(0,-1){8}}
\put(175,5){\arc{5}{3.142}{0}}
\put(195,5){\arc{5}{3.142}{0}}
\put(175,45){\arc{5}{0}{3.142}}
\put(195,45){\arc{5}{0}{3.142}}
\put(162,25){\arc{5}{4.712}{1.571}}
\put(208,25){\arc{5}{1.571}{4.712}}
\Thicklines
\put(175,28.5){\vector(0,-1){0}}
\put(195,32.5){\vector(0,1){0}}
\put(162,25){\line(1,0){11}}
\put(177,25){\line(1,0){31}}
\put(187,25){\vector(1,0){0}}
\put(195,17){\line(0,1){6}}
\put(175,15){\line(0,1){20}}
\put(195,33){\line(0,-1){6}}
\put(158,4){\makebox(0,0){$d$}}
\put(212,4){\makebox(0,0){$d$}}
\put(185,2){\makebox(0,0){$c$}}
\put(185,11){\makebox(0,0){$b$}}
\put(158,46){\makebox(0,0){$a$}}
\put(212,46){\makebox(0,0){$a$}}
\put(185,48){\makebox(0,0){$c$}}
\put(185,39){\makebox(0,0){$b$}}
\put(169,30){\makebox(0,0){\small{$\a^+_\la$}}}
\put(201,30){\makebox(0,0){\small{$\a^-_\mu$}}}
\put(185,20){\makebox(0,0){$\tau$}}
\end{picture}
\end{center}
\caption{Non-degeneracy of the ambichiral braiding}
\label{ndab1}
\end{figure}
We can ``pull in'' the wire $\tau$ since it admits relative
braiding with both $\a^+_\la$ and $\a^-_\mu$, and this way we
obtain the right hand side of Fig.\ \ref{ndab1}.
We can use the expansion of the identity for the
parallel wires $b,c$ on the top and bottom
(cf.\ \cite[Lemma 4.3]{BEK1}) to obtain the left
hand side of Fig.\ \ref{ndab2}.
%
\begin{figure}[htb]
\begin{center}
\unitlength 0.6mm
\begin{picture}(240,60)
\thicklines
\put(18,29){\makebox(0,0){$\displaystyle
\sum_{a,b,c,d,\atop\la,\mu,\tau',\tau''}
\frac{d_b d_c}{w^2}$}}
\put(47,0){\line(0,1){60}}
\put(103,0){\line(0,1){60}}
\put(57,10){\arc{10}{1.571}{3.142}}
\put(57,15){\arc{10}{3.142}{4.712}}
\put(93,15){\arc{10}{4.712}{0}}
\put(93,10){\arc{10}{0}{1.571}}
\put(87,10){\arc{10}{1.571}{3.142}}
\put(87,15){\arc{10}{3.142}{4.712}}
\put(63,15){\arc{10}{4.712}{0}}
\put(63,10){\arc{10}{0}{1.571}}
\put(57,45){\arc{10}{1.571}{3.142}}
\put(57,50){\arc{10}{3.142}{4.712}}
\put(93,50){\arc{10}{4.712}{0}}
\put(93,45){\arc{10}{0}{1.571}}
\put(87,45){\arc{10}{1.571}{3.142}}
\put(87,50){\arc{10}{3.142}{4.712}}
\put(63,50){\arc{10}{4.712}{0}}
\put(63,45){\arc{10}{0}{1.571}}
\put(57,5){\line(1,0){6}}
\put(87,5){\line(1,0){6}}
\put(57,20){\line(1,0){1}}
\put(62,20){\line(1,0){1}}
\put(87,20){\line(1,0){6}}
\put(57,40){\line(1,0){1}}
\put(62,40){\line(1,0){1}}
\put(87,40){\line(1,0){6}}
\put(57,55){\line(1,0){6}}
\put(87,55){\line(1,0){6}}
\put(52,10){\line(0,1){5}}
\put(52,50){\line(0,-1){5}}
\put(68,10){\line(0,1){5}}
\put(68,50){\line(0,-1){5}}
\put(98,10){\line(0,1){5}}
\put(98,50){\line(0,-1){5}}
\put(82,10){\line(0,1){5}}
\put(82,50){\line(0,-1){5}}
\thinlines
\put(60,5){\line(0,1){15}}
\put(90,5){\line(0,1){13}}
\put(60,55){\line(0,-1){15}}
\put(90,55){\line(0,-1){13}}
\put(60,5){\arc{5}{3.142}{0}}
\put(90,5){\arc{5}{3.142}{0}}
\put(60,55){\arc{5}{0}{3.142}}
\put(90,55){\arc{5}{0}{3.142}}
\put(82,12.5){\arc{5}{1.571}{4.712}}
\put(68,12.5){\arc{5}{4.712}{1.571}}
\put(82,47.5){\arc{5}{1.571}{4.712}}
\put(68,47.5){\arc{5}{4.712}{1.571}}
\put(103,30){\arc{5}{1.571}{4.712}}
\put(47,30){\arc{5}{4.712}{1.571}}
\Thicklines
\put(68,12.5){\line(1,0){14}}
\put(68,47.5){\line(1,0){14}}
\put(77,12.5){\vector(1,0){0}}
\put(73,47.5){\vector(-1,0){0}}
\put(90,22){\line(0,1){6}}
\put(60,20){\line(0,1){20}}
\put(90,38){\line(0,-1){6}}
\put(60,33){\vector(0,-1){0}}
\put(90,27){\vector(0,1){0}}
\put(47,30){\line(1,0){11}}
\put(62,30){\line(1,0){41}}
\put(77,30){\vector(1,0){0}}
\put(43,55){\makebox(0,0){$a$}}
\put(43,5){\makebox(0,0){$d$}}
\put(107,55){\makebox(0,0){$a$}}
\put(107,5){\makebox(0,0){$d$}}
\put(65,3){\makebox(0,0){$c$}}
\put(55,3){\makebox(0,0){$b$}}
\put(85,3){\makebox(0,0){$c$}}
\put(95,3){\makebox(0,0){$b$}}
\put(65,57){\makebox(0,0){$c$}}
\put(55,57){\makebox(0,0){$b$}}
\put(85,57){\makebox(0,0){$c$}}
\put(95,57){\makebox(0,0){$b$}}
\put(67,35){\makebox(0,0){$\a^+_\la$}}
\put(83,25){\makebox(0,0){$\a^-_\mu$}}
\put(75,7.5){\makebox(0,0){$\tau''$}}
\put(75,52.5){\makebox(0,0){$\tau'$}}
\put(72,25){\makebox(0,0){$\tau$}}
\thicklines
\put(137,29){\makebox(0,0){$=\,\displaystyle
\sum_{a,b,c,d,\atop\la,\mu,\tau',\tau''}
\frac{d_b d_c}{w^2}$}}
\put(167,0){\line(0,1){60}}
\put(233,0){\line(0,1){60}}
\put(177,20){\arc{10}{1.571}{3.142}}
\put(177,25){\arc{10}{3.142}{4.712}}
\put(223,25){\arc{10}{4.712}{0}}
\put(223,20){\arc{10}{0}{1.571}}
\put(217,20){\arc{10}{1.571}{3.142}}
\put(217,25){\arc{10}{3.142}{4.712}}
\put(183,25){\arc{10}{4.712}{0}}
\put(183,20){\arc{10}{0}{1.571}}
\put(177,45){\arc{10}{1.571}{3.142}}
\put(177,50){\arc{10}{3.142}{4.712}}
\put(223,50){\arc{10}{4.712}{0}}
\put(223,45){\arc{10}{0}{1.571}}
\put(217,45){\arc{10}{1.571}{3.142}}
\put(217,50){\arc{10}{3.142}{4.712}}
\put(183,50){\arc{10}{4.712}{0}}
\put(183,45){\arc{10}{0}{1.571}}
\put(177,15){\line(1,0){6}}
\put(217,15){\line(1,0){6}}
\put(177,30){\line(1,0){1}}
\put(182,30){\line(1,0){1}}
\put(217,30){\line(1,0){6}}
\put(177,40){\line(1,0){1}}
\put(182,40){\line(1,0){1}}
\put(217,40){\line(1,0){6}}
\put(177,55){\line(1,0){6}}
\put(217,55){\line(1,0){6}}
\put(172,20){\line(0,1){5}}
\put(172,50){\line(0,-1){5}}
\put(188,20){\line(0,1){5}}
\put(188,50){\line(0,-1){5}}
\put(228,20){\line(0,1){5}}
\put(228,50){\line(0,-1){5}}
\put(212,20){\line(0,1){5}}
\put(212,50){\line(0,-1){5}}
\thinlines
\put(212,22.5){\arc{5}{1.571}{4.712}}
\put(188,22.5){\arc{5}{4.712}{1.571}}
\put(212,47.5){\arc{5}{1.571}{4.712}}
\put(188,47.5){\arc{5}{4.712}{1.571}}
\put(228,22.5){\arc{5}{1.571}{4.712}}
\put(172,22.5){\arc{5}{4.712}{1.571}}
\put(228,47.5){\arc{5}{1.571}{4.712}}
\put(172,47.5){\arc{5}{4.712}{1.571}}
\put(233,10){\arc{5}{1.571}{4.712}}
\put(167,10){\arc{5}{4.712}{1.571}}
\put(225,27.5){\arc{10}{1.571}{3.142}}
\put(225,42.5){\arc{10}{3.142}{4.712}}
\put(175,42.5){\arc{10}{4.712}{0}}
\put(175,27.5){\arc{10}{0}{1.571}}
\put(172,22.5){\line(1,0){3}}
\put(172,47.5){\line(1,0){3}}
\put(228,22.5){\line(-1,0){3}}
\put(228,47.5){\line(-1,0){3}}
\put(180,30){\line(0,-1){2.5}}
\put(180,40){\line(0,1){2.5}}
\put(220,28){\line(0,-1){0.5}}
\put(220,42){\line(0,1){0.5}}
\Thicklines
\put(167,10){\line(1,0){23}}
\put(210,10){\line(1,0){23}}
\put(190,15){\arc{10}{0}{1.571}}
\put(210,15){\arc{10}{1.571}{3.142}}
\put(200,30){\arc{10}{3.142}{0}}
\put(195,24.5){\line(0,1){5.5}}
\put(195,15){\line(0,1){5.5}}
\put(205,15){\line(0,1){15}}
\put(188,22.5){\line(1,0){15}}
\put(207,22.5){\line(1,0){5}}
\put(188,47.5){\line(1,0){24}}
\put(202,22.5){\vector(1,0){0}}
\put(198,47.5){\vector(-1,0){0}}
\put(220,32){\line(0,1){6}}
\put(217,10){\vector(1,0){0}}
\put(180,30){\line(0,1){10}}
\put(180,33){\vector(0,-1){0}}
\put(220,37){\vector(0,1){0}}
\put(163,55){\makebox(0,0){$a$}}
\put(163,5){\makebox(0,0){$d$}}
\put(238,55){\makebox(0,0){$a$}}
\put(238,5){\makebox(0,0){$d$}}
\put(185,26){\makebox(0,0){$b$}}
\put(175,13){\makebox(0,0){$c$}}
\put(215,26){\makebox(0,0){$b$}}
\put(225,13){\makebox(0,0){$c$}}
\put(185,44){\makebox(0,0){$b$}}
\put(175,57){\makebox(0,0){$c$}}
\put(215,44){\makebox(0,0){$b$}}
\put(225,57){\makebox(0,0){$c$}}
\put(175,35){\makebox(0,0){\small{$\a^+_\la$}}}
\put(226,34){\makebox(0,0){\small{$\a^-_\mu$}}}
\put(215,5){\makebox(0,0){$\tau$}}
\put(208,52.5){\makebox(0,0){$\tau'$}}
\put(200,17.5){\makebox(0,0){$\tau''$}}
\end{picture}
\end{center}
\caption{Non-degeneracy of the ambichiral braiding}
\label{ndab2}
\end{figure}
Here only ambichiral morphisms $\tau,\tau'\in\MXMo$ contribute
in the corresponding sums over $\beta',\beta''\in\MXM$ since
they appear between $\a^+_\la$ and $\a^-_\mu$.
Application of the naturality moves of Figs.\ \ref{natrelui}
and \ref{natreloi} for the relative braiding yields
the right hand side of Fig.\ \ref{ndab2}.
Now we see that intertwiners in $\Hom(\tau',\tau'')$
appear so that we first obtain a factor $\del {\tau'}{\tau''}$.
Then we take the scalar part of the loop separately to obtain
Fig.\ \ref{ndab3},
%
\begin{figure}[htb]
\begin{center}
\unitlength 0.6mm
\begin{picture}(156,60)
\thicklines
\put(24,29){\makebox(0,0){$\displaystyle\sum_{a,b,c,d,\la,\mu,\tau'}
\frac{d_b d_c}{w^2d_{\tau'}}$}}
\put(62,10){\arc{10}{1.571}{3.142}}
\put(62,15){\arc{10}{3.142}{4.712}}
\put(98,15){\arc{10}{4.712}{0}}
\put(98,10){\arc{10}{0}{1.571}}
\put(92,10){\arc{10}{1.571}{3.142}}
\put(92,15){\arc{10}{3.142}{4.712}}
\put(68,15){\arc{10}{4.712}{0}}
\put(68,10){\arc{10}{0}{1.571}}
\put(62,45){\arc{10}{1.571}{3.142}}
\put(62,50){\arc{10}{3.142}{4.712}}
\put(98,50){\arc{10}{4.712}{0}}
\put(98,45){\arc{10}{0}{1.571}}
\put(92,45){\arc{10}{1.571}{3.142}}
\put(92,50){\arc{10}{3.142}{4.712}}
\put(68,50){\arc{10}{4.712}{0}}
\put(68,45){\arc{10}{0}{1.571}}
\put(62,5){\line(1,0){6}}
\put(92,5){\line(1,0){6}}
\put(62,20){\line(1,0){1}}
\put(67,20){\line(1,0){1}}
\put(92,20){\line(1,0){6}}
\put(62,40){\line(1,0){1}}
\put(67,40){\line(1,0){1}}
\put(92,40){\line(1,0){6}}
\put(62,55){\line(1,0){6}}
\put(92,55){\line(1,0){6}}
\put(57,10){\line(0,1){5}}
\put(57,50){\line(0,-1){5}}
\put(73,10){\line(0,1){5}}
\put(73,50){\line(0,-1){5}}
\put(103,10){\line(0,1){5}}
\put(103,50){\line(0,-1){5}}
\put(87,10){\line(0,1){5}}
\put(87,50){\line(0,-1){5}}
\put(120,0){\line(0,1){60}}
\put(150,0){\line(0,1){60}}
\thinlines
\put(103,12.5){\arc{5}{1.571}{4.712}}
\put(57,12.5){\arc{5}{4.712}{1.571}}
\put(103,47.5){\arc{5}{1.571}{4.712}}
\put(57,47.5){\arc{5}{4.712}{1.571}}
\put(87,12.5){\arc{5}{1.571}{4.712}}
\put(73,12.5){\arc{5}{4.712}{1.571}}
\put(87,47.5){\arc{5}{1.571}{4.712}}
\put(73,47.5){\arc{5}{4.712}{1.571}}
\put(120,30){\arc{5}{4.712}{1.571}}
\put(150,30){\arc{5}{1.571}{4.712}}
\put(60,17.5){\arc{10}{0}{1.571}}
\put(60,42.5){\arc{10}{4.712}{0}}
\put(100,17.5){\arc{10}{1.571}{3.142}}
\put(100,42.5){\arc{10}{3.142}{4.712}}
\put(57,12.5){\line(1,0){3}}
\put(57,47.5){\line(1,0){3}}
\put(103,12.5){\line(-1,0){3}}
\put(103,47.5){\line(-1,0){3}}
\put(65,20){\line(0,-1){2.5}}
\put(65,40){\line(0,1){2.5}}
\put(95,20){\line(-1,0){2.5}}
\put(95,40){\line(-1,0){2.5}}
\Thicklines
\put(73,12.5){\line(1,0){14}}
\put(73,47.5){\line(1,0){14}}
\put(82,12.5){\vector(1,0){0}}
\put(78,47.5){\vector(-1,0){0}}
\put(65,20){\line(0,1){20}}
\put(95,22){\line(0,1){16}}
\put(65,28){\vector(0,-1){0}}
\put(95,32){\vector(0,1){0}}
\put(120,30){\line(1,0){6}}
\put(150,30){\line(-1,0){20}}
\put(138,30){\vector(1,0){0}}
\put(128,28){\line(0,1){4}}
\put(133,23){\line(1,0){4}}
\put(133,37){\line(1,0){4}}
\put(133,32){\arc{10}{3.142}{4.712}}
\put(137,32){\arc{10}{4.712}{0}}
\put(133,28){\arc{10}{1.571}{3.142}}
\put(137,28){\arc{10}{0}{1.571}}
\put(138,23){\vector(1,0){0}}
\put(71,23){\makebox(0,0){$b$}}
\put(60,3){\makebox(0,0){$c$}}
\put(89,23){\makebox(0,0){$b$}}
\put(100,3){\makebox(0,0){$c$}}
\put(71,37){\makebox(0,0){$b$}}
\put(60,57){\makebox(0,0){$c$}}
\put(89,37){\makebox(0,0){$b$}}
\put(100,57){\makebox(0,0){$c$}}
\put(116,55){\makebox(0,0){$a$}}
\put(154,55){\makebox(0,0){$a$}}
\put(116,5){\makebox(0,0){$d$}}
\put(154,5){\makebox(0,0){$d$}}
\put(58,30){\makebox(0,0){$\a^+_\la$}}
\put(102,30){\makebox(0,0){$\a^-_\mu$}}
\put(80,7.5){\makebox(0,0){$\tau'$}}
\put(80,52.5){\makebox(0,0){$\tau'$}}
\put(124,25){\makebox(0,0){$\tau$}}
\put(135,18){\makebox(0,0){$\tau'$}}
\end{picture}
\end{center}
\caption{Non-degeneracy of the ambichiral braiding}
\label{ndab3}
\end{figure}
where we need a compensating factor $1/d_{\tau'}$.
By using the ($+$ and $-$ version of the) graphical identity
of Fig.\ \ref{scg=cp} we obtain the left hand side of
Fig.\ \ref{ndab4}.
%
\begin{figure}[htb]
\begin{center}
\unitlength 0.6mm
\begin{picture}(237,60)
\thicklines
\put(20,29){\makebox(0,0){$\displaystyle\sum_{a,b,c,d,\tau'}
\frac{d_b d_c}{w_+^2 d_{\tau'}}$}}
\put(47,10){\arc{10}{1.571}{3.142}}
\put(47,15){\arc{10}{3.142}{4.712}}
\put(83,15){\arc{10}{4.712}{0}}
\put(83,10){\arc{10}{0}{1.571}}
\put(77,10){\arc{10}{1.571}{3.142}}
\put(77,15){\arc{10}{3.142}{4.712}}
\put(53,15){\arc{10}{4.712}{0}}
\put(53,10){\arc{10}{0}{1.571}}
\put(47,45){\arc{10}{1.571}{3.142}}
\put(47,50){\arc{10}{3.142}{4.712}}
\put(83,50){\arc{10}{4.712}{0}}
\put(83,45){\arc{10}{0}{1.571}}
\put(77,45){\arc{10}{1.571}{3.142}}
\put(77,50){\arc{10}{3.142}{4.712}}
\put(53,50){\arc{10}{4.712}{0}}
\put(53,45){\arc{10}{0}{1.571}}
\put(47,5){\line(1,0){6}}
\put(77,5){\line(1,0){6}}
\put(47,20){\line(1,0){6}}
\put(77,20){\line(1,0){6}}
\put(47,40){\line(1,0){6}}
\put(77,40){\line(1,0){6}}
\put(47,55){\line(1,0){6}}
\put(77,55){\line(1,0){6}}
\put(42,10){\line(0,1){5}}
\put(42,50){\line(0,-1){5}}
\put(58,10){\line(0,1){5}}
\put(58,50){\line(0,-1){5}}
\put(88,10){\line(0,1){5}}
\put(88,50){\line(0,-1){5}}
\put(72,10){\line(0,1){5}}
\put(72,50){\line(0,-1){5}}
\put(105,0){\line(0,1){60}}
\put(135,0){\line(0,1){60}}
\thinlines
\put(50,20){\arc{5}{3.142}{0}}
\put(80,20){\arc{5}{3.142}{0}}
\put(50,40){\arc{5}{0}{3.142}}
\put(80,40){\arc{5}{0}{3.142}}
\put(72,12.5){\arc{5}{1.571}{4.712}}
\put(58,12.5){\arc{5}{4.712}{1.571}}
\put(72,47.5){\arc{5}{1.571}{4.712}}
\put(58,47.5){\arc{5}{4.712}{1.571}}
\put(105,30){\arc{5}{4.712}{1.571}}
\put(135,30){\arc{5}{1.571}{4.712}}
\Thicklines
\put(58,12.5){\line(1,0){14}}
\put(58,47.5){\line(1,0){14}}
\put(67,12.5){\vector(1,0){0}}
\put(63,47.5){\vector(-1,0){0}}
\put(50,20){\line(0,1){20}}
\put(80,20){\line(0,1){20}}
\put(50,28){\vector(0,-1){0}}
\put(80,32){\vector(0,1){0}}
\put(105,30){\line(1,0){6}}
\put(135,30){\line(-1,0){20}}
\put(122,30){\vector(1,0){0}}
\put(113,28){\line(0,1){4}}
\put(118,23){\line(1,0){4}}
\put(118,37){\line(1,0){4}}
\put(118,32){\arc{10}{3.142}{4.712}}
\put(122,32){\arc{10}{4.712}{0}}
\put(118,28){\arc{10}{1.571}{3.142}}
\put(122,28){\arc{10}{0}{1.571}}
\put(122,23){\vector(1,0){0}}
\put(56,23){\makebox(0,0){$b$}}
\put(45,3){\makebox(0,0){$c$}}
\put(74,23){\makebox(0,0){$b$}}
\put(85,3){\makebox(0,0){$c$}}
\put(56,37){\makebox(0,0){$b$}}
\put(45,57){\makebox(0,0){$c$}}
\put(74,37){\makebox(0,0){$b$}}
\put(85,57){\makebox(0,0){$c$}}
\put(101,55){\makebox(0,0){$a$}}
\put(139,55){\makebox(0,0){$a$}}
\put(101,5){\makebox(0,0){$d$}}
\put(139,5){\makebox(0,0){$d$}}
\put(45,30){\makebox(0,0){$\tau'$}}
\put(85,30){\makebox(0,0){$\tau'$}}
\put(65,7.5){\makebox(0,0){$\tau'$}}
\put(65,52.5){\makebox(0,0){$\tau'$}}
\put(109,25){\makebox(0,0){$\tau$}}
\put(120,18){\makebox(0,0){$\tau'$}}
\thicklines
\put(170,29){\makebox(0,0){$= \,\displaystyle\sum_{a,b,c,d,\tau'}
\frac{d_bd_c N_{\co c,b}^{\tau'}}{w_+^2}$}}
\put(200,0){\line(0,1){60}}
\put(230,0){\line(0,1){60}}
\put(217,30){\vector(1,0){0}}
\thinlines
\put(200,30){\arc{5}{4.712}{1.571}}
\put(230,30){\arc{5}{1.571}{4.712}}
\Thicklines
\put(200,30){\line(1,0){6}}
\put(230,30){\line(-1,0){20}}
\put(208,28){\line(0,1){4}}
\put(213,23){\line(1,0){4}}
\put(213,37){\line(1,0){4}}
\put(213,32){\arc{10}{3.142}{4.712}}
\put(217,32){\arc{10}{4.712}{0}}
\put(213,28){\arc{10}{1.571}{3.142}}
\put(217,28){\arc{10}{0}{1.571}}
\put(217,23){\vector(1,0){0}}
\put(196,55){\makebox(0,0){$a$}}
\put(234,55){\makebox(0,0){$a$}}
\put(196,5){\makebox(0,0){$d$}}
\put(234,5){\makebox(0,0){$d$}}
\put(204,25){\makebox(0,0){$\tau$}}
\put(215,18){\makebox(0,0){$\tau'$}}
\end{picture}
\end{center}
\caption{Non-degeneracy of the ambichiral braiding}
\label{ndab4}
\end{figure}
Here we used the fact that only the wire $\tau'$ survives
the summations over $\beta\in\MXMpm$ of the chiral horizontal
projectors. The ``bulbs'' give just inner products of
basis isometries. Due to the summation over internal
fusion channels we obtain therefore a multiplicity
$N_{\co c,b}^{\tau'}$ with a closed wire
$\tau'$, evaluated as $d_{\tau'}$. Thus we are left
with the right hand side of Fig.\ \ref{ndab4}.
Note that
$\sum_{b,c} d_bd_c N_{\co c,b}^{\tau'}
=\sum_b d_b^2 d_{\tau'}=w d_{\tau'}$.
Now the $\Hom(a\co a,d\co d)$ part of the right hand side of
Fig.\ \ref{ndab4} must
be equal to the $\Hom(a\co a,d\co d)$ part of $\del \tau 0 e_0$.
Sandwiching this with basis (co-) isometries yields the identity
displayed in Fig. \ref{ndab5}.
%
\begin{figure}[htb]
\begin{center}
\unitlength 0.6mm
\begin{picture}(110,26)
\Thicklines
\put(10,13){\makebox(0,0){$\displaystyle\sum_{\tau'}d_{\tau'}$}}
\put(30,15){\line(1,0){6}}
\put(60,15){\line(-1,0){20}}
\put(47,15){\vector(1,0){0}}
\put(38,13){\line(0,1){4}}
\put(43,8){\line(1,0){4}}
\put(43,22){\line(1,0){4}}
\put(43,17){\arc{10}{3.142}{4.712}}
\put(47,17){\arc{10}{4.712}{0}}
\put(43,13){\arc{10}{1.571}{3.142}}
\put(47,13){\arc{10}{0}{1.571}}
\put(47,8){\vector(1,0){0}}
\put(56,10){\makebox(0,0){$\tau$}}
\put(45,3){\makebox(0,0){$\tau'$}}
\put(93,15){\makebox(0,0){$=\,\del \tau 0\,\displaystyle\frac{w_+^2}w$}}
\end{picture}
\end{center}
\caption{Non-degeneracy of the ambichiral braiding}
\label{ndab5}
\end{figure}
This is the orthogonality relation showing that the braiding
on the ambichiral system is non-degenerate
(cf.\ \cite[Fig.\ 20]{BEK1}). Consequently the number
$w_+^2/w$ must be $w_0$, the ambichiral global index.
\end{proof}

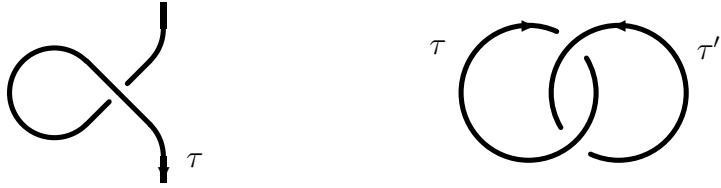
\begin{figure}[htb]
\begin{center}
\unitlength 0.6mm
\begin{picture}(157,40)
\Thicklines
\put(10,20){\arc{20}{0.785}{5.498}}
\put(24.142,20){\line(-1,1){7.071}}
\put(22.142,18){\line(-1,-1){5.071}}
\put(24.142,20){\line(1,-1){7.071}}
\put(26.142,22){\line(1,1){5.071}}
\put(24.142,34.142){\arc{20}{0}{0.785}}
\put(24.142,5.858){\arc{20}{5.498}{6.283}}
\put(34.142,40){\line(0,-1){5.858}}
\put(34.142,5.858){\vector(0,-1){5.858}}
\put(41.142,5){\makebox(0,0){$\tau$}}
\put(115,20){\arc{30}{5.742}{5.142}}
\put(135,20){\arc{30}{2.601}{2.001}}
\put(117,34.8){\vector(1,0){0}}
\put(133,34.8){\vector(-1,0){0}}
\put(95,30){\makebox(0,0){$\tau$}}
\put(155,30){\makebox(0,0){$\tau'$}}
\end{picture}
\end{center}
\caption{Statistics phase $\om_\tau$ and Y-matrix element
$Y_{\tau,\tau'}^\ext$ for the ambichiral system}
\label{omYext}
\end{figure}

Let us define scalars $\omega_\tau,Y_{\tau,\tau'}^\ext\in\bbC$ by
\[ R_\tau^* \Epsr {\co\tau}\tau ^* {\co R}_\tau = \omega_\tau \bfe \,,
\qquad d_\tau d_{\tau'} \phi_\tau (\Epsr {\tau'}\tau \Epsr \tau{\tau'}
)^* = Y_{\tau,\tau'}^\ext\bfe \,, \]
for $\tau,\tau'\in\MXMo$. Note that these numbers can be displayed
graphically as in Fig.\ \ref{omYext}.

%
Putting also
$c_0 = 4 {{\rm{arg}}}(\sum_{\tau\in\MXMo} d_\tau^2 \omega_\tau)/\pi$
we obtain from Theorem \ref{ambinondeg} the following

\begin{corollary}
Matrices $S^\ext$ and $T^\ext$ with matrix elements
$S_{\tau,\tau'}^\ext=w_0^{-1/2} Y_{\tau,\tau'}^\ext$
and $T_{\tau,\tau'}^\ext = \E^{-\pi\I c_0/12} \omega_\tau
\del \tau{\tau'}$, $\tau,\tau'\in\MXMo$, obey the
full Verlinde modular algebra and diagonalize the
fusion rules of the ambichiral system.
\end{corollary}

\subsection{Chiral matrix units}
\begin{figure}[htb]
\begin{center}
\unitlength 0.6mm
\begin{picture}(104,70)
\thicklines
\put(15,34){\makebox(0,0){$\displaystyle\sum_{a,d}$}}
\put(37,0){\line(0,1){20}}
\put(93,0){\line(0,1){20}}
\put(47,20){\arc{20}{3.142}{4.712}}
\put(83,20){\arc{20}{4.712}{0}}
\put(47,15){\arc{10}{1.571}{3.142}}
\put(47,20){\arc{10}{3.142}{4.712}}
\put(83,20){\arc{10}{4.712}{0}}
\put(83,15){\arc{10}{0}{1.571}}
\put(77,15){\arc{10}{1.571}{3.142}}
\put(77,20){\arc{10}{3.142}{4.712}}
\put(53,20){\arc{10}{4.712}{0}}
\put(53,15){\arc{10}{0}{1.571}}
\put(47,50){\arc{20}{1.571}{3.142}}
\put(83,50){\arc{20}{0}{1.571}}
\put(47,50){\arc{10}{1.571}{3.142}}
\put(47,55){\arc{10}{3.142}{4.712}}
\put(83,55){\arc{10}{4.712}{0}}
\put(83,50){\arc{10}{0}{1.571}}
\put(77,50){\arc{10}{1.571}{3.142}}
\put(77,55){\arc{10}{3.142}{4.712}}
\put(53,55){\arc{10}{4.712}{0}}
\put(53,50){\arc{10}{0}{1.571}}
\put(47,10){\line(1,0){6}}
\put(77,10){\line(1,0){6}}
\put(47,25){\line(1,0){1}}
\put(52,25){\line(1,0){1}}
\put(77,25){\line(1,0){6}}
\put(47,30){\line(1,0){1}}
\put(52,30){\line(1,0){31}}
\put(37,70){\line(0,-1){20}}
\put(93,70){\line(0,-1){20}}
\put(47,40){\line(1,0){1}}
\put(52,40){\line(1,0){31}}
\put(47,45){\line(1,0){1}}
\put(52,45){\line(1,0){1}}
\put(77,45){\line(1,0){6}}
\put(47,60){\line(1,0){6}}
\put(77,60){\line(1,0){6}}
\put(42,15){\line(0,1){5}}
\put(42,55){\line(0,-1){5}}
\put(58,15){\line(0,1){5}}
\put(58,55){\line(0,-1){5}}
\put(88,15){\line(0,1){5}}
\put(88,55){\line(0,-1){5}}
\put(72,15){\line(0,1){5}}
\put(72,55){\line(0,-1){5}}
\thinlines
\put(50,10){\line(0,1){15}}
\put(50,30){\line(0,1){10}}
\put(80,10){\line(0,1){13}}
\put(80,32){\line(0,1){6}}
\put(50,33){\vector(0,-1){0}}
\put(80,37){\vector(0,1){0}}
\put(50,60){\line(0,-1){15}}
\put(80,60){\line(0,-1){13}}
\Thicklines
\put(58,17.5){\line(1,0){14}}
\put(58,52.5){\line(1,0){14}}
\put(67,17.5){\vector(1,0){0}}
\put(63,52.5){\vector(-1,0){0}}
\put(80,27){\line(0,1){1}}
\put(50,25){\line(0,1){5}}
\put(80,43){\line(0,-1){1}}
\put(50,45){\line(0,-1){5}}
\put(33,4){\makebox(0,0){$d$}}
\put(57,7){\makebox(0,0){$c_1$}}
\put(42,8){\makebox(0,0){$b_1$}}
\put(73,7){\makebox(0,0){$c_3$}}
\put(88,8){\makebox(0,0){$b_3$}}
\put(33,66){\makebox(0,0){$a$}}
\put(57,62){\makebox(0,0){$c_2$}}
\put(42,62){\makebox(0,0){$b_2$}}
\put(73,62){\makebox(0,0){$c_4$}}
\put(88,62){\makebox(0,0){$b_4$}}
\put(45,35){\makebox(0,0){$\la$}}
\put(85,34){\makebox(0,0){$\mu$}}
\put(65,22.5){\makebox(0,0){$\tau$}}
\put(65,47.5){\makebox(0,0){$\tau$}}
\put(50,4){\makebox(0,0){$t_1$}}
\put(50,66){\makebox(0,0){$t_2^*$}}
\put(80,4){\makebox(0,0){$t_3^*$}}
\put(80,66){\makebox(0,0){$t_4$}}
\put(54.5,17.5){\makebox(0,0){\footnotesize{$X_1^*$}}}
\put(54.5,52.5){\makebox(0,0){\footnotesize{$X_2$}}}
\put(76.5,17.5){\makebox(0,0){\footnotesize{$X_3$}}}
\put(76.5,52.5){\makebox(0,0){\footnotesize{$X_4^*$}}}
\end{picture}
\end{center}
\caption{The element $|\omega^{\tau,\la,+}_{b_1,c_1,t_1,X_1}\rangle
\langle\omega^{\tau,\la,+}_{b_2,c_2,t_2,X_2}| \otimes
|\omega^{\tau,\mu,-}_{b_3,c_3,t_3,X_3}\rangle
\langle\omega^{\tau,\mu,-}_{b_4,c_4,t_4,X_4}|\in\protect\dta$}
\label{omomomom}
\end{figure}
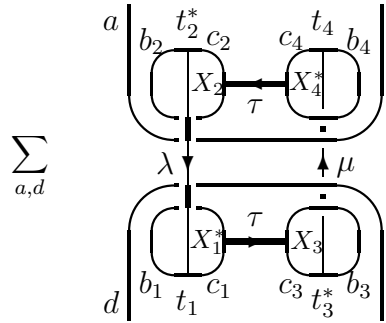

For elements
$|\omega^{\tau,\la,+}_{b_1,c_1,t_1,X_1}\rangle
\langle\omega^{\tau,\la,+}_{b_2,c_2,t_2,X_2}|\in A_{\tau,\la}^+$
and
$|\omega^{\tau,\mu,-}_{b_3,c_3,t_3,X_3}\rangle
\langle\omega^{\tau,\mu,-}_{b_4,c_4,t_4,X_4}|\in A_{\tau,\mu}^-$
we define an element
$|\omega^{\tau,\la,+}_{b_1,c_1,t_1,X_1}\rangle
\langle\omega^{\tau,\la,+}_{b_2,c_2,t_2,X_2}| \otimes
|\omega^{\tau,\mu,-}_{b_3,c_3,t_3,X_3}\rangle
\langle\omega^{\tau,\mu,-}_{b_4,c_4,t_4,X_4}|$
in the double triangle algebra
by the diagram in Fig.\ \ref{omomomom}.
%
Then, for elements
\[ |\varphi^{\tau,\la,+}_1\rangle
\langle\varphi^{\tau,\la,+}_2|
=\sum_{\xi,\xi'} c^\xi_{1,+}(c^{\xi'}_{2,+})^*
|\omega^{\tau,\la,+}_\xi\rangle
\langle\omega^{\tau,\la,+}_{\xi'}| \in A_{\tau,\la}^+ \]
and
\[ |\varphi^{\tau,\mu,-}_3\rangle
\langle\varphi^{\tau,\mu,-}_4|
=\sum_{\xi,\xi'} c^\xi_{3,-}(c^{\xi'}_{4,-})^*
|\omega^{\tau,\mu,-}_\xi\rangle
\langle\omega^{\tau,\mu,-}_{\xi'}| \in A_{\tau,\mu}^- \]
we define an element
$|\varphi^{\tau,\la,+}_1\rangle\langle\varphi^{\tau,\la,+}_2|
\otimes |\varphi^{\tau,\mu,-}_3\rangle
\langle\varphi^{\tau,\mu,-}_4|\in\dta$ by putting
\begin{equation}
\bearl
|\varphi^{\tau,\la,+}_1\rangle\langle\varphi^{\tau,\la,+}_2|
\otimes |\varphi^{\tau,\mu,-}_3\rangle
\langle\varphi^{\tau,\mu,-}_4|=\\[.4em]
\qquad\quad \displaystyle\sum_{\xi,\xi',\xi'',\xi'''}
c^\xi_{1,+}(c^{\xi'}_{2,+})^*
c^{\xi''}_{3,-}(c^{\xi'''}_{4,-})^*
|\omega^{\tau,\la,+}_\xi\rangle
\langle\omega^{\tau,\la,+}_{\xi'}| \otimes
|\omega^{\tau,\mu,-}_{\xi''}\rangle
\langle\omega^{\tau,\mu,-}_{\xi'''}| \,.
\eear
\label{phiphiphiphi}
\end{equation}

\begin{lemma}
\label{linexpppp}
\erf{phiphiphiphi} extends to a bi-linear map
$A_{\tau,\la}^+ \times A_{\tau,\mu}^- \rightarrow \cZ_h$.
\end{lemma}

\begin{proof}
Let $\Phi^+_{a,b}$ denote the
$\Hom(a\tau\co a,b\tau\co b)$ part of
$|\varphi^{\tau,\la,+}_1\rangle\langle\varphi^{\tau,\la,+}_2|$
and similarly $\Phi^-_{a,b}$ the
$\Hom(a\co\tau\co a,b\co\tau\co b)$ part of
$|\varphi^{\tau,\mu,-}_3\rangle\langle\varphi^{\tau,\mu,-}_4|$.
Then the $\Hom(a\co a,b\co b)$ part $\Phi_{a,b}$ of
\[ \Phi=|\varphi^{\tau,\la,+}_1\rangle
\langle\varphi^{\tau,\la,+}_2|
\otimes |\varphi^{\tau,\mu,-}_3\rangle
\langle\varphi^{\tau,\mu,-}_4| \in\dta \]
can be written as
\[ \Phi_{a,b} = d_\tau \sqrt{d_a d_b} \, 
b({\co R}_\tau)^* b \tau (R_b)^*
\Phi^+_{a,b} a\tau \co a (\Phi^-_{a,b}) a \tau(R_a)
a({\co R}_\tau) \,. \]
Thus each component of $\Phi$ is obviously linear
in the components of the vectors in $A_{\tau,\la}^+$
and $A_{\tau,\mu}^-$, proving bi-linearity.
It remains to be shown that $\Phi$ it is
in $\cZ_h$. But this is clear since any element of the form
given in \cite[Fig.\ 33]{BEK1} can be horizontally
``pulled through''. As such elements span the whole
double triangle algebra, the claim is proven.
\end{proof}

We need another graphical identity which refines
\cite[Lemma 6.2]{BEK1}.

\begin{lemma}
\label{keylem1}
We have the identity in Fig.\ \ref{keyfig1} for intertwiners in
$\Hom(\la'\co{\mu'},\la\co\mu)$.
\end{lemma}
%
\begin{figure}[htb]
\begin{center}
\unitlength 0.6mm
\begin{picture}(236,70)
\thicklines
\put(8,34){\makebox(0,0){$\displaystyle\sum_a d_a$}}
\put(27,35){\line(0,1){15}}
\put(83,35){\line(0,1){15}}
\put(37,50){\arc{20}{3.142}{4.712}}
\put(73,50){\arc{20}{4.712}{0}}
\put(37,45){\arc{10}{1.571}{3.142}}
\put(37,50){\arc{10}{3.142}{4.712}}
\put(73,50){\arc{10}{4.712}{0}}
\put(73,45){\arc{10}{0}{1.571}}
\put(37,20){\arc{20}{1.571}{3.142}}
\put(73,20){\arc{20}{0}{1.571}}
\put(37,20){\arc{10}{1.571}{3.142}}
\put(37,25){\arc{10}{3.142}{4.712}}
\put(73,25){\arc{10}{4.712}{0}}
\put(73,20){\arc{10}{0}{1.571}}
\put(67,20){\arc{10}{1.571}{3.142}}
\put(67,25){\arc{10}{3.142}{4.712}}
\put(43,25){\arc{10}{4.712}{0}}
\put(43,20){\arc{10}{0}{1.571}}
\put(37,40){\line(1,0){36}}
\put(37,55){\line(1,0){1}}
\put(42,55){\line(1,0){31}}
\put(37,60){\line(1,0){1}}
\put(42,60){\line(1,0){31}}
\put(27,35){\line(0,-1){15}}
\put(83,35){\line(0,-1){15}}
\put(37,10){\line(1,0){1}}
\put(42,10){\line(1,0){31}}
\put(37,15){\line(1,0){1}}
\put(42,15){\line(1,0){1}}
\put(67,15){\line(1,0){6}}
\put(37,30){\line(1,0){6}}
\put(67,30){\line(1,0){6}}
\put(32,45){\line(0,1){5}}
\put(32,25){\line(0,-1){5}}
\put(48,25){\line(0,-1){5}}
\put(78,45){\line(0,1){5}}
\put(78,25){\line(0,-1){5}}
\put(62,25){\line(0,-1){5}}
\thinlines
\put(40,0){\line(0,1){10}}
\put(40,60){\line(0,1){10}}
\put(40,40){\line(0,1){15}}
\put(70,40){\line(0,1){13}}
\put(70,0){\line(0,1){8}}
\put(70,62){\line(0,1){8}}
\put(40,3){\vector(0,-1){0}}
\put(40,63){\vector(0,-1){0}}
\put(70,7){\vector(0,1){0}}
\put(70,67){\vector(0,1){0}}
\put(40,30){\line(0,-1){15}}
\put(70,30){\line(0,-1){13}}
\Thicklines
\put(48,22.5){\line(1,0){14}}
\put(57,22.5){\vector(1,0){0}}
\put(70,57){\line(0,1){1}}
\put(40,55){\line(0,1){5}}
\put(70,13){\line(0,-1){1}}
\put(40,15){\line(0,-1){5}}
\put(55,36){\makebox(0,0){$c'$}}
\put(55,50){\makebox(0,0){$b'$}}
\put(23,20){\makebox(0,0){$a$}}
\put(48,30){\makebox(0,0){$c$}}
\put(36,22.5){\makebox(0,0){$b$}}
\put(62,30){\makebox(0,0){$c$}}
\put(74,22.5){\makebox(0,0){$b$}}
\put(37,34){\makebox(0,0){$t^*$}}
\put(73,34){\makebox(0,0){$s^*$}}
\put(43,36){\makebox(0,0){$t'$}}
\put(67,36){\makebox(0,0){$s'$}}
\put(44.5,22.5){\makebox(0,0){\footnotesize{$X$}}}
\put(66.5,22.5){\makebox(0,0){\footnotesize{$Y^*$}}}
\put(35,5){\makebox(0,0){$\la$}}
\put(75,4){\makebox(0,0){$\mu$}}
\put(35,65){\makebox(0,0){$\la'$}}
\put(75,65){\makebox(0,0){$\mu'$}}
\put(55,17.5){\makebox(0,0){$\tau$}}
\thicklines
\put(127,34){\makebox(0,0){$=\,\del \la{\la'} \del \mu{\mu'}
\displaystyle\frac w{d_\la d_\mu}$}}
\put(167,15){\arc{10}{1.571}{3.142}}
\put(167,20){\arc{10}{3.142}{4.712}}
\put(203,20){\arc{10}{4.712}{0}}
\put(203,15){\arc{10}{0}{1.571}}
\put(167,50){\arc{10}{1.571}{3.142}}
\put(167,55){\arc{10}{3.142}{4.712}}
\put(203,55){\arc{10}{4.712}{0}}
\put(203,50){\arc{10}{0}{1.571}}
\put(197,50){\arc{10}{1.571}{3.142}}
\put(197,55){\arc{10}{3.142}{4.712}}
\put(173,55){\arc{10}{4.712}{0}}
\put(173,50){\arc{10}{0}{1.571}}
\put(167,10){\line(1,0){36}}
\put(167,25){\line(1,0){1}}
\put(172,25){\line(1,0){31}}
\put(167,45){\line(1,0){1}}
\put(172,45){\line(1,0){1}}
\put(197,45){\line(1,0){6}}
\put(167,60){\line(1,0){6}}
\put(197,60){\line(1,0){6}}
\put(162,15){\line(0,1){5}}
\put(162,55){\line(0,-1){5}}
\put(178,55){\line(0,-1){5}}
\put(208,15){\line(0,1){5}}
\put(208,55){\line(0,-1){5}}
\put(192,55){\line(0,-1){5}}
\thinlines
\put(170,10){\line(0,1){15}}
\put(200,10){\line(0,1){13}}
\put(170,60){\line(0,-1){15}}
\put(200,60){\line(0,-1){13}}
\put(220,0){\line(0,1){70}}
\put(230,0){\line(0,1){70}}
\put(220,33){\vector(0,-1){0}}
\put(230,37){\vector(0,1){0}}
\Thicklines
\put(178,52.5){\line(1,0){14}}
\put(183,52.5){\vector(-1,0){0}}
\put(200,27){\line(0,1){16}}
\put(170,25){\line(0,1){20}}
\put(170,33){\vector(0,-1){0}}
\put(200,37){\vector(0,1){0}}
\put(177,35){\makebox(0,0){$\a^+_\la$}}
\put(193,34){\makebox(0,0){$\a^-_\mu$}}
\put(216,5){\makebox(0,0){$\la$}}
\put(234,5){\makebox(0,0){$\mu$}}
\put(185,5){\makebox(0,0){$c'$}}
\put(185,20){\makebox(0,0){$b'$}}
\put(178,60){\makebox(0,0){$c$}}
\put(166,52.5){\makebox(0,0){$b$}}
\put(192,60){\makebox(0,0){$c$}}
\put(204,52.5){\makebox(0,0){$b$}}
\put(170,64){\makebox(0,0){$t^*$}}
\put(200,64){\makebox(0,0){$s^*$}}
\put(170,6){\makebox(0,0){$t'$}}
\put(200,6){\makebox(0,0){$s'$}}
\put(174.5,52.5){\makebox(0,0){\footnotesize{$X$}}}
\put(196.5,52.5){\makebox(0,0){\footnotesize{$Y^*$}}}
\put(185,47.5){\makebox(0,0){$\tau$}}
\end{picture}
\end{center}
\caption{An identity in $\Hom(\la'\co{\mu'},\la\co\mu)$}
\label{keyfig1}
\end{figure}

\begin{proof}
Using the expansion of the identity (cf.\ \cite[Lemma 4.3]{BEK1})
for the parallel wires $a$ and $b'$ on the top yields the left
hand side of Fig.\ \ref{proofkey1}.
%
\begin{figure}[htb]
\begin{center}
\unitlength 0.6mm
\begin{picture}(208,60)
\thicklines
\put(8,29){\makebox(0,0){$\displaystyle\sum_{a,\nu} d_a$}}
\put(32,35){\arc{10}{3.142}{4.712}}
\put(78,35){\arc{10}{4.712}{0}}
\put(37,20){\arc{20}{1.571}{3.142}}
\put(73,20){\arc{20}{0}{1.571}}
\put(37,20){\arc{10}{1.571}{3.142}}
\put(37,25){\arc{10}{3.142}{4.712}}
\put(73,25){\arc{10}{4.712}{0}}
\put(73,20){\arc{10}{0}{1.571}}
\put(67,20){\arc{10}{1.571}{3.142}}
\put(67,25){\arc{10}{3.142}{4.712}}
\put(43,25){\arc{10}{4.712}{0}}
\put(43,20){\arc{10}{0}{1.571}}
\put(32,40){\line(1,0){46}}
\put(27,35){\line(0,-1){15}}
\put(83,35){\line(0,-1){15}}
\put(37,10){\line(1,0){1}}
\put(42,10){\line(1,0){31}}
\put(37,15){\line(1,0){1}}
\put(42,15){\line(1,0){1}}
\put(67,15){\line(1,0){6}}
\put(37,30){\line(1,0){6}}
\put(67,30){\line(1,0){6}}
\put(32,25){\line(0,-1){5}}
\put(48,25){\line(0,-1){5}}
\put(78,25){\line(0,-1){5}}
\put(62,25){\line(0,-1){5}}
\thinlines
\put(40,0){\line(0,1){10}}
\put(46,40){\line(0,1){20}}
\put(64,40){\line(0,1){8}}
\put(70,0){\line(0,1){8}}
\put(64,52){\line(0,1){8}}
\put(40,3){\vector(0,-1){0}}
\put(46,53){\vector(0,-1){0}}
\put(70,7){\vector(0,1){0}}
\put(64,57){\vector(0,1){0}}
\put(40,30){\line(0,-1){15}}
\put(70,30){\line(0,-1){13}}
\put(34.5,40){\line(0,1){5}}
\put(75.5,40){\line(0,1){5}}
\put(39.5,45){\arc{10}{3.142}{4.712}}
\put(70.5,45){\arc{10}{4.712}{0}}
\put(39.5,50){\line(1,0){4.5}}
\put(48,50){\line(1,0){22.5}}
\put(57,50){\vector(1,0){0}}
\put(34.5,40){\arc{5}{3.142}{0}}
\put(75.5,40){\arc{5}{3.142}{0}}
\Thicklines
\put(48,22.5){\line(1,0){14}}
\put(57,22.5){\vector(1,0){0}}
\put(70,13){\line(0,-1){1}}
\put(40,15){\line(0,-1){5}}
\put(55,44){\makebox(0,0){$c'$}}
\put(41,44){\makebox(0,0){$b'$}}
\put(69,44){\makebox(0,0){$b'$}}
\put(23,20){\makebox(0,0){$a$}}
\put(48,30){\makebox(0,0){$c$}}
\put(36,22.5){\makebox(0,0){$b$}}
\put(62,30){\makebox(0,0){$c$}}
\put(74,22.5){\makebox(0,0){$b$}}
\put(40,34){\makebox(0,0){$t^*$}}
\put(70,34){\makebox(0,0){$s^*$}}
\put(46,36){\makebox(0,0){$t'$}}
\put(64,36){\makebox(0,0){$s'$}}
\put(44.5,22.5){\makebox(0,0){\footnotesize{$X$}}}
\put(66.5,22.5){\makebox(0,0){\footnotesize{$Y^*$}}}
\put(35,5){\makebox(0,0){$\la$}}
\put(75,4){\makebox(0,0){$\mu$}}
\put(41,56){\makebox(0,0){$\la'$}}
\put(69,56){\makebox(0,0){$\mu'$}}
\put(55,55){\makebox(0,0){$\nu$}}
\put(55,17.5){\makebox(0,0){$\tau$}}
\thicklines
\put(119,29){\makebox(0,0){$=\;\displaystyle\sum_\nu \,d_\nu$}}
\put(157,35){\arc{20}{3.142}{4.712}}
\put(193,35){\arc{20}{4.712}{0}}
\put(157,25){\arc{20}{1.571}{3.142}}
\put(193,25){\arc{20}{0}{1.571}}
\put(157,25){\arc{10}{1.571}{3.142}}
\put(157,30){\arc{10}{3.142}{4.712}}
\put(193,30){\arc{10}{4.712}{0}}
\put(193,25){\arc{10}{0}{1.571}}
\put(187,25){\arc{10}{1.571}{3.142}}
\put(187,30){\arc{10}{3.142}{4.712}}
\put(163,30){\arc{10}{4.712}{0}}
\put(163,25){\arc{10}{0}{1.571}}
\put(157,45){\line(1,0){36}}
\put(147,25){\line(0,1){10}}
\put(203,25){\line(0,1){10}}
\put(157,15){\line(1,0){1}}
\put(162,15){\line(1,0){31}}
\put(157,20){\line(1,0){1}}
\put(162,20){\line(1,0){1}}
\put(187,20){\line(1,0){6}}
\put(157,35){\line(1,0){6}}
\put(187,35){\line(1,0){6}}
\put(152,30){\line(0,-1){5}}
\put(168,30){\line(0,-1){5}}
\put(198,30){\line(0,-1){5}}
\put(182,30){\line(0,-1){5}}
\thinlines
\put(160,0){\line(0,1){15}}
\put(160,45){\line(0,1){15}}
\put(190,45){\line(0,1){3}}
\put(190,0){\line(0,1){8}}
\put(190,12){\line(0,1){3}}
\put(190,52){\line(0,1){8}}
\put(160,3){\vector(0,-1){0}}
\put(160,53){\vector(0,-1){0}}
\put(190,7){\vector(0,1){0}}
\put(190,57){\vector(0,1){0}}
\put(160,35){\line(0,-1){15}}
\put(190,35){\line(0,-1){13}}
\put(142,25){\line(0,1){10}}
\put(208,25){\line(0,1){10}}
\put(157,35){\arc{30}{3.142}{4.712}}
\put(193,35){\arc{30}{4.712}{0}}
\put(157,50){\line(1,0){1}}
\put(162,50){\line(1,0){31}}
\put(177,50){\vector(1,0){0}}
\put(157,25){\arc{30}{1.571}{3.142}}
\put(193,25){\arc{30}{0}{1.571}}
\put(157,10){\line(1,0){1}}
\put(162,10){\line(1,0){31}}
\Thicklines
\put(168,27.5){\line(1,0){14}}
\put(177,27.5){\vector(1,0){0}}
\put(190,18){\line(0,-1){1}}
\put(160,20){\line(0,-1){5}}
\put(175,41){\makebox(0,0){$c'$}}
\put(175,19){\makebox(0,0){$b'$}}
\put(168,35){\makebox(0,0){$c$}}
\put(156,27.5){\makebox(0,0){$b$}}
\put(182,35){\makebox(0,0){$c$}}
\put(194,27.5){\makebox(0,0){$b$}}
\put(157,39){\makebox(0,0){$t^*$}}
\put(193,39){\makebox(0,0){$s^*$}}
\put(163,41){\makebox(0,0){$t'$}}
\put(187,41){\makebox(0,0){$s'$}}
\put(164.5,27.5){\makebox(0,0){\footnotesize{$X$}}}
\put(186.5,27.5){\makebox(0,0){\footnotesize{$Y^*$}}}
\put(155,5){\makebox(0,0){$\la$}}
\put(195,4){\makebox(0,0){$\mu$}}
\put(155,56){\makebox(0,0){$\la'$}}
\put(195,56){\makebox(0,0){$\mu'$}}
\put(175,55){\makebox(0,0){$\nu$}}
\put(175,31.5){\makebox(0,0){$\tau$}}
\end{picture}
\end{center}
\caption{The identity in $\Hom(\la'\co{\mu'},\la\co\mu)$}
\label{proofkey1}
\end{figure}
We then slide around the trivalent vertices of the wire $\nu$
along the wire $a$ so that they almost meet at the bottom of
the picture. Turning around their small arcs yields a factor
$d_\nu/d_a$, and we can then see that the summation over
$a$ is just the expansion of the identity
(cf.\ \cite[Lemma 4.3]{BEK1}) which gives us to parallel wires
$b'$ and $\nu$. This way we arrive at the right hand side of
Fig.\ \ref{proofkey1}. Then we apply the expansion of the
identity four times: First twice for the parallel wires
$b$ and $b'$ on the bottom, yielding expansions over $\rho$
and $\rho'$. Next we expand the parallel wires $\tau$ and
$b'$ in the middle lower part of the picture, resulting in
a summation over a wire $a'$. Finally we expand the parallel
wires $c'$ and $a'$ in the center of the picture, yielding
a summation over a wire $\rho''$. This gives us
Fig.\ \ref{proofkey2}.
%
\begin{figure}[htb]
\begin{center}
\unitlength 0.6mm
\begin{picture}(194,70)
\thicklines
\put(20,33){\makebox(0,0){$\displaystyle
\sum_{\nu,a',\rho,\rho',\rho''}\; d_\nu$}}
\put(55,20){\line(0,1){30}}
\put(60,50){\arc{10}{3.142}{4.712}}
\put(60,55){\line(1,0){45}}
\put(105,50){\arc{10}{4.712}{0}}
\put(110,30){\line(0,1){20}}
\put(105,30){\arc{10}{0}{1.571}}
\put(100,25){\line(1,0){5}}
\put(100,20){\line(0,1){5}}
\put(95,20){\arc{10}{0}{1.571}}
\put(90,15){\line(1,0){5}}
\put(90,20){\arc{10}{1.571}{3.142}}
\put(85,20){\line(0,1){20}}
\put(80,40){\arc{10}{4.712}{0}}
\put(70,45){\line(1,0){10}}
\put(70,40){\arc{10}{3.142}{4.712}}
\put(65,20){\line(0,1){20}}
\put(60,20){\arc{10}{0}{3.142}}
\put(185,20){\line(0,1){30}}
\put(180,50){\arc{10}{4.712}{0}}
\put(135,55){\line(1,0){45}}
\put(135,50){\arc{10}{3.142}{4.712}}
\put(130,30){\line(0,1){20}}
\put(135,30){\arc{10}{1.571}{3.142}}
\put(135,25){\line(1,0){5}}
\put(140,20){\line(0,1){5}}
\put(145,20){\arc{10}{1.571}{3.142}}
\put(145,15){\line(1,0){5}}
\put(150,20){\arc{10}{0}{1.571}}
\put(155,20){\line(0,1){20}}
\put(160,40){\arc{10}{3.142}{4.712}}
\put(160,45){\line(1,0){10}}
\put(170,40){\arc{10}{4.712}{0}}
\put(175,20){\line(0,1){20}}
\put(180,20){\arc{10}{0}{3.142}}
\thinlines
\put(75,0){\line(0,1){45}}
\put(75,55){\line(0,1){15}}
\put(165,0){\line(0,1){8}}
\put(165,12){\line(0,1){8.5}}
\put(165,24.5){\line(0,1){20.5}}
\put(165,55){\line(0,1){3}}
\put(165,62){\line(0,1){8}}
\put(75,3){\vector(0,-1){0}}
\put(75,63){\vector(0,-1){0}}
\put(165,7){\vector(0,1){0}}
\put(165,67){\vector(0,1){0}}
\put(50,20){\line(0,1){30}}
\put(190,20){\line(0,1){30}}
\put(60,50){\arc{20}{3.142}{4.712}}
\put(180,50){\arc{20}{4.712}{0}}
\put(60,60){\line(1,0){13}}
\put(77,60){\line(1,0){103}}
\put(60,10){\line(1,0){13}}
\put(77,10){\line(1,0){103}}
\put(122,60){\vector(1,0){0}}
\put(60,20){\arc{20}{1.571}{3.142}}
\put(180,20){\arc{20}{0}{1.571}}
\put(65,22.5){\line(1,0){8}}
\put(85,22.5){\line(-1,0){8}}
\put(72,22.5){\vector(1,0){0}}
\put(155,22.5){\line(1,0){20}}
\put(172,22.5){\vector(1,0){0}}
\put(110,40){\line(1,0){20}}
\put(122,40){\vector(1,0){0}}
\put(110,40){\arc{5}{4.712}{1.571}}
\put(130,40){\arc{5}{1.571}{4.712}}
\put(100,25){\arc{5}{4.712}{1.571}}
\put(140,25){\arc{5}{1.571}{4.712}}
\put(65,22.5){\arc{5}{4.712}{1.571}}
\put(85,22.5){\arc{5}{1.571}{4.712}}
\put(155,22.5){\arc{5}{4.712}{1.571}}
\put(175,22.5){\arc{5}{1.571}{4.712}}
\Thicklines
\put(100,25){\line(0,1){5}}
\put(95,30){\arc{10}{4.712}{0}}
\put(85,35){\line(1,0){10}}
\put(92,35){\vector(1,0){0}}
\put(140,25){\line(0,1){5}}
\put(145,30){\arc{10}{3.142}{4.712}}
\put(145,35){\line(1,0){10}}
\put(152,35){\vector(1,0){0}}
\put(105,50){\makebox(0,0){$c'$}}
\put(135,50){\makebox(0,0){$c'$}}
\put(105,32){\makebox(0,0){$a'$}}
\put(135,32){\makebox(0,0){$a'$}}
\put(85,45){\makebox(0,0){$c$}}
\put(61,30){\makebox(0,0){$b$}}
\put(155,45){\makebox(0,0){$c$}}
\put(179,30){\makebox(0,0){$b$}}
\put(89,30){\makebox(0,0){$b$}}
\put(151,30){\makebox(0,0){$b$}}
\put(92.5,19){\makebox(0,0){$b'$}}
\put(147.5,19){\makebox(0,0){$b'$}}
\put(60,19){\makebox(0,0){$b'$}}
\put(180,19){\makebox(0,0){$b'$}}
\put(72,49){\makebox(0,0){$t^*$}}
\put(168,49){\makebox(0,0){$s^*$}}
\put(78,51){\makebox(0,0){$t'$}}
\put(162,51){\makebox(0,0){$s'$}}
\put(80,35){\makebox(0,0){\footnotesize{$X$}}}
\put(160,35){\makebox(0,0){\footnotesize{$Y^*$}}}
\put(70,5){\makebox(0,0){$\la$}}
\put(170,4){\makebox(0,0){$\mu$}}
\put(70,66){\makebox(0,0){$\la'$}}
\put(170,66){\makebox(0,0){$\mu'$}}
\put(120,65){\makebox(0,0){$\nu$}}
\put(70,17.5){\makebox(0,0){$\rho$}}
\put(170,17.5){\makebox(0,0){$\rho'$}}
\put(120,35){\makebox(0,0){$\rho''$}}
\put(90,40){\makebox(0,0){$\tau$}}
\put(150,40){\makebox(0,0){$\tau$}}
\end{picture}
\end{center}
\caption{The identity in $\Hom(\la'\co{\mu'},\la\co\mu)$}
\label{proofkey2}
\end{figure}
Now we can pull the circle $\nu$ around the middle expansion
$\rho''$, just by virtue of the IBFE moves as well as the
Yang-Baxter relation for thin wires. Due to the prefactor
$d_\nu$, the summation over $\nu$ yields exactly the
orthogonality relation for a non-degenerate braiding
(cf.\ \cite[Fig.\ 20]{BEK1}), the ``killing ring''.
Therefore we obtain zero unless $\rho''=\id$, and our
picture becomes disconnected yielding two intertwiners
in $\Hom(\la',\la)$ and $\Hom(\co{\mu'},\co\mu)$.
Hence we obtain a factor $\del\la{\la'}\del\mu{\mu'}$,
and the whole diagram represents a scalar. To compute the
scalar, we can proceed exactly as in the proof of
\cite[Lemma 6.2]{BEK1}: We go back to the original picture
on the left hand side of Fig.\ \ref{keyfig1} and put now
$\la'=\la$ and $\mu'=\mu$. Then we close the wires
$\la$ and $\mu$ on the right which has to be compensated
by a factor $d_\la^{-1}d_\mu^{-1}$. Next we open the
wire $a$ on the left and close it also on the right.
Then the $a$ loop can be pulled out and the summation over
$a$ gives the global index $w$; we are left with the
right hand side of Fig.\ \ref{keyfig1}.
\end{proof}

Recall from \cite[Thm.\ 6.8]{BEK1} that
$\sum_{\la,\mu} q_{\la,\mu}=e_0$. The $\Hom(a\co a,d\co d)$
part of this relation gives us the graphical identity of
Fig.\ \ref{sumqlamu}.
%
\thinlines
\begin{figure}[htb]
\begin{center}
\unitlength 0.6mm
\begin{picture}(156,50)
\thicklines
\put(15,25){\makebox(0,0){$\displaystyle\sum_{b,c,\la,\mu}
\frac{d_b d_c}{w^2}$}}
\put(37,0){\line(0,1){10}}
\put(83,0){\line(0,1){10}}
\put(47,10){\arc{20}{3.142}{4.712}}
\put(73,10){\arc{20}{4.712}{0}}
\put(47,10){\arc{10}{1.571}{4.712}}
\put(73,10){\arc{10}{4.712}{1.571}}
\put(47,40){\arc{20}{1.571}{3.142}}
\put(73,40){\arc{20}{0}{1.571}}
\put(47,40){\arc{10}{1.571}{4.712}}
\put(73,40){\arc{10}{4.712}{1.571}}
\put(47,5){\line(1,0){26}}
\put(47,15){\line(1,0){1}}
\put(52,15){\line(1,0){21}}
\put(47,20){\line(1,0){1}}
\put(52,20){\line(1,0){21}}
\put(37,50){\line(0,-1){10}}
\put(83,50){\line(0,-1){10}}
\put(47,30){\line(1,0){1}}
\put(52,30){\line(1,0){21}}
\put(47,35){\line(1,0){1}}
\put(52,35){\line(1,0){21}}
\put(47,45){\line(1,0){26}}
\thinlines
\put(50,5){\line(0,1){10}}
\put(50,20){\line(0,1){10}}
\put(70,5){\line(0,1){8}}
\put(70,22){\line(0,1){6}}
\put(50,23){\vector(0,-1){0}}
\put(70,27){\vector(0,1){0}}
\put(50,45){\line(0,-1){10}}
\put(70,45){\line(0,-1){8}}
\put(50,5){\arc{5}{3.142}{0}}
\put(70,5){\arc{5}{3.142}{0}}
\put(50,45){\arc{5}{0}{3.142}}
\put(70,45){\arc{5}{0}{3.142}}
\Thicklines
\put(70,17){\line(0,1){1}}
\put(50,15){\line(0,1){5}}
\put(70,33){\line(0,-1){1}}
\put(50,35){\line(0,-1){5}}
\put(33,4){\makebox(0,0){$d$}}
\put(60,2){\makebox(0,0){$c$}}
\put(60,11){\makebox(0,0){$b$}}
\put(33,46){\makebox(0,0){$a$}}
\put(60,48){\makebox(0,0){$c$}}
\put(60,39){\makebox(0,0){$b$}}
\put(45,25){\makebox(0,0){$\la$}}
\put(75,24){\makebox(0,0){$\mu$}}
\thicklines
\put(110,25){\makebox(0,0){$=\,\del ad \displaystyle\frac 1{d_a}$}}
\put(130,0){\line(0,1){50}}
\put(150,0){\line(0,1){50}}
\put(126,4){\makebox(0,0){$a$}}
\put(154,4){\makebox(0,0){$a$}}
\end{picture}
\end{center}
\caption{A graphical relation from $\sum_{\la,\mu}q_{\la,\mu}=e_0$}
\label{sumqlamu}
\end{figure}
Inserting this in the middle of the left hand side of
Fig.\ \ref{keyfig2},
%
\begin{figure}[htb]
\begin{center}
\unitlength 0.6mm
\begin{picture}(236,70)
\thicklines
\put(8,34){\makebox(0,0){$\displaystyle\sum_a d_a$}}
\put(27,35){\line(0,1){15}}
\put(83,35){\line(0,1){15}}
\put(37,50){\arc{20}{3.142}{4.712}}
\put(73,50){\arc{20}{4.712}{0}}
\put(37,45){\arc{10}{1.571}{3.142}}
\put(37,50){\arc{10}{3.142}{4.712}}
\put(73,50){\arc{10}{4.712}{0}}
\put(73,45){\arc{10}{0}{1.571}}
\put(67,45){\arc{10}{1.571}{3.142}}
\put(67,50){\arc{10}{3.142}{4.712}}
\put(43,50){\arc{10}{4.712}{0}}
\put(43,45){\arc{10}{0}{1.571}}
\put(37,20){\arc{20}{1.571}{3.142}}
\put(73,20){\arc{20}{0}{1.571}}
\put(37,20){\arc{10}{1.571}{3.142}}
\put(37,25){\arc{10}{3.142}{4.712}}
\put(73,25){\arc{10}{4.712}{0}}
\put(73,20){\arc{10}{0}{1.571}}
\put(67,20){\arc{10}{1.571}{3.142}}
\put(67,25){\arc{10}{3.142}{4.712}}
\put(43,25){\arc{10}{4.712}{0}}
\put(43,20){\arc{10}{0}{1.571}}
\put(37,40){\line(1,0){6}}
\put(67,40){\line(1,0){6}}
\put(37,55){\line(1,0){1}}
\put(42,55){\line(1,0){1}}
\put(67,55){\line(1,0){6}}
\put(37,60){\line(1,0){1}}
\put(42,60){\line(1,0){31}}
\put(27,35){\line(0,-1){15}}
\put(83,35){\line(0,-1){15}}
\put(37,10){\line(1,0){1}}
\put(42,10){\line(1,0){31}}
\put(37,15){\line(1,0){1}}
\put(42,15){\line(1,0){1}}
\put(67,15){\line(1,0){6}}
\put(37,30){\line(1,0){6}}
\put(67,30){\line(1,0){6}}
\put(32,45){\line(0,1){5}}
\put(32,25){\line(0,-1){5}}
\put(48,45){\line(0,1){5}}
\put(48,25){\line(0,-1){5}}
\put(78,45){\line(0,1){5}}
\put(78,25){\line(0,-1){5}}
\put(62,45){\line(0,1){5}}
\put(62,25){\line(0,-1){5}}
\thinlines
\put(40,0){\line(0,1){10}}
\put(40,60){\line(0,1){10}}
\put(40,40){\line(0,1){15}}
\put(70,40){\line(0,1){13}}
\put(70,0){\line(0,1){8}}
\put(70,62){\line(0,1){8}}
\put(40,3){\vector(0,-1){0}}
\put(40,63){\vector(0,-1){0}}
\put(70,7){\vector(0,1){0}}
\put(70,67){\vector(0,1){0}}
\put(40,30){\line(0,-1){15}}
\put(70,30){\line(0,-1){13}}
\Thicklines
\put(48,22.5){\line(1,0){14}}
\put(48,47.5){\line(1,0){14}}
\put(57,22.5){\vector(1,0){0}}
\put(53,47.5){\vector(-1,0){0}}
\put(70,57){\line(0,1){1}}
\put(40,55){\line(0,1){5}}
\put(70,13){\line(0,-1){1}}
\put(40,15){\line(0,-1){5}}
\put(49,40){\makebox(0,0){$c'$}}
\put(36,47.5){\makebox(0,0){$b'$}}
\put(61,40){\makebox(0,0){$c'$}}
\put(74,47.5){\makebox(0,0){$b'$}}
\put(23,20){\makebox(0,0){$a$}}
\put(48,30){\makebox(0,0){$c$}}
\put(36,22.5){\makebox(0,0){$b$}}
\put(62,30){\makebox(0,0){$c$}}
\put(74,22.5){\makebox(0,0){$b$}}
\put(37,34){\makebox(0,0){$t^*$}}
\put(73,34){\makebox(0,0){$s^*$}}
\put(43,36){\makebox(0,0){$t'$}}
\put(67,36){\makebox(0,0){$s'$}}
\put(44.5,22.5){\makebox(0,0){\footnotesize{$X$}}}
\put(66.5,22.5){\makebox(0,0){\footnotesize{$Y^*$}}}
\put(44,47.5){\makebox(0,0){\footnotesize{${X'}^*$}}}
\put(66.5,47.5){\makebox(0,0){\footnotesize{$Y'$}}}
\put(35,5){\makebox(0,0){$\la$}}
\put(75,4){\makebox(0,0){$\mu$}}
\put(35,65){\makebox(0,0){$\la'$}}
\put(75,65){\makebox(0,0){$\mu'$}}
\put(55,17.5){\makebox(0,0){$\tau$}}
\put(55,52.5){\makebox(0,0){$\tau'$}}
\thicklines
\put(125,34){\makebox(0,0){$=\,\del \la{\la'} \del \mu{\mu'}
\del \tau{\tau'} \displaystyle\frac w{d_\la d_\mu}$}}
\put(167,15){\arc{10}{1.571}{3.142}}
\put(167,20){\arc{10}{3.142}{4.712}}
\put(203,20){\arc{10}{4.712}{0}}
\put(203,15){\arc{10}{0}{1.571}}
\put(197,15){\arc{10}{1.571}{3.142}}
\put(197,20){\arc{10}{3.142}{4.712}}
\put(173,20){\arc{10}{4.712}{0}}
\put(173,15){\arc{10}{0}{1.571}}
\put(167,50){\arc{10}{1.571}{3.142}}
\put(167,55){\arc{10}{3.142}{4.712}}
\put(203,55){\arc{10}{4.712}{0}}
\put(203,50){\arc{10}{0}{1.571}}
\put(197,50){\arc{10}{1.571}{3.142}}
\put(197,55){\arc{10}{3.142}{4.712}}
\put(173,55){\arc{10}{4.712}{0}}
\put(173,50){\arc{10}{0}{1.571}}
\put(167,10){\line(1,0){6}}
\put(197,10){\line(1,0){6}}
\put(167,25){\line(1,0){1}}
\put(172,25){\line(1,0){1}}
\put(197,25){\line(1,0){6}}
\put(167,45){\line(1,0){1}}
\put(172,45){\line(1,0){1}}
\put(197,45){\line(1,0){6}}
\put(167,60){\line(1,0){6}}
\put(197,60){\line(1,0){6}}
\put(162,15){\line(0,1){5}}
\put(162,55){\line(0,-1){5}}
\put(178,15){\line(0,1){5}}
\put(178,55){\line(0,-1){5}}
\put(208,15){\line(0,1){5}}
\put(208,55){\line(0,-1){5}}
\put(192,15){\line(0,1){5}}
\put(192,55){\line(0,-1){5}}
\thinlines
\put(170,10){\line(0,1){15}}
\put(200,10){\line(0,1){13}}
\put(170,60){\line(0,-1){15}}
\put(200,60){\line(0,-1){13}}
\put(220,0){\line(0,1){70}}
\put(230,0){\line(0,1){70}}
\put(220,33){\vector(0,-1){0}}
\put(230,37){\vector(0,1){0}}
\Thicklines
\put(178,17.5){\line(1,0){14}}
\put(178,52.5){\line(1,0){14}}
\put(187,17.5){\vector(1,0){0}}
\put(183,52.5){\vector(-1,0){0}}
\put(200,27){\line(0,1){16}}
\put(170,25){\line(0,1){20}}
\put(170,33){\vector(0,-1){0}}
\put(200,37){\vector(0,1){0}}
\put(177,35){\makebox(0,0){$\a^+_\la$}}
\put(193,34){\makebox(0,0){$\a^-_\mu$}}
\put(216,5){\makebox(0,0){$\la$}}
\put(234,5){\makebox(0,0){$\mu$}}
\put(177,10){\makebox(0,0){$c'$}}
\put(166,17.5){\makebox(0,0){$b'$}}
\put(191,10){\makebox(0,0){$c'$}}
\put(204,17.5){\makebox(0,0){$b'$}}
\put(178,60){\makebox(0,0){$c$}}
\put(166,52.5){\makebox(0,0){$b$}}
\put(192,60){\makebox(0,0){$c$}}
\put(204,52.5){\makebox(0,0){$b$}}
\put(170,64){\makebox(0,0){$t^*$}}
\put(200,64){\makebox(0,0){$s^*$}}
\put(170,6){\makebox(0,0){$t'$}}
\put(200,6){\makebox(0,0){$s'$}}
\put(174.5,52.5){\makebox(0,0){\footnotesize{$X$}}}
\put(196,52.5){\makebox(0,0){\footnotesize{$Y^*$}}}
\put(174,17.5){\makebox(0,0){\footnotesize{${X'}^*$}}}
\put(196.5,17.5){\makebox(0,0){\footnotesize{$Y'$}}}
\put(185,47.5){\makebox(0,0){$\tau$}}
\put(185,22.5){\makebox(0,0){$\tau$}}
\end{picture}
\end{center}
\caption{An identity in $\Hom(\la'\co{\mu'},\la\co\mu)$}
\label{keyfig2}
\end{figure}
we find that this intertwiner is also a scalar which vanishes
unless $\la=\la'$ and $\mu=\mu'$. It can be evaluated in the
same way, therefore we find a factor $\del \tau{\tau'}$
and thus we arrive at

\begin{corollary}
\label{keylem2}
We have the identity in Fig.\ \ref{keyfig2} for intertwiners in
$\Hom(\la'\co{\mu'},\la\co\mu)$.
\end{corollary}
Using now Fig.\ \ref{<om,om>}, we obtain from
Corollary \ref{keylem2} and Lemma \ref{linexpppp}
the following

\begin{corollary}
We have
\begin{equation}
\bearl
|\varphi^{\tau,\la,+}_1\rangle\langle\varphi^{\tau,\la,+}_2|
\otimes |\varphi^{\tau,\mu,-}_3\rangle
\langle\varphi^{\tau,\mu,-}_4|  *_v
|\varphi^{\tau',\la',+}_5\rangle\langle\varphi^{\tau',\la',+}_6|
\otimes |\varphi^{\tau',\mu',-}_7\rangle
\langle\varphi^{\tau',\mu',-}_8|  \\[.4em]
\qquad = \displaystyle
\frac{\del \tau{\tau'} \del \la{\la'} \del \mu{\mu'}}{wd_\tau}
\langle \varphi^{\tau,\la,+}_2,\varphi^{\tau,\la,+}_5\rangle
\langle \varphi^{\tau,\mu,-}_4,\varphi^{\tau,\mu,-}_7\rangle
|\varphi^{\tau,\la,+}_1\rangle\langle\varphi^{\tau,\la,+}_6|
\otimes |\varphi^{\tau,\mu,-}_3\rangle
\langle\varphi^{\tau,\mu,-}_8| \,.
\eear
\end{equation}
Consequently, defining
\begin{equation}
E_{\tau,\la,\mu;i,k}^{j,l} = w d_\tau
|u^{\tau,\la,+}_i\rangle\langle u^{\tau,\la,+}_j|
\otimes |u^{\tau,\mu,-}_k\rangle
\langle u^{\tau,\mu,-}_l|
\end{equation}
gives a system of matrix units
$\{E_{\tau,\la,\mu;i,k}^{j,l} \}_{\tau,\la,\mu,i,j,k,l}$
in $(\cZ_h,*_v)$, i.e.\ we have
\begin{equation}
E_{\tau,\la,\mu;i,k}^{j,l} *_v E_{\tau',\la',\mu';i',k'}^{j',l'}=
\del \tau{\tau'} \del \la{\la'} \del \mu{\mu'} \del j{i'} \del l{k'}
E_{\tau,\la,\mu;i,k}^{j',l'} \,.
\end{equation}
\end{corollary}

We now define {\sl chiral matrix units} by
\begin{equation}
\bearl
E_{\tau,\la;i}^{+;j} = \sum_\mu \sum_{k=1}^{\dim H_{\tau,\mu}^-}
E_{\tau,\la,\mu;i,k}^{j,k} \,,\\[.4em]
E_{\tau,\mu;k}^{-;l} = \sum_\la \sum_{i=1}^{\dim H_{\tau,\la}^+}
E_{\tau,\la,\mu;i,k}^{i,l} \,.
\eear
\end{equation}
Recall that $\cZ_h^\pm\subset\cZ_h$ are the
chiral vertical subalgebras spanned by elements
$e_{\beta_\pm}$ with $\beta_\pm\in\MXMpm$.

\begin{proposition}
We have $E_{\tau,\la;i}^{\pm;j}\in\cZ_h^\pm$.
\end{proposition}

\begin{proof}
We show $E_{\tau,\la;i}^{+;j}\in\cZ_h^+$. The proof of
$E_{\tau,\la;i}^{-;j}\in\cZ_h^-$ is analogous.
It follows from Lemma \ref{Iexpan} that
$E_{\tau,\la;i}^{+;j} = w d_\tau
|u^{\tau,\la,+}_i\rangle\langle u^{\tau,\la,+}_j|
\otimes I_\tau^-$. Therefore it suffices to show that
$|\om^{\tau,\la,+}_{b',c',t',X'}\rangle
\langle \om^{\tau,\la,+}_{b,c,t,X}|
\otimes I_\tau^-\in\cZ_h^+$.
Such an element is given graphically in
Fig.\ \ref{omomI}.
%
\thinlines
\begin{figure}[htb]
\begin{center}
\unitlength 0.6mm
\begin{picture}(99,65)
\thicklines
\put(20,30){\makebox(0,0){$\displaystyle
\sum_{a,b}\sum_{\beta\in\MXMp}\frac 1{w_+}$}}
\put(47,0){\line(0,1){15}}
\put(93,0){\line(0,1){65}}
\put(57,15){\arc{20}{3.142}{4.712}}
\put(65,35){\arc{10}{4.712}{0}}
\put(57,15){\arc{10}{1.571}{4.712}}
\put(83,15){\arc{10}{4.712}{1.571}}
\put(57,50){\arc{20}{1.571}{3.142}}
\put(65,30){\arc{10}{0}{1.571}}
\put(57,50){\arc{10}{1.571}{4.712}}
\put(83,50){\arc{10}{4.712}{1.571}}
\put(57,10){\line(1,0){26}}
\put(57,20){\line(1,0){1}}
\put(62,20){\line(1,0){21}}
\put(57,25){\line(1,0){1}}
\put(62,25){\line(1,0){3}}
\put(47,65){\line(0,-1){15}}
\put(70,35){\line(0,-1){5}}
\put(57,40){\line(1,0){1}}
\put(62,40){\line(1,0){3}}
\put(57,45){\line(1,0){1}}
\put(62,45){\line(1,0){21}}
\put(57,55){\line(1,0){26}}
\thinlines
\put(60,10){\line(0,1){10}}
\put(60,25){\line(0,1){15}}
\put(60,30.5){\vector(0,-1){0}}
\put(60,55){\line(0,-1){10}}
\put(70,32.5){\arc{5}{4.712}{1.571}}
\put(93,32.5){\arc{5}{1.571}{4.712}}
\Thicklines
\put(70,32.5){\line(1,0){23}}
\put(89,32.5){\vector(1,0){0}}
\put(80,20){\line(0,1){10.5}}
\put(60,20){\line(0,1){5}}
\put(80,45){\line(0,-1){10.5}}
\put(60,45){\line(0,-1){5}}
\put(80,41){\vector(0,1){0}}
\put(43,4){\makebox(0,0){$d$}}
\put(97,4){\makebox(0,0){$d$}}
\put(74,6){\makebox(0,0){$c'$}}
\put(66,16){\makebox(0,0){$b'$}}
\put(43,61){\makebox(0,0){$a$}}
\put(97,61){\makebox(0,0){$a$}}
\put(74,58){\makebox(0,0){$c$}}
\put(66,49){\makebox(0,0){$b$}}
\put(55,32.5){\makebox(0,0){$\la$}}
\put(87.5,27.5){\makebox(0,0){$\beta$}}
\put(84,40){\makebox(0,0){$\tau$}}
\put(60,5){\makebox(0,0){$t'$}}
\put(60,60){\makebox(0,0){$t^*$}}
\put(80,15){\makebox(0,0){\footnotesize{$(X')^*$}}}
\put(80,50){\makebox(0,0){\footnotesize{$X$}}}
\end{picture}
\end{center}
\caption{The element $|\om^{\tau,\la,+}_{b',c',t',X'}\rangle
\langle \om^{\tau,\la,+}_{b,c,t,X}|\otimes I_\tau^-$}
\label{omomI}
\end{figure}
If we multiply horizontally with some $e_{\beta'}$ either
from the left or from the right, then the resulting picture
contains a part which corresponds to an intertwiner in
$\Hom(\beta',\beta)$ or $\Hom(\beta,\beta')$, respectively.
Hence this is zero unless $\beta'\in\MXMp$. But
$\cZ_h$ is spanned by elements $e_\beta$,
$\beta\in\MXM$, and $\cZ_h^+$ is the subspace
spanned by those with $\beta\in\MXMp$. As the
$e_\beta$'s are horizontal projections, the claim follows.
\end{proof}

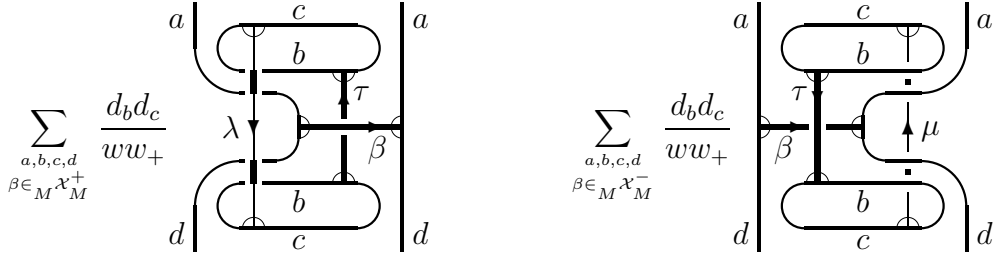
\begin{figure}[htb]
\begin{center}
\unitlength 0.6mm
\begin{picture}(220,55)
\thicklines
\put(18,23){\makebox(0,0){$\displaystyle
\sum_{a,b,c,d\atop\beta\in\MXMp}\frac{d_b d_c}{ww_+}$}}
\put(42,0){\line(0,1){10}}
\put(88,0){\line(0,1){55}}
\put(52,10){\arc{20}{3.142}{4.712}}
\put(60,30){\arc{10}{4.712}{0}}
\put(52,10){\arc{10}{1.571}{4.712}}
\put(78,10){\arc{10}{4.712}{1.571}}
\put(52,45){\arc{20}{1.571}{3.142}}
\put(60,25){\arc{10}{0}{1.571}}
\put(52,45){\arc{10}{1.571}{4.712}}
\put(78,45){\arc{10}{4.712}{1.571}}
\put(52,5){\line(1,0){26}}
\put(52,15){\line(1,0){1}}
\put(57,15){\line(1,0){21}}
\put(52,20){\line(1,0){1}}
\put(57,20){\line(1,0){3}}
\put(42,55){\line(0,-1){10}}
\put(65,30){\line(0,-1){5}}
\put(52,35){\line(1,0){1}}
\put(57,35){\line(1,0){3}}
\put(52,40){\line(1,0){1}}
\put(57,40){\line(1,0){21}}
\put(52,50){\line(1,0){26}}
\thinlines
\put(55,5){\line(0,1){10}}
\put(55,20){\line(0,1){15}}
\put(55,25.5){\vector(0,-1){0}}
\put(55,50){\line(0,-1){10}}
\put(55,5){\arc{5}{3.142}{0}}
\put(75,15){\arc{5}{3.142}{0}}
\put(55,50){\arc{5}{0}{3.142}}
\put(75,40){\arc{5}{0}{3.142}}
\put(65,27.5){\arc{5}{4.712}{1.571}}
\put(88,27.5){\arc{5}{1.571}{4.712}}
\Thicklines
\put(65,27.5){\line(1,0){23}}
\put(84,27.5){\vector(1,0){0}}
\put(75,15){\line(0,1){10.5}}
\put(55,15){\line(0,1){5}}
\put(75,40){\line(0,-1){10.5}}
\put(55,40){\line(0,-1){5}}
\put(75,36){\vector(0,1){0}}
\put(38,4){\makebox(0,0){$d$}}
\put(92,4){\makebox(0,0){$d$}}
\put(65,2){\makebox(0,0){$c$}}
\put(65,11){\makebox(0,0){$b$}}
\put(38,51){\makebox(0,0){$a$}}
\put(92,51){\makebox(0,0){$a$}}
\put(65,53){\makebox(0,0){$c$}}
\put(65,44){\makebox(0,0){$b$}}
\put(50,27.5){\makebox(0,0){$\la$}}
\put(82.5,22.5){\makebox(0,0){$\beta$}}
\put(79,35){\makebox(0,0){$\tau$}}
\thicklines
\put(143,23){\makebox(0,0){$\displaystyle\sum_{a,b,c,d\atop\beta\in\MXMm}
\frac{d_b d_c}{ww_+}$}}
\put(167,0){\line(0,1){55}}
\put(213,0){\line(0,1){10}}
\put(195,30){\arc{10}{3.142}{4.712}}
\put(203,10){\arc{20}{4.712}{0}}
\put(177,10){\arc{10}{1.571}{4.712}}
\put(203,10){\arc{10}{4.712}{1.571}}
\put(195,25){\arc{10}{1.571}{3.142}}
\put(203,45){\arc{20}{0}{1.571}}
\put(177,45){\arc{10}{1.571}{4.712}}
\put(203,45){\arc{10}{4.712}{1.571}}
\put(177,5){\line(1,0){26}}
\put(177,15){\line(1,0){26}}
\put(195,20){\line(1,0){8}}
\put(190,25){\line(0,1){5}}
\put(213,55){\line(0,-1){10}}
\put(195,35){\line(1,0){8}}
\put(177,40){\line(1,0){26}}
\put(177,50){\line(1,0){26}}
\thinlines
\put(200,5){\line(0,1){8}}
\put(200,22){\line(0,1){11}}
\put(200,29.5){\vector(0,1){0}}
\put(200,50){\line(0,-1){8}}
\put(180,15){\arc{5}{3.142}{0}}
\put(200,5){\arc{5}{3.142}{0}}
\put(180,40){\arc{5}{0}{3.142}}
\put(200,50){\arc{5}{0}{3.142}}
\put(167,27.5){\arc{5}{4.712}{1.571}}
\put(190,27.5){\arc{5}{1.571}{4.712}}
\Thicklines
\put(180,32){\vector(0,-1){0}}
\put(200,17){\line(0,1){1}}
\put(180,15){\line(0,1){25}}
\put(200,38){\line(0,-1){1}}
\put(167,27.5){\line(1,0){11}}
\put(182,27.5){\line(1,0){8}}
\put(177.5,27.5){\vector(1,0){0}}
\put(163,4){\makebox(0,0){$d$}}
\put(217,4){\makebox(0,0){$d$}}
\put(190,2){\makebox(0,0){$c$}}
\put(190,11){\makebox(0,0){$b$}}
\put(163,51){\makebox(0,0){$a$}}
\put(217,51){\makebox(0,0){$a$}}
\put(190,53){\makebox(0,0){$c$}}
\put(190,44){\makebox(0,0){$b$}}
\put(205,26.5){\makebox(0,0){$\mu$}}
\put(172.5,22.5){\makebox(0,0){$\beta$}}
\put(176,35){\makebox(0,0){$\tau$}}
\end{picture}
\end{center}
\caption{Chiral vertical projectors
$q_{\tau,\la}^+$ and $q_{\tau,\mu}^-$}
\label{qtlpm}
\end{figure}

Next we define {\sl chiral vertical projectors}
$q_{\tau,\la}^\pm\in\cZ_h^\pm$ by
\[ q_{\tau,\la}^\pm = \sum_{i=1}^{\dim H_{\tau,\la}^\pm}
E_{\tau,\la;i}^{\pm;i} \,. \]
Hence
\[ q_{\tau,\la}^+ = w d_\tau \, I_{\tau,\la}^+ \otimes I_\tau^- =
\frac{\sqrt{d_\tau d_\la}}{w} \sum_\xi |\om^{\tau,\la,+}_\xi\rangle
\langle \om^{\tau,\la,+}_\xi| \otimes I_\tau^- \,, \]
and similarly
\[ q_{\tau,\mu}^- = w d_\tau \, I_\tau^+ \otimes I_{\tau,\mu}^- =
\frac{\sqrt{d_\tau d_\mu}}{w} \sum_\xi I_\tau^+ \otimes
|\om^{\mu,\la,-}_\xi\rangle
\langle \om^{\tau,\mu,-}_\xi| \,. \]
Therefore $q_{\tau,\la}^+$ and $q_{\tau,\mu}^-$ can be displayed
graphically by the left and right hand side of Fig.\ \ref{qtlpm},
respectively.
%
\thinlines

\begin{lemma}
\label{repchir}
Whenever $\beta_\pm\in\MXMpm$ we have
\begin{equation}
\bearll
e_{\beta_+} *_v |\varphi^{\tau,\la,+}_1\rangle
\langle\varphi^{\tau,\la,+}_2| \otimes
|\varphi^{\tau,\mu,-}_3\rangle
\langle\varphi^{\tau,\mu,-}_4|
=|\pi_{\tau,\la}^+ (e_{\beta_+})\varphi^{\tau,\la,+}_1\rangle
\langle\varphi^{\tau,\la,+}_2| \otimes
|\varphi^{\tau,\mu,-}_3\rangle
\langle\varphi^{\tau,\mu,-}_4| ,\\[.4em]
e_{\beta_-} *_v |\varphi^{\tau,\la,+}_1\rangle
\langle\varphi^{\tau,\la,+}_2| \otimes
|\varphi^{\tau,\mu,-}_3\rangle
\langle\varphi^{\tau,\mu,-}_4|
=|\varphi^{\tau,\la,+}_1\rangle
\langle\varphi^{\tau,\la,+}_2| \otimes
|\pi_{\tau,\mu}^- (e_{\beta_-}) \varphi^{\tau,\mu,-}_3\rangle
\langle\varphi^{\tau,\mu,-}_4| .
\eear
\end{equation}
\end{lemma}

\begin{proof}
We only show the first relation; the proof for the second one
is analogous. It suffices to show the relation for vectors
$\omega^{\tau,\la,\pm}_{b,c,t,X}$.
Then the vertical product
\[ e_{\beta_+} *_v |\omega^{\tau,\la,+}_{b_1,c_1,t_1,X_1}\rangle
\langle\omega^{\tau,\la,+}_{b_2,c_2,t_2,X_2}| \otimes
|\omega^{\tau,\mu,-}_{b_3,c_3,t_3,X_3}\rangle
\langle\omega^{\tau,\mu,-}_{b_4,c_4,t_4,X_4}| \]
is given graphically by the left hand side of Fig.\ \ref{ebommm}.
%
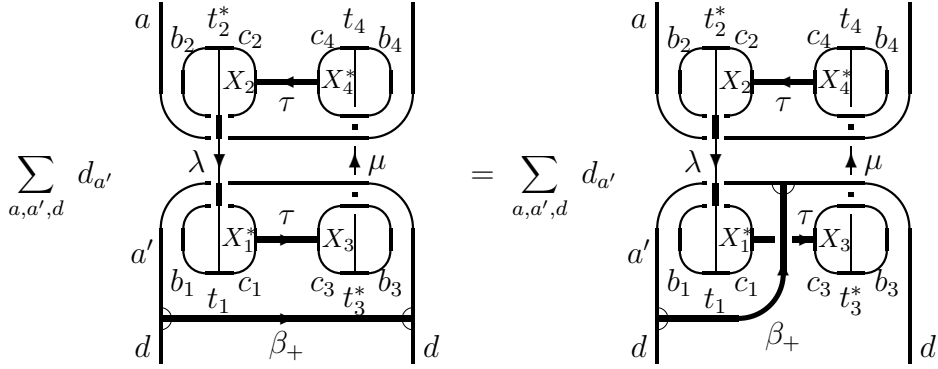
\begin{figure}[htb]
\begin{center}
\unitlength 0.6mm
\begin{picture}(211,80)
\thicklines
\put(15,39){\makebox(0,0){$\displaystyle\sum_{a,a',d}\;d_{a'}$}}
\put(37,0){\line(0,1){30}}
\put(93,0){\line(0,1){30}}
\put(47,30){\arc{20}{3.142}{4.712}}
\put(83,30){\arc{20}{4.712}{0}}
\put(47,25){\arc{10}{1.571}{3.142}}
\put(47,30){\arc{10}{3.142}{4.712}}
\put(83,30){\arc{10}{4.712}{0}}
\put(83,25){\arc{10}{0}{1.571}}
\put(77,25){\arc{10}{1.571}{3.142}}
\put(77,30){\arc{10}{3.142}{4.712}}
\put(53,30){\arc{10}{4.712}{0}}
\put(53,25){\arc{10}{0}{1.571}}
\put(47,60){\arc{20}{1.571}{3.142}}
\put(83,60){\arc{20}{0}{1.571}}
\put(47,60){\arc{10}{1.571}{3.142}}
\put(47,65){\arc{10}{3.142}{4.712}}
\put(83,65){\arc{10}{4.712}{0}}
\put(83,60){\arc{10}{0}{1.571}}
\put(77,60){\arc{10}{1.571}{3.142}}
\put(77,65){\arc{10}{3.142}{4.712}}
\put(53,65){\arc{10}{4.712}{0}}
\put(53,60){\arc{10}{0}{1.571}}
\put(47,20){\line(1,0){6}}
\put(77,20){\line(1,0){6}}
\put(47,35){\line(1,0){1}}
\put(52,35){\line(1,0){1}}
\put(77,35){\line(1,0){6}}
\put(47,40){\line(1,0){1}}
\put(52,40){\line(1,0){31}}
\put(37,80){\line(0,-1){20}}
\put(93,80){\line(0,-1){20}}
\put(47,50){\line(1,0){1}}
\put(52,50){\line(1,0){31}}
\put(47,55){\line(1,0){1}}
\put(52,55){\line(1,0){1}}
\put(77,55){\line(1,0){6}}
\put(47,70){\line(1,0){6}}
\put(77,70){\line(1,0){6}}
\put(42,25){\line(0,1){5}}
\put(42,65){\line(0,-1){5}}
\put(58,25){\line(0,1){5}}
\put(58,65){\line(0,-1){5}}
\put(88,25){\line(0,1){5}}
\put(88,65){\line(0,-1){5}}
\put(72,25){\line(0,1){5}}
\put(72,65){\line(0,-1){5}}
\thinlines
\put(50,20){\line(0,1){15}}
\put(50,40){\line(0,1){10}}
\put(80,20){\line(0,1){13}}
\put(80,42){\line(0,1){6}}
\put(50,43){\vector(0,-1){0}}
\put(80,47){\vector(0,1){0}}
\put(50,70){\line(0,-1){15}}
\put(80,70){\line(0,-1){13}}
\put(37,10){\arc{5}{4.712}{1.571}}
\put(93,10){\arc{5}{1.571}{4.712}}
\Thicklines
\put(37,10){\line(1,0){56}}
\put(67,10){\vector(1,0){0}}
\put(58,27.5){\line(1,0){14}}
\put(58,62.5){\line(1,0){14}}
\put(67,27.5){\vector(1,0){0}}
\put(63,62.5){\vector(-1,0){0}}
\put(80,37){\line(0,1){1}}
\put(50,35){\line(0,1){5}}
\put(80,53){\line(0,-1){1}}
\put(50,55){\line(0,-1){5}}
\put(65,4){\makebox(0,0){$\beta_+$}}
\put(33,4){\makebox(0,0){$d$}}
\put(97,4){\makebox(0,0){$d$}}
\put(33,25){\makebox(0,0){$a'$}}
\put(57,17){\makebox(0,0){$c_1$}}
\put(42,18){\makebox(0,0){$b_1$}}
\put(73,17){\makebox(0,0){$c_3$}}
\put(88,18){\makebox(0,0){$b_3$}}
\put(33,76){\makebox(0,0){$a$}}
\put(57,72){\makebox(0,0){$c_2$}}
\put(42,72){\makebox(0,0){$b_2$}}
\put(73,72){\makebox(0,0){$c_4$}}
\put(88,72){\makebox(0,0){$b_4$}}
\put(45,45){\makebox(0,0){$\la$}}
\put(85,44){\makebox(0,0){$\mu$}}
\put(65,32.5){\makebox(0,0){$\tau$}}
\put(65,57.5){\makebox(0,0){$\tau$}}
\put(50,14){\makebox(0,0){$t_1$}}
\put(50,76){\makebox(0,0){$t_2^*$}}
\put(80,14){\makebox(0,0){$t_3^*$}}
\put(80,76){\makebox(0,0){$t_4$}}
\put(54.5,27.5){\makebox(0,0){\footnotesize{$X_1^*$}}}
\put(54.5,62.5){\makebox(0,0){\footnotesize{$X_2$}}}
\put(76.5,27.5){\makebox(0,0){\footnotesize{$X_3$}}}
\put(76.5,62.5){\makebox(0,0){\footnotesize{$X_4^*$}}}
\thicklines
\put(122,39){\makebox(0,0){$=\;\displaystyle\sum_{a,a',d}\;d_{a'}$}}
\put(147,0){\line(0,1){30}}
\put(203,0){\line(0,1){30}}
\put(157,30){\arc{20}{3.142}{4.712}}
\put(193,30){\arc{20}{4.712}{0}}
\put(157,25){\arc{10}{1.571}{3.142}}
\put(157,30){\arc{10}{3.142}{4.712}}
\put(193,30){\arc{10}{4.712}{0}}
\put(193,25){\arc{10}{0}{1.571}}
\put(187,25){\arc{10}{1.571}{3.142}}
\put(187,30){\arc{10}{3.142}{4.712}}
\put(163,30){\arc{10}{4.712}{0}}
\put(163,25){\arc{10}{0}{1.571}}
\put(157,60){\arc{20}{1.571}{3.142}}
\put(193,60){\arc{20}{0}{1.571}}
\put(157,60){\arc{10}{1.571}{3.142}}
\put(157,65){\arc{10}{3.142}{4.712}}
\put(193,65){\arc{10}{4.712}{0}}
\put(193,60){\arc{10}{0}{1.571}}
\put(187,60){\arc{10}{1.571}{3.142}}
\put(187,65){\arc{10}{3.142}{4.712}}
\put(163,65){\arc{10}{4.712}{0}}
\put(163,60){\arc{10}{0}{1.571}}
\put(157,20){\line(1,0){6}}
\put(187,20){\line(1,0){6}}
\put(157,35){\line(1,0){1}}
\put(162,35){\line(1,0){1}}
\put(187,35){\line(1,0){6}}
\put(157,40){\line(1,0){1}}
\put(162,40){\line(1,0){31}}
\put(147,80){\line(0,-1){20}}
\put(203,80){\line(0,-1){20}}
\put(157,50){\line(1,0){1}}
\put(162,50){\line(1,0){31}}
\put(157,55){\line(1,0){1}}
\put(162,55){\line(1,0){1}}
\put(187,55){\line(1,0){6}}
\put(157,70){\line(1,0){6}}
\put(187,70){\line(1,0){6}}
\put(152,25){\line(0,1){5}}
\put(152,65){\line(0,-1){5}}
\put(168,25){\line(0,1){5}}
\put(168,65){\line(0,-1){5}}
\put(198,25){\line(0,1){5}}
\put(198,65){\line(0,-1){5}}
\put(182,25){\line(0,1){5}}
\put(182,65){\line(0,-1){5}}
\thinlines
\put(160,20){\line(0,1){15}}
\put(160,40){\line(0,1){10}}
\put(190,20){\line(0,1){13}}
\put(190,42){\line(0,1){6}}
\put(160,43){\vector(0,-1){0}}
\put(190,47){\vector(0,1){0}}
\put(160,70){\line(0,-1){15}}
\put(190,70){\line(0,-1){13}}
\put(147,10){\arc{5}{4.712}{1.571}}
\put(175,40){\arc{5}{0}{3.142}}
\Thicklines
\put(147,10){\line(1,0){18}}
\put(175,20){\line(0,1){20}}
\put(175,23){\vector(0,1){0}}
\put(165,20){\arc{20}{0}{1.571}}
\put(168,27.5){\line(1,0){5}}
\put(168,62.5){\line(1,0){14}}
\put(177,27.5){\vector(1,0){5}}
\put(173,62.5){\vector(-1,0){0}}
\put(190,37){\line(0,1){1}}
\put(160,35){\line(0,1){5}}
\put(190,53){\line(0,-1){1}}
\put(160,55){\line(0,-1){5}}
\put(175,6){\makebox(0,0){$\beta_+$}}
\put(143,4){\makebox(0,0){$d$}}
\put(207,4){\makebox(0,0){$d$}}
\put(143,25){\makebox(0,0){$a'$}}
\put(167,17){\makebox(0,0){$c_1$}}
\put(152,18){\makebox(0,0){$b_1$}}
\put(183,17){\makebox(0,0){$c_3$}}
\put(198,18){\makebox(0,0){$b_3$}}
\put(143,76){\makebox(0,0){$a$}}
\put(167,72){\makebox(0,0){$c_2$}}
\put(152,72){\makebox(0,0){$b_2$}}
\put(183,72){\makebox(0,0){$c_4$}}
\put(198,72){\makebox(0,0){$b_4$}}
\put(155,45){\makebox(0,0){$\la$}}
\put(195,44){\makebox(0,0){$\mu$}}
\put(180,32.5){\makebox(0,0){$\tau$}}
\put(175,57.5){\makebox(0,0){$\tau$}}
\put(160,14){\makebox(0,0){$t_1$}}
\put(160,76){\makebox(0,0){$t_2^*$}}
\put(190,14){\makebox(0,0){$t_3^*$}}
\put(190,76){\makebox(0,0){$t_4$}}
\put(164.5,27.5){\makebox(0,0){\footnotesize{$X_1^*$}}}
\put(164.5,62.5){\makebox(0,0){\footnotesize{$X_2$}}}
\put(186.5,27.5){\makebox(0,0){\footnotesize{$X_3$}}}
\put(186.5,62.5){\makebox(0,0){\footnotesize{$X_4^*$}}}
\end{picture}
\end{center}
\caption{The action of $e_{\beta_+}$ on
$A_{\tau,\la}^+\otimes A_{\tau,\mu}^-$}
\label{ebommm}
\end{figure}
Since $\beta_+\in\MXMp$ admits relative braiding with $\a^-_\mu$,
we can slide around the right trivalent vertex of the wire
$\beta_+$ and apply the naturality move of Fig.\ \ref{natreloi}
for the relative braiding to obtain the right hand side of
Fig.\ \ref{ebommm}. In the lower left corner we now recognize
the vector
$\pi_{\tau,\la}^+(e_{\beta_+})
\omega^{\tau,\la,+}_{b_1,c_1,t_1,X_1}$
of Fig.\ \ref{pitlb}, hence the whole diagram
represents the vector
\[ |\pi_{\tau,\la}^+(e_{\beta_+})
\omega^{\tau,\la,+}_{b_1,c_1,t_1,X_1}\rangle
\langle\omega^{\tau,\la,+}_{b_2,c_2,t_2,X_2}| \otimes
|\omega^{\tau,\mu,-}_{b_3,c_3,t_3,X_3}\rangle
\langle\omega^{\tau,\mu,-}_{b_4,c_4,t_4,X_4}| \,,\]
yielding the statement.
\end{proof}

From Lemma \ref{repchir} we now obtain the following

\begin{corollary}
\label{pulloutrep}
We have
\begin{equation}
\label{poreq}
E_{\tau,\la;i}^{\pm;j} *_v e_{\beta_\pm} *_v E_{\tau',\la';k}^{\pm;l}
= \del \tau{\tau'} \del \la{\la'} \langle u^{\tau,\la,\pm}_j ,
\pi_{\tau,\la}^\pm(e_{\beta_\pm}) u^{\tau,\la,\pm}_k \rangle
E_{\tau,\la;i}^{\pm;l} \,.
\end{equation}
\end{corollary}
In the coefficient on the right hand side of \erf{poreq} we
recognize the matrix elements of the chiral representations
$\pi_{\tau,\la}^\pm:\cZ_h^\pm\rightarrow B(H_{\tau,\la}^\pm)$.
We are now ready to prove the main result.

\begin{theorem}
\label{chirdecomthm}
We have completeness
\begin{equation}
\label{complete}
\sum_{\tau\in\MXMo} \sum_{\la,\mu\in\NXN}
\sum_{i=1}^{\dim H_{\tau,\la}^+}\sum_{k=1}^{\dim H_{\tau,\mu}^-}
E_{\tau,\la,\mu;i,k}^{i,k} = e_0 \,.
\end{equation}
Consequently the chiral vertical projectors
$q_{\tau,\la}^\pm$ sum up to the multiplicative unit $e_0$
of $(\cZ_h^\pm,*_v)$. Moreover, $q_{\tau,\la}^\pm=0$ if and only if
$b_{\tau,\la}^\pm=0$, we have mutual orthogonality
$q_{\tau,\la}^\pm *_v q_{\tau',\la'}^\pm
=\del \tau{\tau'} \del \la{\la'} q_{\tau,\la}^\pm$
and $q_{\tau,\la}^\pm$ is a minimal central projection
in $(\cZ_h^\pm,*_v)$ whenever $b_{\tau,\la}^\pm\neq 0$.
Thus the decomposition of the chiral centers into
simple matrix algebras is given as
\begin{equation}
\cZ_h^\pm \simeq \bigoplus_{\tau,\la}
\Mat(b_{\tau,\la}^\pm,\bbC) \,.
\end{equation}
\end{theorem}

\begin{proof}
All we have to show is the completeness relation
\erf{complete}; the rest is clear since then each
$e_\beta$, $\beta\in\MXMpm$ can be expanded in
the chiral matrix units. We have
\[ \sum_{\tau,\la,\mu,i,k} E_{\tau,\la,\mu;i,k}^{i,k}
= \sum_{\tau,\la,\mu,i,k} w d_\tau
|u^{\tau,\la,+}_i\rangle\langle u^{\tau,\la,+}_u| \otimes
|u^{\tau,\mu,-}_k\rangle\langle u^{\tau,\mu,-}_k|
= \sum_\tau w d_\tau I_\tau^+ \otimes I_\tau^- \,, \]
and this is given graphically by the left hand side
of Fig.\ \ref{complpic}.
%
\begin{figure}[htb]
\begin{center}
\unitlength 0.6mm
\begin{picture}(231,40)
\thicklines
\put(25,15){\makebox(0,0){$\displaystyle
\sum_{\beta_+\in\MXMp\atop\beta_-\in\MXMm}\sum_{\tau,a,b}\;
\frac{wd_\tau}{w_+^2}$}}
\put(60,0){\line(0,1){40}}
\put(135,0){\line(0,1){40}}
\put(90,15){\line(0,1){10}}
\put(95,25){\arc{10}{3.142}{4.712}}
\put(95,30){\line(1,0){5}}
\put(100,25){\arc{10}{4.712}{0}}
\put(105,15){\line(0,1){10}}
\put(95,15){\arc{10}{1.571}{3.142}}
\put(95,10){\line(1,0){5}}
\put(100,15){\arc{10}{0}{1.571}}
\Thicklines
\put(60,20){\line(1,0){13}}
\put(90,20){\line(-1,0){13}}
\put(105,20){\line(1,0){30}}
\put(75,10){\line(0,1){20}}
\put(85,30){\arc{20}{3.142}{4.712}}
\put(85,40){\line(1,0){25}}
\put(110,30){\arc{20}{4.712}{0}}
\put(120,22){\line(0,1){8}}
\put(120,18){\line(0,-1){8}}
\put(85,10){\arc{20}{1.571}{3.142}}
\put(85,0){\line(1,0){25}}
\put(110,10){\arc{20}{0}{1.571}}
\put(99.5,0){\vector(1,0){0}}
\put(70,20){\vector(1,0){0}}
\put(132,20){\vector(1,0){0}}
\thinlines
\put(60,20){\arc{5}{4.712}{1.571}}
\put(135,20){\arc{5}{1.571}{4.712}}
\put(105,20){\arc{5}{4.712}{1.571}}
\put(90,20){\arc{5}{1.571}{4.712}}
\put(55,35){\makebox(0,0){$a$}}
\put(55,5){\makebox(0,0){$b$}}
\put(140,35){\makebox(0,0){$a$}}
\put(140,5){\makebox(0,0){$b$}}
\put(97.5,25){\makebox(0,0){$a$}}
\put(97.5,15){\makebox(0,0){$b$}}
\put(110,5){\makebox(0,0){$\tau$}}
\put(67.5,15){\makebox(0,0){$\beta_-$}}
\put(127.5,15){\makebox(0,0){$\beta_+$}}
\thicklines
\put(172,19){\makebox(0,0){$=\displaystyle
\sum_{\tau,\tau',a,b}\;\frac{wd_\tau}{w_+^2}$}}
\put(200,0){\line(0,1){40}}
\put(224,0){\line(0,1){40}}
\Thicklines
\put(200,20){\line(1,0){5}}
\put(209,20){\line(1,0){16}}
\put(207,5){\line(0,1){30}}
\put(217,22){\line(0,1){13}}
\put(217,18){\line(0,-1){13}}
\put(212,35){\arc{10}{3.142}{0}}
\put(212,5){\arc{10}{0}{3.142}}
\put(214,20){\vector(1,0){0}}
\put(217,15){\vector(0,1){0}}
\thinlines
\put(200,20){\arc{5}{4.712}{1.571}}
\put(224,20){\arc{5}{1.571}{4.712}}
\put(212,6){\makebox(0,0){$\tau$}}
\put(212,25){\makebox(0,0){$\tau'$}}
\put(195,35){\makebox(0,0){$a$}}
\put(195,5){\makebox(0,0){$b$}}
\put(229,35){\makebox(0,0){$a$}}
\put(229,5){\makebox(0,0){$b$}}
\end{picture}
\end{center}
\caption{Completeness
$\sum_{\tau,\la,\mu,i,k} E_{\tau,\la,\mu;i,k}^{i,k}=e_0$}
\label{complpic}
\end{figure}
Looking at the middle part we observe that we obtain
a factor $\del {\beta_+}{\beta_-}$, and therefore we
only have a summation over $\tau'\in\MXMo$.
Then the middle bulb gives just the inner product of basis
isometries, so that only one summation over internal
fusion channels remains and we are left with the right
hand side of Fig.\ \ref{complpic}. But now we obtain
a factor $\del {\tau'}0$ and this yields exactly $e_0$
by virtue of the non-degeneracy of the ambichiral braiding,
Theorem \ref{ambinondeg}.
\end{proof}

\begin{corollary}
The total numbers of morphisms in the chiral systems
$\MXMpm$ are given by 
$\tr({}^{{\rm{t}}}\!b^\pm b^\pm)= \tr (b^\pm {}^{{\rm{t}}}\!b^\pm)=
\sum_{\tau,\la} (b_{\tau,\la}^\pm)^2$.
\end{corollary}

From Lemma \ref{keylem1} we conclude that
$q_{\la,\mu} *_v E_{\tau,\la',\mu';i,k}^{j,l}=0$ unless
$\la=\la'$ and $\mu=\mu'$. Since
$\sum_{\la,\mu} q_{\la,\mu} =e_0$ by \cite[Thm.\ 6.8]{BEK1}
we therefore obtain
$E_{\tau,\la,\mu;i,k}^{j,l}= q_{\la,\mu} *_v
E_{\tau,\la,\mu;i,k}^{j,l}$. On the other hand the
completeness relation \erf{complete} yields similarly
$q_{\la,\mu}=\sum_{\tau,i,k} q_{\la,\mu} *_v 
E_{\tau,\la,\mu;i,k}^{i,k}$. Hence we arrive at

\begin{corollary}
\label{qlm=sumE}
The vertical projector $q_{\la,\mu}$ can be expanded as
\begin{equation}
q_{\la,\mu}=\sum_{\tau\in\MXMo} \sum_{i=1}^{\dim H_{\tau,\la}^+}
\sum_{k=1}^{\dim H_{\tau,\mu}^-}E_{\tau,\la,\mu;i,k}^{i,k}
\end{equation}
for any $\la,\mu\in\NXN$.
\end{corollary}

Note that this expansion corresponds exactly to the
expansion of the modular invariant mass matrix in
chiral branching coefficients in \erf{Z=bb}.

\subsection{Representations of fusion rules and exponents}
\label{refurex}

Recall that
$\chi_\la(\nu)=Y_{\la,\nu}/d_\la=S_{\la,\nu}/S_{\la,0}$
are the evaluations of the statistics characters,
$\la,\nu\in\NXN$. Similarly we have statistics characters
for the ambichiral system:
$\chi^\ext_\tau(\tau')=Y^\ext_{\tau,\tau'}/d_\tau
=S^\ext_{\tau,\tau'}/S^\ext_{\tau,0}$. As derived in the
general theory of $\a$-induction \cite{BE1,BE3},
sectors $[\a^\pm_\nu]$ commute with all subsectors
of $[\a^+_\la][\a^-_\mu]$, thus with all sectors
arising from $\MXM$ and in particular from $\MXMpm$.
Consequently they must be scalar multiples
of the identity in the irreducible representations
of the corresponding fusion rules. In fact these
scalars must be given by the evaluations of the
chiral characters of the system $\NXN$ by virtue of
the homomorphism property of $\a$-induction
(cf.\ \cite{BE2}). We will now precisely determine
the multiplicities of the occurring characters i.e.\
the multiplicities of the eigenvalues of the
representation matrices.

\begin{lemma}
\label{eigen}
For $\la,\mu,\nu,\rho\in\NXN$ and $\tau,\tau'\in\MXMo$
we have vertical multiplication rules
\begin{equation}
e_{\tau'} *_v  p_\nu^+ *_v p_\rho^- *_v
E_{\tau,\la,\mu;i,k}^{j,l}  = d_{\tau'} \chi^\ext_\tau(\tau')
\, d_\nu \chi_\la(\nu) \, d_\rho \chi_\mu(\rho) \,
E_{\tau,\la,\mu;i,k}^{j,l} \,.
\end{equation}
\end{lemma}

\begin{proof}
It suffices to show the relation using elements
given in Fig.\ \ref{omomomom}
instead of matrix units $E_{\tau,\la,\mu;i,k}^{j.l}$.
The product
\[ p^+_\nu *_v |\omega^{\tau,\la,+}_{b_1,c_1,t_1,X_1}\rangle
\langle\omega^{\tau,\la,+}_{b_2,c_2,t_2,X_2}| \otimes
|\omega^{\tau,\mu,-}_{b_3,c_3,t_3,X_3}\rangle
\langle\omega^{\tau,\mu,-}_{b_4,c_4,t_4,X_4}| \]
is given graphically by the left hand side
of Fig.\ \ref{p+ommm1}.
%
\begin{figure}[htb]
\begin{center}
\unitlength 0.6mm
\begin{picture}(211,90)
\thicklines
\put(15,43){\makebox(0,0){$\displaystyle
\sum_{a,a',d,\atop\rho,\rho'}\;d_{a'}$}}
\put(37,0){\line(0,1){18}}
\put(93,0){\line(0,1){18}}
\put(37,22){\line(0,1){25.5}}
\put(93,22){\line(0,1){25.5}}
\put(39.5,47.5){\arc{5}{3.142}{0}}
\put(90.5,47.5){\arc{5}{3.142}{0}}
\put(47,35){\arc{10}{1.571}{3.142}}
\put(60.5,47.5){\arc{5}{3.142}{4.712}}
\put(69.5,47.5){\arc{5}{4.712}{0}}
\put(83,35){\arc{10}{0}{1.571}}
\put(77,35){\arc{10}{1.571}{3.142}}
\put(53,35){\arc{10}{0}{1.571}}
\put(47,70){\arc{20}{1.571}{3.142}}
\put(83,70){\arc{20}{0}{1.571}}
\put(47,70){\arc{10}{1.571}{3.142}}
\put(47,75){\arc{10}{3.142}{4.712}}
\put(83,75){\arc{10}{4.712}{0}}
\put(83,70){\arc{10}{0}{1.571}}
\put(77,70){\arc{10}{1.571}{3.142}}
\put(77,75){\arc{10}{3.142}{4.712}}
\put(53,75){\arc{10}{4.712}{0}}
\put(53,70){\arc{10}{0}{1.571}}
\put(47,30){\line(1,0){6}}
\put(77,30){\line(1,0){6}}
\put(60.5,50){\line(1,0){9}}
\put(37,90){\line(0,-1){20}}
\put(93,90){\line(0,-1){20}}
\put(47,60){\line(1,0){1}}
\put(52,60){\line(1,0){31}}
\put(47,65){\line(1,0){1}}
\put(52,65){\line(1,0){1}}
\put(77,65){\line(1,0){6}}
\put(47,80){\line(1,0){6}}
\put(77,80){\line(1,0){6}}
\put(42,35){\line(0,1){12.5}}
\put(42,75){\line(0,-1){5}}
\put(58,35){\line(0,1){12.5}}
\put(58,75){\line(0,-1){5}}
\put(88,35){\line(0,1){12.5}}
\put(88,75){\line(0,-1){5}}
\put(72,35){\line(0,1){12.5}}
\put(72,75){\line(0,-1){5}}
\thinlines
\put(50,30){\line(0,1){30}}
\put(80,30){\line(0,1){13}}
\put(80,47){\line(0,1){11}}
\put(50,53){\vector(0,-1){0}}
\put(80,57){\vector(0,1){0}}
\put(50,80){\line(0,-1){15}}
\put(80,80){\line(0,-1){13}}
\put(42,45){\line(1,0){6}}
\put(58,45){\line(-1,0){6}}
\put(72,45){\line(1,0){16}}
\put(56,45){\vector(1,0){0}}
\put(79,45){\vector(1,0){0}}
\put(93,15){\arc{10}{4.712}{1.571}}
\put(37,15){\arc{10}{1.571}{4.712}}
\put(93,10){\arc{5}{4.712}{1.571}}
\put(37,10){\arc{5}{1.571}{4.712}}
\put(42,45){\arc{5}{4.712}{1.571}}
\put(58,45){\arc{5}{1.571}{4.712}}
\put(72,45){\arc{5}{4.712}{1.571}}
\put(88,45){\arc{5}{1.571}{4.712}}
\Thicklines
\put(37,20){\line(1,0){56}}
\put(67,20){\vector(1,0){0}}
\put(58,37.5){\line(1,0){14}}
\put(58,72.5){\line(1,0){14}}
\put(67,37.5){\vector(1,0){0}}
\put(63,72.5){\vector(-1,0){0}}
\put(80,63){\line(0,-1){1}}
\put(50,65){\line(0,-1){5}}
\put(65,14){\makebox(0,0){$\a^+_\nu$}}
\put(33,4){\makebox(0,0){$d$}}
\put(97,4){\makebox(0,0){$d$}}
\put(33,35){\makebox(0,0){$a'$}}
\put(97,35){\makebox(0,0){$a'$}}
\put(65,47){\makebox(0,0){\footnotesize{$a'$}}}
\put(57,27){\makebox(0,0){$c_1$}}
\put(42,28){\makebox(0,0){$b_1$}}
\put(73,27){\makebox(0,0){$c_3$}}
\put(88,28){\makebox(0,0){$b_3$}}
\put(33,86){\makebox(0,0){$a$}}
\put(57,82){\makebox(0,0){$c_2$}}
\put(42,82){\makebox(0,0){$b_2$}}
\put(73,82){\makebox(0,0){$c_4$}}
\put(88,82){\makebox(0,0){$b_4$}}
\put(45,55){\makebox(0,0){$\la$}}
\put(85,54){\makebox(0,0){$\mu$}}
\put(54.5,49){\makebox(0,0){\footnotesize{$\rho$}}}
\put(76.5,50){\makebox(0,0){\footnotesize{$\rho'$}}}
\put(65,33.5){\makebox(0,0){$\tau$}}
\put(65,67.5){\makebox(0,0){$\tau$}}
\put(50,24.5){\makebox(0,0){$t_1$}}
\put(50,86){\makebox(0,0){$t_2^*$}}
\put(80,25){\makebox(0,0){$t_3^*$}}
\put(80,86){\makebox(0,0){$t_4$}}
\put(54.5,37.5){\makebox(0,0){\footnotesize{$X_1^*$}}}
\put(54.5,72.5){\makebox(0,0){\footnotesize{$X_2$}}}
\put(76.5,37.5){\makebox(0,0){\footnotesize{$X_3$}}}
\put(76.5,72.5){\makebox(0,0){\footnotesize{$X_4^*$}}}
\thicklines
\put(122,44){\makebox(0,0){$=\;\displaystyle\sum_{a,a',d}\;d_{a'}$}}
\put(147,0){\line(0,1){20}}
\put(203,0){\line(0,1){20}}
\put(157,20){\arc{20}{3.142}{4.712}}
\put(193,20){\arc{20}{4.712}{0}}
\put(157,15){\arc{10}{1.571}{3.142}}
\put(157,20){\arc{10}{3.142}{4.712}}
\put(193,20){\arc{10}{4.712}{0}}
\put(193,15){\arc{10}{0}{1.571}}
\put(187,15){\arc{10}{1.571}{3.142}}
\put(187,20){\arc{10}{3.142}{4.712}}
\put(163,20){\arc{10}{4.712}{0}}
\put(163,15){\arc{10}{0}{1.571}}
\put(157,70){\arc{20}{1.571}{3.142}}
\put(193,70){\arc{20}{0}{1.571}}
\put(157,70){\arc{10}{1.571}{3.142}}
\put(157,75){\arc{10}{3.142}{4.712}}
\put(193,75){\arc{10}{4.712}{0}}
\put(193,70){\arc{10}{0}{1.571}}
\put(187,70){\arc{10}{1.571}{3.142}}
\put(187,75){\arc{10}{3.142}{4.712}}
\put(163,75){\arc{10}{4.712}{0}}
\put(163,70){\arc{10}{0}{1.571}}
\put(157,10){\line(1,0){6}}
\put(187,10){\line(1,0){6}}
\put(157,25){\line(1,0){1}}
\put(162,25){\line(1,0){1}}
\put(187,25){\line(1,0){6}}
\put(157,30){\line(1,0){1}}
\put(162,30){\line(1,0){31}}
\put(147,90){\line(0,-1){20}}
\put(203,90){\line(0,-1){20}}
\put(157,60){\line(1,0){1}}
\put(162,60){\line(1,0){31}}
\put(157,65){\line(1,0){1}}
\put(162,65){\line(1,0){1}}
\put(187,65){\line(1,0){6}}
\put(157,80){\line(1,0){6}}
\put(187,80){\line(1,0){6}}
\put(152,15){\line(0,1){5}}
\put(152,75){\line(0,-1){5}}
\put(168,15){\line(0,1){5}}
\put(168,75){\line(0,-1){5}}
\put(198,15){\line(0,1){5}}
\put(198,75){\line(0,-1){5}}
\put(182,15){\line(0,1){5}}
\put(182,75){\line(0,-1){5}}
\thinlines
\put(160,10){\line(0,1){15}}
\put(160,30){\line(0,1){13}}
\put(190,10){\line(0,1){13}}
\put(190,32){\line(0,1){11}}
\put(160,60){\line(0,-1){13}}
\put(190,58){\line(0,-1){11}}
\put(160,53){\vector(0,-1){0}}
\put(190,57){\vector(0,1){0}}
\put(160,80){\line(0,-1){15}}
\put(190,80){\line(0,-1){13}}
\put(157,45){\line(1,0){36}}
\put(177,45){\vector(1,0){0}}
\put(157,30){\arc{30}{3.142}{4.712}}
\put(193,30){\arc{30}{4.712}{0}}
\put(142,15){\line(0,1){15}}
\put(208,15){\line(0,1){15}}
\put(147,15){\arc{10}{1.571}{3.142}}
\put(203,15){\arc{10}{0}{1.571}}
\put(147,10){\arc{5}{1.571}{4.712}}
\put(203,10){\arc{5}{4.712}{1.571}}
\Thicklines
\put(168,17.5){\line(1,0){14}}
\put(168,72.5){\line(1,0){14}}
\put(177,17.5){\vector(1,0){0}}
\put(173,72.5){\vector(-1,0){0}}
\put(190,27){\line(0,1){1}}
\put(160,25){\line(0,1){5}}
\put(190,63){\line(0,-1){1}}
\put(160,65){\line(0,-1){5}}
\put(175,50){\makebox(0,0){$\nu$}}
\put(143,4){\makebox(0,0){$d$}}
\put(207,4){\makebox(0,0){$d$}}
\put(175,35){\makebox(0,0){$a'$}}
\put(167,7){\makebox(0,0){$c_1$}}
\put(152,8){\makebox(0,0){$b_1$}}
\put(183,7){\makebox(0,0){$c_3$}}
\put(198,8){\makebox(0,0){$b_3$}}
\put(143,86){\makebox(0,0){$a$}}
\put(167,82){\makebox(0,0){$c_2$}}
\put(152,82){\makebox(0,0){$b_2$}}
\put(183,82){\makebox(0,0){$c_4$}}
\put(198,82){\makebox(0,0){$b_4$}}
\put(155,55){\makebox(0,0){$\la$}}
\put(195,54){\makebox(0,0){$\mu$}}
\put(175,22.5){\makebox(0,0){$\tau$}}
\put(175,67.5){\makebox(0,0){$\tau$}}
\put(160,4){\makebox(0,0){$t_1$}}
\put(160,86){\makebox(0,0){$t_2^*$}}
\put(190,4){\makebox(0,0){$t_3^*$}}
\put(190,86){\makebox(0,0){$t_4$}}
\put(164.5,17.5){\makebox(0,0){\footnotesize{$X_1^*$}}}
\put(164.5,72.5){\makebox(0,0){\footnotesize{$X_2$}}}
\put(186.5,17.5){\makebox(0,0){\footnotesize{$X_3$}}}
\put(186.5,72.5){\makebox(0,0){\footnotesize{$X_4^*$}}}
\end{picture}
\end{center}
\caption{The action of $p^+_\nu$ on
$A_{\tau,\la}^+\otimes A_{\tau,\mu}^-$}
\label{p+ommm1}
\end{figure}
Here we have used the expansion of the identity to replace
parallel wires $a',b_1$ and $a',b_3$ by summations over
wires $\rho$ and $\rho'$. By virtue of the unitarity of
braiding operators, the IBFE symmetries and the Yang-Baxter
relation for thin wires, the wire $\a^+_\nu$ can now be
pulled over the trivalent vertices and crossings to obtain
the right hand side of Fig.\ \ref{p+ommm1}.
Here we have already resolved the summations over
$\rho,\rho'$ back to parallel wires $a',b_1$ and
$a',b_3$, respectively.
Then we slide the trivalent vertices of the wire
$\nu$ along the wire $a'$ so that we obtain the left
hand side of Fig.\ \ref{p+ommm2}.
%
\begin{figure}[htb]
\begin{center}
\unitlength 0.6mm
\begin{picture}(206,90)
\thicklines
\put(15,42){\makebox(0,0){$\displaystyle\sum_{a,a',d}\;d_{a'}$}}
\put(37,0){\line(0,1){20}}
\put(93,0){\line(0,1){20}}
\put(47,20){\arc{20}{3.142}{4.712}}
\put(83,20){\arc{20}{4.712}{0}}
\put(47,15){\arc{10}{1.571}{3.142}}
\put(47,20){\arc{10}{3.142}{4.712}}
\put(83,20){\arc{10}{4.712}{0}}
\put(83,15){\arc{10}{0}{1.571}}
\put(77,15){\arc{10}{1.571}{3.142}}
\put(77,20){\arc{10}{3.142}{4.712}}
\put(53,20){\arc{10}{4.712}{0}}
\put(53,15){\arc{10}{0}{1.571}}
\put(47,70){\arc{20}{1.571}{3.142}}
\put(83,70){\arc{20}{0}{1.571}}
\put(47,70){\arc{10}{1.571}{3.142}}
\put(47,75){\arc{10}{3.142}{4.712}}
\put(83,75){\arc{10}{4.712}{0}}
\put(83,70){\arc{10}{0}{1.571}}
\put(77,70){\arc{10}{1.571}{3.142}}
\put(77,75){\arc{10}{3.142}{4.712}}
\put(53,75){\arc{10}{4.712}{0}}
\put(53,70){\arc{10}{0}{1.571}}
\put(47,10){\line(1,0){6}}
\put(77,10){\line(1,0){6}}
\put(47,25){\line(1,0){1}}
\put(52,25){\line(1,0){1}}
\put(77,25){\line(1,0){6}}
\put(47,30){\line(1,0){1}}
\put(52,30){\line(1,0){31}}
\put(37,90){\line(0,-1){20}}
\put(93,90){\line(0,-1){20}}
\put(47,60){\line(1,0){1}}
\put(52,60){\line(1,0){31}}
\put(47,65){\line(1,0){1}}
\put(52,65){\line(1,0){1}}
\put(77,65){\line(1,0){6}}
\put(47,80){\line(1,0){6}}
\put(77,80){\line(1,0){6}}
\put(42,15){\line(0,1){5}}
\put(42,75){\line(0,-1){5}}
\put(58,15){\line(0,1){5}}
\put(58,75){\line(0,-1){5}}
\put(88,15){\line(0,1){5}}
\put(88,75){\line(0,-1){5}}
\put(72,15){\line(0,1){5}}
\put(72,75){\line(0,-1){5}}
\thinlines
\put(50,10){\line(0,1){15}}
\put(50,30){\line(0,1){18}}
\put(80,10){\line(0,1){13}}
\put(80,32){\line(0,1){26}}
\put(50,60){\line(0,-1){8}}
\put(50,53){\vector(0,-1){0}}
\put(80,47){\vector(0,1){0}}
\put(50,80){\line(0,-1){15}}
\put(80,80){\line(0,-1){13}}
\put(47,50){\line(1,0){13}}
\put(59,50){\vector(1,0){0}}
\put(47,45){\arc{10}{1.571}{4.712}}
\put(60,40){\arc{20}{4.712}{0}}
\put(60,30){\line(0,1){5}}
\put(70,30){\line(0,1){10}}
\put(47,40){\line(1,0){1}}
\put(52,40){\line(1,0){3}}
\put(55,35){\arc{10}{4.712}{0}}
\put(60,30){\arc{5}{3.142}{0}}
\put(70,30){\arc{5}{3.142}{0}}
\Thicklines
\put(58,17.5){\line(1,0){14}}
\put(58,72.5){\line(1,0){14}}
\put(67,17.5){\vector(1,0){0}}
\put(63,72.5){\vector(-1,0){0}}
\put(80,27){\line(0,1){1}}
\put(50,25){\line(0,1){5}}
\put(80,63){\line(0,-1){1}}
\put(50,65){\line(0,-1){5}}
\put(66,52){\makebox(0,0){$\nu$}}
\put(33,4){\makebox(0,0){$d$}}
\put(97,4){\makebox(0,0){$d$}}
\put(65,35){\makebox(0,0){$a'$}}
\put(57,7){\makebox(0,0){$c_1$}}
\put(42,8){\makebox(0,0){$b_1$}}
\put(73,7){\makebox(0,0){$c_3$}}
\put(88,8){\makebox(0,0){$b_3$}}
\put(33,86){\makebox(0,0){$a$}}
\put(57,82){\makebox(0,0){$c_2$}}
\put(42,82){\makebox(0,0){$b_2$}}
\put(73,82){\makebox(0,0){$c_4$}}
\put(88,82){\makebox(0,0){$b_4$}}
\put(45,55){\makebox(0,0){$\la$}}
\put(85,45){\makebox(0,0){$\mu$}}
\put(65,22.5){\makebox(0,0){$\tau$}}
\put(65,67.5){\makebox(0,0){$\tau$}}
\put(50,4){\makebox(0,0){$t_1$}}
\put(50,86){\makebox(0,0){$t_2^*$}}
\put(80,4){\makebox(0,0){$t_3^*$}}
\put(80,86){\makebox(0,0){$t_4$}}
\put(54.5,17.5){\makebox(0,0){\footnotesize{$X_1^*$}}}
\put(54.5,72.5){\makebox(0,0){\footnotesize{$X_2$}}}
\put(76.5,17.5){\makebox(0,0){\footnotesize{$X_3$}}}
\put(76.5,72.5){\makebox(0,0){\footnotesize{$X_4^*$}}}
\thicklines
\put(126,42){\makebox(0,0){$=\;\displaystyle\sum_{a,d}\;d_\nu$}}
\put(147,0){\line(0,1){20}}
\put(203,0){\line(0,1){20}}
\put(157,20){\arc{20}{3.142}{4.712}}
\put(193,20){\arc{20}{4.712}{0}}
\put(157,15){\arc{10}{1.571}{3.142}}
\put(157,20){\arc{10}{3.142}{4.712}}
\put(193,20){\arc{10}{4.712}{0}}
\put(193,15){\arc{10}{0}{1.571}}
\put(187,15){\arc{10}{1.571}{3.142}}
\put(187,20){\arc{10}{3.142}{4.712}}
\put(163,20){\arc{10}{4.712}{0}}
\put(163,15){\arc{10}{0}{1.571}}
\put(157,70){\arc{20}{1.571}{3.142}}
\put(193,70){\arc{20}{0}{1.571}}
\put(157,70){\arc{10}{1.571}{3.142}}
\put(157,75){\arc{10}{3.142}{4.712}}
\put(193,75){\arc{10}{4.712}{0}}
\put(193,70){\arc{10}{0}{1.571}}
\put(187,70){\arc{10}{1.571}{3.142}}
\put(187,75){\arc{10}{3.142}{4.712}}
\put(163,75){\arc{10}{4.712}{0}}
\put(163,70){\arc{10}{0}{1.571}}
\put(157,10){\line(1,0){6}}
\put(187,10){\line(1,0){6}}
\put(157,25){\line(1,0){1}}
\put(162,25){\line(1,0){1}}
\put(187,25){\line(1,0){6}}
\put(157,30){\line(1,0){1}}
\put(162,30){\line(1,0){31}}
\put(147,90){\line(0,-1){20}}
\put(203,90){\line(0,-1){20}}
\put(157,60){\line(1,0){1}}
\put(162,60){\line(1,0){31}}
\put(157,65){\line(1,0){1}}
\put(162,65){\line(1,0){1}}
\put(187,65){\line(1,0){6}}
\put(157,80){\line(1,0){6}}
\put(187,80){\line(1,0){6}}
\put(152,15){\line(0,1){5}}
\put(152,75){\line(0,-1){5}}
\put(168,15){\line(0,1){5}}
\put(168,75){\line(0,-1){5}}
\put(198,15){\line(0,1){5}}
\put(198,75){\line(0,-1){5}}
\put(182,15){\line(0,1){5}}
\put(182,75){\line(0,-1){5}}
\thinlines
\put(160,10){\line(0,1){15}}
\put(160,30){\line(0,1){18}}
\put(190,10){\line(0,1){13}}
\put(190,32){\line(0,1){26}}
\put(160,60){\line(0,-1){8}}
\put(160,33){\vector(0,-1){0}}
\put(190,47){\vector(0,1){0}}
\put(160,80){\line(0,-1){15}}
\put(190,80){\line(0,-1){13}}
\put(154,50){\line(1,0){12}}
\put(159,50){\vector(1,0){0}}
\put(154,45){\arc{10}{1.571}{4.712}}
\put(166,45){\arc{10}{4.712}{1.571}}
\put(154,40){\line(1,0){4}}
\put(162,40){\line(1,0){4}}
\Thicklines
\put(168,17.5){\line(1,0){14}}
\put(168,72.5){\line(1,0){14}}
\put(177,17.5){\vector(1,0){0}}
\put(173,72.5){\vector(-1,0){0}}
\put(190,27){\line(0,1){1}}
\put(160,25){\line(0,1){5}}
\put(190,63){\line(0,-1){1}}
\put(160,65){\line(0,-1){5}}
\put(176,45){\makebox(0,0){$\nu$}}
\put(143,4){\makebox(0,0){$d$}}
\put(167,7){\makebox(0,0){$c_1$}}
\put(152,8){\makebox(0,0){$b_1$}}
\put(183,7){\makebox(0,0){$c_3$}}
\put(198,8){\makebox(0,0){$b_3$}}
\put(143,86){\makebox(0,0){$a$}}
\put(167,82){\makebox(0,0){$c_2$}}
\put(152,82){\makebox(0,0){$b_2$}}
\put(183,82){\makebox(0,0){$c_4$}}
\put(198,82){\makebox(0,0){$b_4$}}
\put(155,35){\makebox(0,0){$\la$}}
\put(195,45){\makebox(0,0){$\mu$}}
\put(175,22.5){\makebox(0,0){$\tau$}}
\put(175,67.5){\makebox(0,0){$\tau$}}
\put(160,4){\makebox(0,0){$t_1$}}
\put(160,86){\makebox(0,0){$t_2^*$}}
\put(190,4){\makebox(0,0){$t_3^*$}}
\put(190,86){\makebox(0,0){$t_4$}}
\put(164.5,17.5){\makebox(0,0){\footnotesize{$X_1^*$}}}
\put(164.5,72.5){\makebox(0,0){\footnotesize{$X_2$}}}
\put(186.5,17.5){\makebox(0,0){\footnotesize{$X_3$}}}
\put(186.5,72.5){\makebox(0,0){\footnotesize{$X_4^*$}}}
\end{picture}
\end{center}
\caption{The action of $p^+_\nu$ on
$A_{\tau,\la}^+\otimes A_{\tau,\mu}^-$}
\label{p+ommm2}
\end{figure}
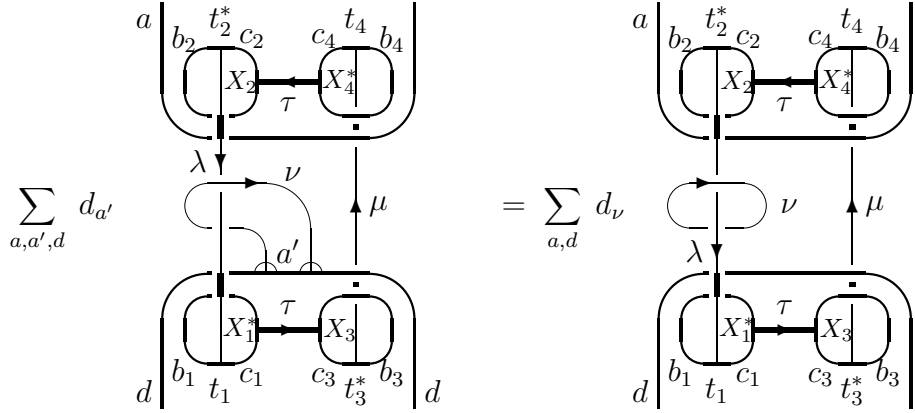
Next we turn around the small arcs at the trivalent vertices
of the wire $\nu$, yielding a factor $d_\nu/d_{a'}$, so that
the summation over $a'$ is just identified as another expansion
of the identity. Thus we arrive at the right hand side
of Fig.\ \ref{p+ommm2}. The circle $\nu$ around the wire $\la$
is evaluated as the statistics character $\chi_\la(\nu)$
(cf.\ \cite[Fig.\ 18]{BEK1}). Therefore the resulting diagram
represents
\[ d_\nu \chi_\la(\nu) \,
|\omega^{\tau,\la,+}_{b_1,c_1,t_1,X_1}\rangle
\langle\omega^{\tau,\la,+}_{b_2,c_2,t_2,X_2}| \otimes
|\omega^{\tau,\mu,-}_{b_3,c_3,t_3,X_3}\rangle
\langle\omega^{\tau,\mu,-}_{b_4,c_4,t_4,X_4}| \,. \]
The proof for $p^-_\rho$ is analogous.
Finally we consider
\[ e_{\tau'} *_v
|\omega^{\tau,\la,+}_{b_1,c_1,t_1,X_1}\rangle
\langle\omega^{\tau,\la,+}_{b_2,c_2,t_2,X_2}| \otimes
|\omega^{\tau,\mu,-}_{b_3,c_3,t_3,X_3}\rangle
\langle\omega^{\tau,\mu,-}_{b_4,c_4,t_4,X_4}|  \]
for $\tau'\in\MXMo$. We proceed graphically as
in the proof of Lemma \ref{repchir},
Fig.\ \ref{ebommm}. But now we can slide around
the trivalent vertices of the wire $\tau'$ and apply
the naturality moves of Figs.\ \ref{natrelui}
and \ref{natreloi} on both sides as $\tau'$ is
ambichiral. Therefore we obtain Fig.\ \ref{etommm}.
%
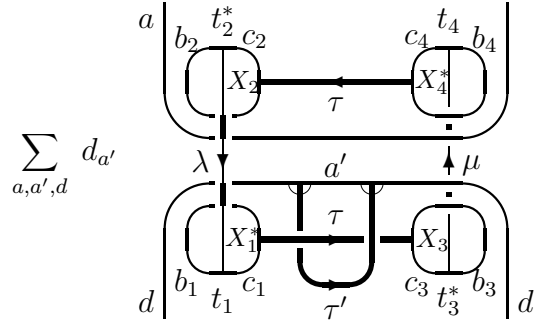
\begin{figure}[htb]
\begin{center}
\unitlength 0.6mm
\begin{picture}(122,70)
\thicklines
\put(15,34){\makebox(0,0){$\displaystyle\sum_{a,a',d}\; d_{a'}$}}
\put(37,0){\line(0,1){20}}
\put(113,0){\line(0,1){20}}
\put(47,20){\arc{20}{3.142}{4.712}}
\put(103,20){\arc{20}{4.712}{0}}
\put(47,15){\arc{10}{1.571}{3.142}}
\put(47,20){\arc{10}{3.142}{4.712}}
\put(103,20){\arc{10}{4.712}{0}}
\put(103,15){\arc{10}{0}{1.571}}
\put(97,15){\arc{10}{1.571}{3.142}}
\put(97,20){\arc{10}{3.142}{4.712}}
\put(53,20){\arc{10}{4.712}{0}}
\put(53,15){\arc{10}{0}{1.571}}
\put(47,50){\arc{20}{1.571}{3.142}}
\put(103,50){\arc{20}{0}{1.571}}
\put(47,50){\arc{10}{1.571}{3.142}}
\put(47,55){\arc{10}{3.142}{4.712}}
\put(103,55){\arc{10}{4.712}{0}}
\put(103,50){\arc{10}{0}{1.571}}
\put(97,50){\arc{10}{1.571}{3.142}}
\put(97,55){\arc{10}{3.142}{4.712}}
\put(53,55){\arc{10}{4.712}{0}}
\put(53,50){\arc{10}{0}{1.571}}
\put(47,10){\line(1,0){6}}
\put(97,10){\line(1,0){6}}
\put(47,25){\line(1,0){1}}
\put(52,25){\line(1,0){1}}
\put(97,25){\line(1,0){6}}
\put(47,30){\line(1,0){1}}
\put(52,30){\line(1,0){51}}
\put(37,70){\line(0,-1){20}}
\put(113,70){\line(0,-1){20}}
\put(47,40){\line(1,0){1}}
\put(52,40){\line(1,0){51}}
\put(47,45){\line(1,0){1}}
\put(52,45){\line(1,0){1}}
\put(97,45){\line(1,0){6}}
\put(47,60){\line(1,0){6}}
\put(97,60){\line(1,0){6}}
\put(42,15){\line(0,1){5}}
\put(42,55){\line(0,-1){5}}
\put(58,15){\line(0,1){5}}
\put(58,55){\line(0,-1){5}}
\put(108,15){\line(0,1){5}}
\put(108,55){\line(0,-1){5}}
\put(92,15){\line(0,1){5}}
\put(92,55){\line(0,-1){5}}
\thinlines
\put(50,10){\line(0,1){15}}
\put(50,30){\line(0,1){10}}
\put(100,10){\line(0,1){13}}
\put(100,32){\line(0,1){6}}
\put(50,33){\vector(0,-1){0}}
\put(100,37){\vector(0,1){0}}
\put(50,60){\line(0,-1){15}}
\put(100,60){\line(0,-1){13}}
\put(67,30){\arc{5}{0}{3.142}}
\put(83,30){\arc{5}{0}{3.142}}
\Thicklines
\put(58,17.5){\line(1,0){23}}
\put(92,17.5){\line(-1,0){7}}
\put(58,52.5){\line(1,0){34}}
\put(77,17.5){\vector(1,0){0}}
\put(73,52.5){\vector(-1,0){0}}
\put(100,27){\line(0,1){1}}
\put(50,25){\line(0,1){5}}
\put(100,43){\line(0,-1){1}}
\put(50,45){\line(0,-1){5}}
\put(83,30){\line(0,-1){17.5}}
\put(67,30){\line(0,-1){10.5}}
\put(67,15.5){\line(0,-1){3}}
\put(72,7.5){\line(1,0){6}}
\put(77,7.5){\vector(1,0){0}}
\put(78,12.5){\arc{10}{0}{1.571}}
\put(72,12.5){\arc{10}{1.571}{3.142}}
\put(33,4){\makebox(0,0){$d$}}
\put(117,4){\makebox(0,0){$d$}}
\put(57,7){\makebox(0,0){$c_1$}}
\put(42,8){\makebox(0,0){$b_1$}}
\put(93,7){\makebox(0,0){$c_3$}}
\put(108,8){\makebox(0,0){$b_3$}}
\put(33,66){\makebox(0,0){$a$}}
\put(57,62){\makebox(0,0){$c_2$}}
\put(42,62){\makebox(0,0){$b_2$}}
\put(93,62){\makebox(0,0){$c_4$}}
\put(108,62){\makebox(0,0){$b_4$}}
\put(45,35){\makebox(0,0){$\la$}}
\put(105,34){\makebox(0,0){$\mu$}}
\put(75,22.5){\makebox(0,0){$\tau$}}
\put(75,47.5){\makebox(0,0){$\tau$}}
\put(75,2.5){\makebox(0,0){$\tau'$}}
\put(75,34){\makebox(0,0){$a'$}}
\put(50,4){\makebox(0,0){$t_1$}}
\put(50,66){\makebox(0,0){$t_2^*$}}
\put(100,4){\makebox(0,0){$t_3^*$}}
\put(100,66){\makebox(0,0){$t_4$}}
\put(54.5,17.5){\makebox(0,0){\footnotesize{$X_1^*$}}}
\put(54.5,52.5){\makebox(0,0){\footnotesize{$X_2$}}}
\put(96.5,17.5){\makebox(0,0){\footnotesize{$X_3$}}}
\put(96.5,52.5){\makebox(0,0){\footnotesize{$X_4^*$}}}
\end{picture}
\end{center}
\caption{The action of $e_{\tau'}$ on
$A_{\tau,\la}^+\otimes A_{\tau,\mu}^-$}
\label{etommm}
\end{figure}
Then the small arcs of the trivalent vertices of the
wire $\tau'$ can again be turned around so that we
obtain a factor $d_{\tau'}/d_{a'}$ and that the
summation over $a'$ yields just the expansion of
the identity leaving us with parallel wires
$d$ and $\tau'$. We conclude that the resulting
diagram represents
\[ d_{\tau'} \chi_\tau^\ext(\tau') \,
|\omega^{\tau,\la,+}_{b_1,c_1,t_1,X_1}\rangle
\langle\omega^{\tau,\la,+}_{b_2,c_2,t_2,X_2}| \otimes
|\omega^{\tau,\mu,-}_{b_3,c_3,t_3,X_3}\rangle
\langle\omega^{\tau,\mu,-}_{b_4,c_4,t_4,X_4}| \,, \]
completing the proof.
\end{proof}

Recall from \cite[Sect.\ 6]{BEK1} that the irreducible
representations $\pi_{\la,\mu}$ of the full center
$(\cZ_h,*_v)$ are labelled by pairs $\la,\mu\in\NXN$ with
$Z_{\la,\mu}\neq 0$, and that they act on
$Z_{\la,\mu}$-dimensional representation spaces $H_{\la,\mu}$.
From Corollary \ref{pulloutrep} and Corollary \ref{qlm=sumE}
we now obtain the following

\begin{corollary}
\label{scalarrep}
For $\la,\mu,\nu,\rho\in\NXN$ and $\tau,\tau'\in\MXMo$
we have
\begin{equation}
\bearll
\pi_{\tau,\la}^\pm (e_{\tau'} *_v  p_\nu^\pm)
&= \; d_{\tau'} \chi^\ext_\tau(\tau') \,
d_\nu \chi_\la(\nu) \,\bfe_{H_{\tau,\la}^\pm} \,, \\[.4em]
\pi_{\la,\mu} (p_\nu^+ *_v p_\rho^-)
&= \; d_\nu \chi_\la(\nu) \, d_\rho \chi_\mu(\rho)
\,\bfe_{H_{\la,\mu}} \,.
\eear
\end{equation}
\end{corollary}

Let $\Gamma_{\nu,\rho}$, $\nu,\rho\in\NXN$, denote the
representation matrix of $[\a^+_\nu\a^-_\rho]$ in the
regular representation, i.e.\ the matrix elements
are given by
\[ \Gamma_{\nu,\rho;\beta}^{\beta'} = \langle
\beta \a^+_\nu\a^-_\rho , \beta' \rangle \,,
\qquad \beta,\beta'\in\MXM \,. \]
We can consider $\Gamma_{\nu,\rho}$ as the adjacency
matrix of the simultaneous fusion graph of
$[\a^+_\nu]$ and $[\a^-_\rho]$ on the $M$-$M$ sectors.
Similarly, let $G_\nu$, $\nu\in\NXN$, denote the
representation matrix of $[\a^\pm_\nu]$ in the
representation $\varrho\circ\Phi$ on the $M$-$N$
sectors (cf.\ \cite[Thm.\ 6.12]{BEK1}),
i.e.\ the matrix elements are given by
\[ G_{\nu;a}^b = \langle a \a^\pm_\nu , b \rangle
= \langle \nu a,b \rangle \,, \qquad a,b\in\NXM \,,\]
where the second equality is due to
\cite[Prop.\ 3.1]{BEK1}, and hence there is no
distinction between $+$ and $-$. We can consider $G_\nu$
as the adjacency matrix of the fusion graph of $[\a_\nu^\pm]$
on the $M$-$N$ sectors via left multiplication.
Finally, let $\Gamma_{\tau',\nu}^\pm$, $\tau'\in\MXMo$,
$\nu\in\NXN$, denote the representation matrices of
$[\tau' \a_\nu^\pm]$ in the chiral regular representations,
i.e.\ the matrix elements are given by
\[ \Gamma_{\tau',\nu;\beta}^{\pm;\beta'} = \langle
\beta \tau' \a_\nu^\pm, \beta' \rangle \,,
\qquad \beta,\beta'\in\MXMpm \,. \]
We now arrive at our classification result.

\pagebreak[4]

\begin{theorem}
\label{adjacency}
The eigenvalues (``exponents'') of $\Gamma_{\nu,\rho}$,
$G_\nu$ and $\Gamma_{\tau,\nu}^\pm$ for
$\nu,\rho\in\NXN$, $\tau'\in\MXMo$ are given by
$\chi_\la(\nu)\chi_\mu(\rho)$, $\chi_\la(\nu)$, and
$\chi_\tau^\ext(\tau')\chi_\la(\nu)$, respectively,
where $\la,\mu\in\NXN$ and $\tau\in\MXMo$. They
occur with the following multiplicities:
\begin{enumerate}
\item $\mult (\chi_\la(\nu)\chi_\mu(\rho)) = Z_{\la,\mu}^2$
for $\Gamma_{\nu,\rho}$,
\item $\mult (\chi_\la(\nu)) = Z_{\la,\la}$
for $G_\nu$,
\item $\mult (\chi_\tau^\ext(\tau')\chi_\la(\nu))
= (b_{\tau,\la}^\pm)^2$ for $\Gamma_{\tau',\nu}^\pm$.
\end{enumerate}
\end{theorem}

\begin{proof}
From the decomposition of the chiral centers in
Theorem \ref{chirdecomthm} it follows that the (left)
regular representations $\pi_{{\rm{reg}}}^\pm$ of
$(\cZ_h^\pm,*_v)$ decompose into irreducibles as
$\pi_{{\rm{reg}}}^\pm=\bigoplus_{\tau,\la} b_{\tau,\la}^\pm
\pi_{\tau,\la}^\pm$. It follows similarly from
\cite[Thm.\ 6.8]{BEK1} that the (left) regular
representation $\pi_{{\rm{reg}}}$ of $(\cZ_h,*_v)$
decomposes into irreducibles as 
$\pi_{{\rm{reg}}}=\bigoplus_{\la,\mu} Z_{\la,\mu}\pi_{\la,\mu}$.
Representations of the corresponding fusion rule algebras
of $M$-$M$ sectors are obtained by composition with the
isomorphisms $\Phi$ mapping the $M$-$M$ fusion rule algebra
to $(\cZ_h,*_v)$.
It was established in \cite[Thm.\ 6.12]{BEK1} that the
representation $\varrho\circ\Phi$ of the full $M$-$M$ fusion
rule algebra obtained by left action multiplication on the
$M$-$N$ sectors decomposes into irreducibles as
$\varrho\circ\Phi=\bigoplus_\la \pi_{\la,\la}\circ\Phi$.
The claim follows now since $\Phi$ fulfills
$\Phi([\beta])=d_\beta^{-1} e_\beta$ by definition
(cf.\ \cite[Def.\ 4.5]{BEK1}) and
$\Phi([\a^\pm_\nu])=d_\nu^{-1} p_\nu^\pm$
by the identification theorem \cite[Thm.\ 5.3]{BEK1}.
\end{proof}

Recall that chiral locality implies for the branching
coefficients $b_{\tau,\la}=b^\pm_{\tau,\la}$.
The third statement of Theorem \ref{adjacency} was
actually conjectured in \cite[Subsect.\ 4.2]{BE2} for
conformal inclusions and (local) simple current extensions
as a refinement of \cite[Thm.\ 4.10]{BE2},
and such a connection between branching coefficients
and dimensions of eigenspaces was first raised as a
question in \cite[Page 21]{X2} in the context
of conformal inclusions.

\section{The A-D-E classification of $\SUz$ modular invariants}
\label{ADESUz}

We now consider $\SUz_k$ braided subfactors,
i.e.\ we are dealing with subfactors $N\subset M$
where the system $\NXN$ is given by morphisms
$\la_j$, $j=0,1,2,...,k$, $\la_0=\id$,
such that we have fusion rules
$[\la_j][\la_{j'}]=\bigoplus_{j''}N_{j,j'}^{j''}[\la_{j''}]$
with
\begin{equation}
\label{SU2kfuru}
N_{j,j'}^{j''} = \left\{ \begin{array}{lc}
1 \qquad & |j-j'| \le j'' \le \min(j+j',2k-j-j') \,,
\quad j+j'+j''\in2\bbZ \,, \\ 0 & {{\rm{otherwise}}},
\end{array} \right.
\end{equation}
and that the statistics phases are given by
\[ \omega_j = \E^{2\pi\I h_j} \,, \qquad
h_j = \frac{j(j+2)}{4k+8} \, \]
where $k=1,2,3,...$ is the level. Therefore we
are constructing modular invariants of the well-known
representations of $\SLZ$ arising from the $\SUz$
level $k$ WZW models.

\subsection{The local inclusions: $\rmA_\ell$,
$\rmD_{2\ell}$, $\rmE_6$ and $\rmE_8$}

We first recall the treatment of the local extensions,
i.e.\ inclusions where the chiral locality condition
is met. Namely, we consider ``quantum field theoretical
nets of subfactors'' \cite{LR} $N(I)\subset M(I)$ on the
punctured circle along the lines of \cite{BE1,BE2,BE3}.
Here these algebras live on a Hilbert space $\cH$, and
the restriction of the algebras $N(I)$ to the vacuum
subspace $\cH_0$ is of the form $\pi_0(\LISUz)''$
with $\pi_0$ being the level $k$ vacuum representation
of $\LSUz$. We choose some interval $I_\circ$ to obtain
a single subfactor $N=N(I_\circ)\subset M(I_\circ)=M$.
Then the system $\NXN=\{\la_j\}$ is given by the
restrictions of DHR endomorphisms to the local algebras
which arise from Wassermann's \cite{W2} bimodule
construction (see \cite{BE2} for more explanation).
The braiding is then given by the DHR statistics
operators.

A rather trivial situation is
clearly given by the trivial inclusion
$N(I)=M(I)=\pi_0(\LISUz)''$ corresponding to
$[\canr]=[\id]$. We then obviously have
$[\a^\pm_j]=[\la_j]$ for all $j$. (We denote
$[\a^\pm_j]\equiv[\a^\pm_{\la_j}]$.) Therefore
we just produce the trivial modular invariant
$Z_{j,j'}=\del j{j'}$, and the simultaneous fusion
graph of $[\a^+_1]$ and $[\a^-_1]$ is nothing but
one and the same graph $\rmA_{k+1}$.

More interesting are the local simple current
extensions (or ``orbifold inclusions'') considered
in \cite{BE2,BE3}. They occur at levels $k=4\ell-4$,
$\ell=2,3,4,...$, and are constructed by means of the
simple current $\la_k$ which satisfies $\la_k^2=\id$
and so that $[\canr]=[\id]\oplus[\la_k]$.
The structure of the full system $\MXM$, producing the
$\rmD_{2\ell}$ modular invariant, has been determined
in \cite[Subsect.\ 6.2]{BE3}. The fusion graphs of
$[\a^\pm_1]$ in the chiral systems were already identified
in \cite{BE2} as $\rmD_{2\ell}$. Note that these are also
the graphs with adjacency matrix $G_1$, arising from
the multiplication on $M$-$N$ sectors. This is actually
a general fact rather than a coincidence: Whenever
the chiral locality condition
$\epsp\canr\canr\can(v)=\can(v)$ holds, then the set
$\MXN$ consists of morphisms $\beta\iota$ where $\beta$
varies in either $\MXMp$ or equivalently in $\MXMm$
due to \cite[Lemma 4.1]{BE3}.

The exceptional invariants labelled by $\rmE_6$ and $\rmE_8$
arise form conformal inclusions $\SUz_{10}\subset\SOf_1$
and $\SUz_{28}\subset(\Gtwo)_1$, respectively, and have
been treated in the nets of subfactors setting in
\cite{X1,BE2,BE3}. The structure of the full systems has
been completely determined in \cite[Subsect.\ 6.1]{BE3}.
Note that in all these $\SUz$ cases the simultaneous fusion
graphs of $[\a^+_1]$ and $[\a^-_1]$ turn out
\cite[Figs.\ 2,5,8,9]{BE3} (and similarly for the
non-local examples Figs.\ \ref{Dodd} and \ref{E7} below)
to coincide with Ocneanu's
diagrams for his ``quantum symmetry on Coxeter graphs''
\cite{O7}. The reason for this coincidence reflects
the relation between $\a$-induction and
chiral generators for double triangle algebras
\cite[Thm.\ 5.3]{BEK1}. (See also the appendix of
this paper for relations between our subfactors
specified by canonical endomorphisms in a $\SUz_k$
sector system and GHJ subfactors used in \cite{O7}.)

\subsection{The non-local simple current extensions:
$\rmD_{2\ell+1}$}

We are now passing to the non-local examples which were
not treated in \cite{BE2,BE3}. Without chiral locality
we only have the inequality
\begin{equation}
\label{mainineq}
\langle\a^\pm_\la,\a^\pm_\mu\rangle \le \langle
\canr\la,\mu\rangle
\end{equation}
rather then the ``main formula'' \cite[Thm.\ 3.9]{BE1}
because the ``$\ge$'' part of the proof of the main formula
relies on the chiral locality condition. We remark that
\erf{mainineq} is the analogue of Ocneanu's ``gap''
argument used in his A-D-E setup of \cite{O7} but
\erf{mainineq} is the suitable formulation for our
more general setting which can in particular be used
for non-local simple current extensions and other
non-local inclusions of $\LSUn$ theories. Moreover, we
know that for the local cases, e.g.\ conformal inclusions
and local simple current extensions of $\LSUn$ as
treated in \cite{BE2,BE3}, we have exact equality
and this makes concrete computations much easier.

\begin{figure}[htb]
\unitlength 0.6mm
\begin{center}
\begin{picture}(210,110)
\thinlines 
\multiput(40,25)(0,15){5}{\circle*{2}}
\multiput(10,55)(60,0){2}{\circle*{2}}
\multiput(40,70)(0,15){2}{\circle{4}}
\multiput(40,25)(0,15){2}{\circle{4}}
\multiput(40,25)(0,15){5}{\circle{6}}
\multiput(10,55)(60,0){2}{\circle{6}}
\thicklines 
\path(40,85)(10,55)(40,70)
\path(39.5,70)(39.5,40)
\path(40,40)(70,55)(40,25)
\dottedline{1.8}(40,85)(70,55)(40,70)
\dottedline{1.8}(40.5,70)(40.5,40)
\dottedline{1.8}(40,40)(10,55)(40,25)
\put(40,92){\makebox(0,0){\footnotesize{$[\id]$}}}
\put(5,62){\makebox(0,0){\footnotesize{$[\a^+_1]$}}}
\put(40,76){\makebox(0,0){\footnotesize{$[\a^+_2]$}}}
\put(50,55){\makebox(0,0){\footnotesize{$[\a^+_3]$}}}
\put(40,34){\makebox(0,0){\footnotesize{$[\a^+_4]$}}}
\put(75,48){\makebox(0,0){\footnotesize{$[\a^+_5]$}}}
\put(40,18){\makebox(0,0){\footnotesize{$[\a^+_6]$}}}
\thinlines 
\multiput(155,10)(0,15){7}{\circle*{2}}
\multiput(110,55)(90,0){2}{\circle*{2}}
\multiput(125,55)(60,0){2}{\circle*{2}}
\multiput(155,70)(0,15){3}{\circle{4}}
\multiput(155,10)(0,15){3}{\circle{4}}
\multiput(155,10)(0,15){7}{\circle{6}}
\multiput(110,55)(90,0){2}{\circle{6}}
\multiput(125,55)(60,0){2}{\circle{6}}
\thicklines 
\path(155,100)(110,55)(155,85)(125,55)(155,70)
\path(154.5,70)(154.5,40)
\path(155,40)(185,55)(155,25)(200,55)(155,10)
\dottedline{1.8}(155,100)(200,55)(155,85)(185,55)(155,70)
\dottedline{1.8}(155.5,70)(155.5,40)
\dottedline{1.8}(155,40)(125,55)(155,25)(110,55)(155,10)
\put(155,107){\makebox(0,0){\footnotesize{$[\id]$}}}
\put(105,62){\makebox(0,0){\footnotesize{$[\a^+_1]$}}}
\put(155,91){\makebox(0,0){\footnotesize{$[\a^+_2]$}}}
\put(137,55){\makebox(0,0){\footnotesize{$[\a^+_3]$}}}
\put(155,76){\makebox(0,0){\footnotesize{$[\a^+_4]$}}}
\put(162,49){\makebox(0,0){\footnotesize{$[\a^+_5]$}}}
\put(155,34){\makebox(0,0){\footnotesize{$[\a^+_6]$}}}
\put(173,55){\makebox(0,0){\footnotesize{$[\a^+_7]$}}}
\put(155,19){\makebox(0,0){\footnotesize{$[\a^+_8]$}}}
\put(205,48){\makebox(0,0){\footnotesize{$[\a^+_9]$}}}
\put(155,3){\makebox(0,0){\footnotesize{$[\a^+_{10}]$}}}
\end{picture}
\caption{$\rmD_5$ and $\rmD_7$: Fusion graphs
of $[\a^+_1]$ and $[\a^-_1]$}
\label{Dodd}
\end{center}
\end{figure}
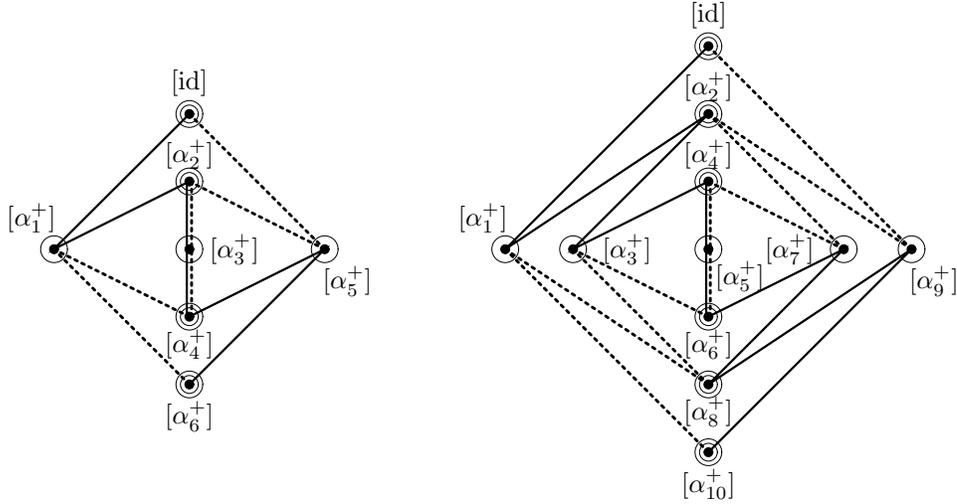

As our first non-local example we consider the simple
current extensions of $\LSUz$ which, as we will see,
produce the $\rmD_{{\rm{odd}}}$ modular invariants.
We start again with a net of local algebras for the
$\LSUz$ theories and construct nets of subfactors by
simple current extensions along the lines of
\cite[Sect.\ 3]{BE2} and \cite[Subsect.\ 6.2]{BE3}.
Using the simple current $[\la_k]$ at level $k$
satisfying the fusion rule $[\la_k^2]=[\id]$,
it was found in \cite{BE2} that a local extension
is only possible for $k\in 4\bbZ$. However, to proceed
with the crossed product construction we only need
the existence of a representative morphism $\la_k$
of the sector $[\la_k]$ which satisfies $\la_k^2=\id$
as an endomorphism. By \cite[Lemma 4.4]{R2},
such a choice is possible if and only if the
statistics phase $\om_k$ of $[\la_k]$ fulfills
$\om_k^2=1$. As $\om_k=\E^{2\pi\I h_k}$ by the
conformal spin and statistics theorem \cite{GL}
(see also \cite{FG,FRS2})
and since this conformal dimension is given by
$h_k=k/4$, an extension can be constructed
whenever the level is even. Now $k=4\ell-4$ is
the local case producing $\rmD_{2\ell}$, so here
we are looking at $k=4\ell-2$ where $\ell=2,3,4,...$.
Because of \erf{mainineq}
we find with $[\canr]=[\id]\oplus[\la_k]$ that
$\langle\a^\pm_j,\a^\pm_{j'}\rangle
\le \del j{j'}+\del j{k-j'}$ and hence all
$[\a^\pm_j]$'s are forced to be irreducible
except $[\a^\pm_{2\ell-1}]$ which may either be
irreducible or decompose into two irreducibles.
Moreover, we conclude
$Z_{0,j}=\langle\id,\a^-_j\rangle
=\langle\canr,\la_j\rangle=0$ for $j=1,2,...,k-1$.
But we also obtain $Z_{0,k}=0$ from $\omega_k=-1$
and $[T,Z]=0$. Thus we have $Z_{0,j}=\del 0j$,
and this forces a modular invariant mass matrix
already to be a permutation matrix by
Proposition \ref{Zperm}. Now let us look
at the $M$-$N$ sectors which are
subsectors of the $[\iota\la]$'s. By Frobenius
reciprocity, we have in general
\begin{equation}
\label{NMmain}
\langle\iota\la,\iota\mu\rangle
=\langle\canr\la,\mu\rangle \,,
\qquad \la,\mu\in\NXN \,.
\end{equation}
Therefore we find here
$\langle\iota\la_j,\iota\la_{j'}\rangle
=\del j{j'} + \del j{k-j'}$. This is enough
to conclude that we have $2\ell+1$ irreducible
$M$-$N$ morphisms which can be given by $\iota\la_j$,
$j=0,1,2,...,2\ell-2$, and $\co b, \co{b'}$ with
$[\iota\la_{2\ell-1}]=[\co b]\oplus[\co{b'}]$.
As a consequence, the matrix $G_1$ (i.e.\ the
matrix $G_\nu$ for $\nu=\la_1$) is determined
to be the adjacency matrix of $\rmD_{2\ell+1}$.
The exponents of $\rmD_{2\ell+1}$
are known to be (see e.g.\ \cite{GHJ})
\[ \Exp(\rmD_{2\ell+1})= \{ 0,2,4,\ldots,
4\ell-2,2\ell-1 \} \]
and all occur with multiplicity one.
Theorem \ref{adjacency} forces
the diagonal part of $Z$ to be
\[ Z_{j,j} = \left\{ \bearll
1 \qquad & j\in\Exp(\rmD_{2\ell+1}) \\
0 & j\notin\Exp(\rmD_{2\ell+1})
\eear \right. . \]
By virtue of the classification of $\SUz$
modular invariants \cite{CIZ2,Kt} we could now
argue that $Z$ must be the mass matrix labelled
by $\rmD_{2\ell+1}$, however, this is not
necessary since simple and general arguments
already allow to construct $Z$ directly.
Namely, as $Z$ is a permutation matrix we have
$Z_{j,j'}=\del j{\pi(j')}$ with $\pi$ a permutation
such that $\pi(j)=j$ for $j\in\Exp(\rmD_{2\ell+1})$
and $\pi(j)\neq j$ for $j\notin\Exp(\rmD_{2\ell+1})$.
But since $\pi$ defines a fusion rule automorphism
we necessarily have
$d_{\pi(j)}=d_j$. The values of the statistical
dimensions for $\SUz_k$ then allow only $\pi(j)=j$
or $\pi(j)=k-j$. We therefore have derived
\[ Z_{j,j'} = \left\{ \bearll
\del j{j'} \qquad & j\in\Exp(\rmD_{2\ell+1}) \\
\del j{k-j'} & j\notin\Exp(\rmD_{2\ell+1})
\eear \right. , \qquad j,j'=0,1,2,...,k \,.\]
This is the well-known mass matrix which was
labelled by $\rmD_{2\ell+1}$ in \cite{CIZ1}.
Note that we have $\MXMpm=\MXM$ here. We can now
easily draw the simultaneous fusion graph of
$[\a^+_1]$ and $[\a^-_1]$ which we display
in Fig.\ \ref{Dodd} for $\rmD_5$ and $\rmD_7$.
%
%
As in \cite{BE3}, we draw straight lines for the fusion
with $[\a^+_1]$ and dotted lines for the fusion with
$[\a^-_1]$. (Note that $[\a^-_1]=[\a^+_{k-1}]$ here.)
We also encircle even vertices by small
circles and ambichiral (i.e.\ ``marked'') vertices by
large circles.

Note that we have
$\langle\a^\pm_k,\can\rangle
=\langle\a^\pm_k\iota,\iota\rangle
=\langle\iota\la_k,\iota\rangle
=\langle\la_k,\canr\rangle=1$ by Frobenius
reciprocity. Since $d_\can=d_\canr=2$ we
conclude $[\can]=[\id]\oplus[\a^+_k]$.
This shows that \cite[Lemma 3.17]{BE3}
(and in turn \cite[Cor.\ 3.18]{BE3}) does not
hold true without chiral locality.

\subsection{$\rm E_7$}

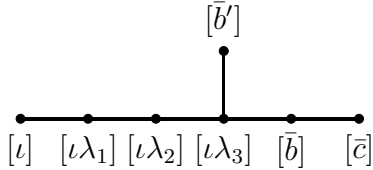
\begin{figure}[htb]
\unitlength 0.6mm
\begin{center}
\begin{picture}(95,35)
\thinlines 
\multiput(10,10)(15,0){6}{\circle*{2}}
\put(55,25){\circle*{2}}
\thicklines
\put(10,10){\line(1,0){75}}
\put(55,10){\line(0,1){15}}
\put(10,3){\makebox(0,0){$[\iota]$}}
\put(25,3){\makebox(0,0){$[\iota\la_1]$}}
\put(40,3){\makebox(0,0){$[\iota\la_2]$}}
\put(55,3){\makebox(0,0){$[\iota\la_3]$}}
\put(70,3){\makebox(0,0){$[\co b]$}}
\put(85,3){\makebox(0,0){$[\co c]$}}
\put(55,32){\makebox(0,0){$[\co b']$}}
\end{picture}
\caption{$G_1$ is the adjacency matrix of $\rmE_7$}
\label{MNE7}
\end{center}
\end{figure}

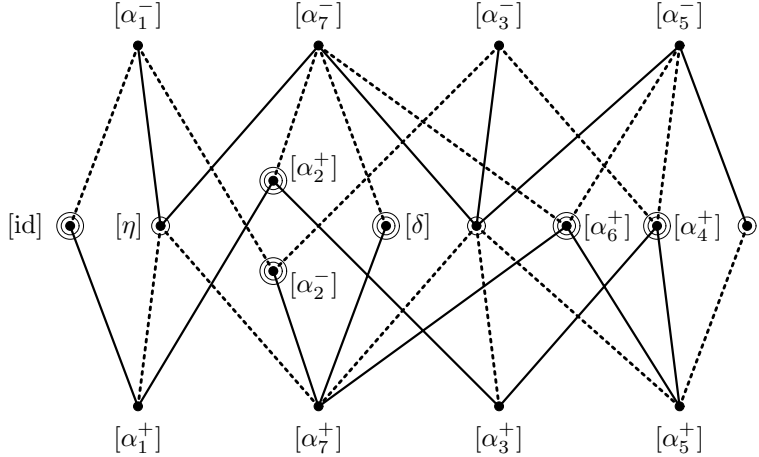
\begin{figure}[htb]
\unitlength 0.6mm
\begin{center}
\begin{picture}(180,100)
\thinlines 
\multiput(30,10)(40,0){4}{\circle*{2}}
\multiput(30,90)(40,0){4}{\circle*{2}}
\multiput(15,50)(20,0){2}{\circle*{2}}
\multiput(60,40)(0,20){2}{\circle*{2}}
\multiput(85,50)(20,0){5}{\circle*{2}}
\multiput(15,50)(20,0){2}{\circle{4}}
\multiput(60,40)(0,20){2}{\circle{4}}
\multiput(85,50)(20,0){5}{\circle{4}}
\multiput(15,50)(70,0){2}{\circle{6}}
\multiput(60,40)(0,20){2}{\circle{6}}
\multiput(125,50)(20,0){2}{\circle{6}}
\thicklines 
\path(15,50)(30,10)(60,60)(110,10)(145,50)(150,10)
(125,50)(70,10)(85,50)
\path(70,10)(60,40)
\path(30,90)(35,50)(70,90)(105,50)(150,90)(165,50)
\path(105,50)(110,90)
\dottedline{1.8}(15,50)(30,90)(60,40)(110,90)(145,50)
(150,90)(125,50)(70,90)(85,50)
\dottedline{1.8}(70,90)(60,60)
\dottedline{1.8}(30,10)(35,50)(70,10)(105,50)(150,10)(165,50)
\dottedline{1.8}(105,50)(110,10)
\put(5,50){\makebox(0,0){\footnotesize{$[\id]$}}}
\put(28,50){\makebox(0,0){\footnotesize{$[\eta]$}}}
\put(30,3){\makebox(0,0){\footnotesize{$[\a^+_1]$}}}
\put(70,3){\makebox(0,0){\footnotesize{$[\a^+_7]$}}}
\put(110,3){\makebox(0,0){\footnotesize{$[\a^+_3]$}}}
\put(150,3){\makebox(0,0){\footnotesize{$[\a^+_5]$}}}
\put(30,97){\makebox(0,0){\footnotesize{$[\a^-_1]$}}}
\put(70,97){\makebox(0,0){\footnotesize{$[\a^-_7]$}}}
\put(110,97){\makebox(0,0){\footnotesize{$[\a^-_3]$}}}
\put(150,97){\makebox(0,0){\footnotesize{$[\a^-_5]$}}}
\put(69,37){\makebox(0,0){\footnotesize{$[\a^-_2]$}}}
\put(69,63){\makebox(0,0){\footnotesize{$[\a^+_2]$}}}
\put(92,50){\makebox(0,0){\footnotesize{$[\delta]$}}}
\put(134,50){\makebox(0,0){\footnotesize{$[\a^+_6]$}}}
\put(154,50){\makebox(0,0){\footnotesize{$[\a^+_4]$}}}
\end{picture}
\caption{$\rmE_7$: Fusion graphs
of $[\a^+_1]$ and $[\a^-_1]$}
\label{E7}
\end{center}
\end{figure}

We put $N=\pi_0(\LISUz)''$ where $\pi_0$ here denotes
the level $16$ vacuum representation of $\LSUz$.
We will show in the appendix (Lemma \ref{E7can})
that there is an endomorphism $\canr\in\Mor(N,N)$
at level $k=16$ such that
$[\canr]=[\id]\oplus[\la_8]\oplus[\la_{16}]$ and
which is the dual canonical endomorphism of some
subfactor $N\subset M$.
We will now show that this dual canonical endomorphism
produces the $\rmE_7$ modular invariant.
From \erf{NMmain} we obtain
$\langle\iota\la_j,\iota\la_{j'}\rangle=\del j{j'}+
N_{8,j}^{j'}+\del j{k-j'}$ where the fusion rules come
from \erf{SU2kfuru} with $k=16$. With this it is
straightforward to check that $[\iota\la_j]$, $j=0,1,2,3$,
are irreducible and distinct $M$-$N$ sectors.
As $\langle\iota\la_4,\iota\la_4\rangle=2$ but
$\langle\iota\la_4,\iota\la_j\rangle=0$ for $j=0,1,2,3$
we conclude that $[\iota\la_4]$ decomposes into two
new different sectors, $[\iota\la_4]=[\co b]\oplus[\co b']$.
Similarly, $[\iota\la_5]$ decomposes into two
sectors but here we have
$[\iota\la_5]=[\iota\la_3]\oplus[\co c]$ with only one
new $M$-$N$ sector $[\co c]$ because
$\langle\iota\la_5,\iota\la_3\rangle=1$. We have
$\langle\iota\la_6,\iota\la_j\rangle=1$ for $j=2$
and $j=4$, so $[\iota\la_6]$ has one subsector in
common with $[\iota\la_4]$, say $[\co b]$:
$[\iota\la_6]=[\iota\la_2]\oplus[\co b]$. We similarly
find that the other $[\iota\la_j]$'s do not produce
new $M$-$N$ sectors. From
$[\iota\la_5][\la_1]=[\iota\la_4]\oplus[\iota\la_6]$
and 
$[\iota\la_3][\la_1]=[\iota\la_2]\oplus[\iota\la_4]$
we now obtain $[\co c][\la_1]=[\co b]$. Thanks to
Frobenius reciprocity we find also that $[\co c]$
appears in the decomposition of $[\co b][\la_1]$.
This forces $[\co b][\la_1]=[\iota\la_3]\oplus[\co c]$
and $[\co b'][\la_1]=[\iota\la_3]$. We therefore
have determined the matrix $G_1$ to be the adjacency
matrix of $\rmE_7$, see Fig.\ \ref{MNE7}.
%
%
The exponents of $\rmE_7$ are given by
$\Exp(\rmE_7)=\{0,4,6,8,10,12,16\}$
and all occur with multiplicity one.
Theorem \ref{adjacency} forces
the diagonal part of $Z$ to be
\[ Z_{j,j} = \left\{ \bearll
1 \qquad & j\in\Exp(\rmE_7) \\
0 & j\notin\Exp(\rmE_7)
\eear \right. . \]
By virtue of the classification of $\SUz$
modular invariants \cite{CIZ2,Kt} we could now
argue that $Z$ must be the mass matrix labelled
by $\rmE_7$ but, as it is quite instructive, we
prefer again to construct $Z$ directly. From
\erf{mainineq} we conclude that among the
zero-co\-lumn/row only $Z_{0,0}$, $Z_{0,8}$,
$Z_{8,0}$, $Z_{0,16}$ and $Z_{16,0}$ can at most
be one. But $[T,Z]=0$ and $h_8=10/9$ forces
$Z_{0,8}=Z_{8,0}=0$. Now assume for contradiction
that $Z_{0,16}$ (and hence $Z_{16,0}$) is zero.
Then $Z$ would be a permutation matrix by
Proposition \ref{Zperm}. As $Z_{1,1}=0$ this
would imply that $Z_{1,j}\neq 0$ for some $j\neq0$,
but this contradicts $[T,Z]=0$ because $h_1=1/24$
and there is no other $j$ with
$h_j=1/24\,{{\rm{mod}}}\,\bbZ$. Consequently
$Z_{0,16}=Z_{16,0}=1$. But the zero-column determines
$\langle\a^+_j,\a^+_{j'}\rangle$ since
\[ \langle\a^+_j,\a^+_{j'}\rangle =
\langle\a^+_j\a^+_{j'},\id\rangle=\sum_{j''}
N_{j,j'}^{j''}Z_{j'',0}=\del j{j'}+\del j{16-j'}\,,\]
and similarly the zero row determines
$\langle\a^-_j,\a^-_{j'}\rangle=\del j{j'}+\del j{16-j'}$.
This forces the fusion graphs of $[\a^\pm_1]$ in
the chiral sector systems to be $\rmD_{10}$, and then
the whole fusion tables for the systems $\MXMpm$
are determined completely \cite{I0}. Moreover,
we learn $w_+=w/2$ from Proposition \ref{sumcg} and
$w_0=w_+/2$ from Theorem \ref{ambinondeg}.
This forces the subsystem $\MXMo\subset\MXMpm$
to correspond to the even vertices of the
$\rmD_{10}$ graph so that it can be given by
$\MXMo=\{\id,\a^+_2,\a^+_4,\a^+_6,\delta,\delta'\}$
with $\delta,\delta'\in\Mor(M,M)$ such that
$[\a^+_8]=[\delta]\oplus[\delta']$. The well-known 
Perron-Frobenius eigenvector of $\rmD_{10}$ tells us
$d_\delta=d_{\delta'}=d_8/2$. Note that
$[\a_8^+]$ and $[\a_8^-]$ have only one sector
in common, say $[\delta]$, since $Z_{8,8}=1$.
On the other hand, $[\a^-_8]$ decomposes into two
sectors, $[\a^-_8]=[\delta]\oplus[\delta'']$,
which correspond to even vertices on the
fusion graph $\rmD_{10}$ of $[\a^-_1]$, hence they
are both ambichiral. The statistical dimensions then
allow only $[\delta'']=[\a^+_2]$ and similarly
$[\delta']=[\a^-_2]$. Having now determined
$b^\pm_{\tau,j}=\langle\tau,\a^\pm_j\rangle$
for each $j$ and $\tau\in\MXMo$ we can now read off
the mass matrix $Z$ from \erf{Z=bb} and find that
it is the $\rmE_7$ one of \cite{CIZ1}. We can also
easily draw the simultaneous fusion graph of
$[\a^+_1]$ and $[\a^-_1]$ in the entire $M$-$M$
fusion rule algebra and we present it in
Fig.\ \ref{E7}.
%
%
Again, we encircled even vertices by small and ambichiral
(``marked'') vertices by large circles.

It is instructive to determine the canonical endomorphism
sector $[\can]$. From
\[ \langle \a^+_1 \a^-_1, \a^+_1 \a^-_1 \rangle =
\langle \a^+_1 \a^+_1, \a^-_1 \a^-_1 \rangle =
Z_{0,0}+Z_{0,2}+Z_{2,0}+Z_{2,2}=1 \]
we conclude that $[\eta]=[\a^+_1 \a^-_1]$ is an
irreducible sector which is a subsector of $[\can]$ since
$\langle\a^+_1\a^-_1,\can\rangle=
\langle\la_1\la_1,\canr\rangle=1$
by Frobenius reciprocity. Similarly we find
$\langle\a^\pm_8,\can\rangle=\langle\la_8,\canr\rangle=1$
which implies that $[\delta]$ is a subsector of $[\can]$.
As $\langle\can,\can\rangle=\langle\canr,\canr\rangle=3$
by \cite[Lemma 3.16]{BE3}, we conclude
$[\can]=[\id]\oplus[\eta]\oplus[\delta]$.

\subsection{A-D-E and representations of the Verlinde
fusion rules}

We have realized all $\SUz$ modular invariants from
subfactors. All canonical endomorphisms of these
subfactors have only subsectors $[\la_j]$ with $j$ even.
Therefore \erf{NMmain} transfers the two-coloring of
the $\SUz$ sectors to the $M$-$N$ sectors: Set
the colour of an $M$-$N$ sector $[\co a]$ to be
0 (respectively 1) whenever it is a subsector of
$[\iota\la_j]$ with $j$ even (respectively odd).
Consequently the matrix $G_1$ is the adjacency matrix of
a bi-colourable graph. Moreover, $G_1$ is irreducible
(i.e.\ the graph is connected) since $\la_1$
generates the whole $N$-$N$ system. We also have
$\|G_1\|=d_1<2$. Hence $G_1$ must be one of
the A-D-E cases (see e.g.\ \cite{GHJ}).
As Theorem \ref{adjacency} forces the diagonal
entries $Z_{j,j}$ of the modular invariant mass matrix
to be given as the multiplicities of the eigenvalues
$\chi_j(1)$ of $G_1$, our results explain why they
happen to be the multiplicities of the Coxeter exponents
of A-D-E Dynkin diagrams. We summarize several
data about the sector systems for the $\SUz$ modular
invariants in Table \ref{ADEtable}.
\begin{table}[htb]
\begin{center}
  \begin{tabular}{|c|c|c|c|c|c|c|c|} \hline &&&&&&&\\[-.9em]
Invariant $\leftrightarrow G_1$ & Level $k$ & $\# \MXM$
 & $\#\MXN$ & $\#\MXMpm$ & $\#\MXMo$ & $\Gamma^\pm_{0,1}$
 & $\Gamma^\pm_{\tau_{{\rm{gen}}},0}$
  \\[-.9em]&&&&&&&\\ \hline\hline &&&&&&&\\[-.9em]
$\rmA_\ell \,,\,\,\, \ell \ge 2$ & $\ell -1$ & $\ell $
 & $\ell$ & $\ell$ & $\ell$ & $\rmA_\ell$ & $\rmA_\ell$
 \\ &&&&&&&\\[-.9em]
$\rmD_{2\ell}\,,\,\,\,\ell \ge 2$ & $4\ell-4$ & $4\ell$
 & $2\ell$ & $2\ell$ & $\ell+1$ & $\rmD_{2\ell}$
 & $\rmD_{2\ell}^{{\rm{even}}}$
 \\ &&&&&&&\\[-.9em]
$\rmD_{2\ell+1}\,,\,\,\,\ell \ge 2$ & $4\ell-2$ & $4\ell-1$
 & $2\ell+1$ & $4\ell-1$ & $4\ell-1$ & $\rmA_{4\ell-1}$
 & $\rmA_{4\ell-1}$
 \\ &&&&&&&\\[-.9em]
$\rmE_6$ & $10$ & $12$ & $6$ & $6$ & $3$ & $\rmE_6$ & $\rmA_3$
 \\ &&&&&&&\\[-.9em]
$\rmE_7$ & $16$ & $17$ & $7$ & $10$ & $6$ & $\rmD_{10}$
 & $\rmD_{10}^{{\rm{even}}}$
 \\ &&&&&&&\\[-.9em]
$\rmE_8$ & $28$ & $32$ & $8$ & $8$ & $2$ & $\rmE_8$
 & $\rmA_4^{{\rm{even}}}$
 \\[-.8em] &&&&&&&\\
\hline \multicolumn3c {} \\[.05em] \end{tabular}
\end{center}
\caption{The A-D-E classification of $\SUz$ modular invariants}
\label{ADEtable}
\end{table}
The last column has the following meaning. We chose
an element $\tau_{{\rm{gen}}}\in\MXMo$ such that
$[\tau_{{\rm{gen}}}]$ is a subsector of $[\a^+_j]$ for
the smallest possible $j\ge1$. This element turns out
to generate the whole ambichiral system. For example,
in the $\rmE_7$ case we take $\tau_{{\rm{gen}}}=\a^+_2$.
The (adjacency matrix of the) fusion graph of
$[\tau_{{\rm{gen}}}]$ in the ambichiral system is
given in the last column.

Let us finally explain how the representation
$\varrho\circ\Phi$ which arises from left multiplication
of $M$-$M$ sectors on the $M$-$N$ sectors is
related to a fusion rule algebra
for (some) type \nolinebreak I invariants.
Let $V_1$ be the adjacency matrix of
one of the Dynkin diagrams. Then there is a
unitary matrix which dia\-gona\-lizes $V_1$,
i.e.\ $\psi^* V_1 \psi$ is the diagonal matrix
giving the eigenvalues corresponding to the
Coxeter exponents. In fact, Di Francesco and Zuber
\cite{DZ1,DZ2} built up a whole family of matrices
$V_\la$ with non-negative integer entries
($\la$ running over the spins for the time being), 
dia\-go\-na\-lized simultaneously by $\psi$ and providing
a representation of the Verlinde fusion rules,
$V_\la V_\mu = \sum_\nu N_{\la,\mu}^\nu V_\nu$.
Among the column vectors $\psi_m$,
$m$ labelling the eigenvalues including multiplicities
of the diagram at hand,
there is necessarily a Perron-Frobenius eigenvector
$\psi_0$ of $V_1$ with only strictly positive entries:
$\psi_{a,0}>0$ for all vertices $a$ of the diagram.
It turned out, actually first noticed in \cite{Pa},
that for $\rmD_{{\rm{even}}}$, $\rmE_6$ and $\rmE_8$,
which label the type \nolinebreak I modular invariants,
it was possible to choose\footnote{The matrix
$\psi$ is determined up to a rotation in each
multiplicity space of the eigenvalues (exponents).
So it is only $\rmD_{{\rm{even}}}$ where one needs
to make a choice to produce non-negative integers.}
$\psi$ such that also all $\psi_{0,m}>0$, here $a=0$ refers
to the extremal vertex, and that it has a
remarkable property: Plugged in a Verlinde type formula,
\begin{equation}
\label{verlinde'}
N_{a,b}^c=\sum_m \frac{\psi_{a,m}}{\psi_{0,m}}
\psi_{b,m} \psi_{c,m}^* \,,
\end{equation}
it yields non-negative integers $N_{a,b}^c$ which
could be interpreted as structure constants of a
fusion algebra, the ``graph algebra''.
This procedure worked analogously for the graphs
Di Francesco and Zuber \cite{DZ1,DZ2} associated to
some $\SUn$ type \nolinebreak I modular invariants
essentially by matching the spectra with the diagonal
entries of the mass matrices, whereas for
type \nolinebreak II invariants, in particular
$\rmD_{{\rm{odd}}}$ and $\rmE_7$ for $\SUz$,
it did not work. For instance, for $\rmE_7$ there
appeared some negative structure constants.

These observations find a natural explanation
in our setting. The graphs Di Fran\-cesco and Zuber
associated empirically to modular invariants are
recognized as the fusion graphs of $[\a^+_\la]$
obtained by multiplication from the left on the
$M$-$N$ sectors (or, equivalently, from the right
on $N$-$M$ sectors), i.e.\ $V_\la=G_\la$.
A priori, there is no reason why a matrix $\psi$
which diagonalizes the adjacency matrix of the graph(s)
should produce non-negative integer structure constants
because the $N$-$M$ morphisms alone do not form a fusion
algebra on their own:
You cannot multiply two $N$-$M$ morphisms,
and there is no identity. However, whenever the
chiral locality condition holds, then there is a
canonical bijection between the $N$-$M$ system and
either chiral induced system \cite[Lemma 4.1]{BE3}:
Any $N$-$M$ sector $[a]$, $a\in\NXM$, is of the form
$[a]=[\co\iota\beta]$, where either $\beta\in\MXMp$
or $\beta\in\MXMm$. This implies that, in the
notation of Subsect.\ \ref{refurex}, we have equality of
matrices $V_\nu=G_\nu=\Gamma^+_{0,\nu}=\Gamma^-_{0,\nu}$.
Recall that chiral locality implies by
Proposition \ref{asreci} that
$b^+_{\tau\la}=b^-_{\tau,\la}=b_{\tau,\la}$,
with restriction coefficients 
$b_{\tau,\la}=\langle\co\iota\tau\iota,\la\rangle$,
and that then the modular invariant is of
type \nolinebreak I:
$Z_{\la,\mu}=\sum_\tau b_{\tau,\la} b_{\tau,\mu}.$
In fact, we read off from Theorem \ref{adjacency}
that the eigenvalue $\chi_\la(\nu)$ of
$G_\nu=\Gamma^\pm_{0,\nu}$ appears with multiplicity
$Z_{\la,\la}=\sum_\tau b_{\tau,\la}^2$.
Now let $N_\beta$ be the fusion matrix of
$\beta\in\MXMp$ in the chiral system, i.e.\
$(N_\beta)_{\beta',\beta''}=N_{\beta',\beta}^{\beta''}=
\langle \beta'\beta,\beta'' \rangle$,
$\beta',\beta''\in\MXMp$. Then we have
$\Gamma^+_{0,\nu}=\sum_\beta
\langle \beta,\a^+_\nu \rangle N_\beta$.
Consequently, as long as the chiral system is
commutative\footnote{The ``first'' example of a
non-commutative chiral system is the
type \nolinebreak I invariant
coming from the conformal inclusion
${{\it{SU}}}(4)_4\in{{\it{SO}}}(15)_1$
\cite{X1,BE2}. In fact, that there are difficulties
to obtain non-negativity of structure constants from
a Verlinde type formula was noticed in \cite{PZ}.
A general analysis taking care of non-commutative
chiral systems as well as a discussion of
``marked vertices'' can be found in \cite{BE2}.},
there is always a unitary matrix $\psi$ which
diagonalizes the fusion matrices $N_\beta$ simultaneously,
and in turn their linear combinations $\Gamma^+_{0,\nu}$.
Evaluation of the zero-component of
$N_\beta \psi_m = \gamma_m(\beta) \psi_m$,
with $\gamma_m(\beta)$ some eigenvalue, yields
$\psi_{\beta,m}=\gamma_m(\beta) \psi_{0,m}$,
hence vanishing $\psi_{0,m}$ would contradict
unitarity of $\psi$, and thus one can choose
$\psi_{0,m}>0$. (See e.g.\ \cite{Kw} or
\cite[Sect.\ 8.7]{EK3} for such computations.)
Consequently the eigenvalues are given as
$\gamma_m(\beta)=\psi_{\beta,m}/ \psi_{0,m}$,
so that the structure constants are in fact
given by \erf{verlinde'}, using the bijection
$\NXM\ni a\leftrightarrow\beta\in\MXMp$.

Type \nolinebreak II modular invariants necessarily violate
the chiral locality condition, and without chiral locality
the bijection between $N$-$M$ system and the chiral
systems in general breaks down. For $\SUz$ this can
nicely be seen in Table \ref{ADEtable}: For the
$\rmD_{{\rm{odd}}}$ invariants, $G_1=\rmD_{2\ell+1}$,
we see that $G_1$ is in fact different from
$\Gamma^\pm_{0,1}=\rmA_{4\ell-1}$. Similarly we have
$\Gamma^\pm_{0,1}=\rmD_{10}$ for $G_1=\rmE_7$.

\section{More examples}
\label{morex}

\subsection{Conformal inclusions of $\SUd$}

We discuss two more examples arising from conformal
inclusions of $\SUd$.
Combining the methods and results in \cite{BE1,BE2,BE3},
and \cite{BEK1}, we can compute examples along the lines of
\cite[Sect.\ 6]{BE3}.

The first example is the conformal inclusion
$\SUd_3\subset{{\it{SO}}}(8)_1$. The associated
modular invariant is
$$ Z_{{\cal D}^{(6)}}=
|\chi_{(0,0)} + \chi_{(3,0)} + \chi_{(3,3)}|^2 +
3|\chi_{(2,1)}|^2,$$
and was labelled by the orbifold graph ${\cal D}^{(6)}$.
In fact, this conformal inclusion can also be treated
as a local simple current extension, similar to the
$\rmD_4$ case for $\SUz$. The chiral systems $\MXMpm$,
i.e.\ the images of $\a^\pm$-induction
were determined in \cite[Fig.\ 10]{BE2}. 
Here we describe the structure of the full system $\MXM$.
By \cite[Prop.\ 5.1]{BE3} and \cite[Thm.\ 5.10]{BEK1}, we
know that the intersection $\MXMo$ are the ``marked vertices''
of \cite[Fig.\ 10]{BE2}. From \cite[Cor.\ 6.10]{BEK1},
we learn that the number of the $M$-$M$ morphisms is 18.
We next note that the dual canonical endomorphism
$\theta$ decomposes into three mutually
inequivalent irreducible $N$-$N$ morphisms of dimension 1 
and the fusion rules of these three morphisms are given by
the group $\bbZ_3$.  This implies that the
canonical endomorphism $\gamma$ also decomposes into three mutually
inequivalent irreducible $M$-$M$ morphisms of dimension 1 
and the fusion rule of these three morphisms is again given by
the group $\bbZ_3$.
We compute
\[ \bearll
\lan \a^+_{(1,0)}\a^-_{(1,0)}, \a^+_{(1,0)}\a^-_{(1,0)} \ran
&= \lan \bar \a^+_{(1,0)} \a^+_{(1,0)}, \a^-_{(1,0)}
   \bar \a^-_{(1,0)} \ran
   =\lan \a^+_{(1,1)} \a^+_{(1,0)},
   \a^-_{(1,0)} \a^-_{(1,1)}\ran \\[.4em]
&=Z_{(0,0),(0,0)}+Z_{(2,1),(0,0)}+Z_{(0,0),(2,1)}+Z_{(2,1),(2,1)}=4\,,
\eear \]
and we similarly have
\[
\lan \a^+_{(1,1)}\a^-_{(1,1)}, \a^+_{(1,1)}\a^-_{(1,1)} \ran=
\lan \a^+_{(1,1)}\a^-_{(1,0)}, \a^+_{(1,1)}\a^-_{(1,0)} \ran=
\lan \a^+_{(1,0)}\a^-_{(1,1)}, \a^+_{(1,0)}\a^-_{(1,1)} \ran=4 \,.
\]
We have two ways of expressing the number four as a sum of
squares, $4=1+1+1+1=2^2$.  This means that
$\a^+_{(1,0)}\a^-_{(1,0)}$ decomposes either into
four mutually inequivalent irreducible $M$-$M$ morphisms
or into two copies of one irreducible $M$-$M$ morphism.
The statistical dimension of $\a^+_{(1,0)}\a^-_{(1,0)}$ is 4, and this
means that $\a^+_{(1,0)}\a^-_{(1,0)}$
decomposes into four mutually
inequivalent irreducible $M$-$M$ morphisms of dimension 1 each or
into two copies of one irreducible $M$-$M$ morphism of dimension 2.
In either case, the square sum of the statistical dimensions of the
irreducible morphisms appearing in the decomposition is 4.
The same holds for $\a^+_{(1,1)}\a^-_{(1,1)}$, 
$\a^+_{(1,1)}\a^-_{(1,0)}$ and $\a^+_{(1,0)}\a^-_{(1,1)}$.
We have
$$\lan \a^+_{(1,0)}\a^-_{(1,0)}, \a^+_{(1,1)}\a^-_{(1,1)} \ran=0,$$
which implies the morphisms appearing in the irreducible decomposition
of $\a^+_{(1,0)}\a^-_{(1,0)}$ and those in the decomposition of
$\a^+_{(1,1)}\a^-_{(1,1)}$ are disjoint, and the disjointness
holds for any distinct two of the four endomorphisms
$\a^+_{(1,0)}\a^-_{(1,0)}$, $\a^+_{(1,1)}\a^-_{(1,1)}$,
$\a^+_{(1,1)}\a^-_{(1,0)}$ and $\a^+_{(1,0)}\a^-_{(1,1)}$.
The morphisms appearing in the decompositions of these four $M$-$M$
morphisms are also disjoint with those in $\MXMp \cup \MXMm$.
The contribution of the morphisms in $\MXMp \cup \MXMm$ to
the global index $w=36$ is $12+12-4=20$.  The contribution of
the morphisms appearing in the irreducible decompositions of
$\a^+_{(1,0)}\a^-_{(1,0)}$, $\a^+_{(1,1)}\a^-_{(1,1)}$,
$\a^+_{(1,1)}\a^-_{(1,0)}$ and $\a^+_{(1,0)}\a^-_{(1,1)}$
is $4\times 4=16$ and this means that these morphisms,
together with those in $\MXMp \cup \MXMm$, give the
entire system $\MXM$.  Since the total number of the
$M$-$M$ morphisms in $\MXM$ is 18, we conclude that two of
$\a^+_{(1,0)}\a^-_{(1,0)}$, $\a^+_{(1,1)}\a^-_{(1,1)}$,
$\a^+_{(1,1)}\a^-_{(1,0)}$ and $\a^+_{(1,0)}\a^-_{(1,1)}$
decompose into four mutually
inequivalent irreducible $M$-$M$ morphisms of dimension 1 each 
and the other two decompose
into two copies of one irreducible $M$-$M$ morphism of dimension 2,
respectively.

We next compute
$$
\lan \a^+_{(1,0)}\a^-_{(1,1)}, \gamma \ran=
\lan \a^+_{(1,0)}\a^+_{(1,1)}, \gamma \ran=
\lan \a^+_{(0,0)}, \gamma \ran+
\lan \a^+_{(2,1)}, \gamma \ran
= 1 \,,
$$
since $\lan \la_{(2,1)}, \canr \ran=0$ and we have
in general
$\lan\a^\pm_\la,\can\ran=\lan\la,\canr\ran$
by Frobenius reciprocity. Hence we conclude that
$\lan \a^+_{(1,0)}\a^-_{(1,1)}, \gamma \ran=1$.
This is impossible if $\a^+_{(1,0)}\a^-_{(1,1)}$
decomposes into two copies of one irreducible
$M$-$M$ morphism of dimension 2, thus
$\a^+_{(1,0)}\a^-_{(1,1)}$ must
decompose into four mutually inequivalent
irreducible $M$-$M$ morphisms of dimension 1 each.
The same conclusion holds for $\a^+_{(1,1)}\a^-_{(1,0)}$.
Hence we know that both $\a^+_{(1,0)}\a^-_{(1,0)}$ and
$\a^+_{(1,1)}\a^-_{(1,1)}$ decompose into two copies of
one irreducible $M$-$M$ morphism of dimension 2,
respectively.

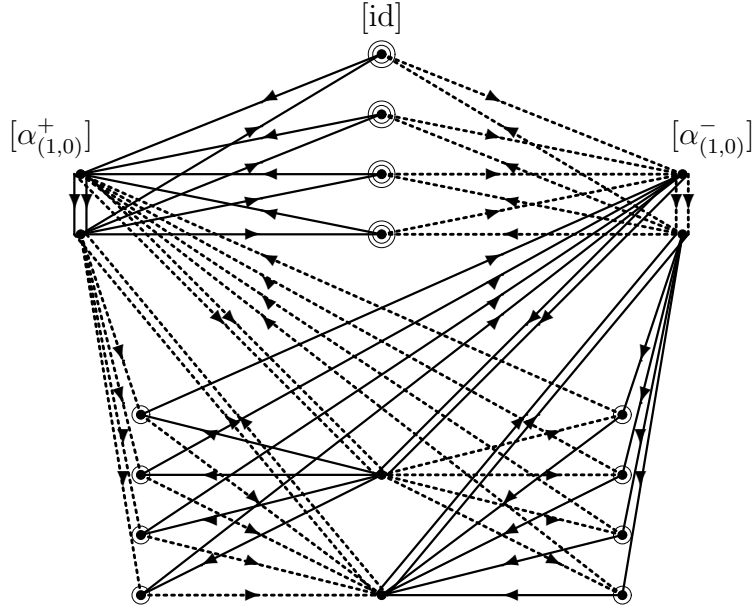
\begin{figure}[htb]
\unitlength 0.8mm
\begin{center}
\begin{picture}(120,110)
\thinlines 
\multiput(10,70)(0,10){2}{\circle*{1.5}}
\multiput(20,10)(0,10){4}{\circle*{1.5}}
\multiput(60,10)(0,20){2}{\circle*{1.5}}
\multiput(60,70)(0,10){4}{\circle*{1.5}}
\multiput(100,10)(0,10){4}{\circle*{1.5}}
\multiput(110,70)(0,10){2}{\circle*{1.5}}
\multiput(20,10)(0,10){4}{\circle{3}} 
\multiput(60,70)(0,10){4}{\circle{3}} 
\multiput(100,10)(0,10){4}{\circle{3}} 
\multiput(60,70)(0,10){4}{\circle{4.5}} 
\thicklines 
\path(60,100)(10,80)
\path(60,90)(10,80)
\path(60,80)(10,80)
\path(60,70)(10,80)
\path(10,70)(60,100)
\path(10,70)(60,90)
\path(10,70)(60,80)
\path(10,70)(60,70)
\path(11,80)(11,70)
\path(9,80)(9,70)
\path(20,40)(110,80)
\path(20,30)(110,80)
\path(20,20)(110,80)
\path(20,10)(110,80)
\path(60,30)(20,40)
\path(60,30)(20,30)
\path(60,30)(20,20)
\path(60,30)(20,10)
\path(111,80)(61,30)
\path(109,80)(59,30)
\path(100,40)(60,10)
\path(100,30)(60,10)
\path(100,20)(60,10)
\path(100,10)(60,10)
\path(110,70)(100,40)
\path(110,70)(100,30)
\path(110,70)(100,20)
\path(110,70)(100,10)
\path(111,70)(61,10)
\path(109,70)(59,10)
\dottedline{1.3}(60,100)(110,80)
\dottedline{1.3}(60,90)(110,80)
\dottedline{1.3}(60,80)(110,80)
\dottedline{1.3}(60,70)(110,80)
\dottedline{1.3}(110,70)(60,100)
\dottedline{1.3}(110,70)(60,90)
\dottedline{1.3}(110,70)(60,80)
\dottedline{1.3}(110,70)(60,70)
\dottedline{1.3}(109,80)(109,70)
\dottedline{1.3}(111,80)(111,70)
\dottedline{1.3}(100,40)(10,80)
\dottedline{1.3}(100,30)(10,80)
\dottedline{1.3}(100,20)(10,80)
\dottedline{1.3}(100,10)(10,80)
\dottedline{1.3}(60,30)(100,40)
\dottedline{1.3}(60,30)(100,30)
\dottedline{1.3}(60,30)(100,20)
\dottedline{1.3}(60,30)(100,10)
\dottedline{1.3}(9,80)(59,30)
\dottedline{1.3}(11,80)(61,30)
\dottedline{1.3}(20,40)(60,10)
\dottedline{1.3}(20,30)(60,10)
\dottedline{1.3}(20,20)(60,10)
\dottedline{1.3}(20,10)(60,10)
\dottedline{1.3}(10,70)(20,40)
\dottedline{1.3}(10,70)(20,30)
\dottedline{1.3}(10,70)(20,20)
\dottedline{1.3}(10,70)(20,10)
\dottedline{1.3}(9,70)(59,10)
\dottedline{1.3}(11,70)(61,10)
\put(60,106){\makebox(0,0){$[\id]$}}
\put(5,86){\makebox(0,0){$[\alpha^+_{(1,0)}]$}}
\put(115,86){\makebox(0,0){$[\alpha^-_{(1,0)}]$}}
\put(40,92){\vector(-2,-1){0}}
\put(40,86){\vector(-4,-1){0}}
\put(40,80){\vector(-1,0){0}}
\put(40,74){\vector(-4,1){0}}
\put(40,88){\vector(2,1){0}}
\put(40,82){\vector(2,1){0}}
\put(40,76){\vector(4,1){0}}
\put(40,70){\vector(1,0){0}}
\put(9,74){\vector(0,-1){0}}
\put(11,74){\vector(0,-1){0}}
\put(80,66.7){\vector(2,1){0}}
\put(80,63.3){\vector(1,1){0}}
\put(80,60){\vector(1,1){0}}
\put(80,56.7){\vector(1,1){0}}
\put(30,37.5){\vector(-3,1){0}}
\put(30,30){\vector(-1,0){0}}
\put(30,22.5){\vector(-1,0){0}}
\put(30,15){\vector(-2,-1){0}}
\put(86,55){\vector(-1,-1){0}}
\put(84,55){\vector(-1,-1){0}}
\put(80,25){\vector(-1,-1){0}}
\put(80,20){\vector(-2,-1){0}}
\put(80,15){\vector(-3,-1){0}}
\put(80,10){\vector(-1,0){0}}
\put(103,49){\vector(-1,-3){0}}
\put(103,42){\vector(-1,-3){0}}
\put(103,35){\vector(0,-1){0}}
\put(103,28){\vector(0,-1){0}}
\put(86,40){\vector(1,1){0}}
\put(84,40){\vector(1,1){0}}
\put(80,92){\vector(2,-1){0}}
\put(80,86){\vector(4,-1){0}}
\put(80,80){\vector(1,0){0}}
\put(80,74){\vector(4,1){0}}
\put(80,88){\vector(-2,1){0}}
\put(80,82){\vector(-2,1){0}}
\put(80,76){\vector(-4,1){0}}
\put(80,70){\vector(-1,0){0}}
\put(111,74){\vector(0,-1){0}}
\put(109,74){\vector(0,-1){0}}
\put(40,66.7){\vector(-2,1){0}}
\put(40,63.3){\vector(-1,1){0}}
\put(40,60){\vector(-1,1){0}}
\put(40,56.7){\vector(-1,1){0}}
\put(90,37.5){\vector(3,1){0}}
\put(90,30){\vector(1,0){0}}
\put(90,22.5){\vector(1,0){0}}
\put(90,15){\vector(2,-1){0}}
\put(34,55){\vector(1,-1){0}}
\put(36,55){\vector(1,-1){0}}
\put(40,25){\vector(1,-1){0}}
\put(40,20){\vector(2,-1){0}}
\put(40,15){\vector(3,-1){0}}
\put(40,10){\vector(1,0){0}}
\put(17,49){\vector(1,-3){0}}
\put(17,42){\vector(1,-3){0}}
\put(17,35){\vector(0,-1){0}}
\put(17,28){\vector(0,-1){0}}
\put(34,40){\vector(-1,1){0}}
\put(36,40){\vector(-1,1){0}}
\end{picture}
\caption{$SU(3)_3\subset SO(8)_1$, ${\cal D}^{(6)}$:
Fusion graph of $[\a^+_{(1,0)}]$ and $[\a^-_{(1,0)}]$}
\label{D(6)}
\end{center}
\end{figure}

With these informations, it is easy to determine the
simultaneous fusion graph of $[\a^+_{(1,0)}]$ (straight lines)
and $[\alpha^-_{(1,0)}]$ (dotted lines) as in Fig.\  \ref{D(6)}.
We have encircled the marked vertices by big circles
and the colour zero vertices by small vertices.
(Because the vacuum block has only colour zero
contributions, the full system inherits the three
colouring of the $\SUd_3$ system $\NXN$ here.)
As the modular invariant contains an entry 3 we conclude
by \cite[Cor.\ 6.9]{BEK1}  that the entire $M$-$M$ fusion
rule algebra is non-commutative.  The colour zero part has
12 vertices which all correspond to simple sectors. Therefore
they form a closed subsystem corresponding to a group.
This group must contain a $\bbZ_2\times\bbZ_2$ subgroup
corresponding to the ${{\it{SO}}}(8)_1$ fusion rules
of the marked vertices. Note that any $M$-$M$ sector of
non-zero colour is a product $[\a^\pm_{(1,0)}][\beta]$
or $[\a^\pm_{(1,1)}][\beta]$
with $\beta\in\MXM$ a colour zero morphism. Since
$[\a^\pm_{(1,0)}]$ and $[\a^\pm_{(1,1)}]$ commute with
each $M$-$M$ sector by \cite[Lemma 3.20]{BE3}, they will
be scalars in any irreducible representation of the
$M$-$M$ fusion rules. Consequently, the representation
$\pi_{(2,1),(2,1)}$ of dimension $Z_{(2,1),(2,1)}=3$
will remain irreducible upon restriction to the
group of colour zero sectors. Therefore its group
dual is forced to consist of one 3-dimensional and
3 scalar representations, and in turn we identify
the group of colour zero sectors to be the tetrahedral
group $A_4=(\bbZ_2\times\bbZ_2)\rtimes\bbZ_3$.

The next example is a conformal inclusion
$\SUd_5\subset{{\it{SU}}}(6)_1$. The associated modular
invariant, labelled as ${\cal E}^{(8)}$, is given by
\begin{eqnarray*}
Z_{{\cal E}^{(8)}}&=& |\chi_{(0,0)} + \chi_{(4,2)}|^2 +
|\chi_{(2,0)} + \chi_{(5,3)}|^2
+ |\chi_{(2,2)} + \chi_{(5,2)}|^2 \\
&&\qquad\qquad + \,|\chi_{(3,0)} + \chi_{(3,3)}|^2 +
|\chi_{(3,1)} + \chi_{(5,5)} |^2
+ |\chi_{(3,2)} + \chi_{(5,0)}|^2\,.
\end{eqnarray*}
The chiral systems $\MXMpm$ were also determined in
\cite[Subsect.\ 2.3]{BE2} and here we describe the
structure of the full system $\MXM$.

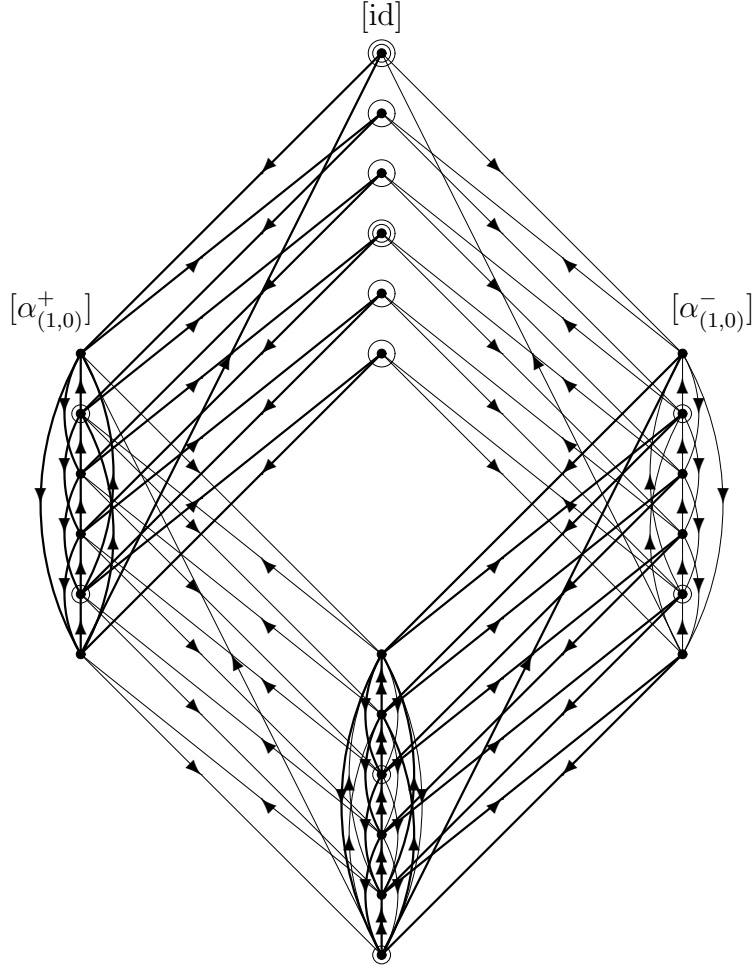
\begin{figure}[htb]
\unitlength 0.8mm
\begin{center}
\begin{picture}(120,170)
\thinlines 
\multiput(10,60)(0,10){6}{\circle*{1.5}}
\multiput(60,10)(0,10){6}{\circle*{1.5}}
\multiput(60,110)(0,10){6}{\circle*{1.5}}
\multiput(110,60)(0,10){6}{\circle*{1.5}}
\multiput(10,70)(0,30){2}{\circle{3}} 
\multiput(60,10)(0,30){2}{\circle{3}} 
\multiput(60,130)(0,30){2}{\circle{3}} 
\multiput(110,70)(0,30){2}{\circle{3}} 
\multiput(60,110)(0,10){6}{\circle{4.5}} 
\thicklines
\path(60,160)(10,110)
\path(60,150)(10,100)
\path(60,140)(10,90)
\path(60,130)(10,80)
\path(60,120)(10,70)
\path(60,110)(10,60)
\path(10,110)(60,150)
\path(10,100)(60,140)
\path(10,90)(60,130)
\path(10,80)(60,120)
\path(10,70)(60,110)
\path(10,60)(60,160)
\path(10,60)(10,110)
\put(27.3,100){\arc{40}{2.62}{3.66}}
\put(27.3,90){\arc{40}{2.62}{3.66}}
\put(27.3,80){\arc{40}{2.62}{3.66}}
\put(27.3,70){\arc{40}{2.62}{3.66}}
\put(-24.6,80){\arc{80}{5.76}{6.80}}
\put(-24.6,90){\arc{80}{5.76}{6.80}}
\put(53.3,85){\arc{100}{2.62}{3.66}}
\put(40,140){\vector(-1,-1){0}}
\put(40,130){\vector(-1,-1){0}}
\put(40,120){\vector(-1,-1){0}}
\put(40,110){\vector(-1,-1){0}}
\put(40,100){\vector(-1,-1){0}}
\put(40,90){\vector(-1,-1){0}}
\put(30,126){\vector(1,1){0}}
\put(30,116){\vector(1,1){0}}
\put(30,106){\vector(1,1){0}}
\put(30,96){\vector(1,1){0}}
\put(30,86){\vector(1,1){0}}
\put(35,110){\vector(1,2){0}}
\put(10,66){\vector(0,1){0}}
\put(10,76){\vector(0,1){0}}
\put(10,86){\vector(0,1){0}}
\put(10,96){\vector(0,1){0}}
\put(10,106){\vector(0,1){0}}
\put(7.3,100){\vector(0,-1){0}}
\put(7.3,90){\vector(0,-1){0}}
\put(7.3,80){\vector(0,-1){0}}
\put(7.3,70){\vector(0,-1){0}}
\put(3.3,85){\vector(0,-1){0}}
\put(15.4,80){\vector(0,1){0}}
\put(15.4,90){\vector(0,1){0}}
\put(60,166){\makebox(0,0){$[\id]$}}
\put(5,117){\makebox(0,0){$[\alpha^+_{(1,0)}]$}}
\put(115,117){\makebox(0,0){$[\alpha^-_{(1,0)}]$}}
\path(110,110)(60,60)
\path(110,100)(60,50)
\path(110,90)(60,40)
\path(110,80)(60,30)
\path(110,70)(60,20)
\path(110,60)(60,10)
\path(60,60)(110,100)
\path(60,50)(110,90)
\path(60,40)(110,80)
\path(60,30)(110,70)
\path(60,20)(110,60)
\path(60,10)(110,110)
\path(60,10)(60,60)
\put(77.3,50){\arc{40}{2.62}{3.66}}
\put(77.3,40){\arc{40}{2.62}{3.66}}
\put(77.3,30){\arc{40}{2.62}{3.66}}
\put(77.3,20){\arc{40}{2.62}{3.66}}
\put(25.4,30){\arc{80}{5.76}{6.80}}
\put(25.4,40){\arc{80}{5.76}{6.80}}
\put(103.3,35){\arc{100}{2.62}{3.66}}
\put(90,90){\vector(-1,-1){0}}
\put(90,80){\vector(-1,-1){0}}
\put(90,70){\vector(-1,-1){0}}
\put(90,60){\vector(-1,-1){0}}
\put(90,50){\vector(-1,-1){0}}
\put(90,40){\vector(-1,-1){0}}
\put(80,76){\vector(1,1){0}}
\put(80,66){\vector(1,1){0}}
\put(80,56){\vector(1,1){0}}
\put(80,46){\vector(1,1){0}}
\put(80,36){\vector(1,1){0}}
\put(85,60){\vector(1,2){0}}
\put(60,16){\vector(0,1){0}}
\put(60,26){\vector(0,1){0}}
\put(60,36){\vector(0,1){0}}
\put(60,46){\vector(0,1){0}}
\put(60,56){\vector(0,1){0}}
\put(57.3,50){\vector(0,-1){0}}
\put(57.3,40){\vector(0,-1){0}}
\put(57.3,30){\vector(0,-1){0}}
\put(57.3,20){\vector(0,-1){0}}
\put(53.3,35){\vector(0,-1){0}}
\put(65.4,30){\vector(0,1){0}}
\put(65.4,40){\vector(0,1){0}}
\thinlines
\path(60,160)(110,110)
\path(60,150)(110,100)
\path(60,140)(110,90)
\path(60,130)(110,80)
\path(60,120)(110,70)
\path(60,110)(110,60)
\path(110,110)(60,150)
\path(110,100)(60,140)
\path(110,90)(60,130)
\path(110,80)(60,120)
\path(110,70)(60,110)
\path(110,60)(60,160)
\path(110,60)(110,110)
\put(92.7,100){\arc{40}{5.76}{6.80}}
\put(92.7,90){\arc{40}{5.76}{6.80}}
\put(92.7,80){\arc{40}{5.76}{6.80}}
\put(92.7,70){\arc{40}{5.76}{6.80}}
\put(144.6,80){\arc{80}{2.62}{3.66}}
\put(144.6,90){\arc{80}{2.62}{3.66}}
\put(66.7,85){\arc{100}{5.76}{6.80}}
\put(80,140){\vector(1,-1){0}}
\put(80,130){\vector(1,-1){0}}
\put(80,120){\vector(1,-1){0}}
\put(80,110){\vector(1,-1){0}}
\put(80,100){\vector(1,-1){0}}
\put(80,90){\vector(1,-1){0}}
\put(90,126){\vector(-1,1){0}}
\put(90,116){\vector(-1,1){0}}
\put(90,106){\vector(-1,1){0}}
\put(90,96){\vector(-1,1){0}}
\put(90,86){\vector(-1,1){0}}
\put(85,110){\vector(-1,2){0}}
\put(110,66){\vector(0,1){0}}
\put(110,76){\vector(0,1){0}}
\put(110,86){\vector(0,1){0}}
\put(110,96){\vector(0,1){0}}
\put(110,106){\vector(0,1){0}}
\put(112.7,100){\vector(0,-1){0}}
\put(112.7,90){\vector(0,-1){0}}
\put(112.7,80){\vector(0,-1){0}}
\put(112.7,70){\vector(0,-1){0}}
\put(116.7,85){\vector(0,-1){0}}
\put(104.6,80){\vector(0,1){0}}
\put(104.6,90){\vector(0,1){0}}
\path(10,110)(60,60)
\path(10,100)(60,50)
\path(10,90)(60,40)
\path(10,80)(60,30)
\path(10,70)(60,20)
\path(10,60)(60,10)
\path(60,60)(10,100)
\path(60,50)(10,90)
\path(60,40)(10,80)
\path(60,30)(10,70)
\path(60,20)(10,60)
\path(60,10)(10,110)
\put(42.7,50){\arc{40}{5.76}{6.80}}
\put(42.7,40){\arc{40}{5.76}{6.80}}
\put(42.7,30){\arc{40}{5.76}{6.80}}
\put(42.7,20){\arc{40}{5.76}{6.80}}
\put(94.6,30){\arc{80}{2.62}{3.66}}
\put(94.6,40){\arc{80}{2.62}{3.66}}
\put(16.7,35){\arc{100}{5.76}{6.80}}
\put(30,90){\vector(1,-1){0}}
\put(30,80){\vector(1,-1){0}}
\put(30,70){\vector(1,-1){0}}
\put(30,60){\vector(1,-1){0}}
\put(30,50){\vector(1,-1){0}}
\put(30,40){\vector(1,-1){0}}
\put(40,76){\vector(-1,1){0}}
\put(40,66){\vector(-1,1){0}}
\put(40,56){\vector(-1,1){0}}
\put(40,46){\vector(-1,1){0}}
\put(40,36){\vector(-1,1){0}}
\put(35,60){\vector(-1,2){0}}
\put(60,18){\vector(0,1){0}}
\put(60,28){\vector(0,1){0}}
\put(60,38){\vector(0,1){0}}
\put(60,48){\vector(0,1){0}}
\put(60,58){\vector(0,1){0}}
\put(62.7,50){\vector(0,-1){0}}
\put(62.7,40){\vector(0,-1){0}}
\put(62.7,30){\vector(0,-1){0}}
\put(62.7,20){\vector(0,-1){0}}
\put(66.7,35){\vector(0,-1){0}}
\put(54.6,30){\vector(0,1){0}}
\put(54.6,40){\vector(0,1){0}}
\end{picture}
\caption{$SU(3)_5\subset SU(6)_1$, ${\cal E}^{(8)}$:
Fusion graph of $[\a^+_{(1,0)}]$ and $[\a^-_{(1,0)}]$}
\label{E(8)}
\end{center}
\end{figure}

Again by \cite[Prop.\ 5.1]{BE3} and \cite[Thm.\ 5.10]{BEK1}, we
know that the intersection $\MXMo$ are the ``marked vertices''
of \cite[Fig.\ 11]{BE2} and this consists of six morphisms of
dimension 1.  From \cite[Cor.\ 6.10]{BEK1},
we learn that the number of the $M$-$M$ morphisms is 24.
Since $\MXMp\cup \MXMm$ has 18 morphisms, we need to
find 6 more morphisms.  We compute
\[ \bearll
\lan \a^+_{(1,0)}\a^-_{(1,0)},\a^+_{(1,0)}\a^-_{(1,0)}\ran
&=\lan \a^+_{(1,1)}\a^+_{(1,0)},\a^-_{(1,0)}\a^-_{(1,1)}\ran \\[.4em]
&=Z_{(0,0),(0,0)}+Z_{(0,0),(2,1)}+Z_{(2,1),(0,0)}+Z_{(2,1),(2,1)}=1\,,
\eear \]
which shows $\a^+_{(1,0)}\a^-_{(1,0)}$ is an irreducible $M$-$M$
morphism outside of $\MXMp\cup \MXMm$.  Similarly,
$\a^+_{(5,4)}\a^-_{(1,0)}$,
$\a^+_{(4,4)}\a^-_{(1,0)}$,
$\a^+_{(4,0)}\a^-_{(1,0)}$,
$\a^+_{(5,1)}\a^-_{(1,0)}$ and
$\a^+_{(1,1)}\a^-_{(1,0)}$ are also irreducible $M$-$M$
morphisms outside of $\MXMp\cup \MXMm$. Similarly we compute
\[ \bearll
\lan \a^+_{(1,0)}\a^-_{(1,0)},\a^+_{(5,4)}\a^-_{(1,0)}\ran
&=\lan \a^+_{(5,1)}\a^+_{(1,0)},\a^-_{(1,0)}\a^-_{(1,1)}\ran\\[.4em]
&=Z_{(5,2),(0,0)}+Z_{(5,2),(2,1)}+Z_{(4,0),(0,0)}+Z_{(4,0),(2,1)}=0\,,
\eear \]
which shows that the irreducible morphisms $\a^+_{(5,4)}\a^-_{(1,0)}$ and
$\a^+_{(5,4)}\a^-_{(1,0)}$ are not equivalent.  Similarly,
we know that all of $\a^+_{(1,0)}\a^-_{(1,0)}$,
$\a^+_{(5,4)}\a^-_{(1,0)}$,
$\a^+_{(4,4)}\a^-_{(1,0)}$,
$\a^+_{(4,0)}\a^-_{(1,0)}$,
$\a^+_{(5,1)}\a^-_{(1,0)}$ and
$\a^+_{(1,1)}\a^-_{(1,0)}$ are mutually inequivalent, and these
thus give the missing six irreducible morphisms in $\MXM$.
We also have 
\[ \bearll
\lan \a^+_{(5,4)}\a^-_{(1,0)},\a^+_{(1,0)}\a^-_{(5,4)}\ran
&=\lan \a^-_{(1,0)}\a^-_{(5,1)},\a^+_{(5,1)}\a^+_{(1,0)}\ran\\[.4em]
&=Z_{(5,2),(5,2)}+Z_{(5,2),(4,0)}+Z_{(4,0),(5,2)}+Z_{(4,0),(4,0)}=1\,,
\eear \]
which implies that the morphisms $\a^+_{(5,4)}\a^-_{(1,0)}$,
$\a^+_{(1,0)}\a^-_{(5,4)}$ are equivalent.  We similarly have
$[\a^+_{(4,4)}\a^-_{(1,0)}]=[\a^-_{(4,4)}\a^+_{(1,0)}]$,
$[\a^+_{(4,0)}\a^-_{(1,0)}]=[\a^-_{(4,0)}\a^+_{(1,0)}]$,
$[\a^+_{(5,1)}\a^-_{(1,0)}]=[\a^-_{(5,1)}\a^+_{(1,0)}]$,
$[\a^+_{(1,1)}\a^-_{(1,0)}]=[\a^-_{(1,1)}\a^+_{(1,0)}]$.
We then can compute the
simultaneous fusion graph of $[\alpha^+_{(1,0)}]$ (thick lines)
and $[\alpha^-_{(1,0)}]$ (thin lines) easily as in Fig.\ \ref{E(8)}.

\subsection{The trivial invariant from a non-trivial inclusion}
\label{Ising}

Here we give an example of a non-trivial inclusion
$N\subset M$ which however produces the trivial
modular invariant
$Z_{\la,\mu}=\langle\a^+_\la,\a^-_\mu\rangle=\del\la\mu$.
This is clearly only possible if the chiral locality
condition is violated because chiral locality implies
the formula $\langle\a^\pm_\la,\a^\pm_\mu\rangle
=\langle\canr\la,\mu\rangle$, as derived in
\cite[Thm.\ 3.9]{BE1}; hence, if $[\mu]$ is a non-trivial
subsector of $[\canr]$, then
$Z_{0,\mu}=\langle\id,\a^-_\mu\rangle
=\langle\a^-_\id,\a^-_\mu\rangle$ must be non-zero.
Consequently a ``local extension'' can only exist
if there exists a non-trivial mass matrix $Z$ commuting
with the S- and T-matrices arising from the braiding.

An example for such a non-local extension
is provided by the chiral conformal
Ising model in the algebraic formulation of
\cite{MaS,Bo1} where the local observable algebras
are realized as $\bbZ_2$ gauge invariant subalgebras
of fermionic algebras.
Let $\Gamma$ denote the complex conjugation
on $L^2(S^1;\bbC)$. The fermion algebra (in Araki's
self-dual CAR formalism \cite{A1}) is the unital
$C^*$-algebra generated by the image of a linear
map $f\mapsto\psi(f)$ subject to relations
\[ \psi(f)^*=\psi(\Gamma f) \,, \qquad
\psi(f)^*\psi(g)+\psi(g)\psi(f)^*
=\langle f,g \rangle \bfe \,, \qquad
f,g\in L^2(S^1) \,. \]
Any real isometry
$V=\Gamma V \Gamma\in B(L^2(S^1;\bbC))$ induces
a unital Bogoliubov endomorphism $\varrho_V$
of the fermion algebra via
$\varrho_V(\psi(f))=\psi(Vf)$. A simple example
is the outer $\bbZ_2$ gauge automorphism
$\varrho_{-\bfe}$ sending $\psi(f)$ to $-\psi(f)$.
The fermion algebra possesses a faithful irreducible
representation $\pi_{{\rm{NS}}}$ on the anti-symmetric
Fock space $\cF_-(P_{{\rm{NS}}}L^2(S^1;\bbC))$, where
$P_{{\rm{NS}}}$ is the Neveu-Schwarz polarization
used in \cite{Bo1}. We may and do now identify
$\pi_{{\rm{NS}}}(\psi(f))$ with $\psi(f)$.
After removing a ``point at
infinity'' from the circle $S^1$, a directed net
of factors $\{M(I)\}$ on the Fock space is obtained
by defining $M(I)$ to be the von Neumann algebra
generated by $\psi(f)$'s where
${{\rm{supp}}}(f)\subset I$ with $I$ varying
over the intervals on the punctured circle.
A net of subfactors is obtained by putting
$N(I)=M(I)^{\bbZ_2}$ where the $\bbZ_2$
action comes from the gauge automorphism
$\varrho_{-\bfe}$ which is still outer
on each $M(I)$. As usual, we denote the
associated $C^*$-algebras i.e.\ the norm closures
of the unions of all $N(I)$ respectively $M(I)$
by $\cN$ and $\cM$. The net $\{N(I)\}$ is local
and its restriction to the even Fock space,
defining to the vacuum representation $\pi_0$,
is Haag dual. The superselection sectors $[\id]$,
$[\eta]$, $[\sigma]$ of $\cN$ are those of the
Ising model can be realized by localized
endomorphisms $\id$, $\eta$,
$\sigma$ of $\cN$ which are restrictions of
Bogoliubov endomorphisms of $\cM$, and they
satisfy in fact the Ising fusion rules
$[\eta^2]=[\id]$, $[\eta][\sigma]=[\sigma]$
and $[\sigma^2]=[\id]\oplus[\eta]$, as
shown in \cite{Bo1}.
Fix an interval $I_\circ$ in which these
morphisms are localized. Note that then
the restrictions of these morphisms to the local
algebra $N(I_\circ)$ obey the same fusion
rules as sectors of this factor since
$\cN$ has a Haag dual subrepresentation.
The canonical endomorphism sector of the inclusion
$N=N(I_\circ)\subset M(I_\circ)=M$
is given by $[\canr]=[\id]\oplus[\eta]$
since the Fock representation $\pi_{{\rm{NS}}}$
decomposes into these superselection sectors
upon restriction to the gauge invariant
fermionic algebra $\cN$.

Now let us compute the $\alpha$-induced
endomorphisms explicitly. We here consider them
as endomorphisms of the entire algebra $\cM$
as in \cite{BE1} and we denote there restrictions
to a local algebra $M(I_\circ)$ by the same symbols
e.g.\ when we consider sectors of $M(I_\circ)$.
Choose a real function
$h=\Gamma h\in L^2(S^1;\bbC)$ with support in
$I_\circ$ and such that $\|h\|^2=2$. Recall from
\cite[Def.\ 3.7]{Bo1} that then
the localized endomorphism $\eta$ is the
restriction of the Bogoliubov endomorphism
$\varrho_U$ to $\cN$, where
$U=|h\rangle\langle h|-\bfe$. Moreover,
on the entire fermionic algebra $\cM$, $\varrho_U$
is in fact inner: $\varrho_U=\Ad(\psi(h))$.
(Note the $\psi(h)$ is unitary and self-adjoint.)
Thus $\psi(h)\in\Hom(\iota,\iota\eta)$ and hence
$\a_\eta^\pm(\psi(h))=\epsmp\eta\eta \psi(h)$
by \cite[Lemma 3.25]{BE1}. It was checked
in \cite{Bo1} that $\epspm\eta\eta=-\bfe$,
hence
$\a_\eta^\pm(\psi(h))=-\psi(h)$.
As $\a^\pm_\eta$ extends $\eta$ and since
$\psi(f)=\psi(f)\psi(h)^2$ we conclude that
$\a^\pm_\eta(\psi(f))=
-\varrho_U(\psi(f)\psi(h))\psi(h)
=-\varrho_U(\psi(f))$
for any $f\in L^2(S^1;\bbC)$.
We have shown
\[ \a^+_\eta =\a^-_\eta =\varrho_{-U}
=\varrho_{-\bfe}\circ \Ad(\psi(h)) \]
which we in fact recognize as a localized
automorphism. Because $\varrho_{-\bfe}$
is outer on $M(I_\circ)$ we find that
$[\a^\pm_\eta]$ and $[\id]$ are different sectors
of $M(I_\circ)$. We remark that, however,
$\a^\pm_\eta$ does not provide a new superselection
sector of the entire fermionic algebra $\cM$ because
$\varrho_{-\bfe}$ is unitarily implemented on the
Fock space by the parity operator.

Next we compute $\a^\pm_\sigma$. We obtain
$\a^\pm_\sigma(\psi(h))=\epsmp\eta\sigma \psi(h)
= \epspm\sigma\eta ^* \psi(h)$ from
\cite[Lemma 3.25]{BE1}. The statistics operators
are given by
$\epspm\sigma\eta = u_{\eta,\pm}^*\sigma(u_{\eta,\pm})$
where $u_{\eta,\pm}\in\cN$ are unitary charge transporters
for $\eta$. (See \cite[Subsect.\ 2.2]{BE1} for a
review on these matters using our notation.)
They can in fact be given explicitly as
$u_{\eta,\pm}=\psi(h_\pm)\psi(h)$ where $h_-$ and
$h_+$ are real functions with support in the left
respectively right complement of $I_\circ$ and
$\|h_\pm\|^2=2$ (cf.\ \cite[Subsect.\ 4.4]{Bo1}).
Now recall from \cite[Def.\ 3.9]{Bo1}
that $\sigma$ is the restriction of a Bogoliubov
endomorphism $\varrho_V$ where $V$ is a
``pseudo-localized'' isometry which
acts as the identity on the left and as
minus the identity on the right of $I_\circ$.
Consequently we have
$\varrho_V(\psi(h_\pm))=\mp \psi(h_\pm)$
and therefore
\[ \epspm\sigma\eta = \psi(h)\psi(h_\pm)
\varrho_V(\psi(h_\pm)\psi(h))=
\mp \psi(h) \varrho_V(\psi(h)) \,.\]
Hence we find
$\a^\pm_\sigma(\psi(h))=\mp\varrho_V(\psi(h))
\psi(h)^2=\mp\varrho_V(\psi(h))$.
As $\a^\pm_\sigma$ extends $\sigma$ and since
$\psi(f)=\psi(f)\psi(h)^2$ we conclude that
$\a^\pm_\sigma(\psi(f))=\mp\varrho_V(\psi(f))$
for any $f\in L^2(S^1;\bbC)$. We have shown
\[ \a^\pm_\sigma = \varrho_{\mp V} \,. \]
Hence
$\a^\pm_\sigma=\varrho_{-\bfe}\circ\a^\mp_\sigma$
but, though $\varrho_{-\bfe}$ is outer,
it does not mean that they produce distinct
sectors of $M(I_\circ)$: Let $v_0$ be the real
function which spans the co-kernel of $V$ and
$\|v_0\|^2=2$. Then we see that we can in fact write
$\a^\pm_\sigma=\Ad(\psi(v_0))\circ\a^\mp_\sigma$,
i.e.\ $\a^+_\sigma$ and $\a^-_\sigma$ are connected
by an automorphism which is inner since
${{\rm{supp}}}(v_0)\subset I_\circ$. Note that this
shows, maybe not surprisingly, that the implication
$[\a^+_\la]=[\a^-_\la]\Rightarrow\a^+_\la=\a^-_\la$
of \cite[Prop.\ 3.23]{BE1} does not hold without the
chiral locality assumption.
Also note that $\a^+_\sigma$ and $\a^-_\sigma$ are
solitons with different chirality which however
produce the same sector of $M(I_\circ)$.
Finally we argue that the sector $[\a^\pm_\sigma]$
must be different from $[\id]$ and $[\a^\pm_\eta]$
since $d_\sigma=\sqrt 2$. Summarizing we have found
three different irreducible sectors, $[\id]$,
$[\a^+_\eta]=[\a^-_\eta]$ and
$[\a^+_\sigma]=[\a^-_\sigma]$ and consequently
$\langle\a^+_\la,\a^-_\mu\rangle=\del\la\mu$.

Though explicit and instructive, this was the
pedestrians method to conclude
$Z_{\la,\mu}=\del\la\mu$! The conformal weights of the
Ising model and consequently the lowest eigenvalues
of the generator $L_0$ of rotations on the circle
in the vacuum representation $\pi_0$ and the
representations $\pi_0\circ\eta$ and
$\pi_0\circ\sigma$ of $\cN$ are given by
$h_0=0$, $h_\eta=1/2$ and $h_\sigma=1/16$, respectively.
Hence $\omega_0=1$, $\omega_\eta=-1$ and
$\omega_\sigma=\E^{\pi\I/8}$ by the conformal spin
and statistics theorem \cite{FG,FRS2,GL}. Since
the sectors obey the Ising fusion rules, this
determines Rehren's monodromy matrix $Y$ completely
by \erf{YT} so that it coincides with
the modular S-matrix of the Ising model
up to a normalization factor $1/\sqrt w=1/2$.
(And the braiding from the statistics operators
is in particular non-degenerate.)
We then apply $\a$-induction and find that
the matrix $Z$ with entries
$Z_{\la,\mu}=\langle\a^+_\la,\a^-_\mu\rangle$
gives us a modular invariant by
\cite[Thm.\ 5.7]{BEK1}. But the representation
of the modular group arising from the S- and
T-matrices of the Ising model does not
possess a non-trivial modular invariant!
Therefore we must have $Z_{\la,\mu}=\del\la\mu$.

\subsection{Degenerate braidings}

We first consider a completely degenerate example,
arising from the classical DHR theory \cite{DHR1}.
The subfactor $N\subset M$ is given by a local
subfactor $A(\cO)\subset F(\cO)$, arising from
a net of inclusions of observable algebras in
field algebras over the Minkowski space, arising
from a compact gauge group $G$. Then $A(\cO)$
is given as the fixed point algebra under the
outer action of the gauge group,
$A(\cO)=F(\cO)^G$. The canonical endomorphism
sector $[\canr]$ decomposes as
$[\canr]=\bigoplus d_\la [\la]$, where the sum
runs over DHR endomorphisms $\la$ labelled by the
irreducible representations of $G$. (By abuse
of notation we use the same symbol $\la$ for
the morphisms as for the the elements of the
group dual $\hat G$.) These DHR morphisms obey
the fusion rules of $\hat G$ so that the statistical
dimension $d_\la$ is in particular the dimension of
the group representation. We assume that $G$ is finite
and choose the system $\NXN$ to be given by all the
$\la$'s. Moreover, we assume that the field net is
purely bosonic, i.e.\ local, so that we have
$\om_\la=1$ for all $\la$. It is straightforward
to check that then (see \cite[Subsect.\ 2.2]{BEK1})
$w=\sum_\la d_\la^2=\# G$,
$S_{\la,\mu}=(\#G)^{-1} d_\la d_\mu$ and
$T_{\la,\mu}=\del\la\mu$. Note that the S-matrix
is a rank one projection here. Due to locality
of the field net, the chiral locality
condition\footnote{We admit that the name
``{\sl chiral} locality condition'' does not make
much sense when using the Minkowski space instead
of a compactified light cone axis $S^1$.}
holds here, and consequently
$\langle \a^\pm_\la , \a^\pm_\mu \rangle
= \langle \canr \la, \mu \rangle = d_\la d_\mu$,
which forces $[\a^\pm_\la]=d_\la[\id]$. Hence we
find $Z_{\la,\mu}=d_\la d_\mu$, i.e.\ $Z=wS$.
Note that $\tr(Z)=\#G$ and $\tr(Z^*Z)=(\#G)^2$.
However, we have $\# \NXM=1$ as
$[\iota\la]=[\a_\la^\pm\iota]=d_\la[\iota]$
and $\#\MXM=\# G$ since similarly
$[\iota\la\co\iota]=d_\la[\can]$,
and since it is known \cite{L3} that $\can$ decomposes
into automorphisms corresponding to the group elements.
So we observe that, due to the
degeneracy, the generating property of $\a$-induction
\cite[Thm.\ 5.10]{BEK1} does not hold, neither
the countings of \cite[Cors.\ 6.10 and 6.13]{BEK1}
are true here; we have an over-counting by $\# G$.

\begin{figure}[htb]
\unitlength 0.7mm
\begin{center}
\begin{picture}(155,40)
\thicklines
\multiput(10,10)(20,0){5}{\circle*{2}}
\multiput(10,30)(20,0){5}{\circle*{2}}
\multiput(20,20)(20,0){4}{\circle*{2}}
\multiput(10,10)(20,0){4}{\line(1,1){20}}
\multiput(10,30)(20,0){4}{\line(1,-1){20}}
\put(10,36){\makebox(0,0){$\la_0$}}
\put(30,36){\makebox(0,0){$\la_2$}}
\put(50,36){\makebox(0,0){$\la_4$}}
\put(70,36){\makebox(0,0){$\la_6$}}
\put(90,36){\makebox(0,0){$\la_8$}}
\put(10,3){\makebox(0,0){$\beta_0$}}
\put(30,3){\makebox(0,0){$\beta_2$}}
\put(50,3){\makebox(0,0){$\beta_4$}}
\put(70,3){\makebox(0,0){$\beta_6$}}
\put(90,3){\makebox(0,0){$\beta_8$}}
\put(122,10){\makebox(0,0){$M$-$M$ morphisms}}
\put(122,20){\makebox(0,0){$M$-$N$ morphisms}}
\put(122,30){\makebox(0,0){$N$-$N$ morphisms}}
\end{picture}
\caption{$\a$-induction for the Jones-Wassermann
subfactor of type $\rmA_9$}
\label{A9}
\end{center}
\end{figure}
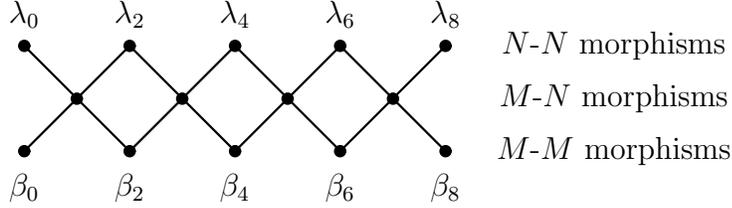

Maybe a more interesting and only partially degenerate
example is given by the following.
Instead of a conformal inclusion, we now consider the
Jones-Wassermann subfactor
\[ N=\pi_1(\LISUz)''\subset\pi_1(\LIcSUz)'=M \,, \]
where $\pi_1$ is the spin $j=1$ level $k$
positive energy representation of $\LSUz$,
$I\subset S^1$ a proper interval and $I'$ its
complement.
By \cite{W2} and \cite[Cor.\ 6.4]{P}, this is a
(type \nolinebreak III$_1$)
Jones subfactor \cite{J} with principal graph $\rmA_{k+1}$.
We label $N$-$N$ morphisms and $M$-$M$ morphisms
as exemplified for level $k=8$ in Fig.\ \ref{A9}.
We denote $d_j=d_{\la_j}=d_{\beta_j}$
for even spins $j$ and make the ``minimal choice''
$\NXN=\{\la_j:j \,{{\rm{even}}}\}$.
Note that $[\canr]=[\la_0]\oplus[\la_2]$ here,
and we have
$\langle\a^\pm_{j'},\iota\la_j\co\iota\rangle
=\langle\canr\la_{j'},\la_j\rangle$.
But $[\iota\la_j\co\iota]$ can be read off from
the Bratteli diagram, we have in particular
$[\iota\la_j\co\iota]=[\beta_{j-2}]\oplus
2[\beta_j]\oplus[\beta_{j+2}]$ for even spins
$0<j<k-1$. This forces in fact
$[\a^\pm_j]=[\beta_j]$.
We found in particular that $[\a^+_j]=[\a^-_j]$ for
all $\la\in\NXN$, though the system $\NXN$ is not completely
degenerate.  (The complete degeneracy means that any
monodromy operator is trivial. However, for $k>2$
the self-monodromy of the morphism $\la_2$ has always
the non-trivial eigenvalue $\E^{-8\pi\I/k+2}$ corresponding
to the fusion rule $N_{2,2}^0=1$ and due to
$h_2=2/k+2$, cf.\ \cite[Eq.\ (11)]{BEK1}.)
But since the non-degeneracy condition does not hold either,
this example shows that the generating property of
$\a$-induction \cite[Thm.\ 5.10]{BEK1} can even hold
without non-degeneracy in particular cases.
Note that we have also $Z_{\la,\mu}=\del\la\mu$ here,
but there is no representation of the modular group
arising from the braiding because of the degeneracy.
The degeneracy can be removed by extending the
$N$-$N$ system to all spins, $\NXN=\{\la_j:j=0,1,...,k\}$.
Still, the choice of the canonical endomorphism,
$[\canr]=[\la_0]\oplus[\la_2]$ produces the trivial
modular invariant. Namely, it is easy to check that
the formula
$\langle\iota\la_j,\iota\la_{j'}\rangle
=\del j{j'}+N_{2,j}^{j'}$
yields $k+1$ different $M$-$N$ sectors, and the
graph corresponding to left multiplication by
$[\a^\pm_1]$ gives just the $\rmA_{k+1}$ graph.
(Though it appears ``creased'', i.e.\ the sector $[\iota]$
is not an external vertex.) Consequently, we have
$Z_{j,j}=1$ for all $j$, implying that $Z$ is trivial.

\begin{appendix}
\section{The dual canonical endomorphism for $\rmE_7$}

\begin{lemma}
\label{E7can}
For $SU(2)$ at level $k=16$, there is an endomorphism
$\canr\in\Mor(N,N)$ such that
$[\canr]=[\id]\oplus[\la_8]\oplus[\la_{16}]$ and
which is the dual canonical endomorphism of a
subfactor $N\subset M$.
\end{lemma}

\begin{proof}
First note that the subfactor $\la_1(N)\subset N$
arising from the loop group construction for $\SUz_{16}$
in \cite{W2} is isomorphic to 
$P\otimes R \subset Q\otimes R$, where $Q$ is a hyperfinite
II$_1$ factor, $P\subset Q$ is the
Jones subfactor \cite{J} with principal graph $A_{17}$, and
$R$ is an injective III$_1$ factor, by 
\cite[Cor.\ 6.4]{P}.  This shows that the subfactor
$\canr(N)\subset N$ for
$[\canr]=[\id]\oplus[\la_8]\oplus[\la_{16}]$ is isomorphic
to $pP \subset p(Q_{15})p$, where $P\subset Q\subset Q_1
\subset Q_2\subset\cdots$ is the Jones tower of $P\subset Q$
and $p$ is a sum of three minimal projections in 
$P'\cap Q_{15}$ corresponding to $\id, \la_8, \la_{16}$.
It is thus enough to prove that the subfactor
$pP \subset p(Q_{15})p$ is a basic construction of some subfactor.

We recall a construction in \cite[Sect.\ 4.5]{GHJ}.
Let $\Gamma$ be one of the Dynkin diagrams of type A, D, E.
Let $A_0$ be an abelian von Neumann algebra ${\bf C}^n$ and
$A_1$ be a finite dimensional von Neumann algebra containing
$A_0$ such that the Bratteli diagram for $A_0\subset A_1$ is
$\Gamma$.  Using the unique normalized Markov trace on $A_1$,
we repeat basic constructions to get a tower
$A_0\subset A_1\subset A_2\subset\cdots$ with the Jones
projections $e_1, e_2, e_3,\cdots$.  Let $\tilde C$ be the
GNS-completion of $\bigcup_{m\ge 0} A_m$
with respect to the trace and $\tilde B$
its von Neumann subalgebra generated by $\{e_m\}_{m\ge 1}$.
We have $\tilde B'\cap \tilde C=A_0$ by Skau's lemma.
For a projection $q\in A_0$, we have a subfactor
$B=q \tilde B \subset q\tilde C q=C$, which is called a
Goodman-de la Harpe-Jones (GHJ) subfactor.

Let $\Gamma$ be E$_7$ and $q$ be the projection corresponding
to the vertex of E$_7$ with minimum Perron-Frobenius 
eigenvector entry.  We study the subfactor $B\subset C$
in this setting. 
Set $B_m=q \langle e_1, e_2,\dots,e_{m-1}\rangle$,
$C_m=q A_m q$.  The sequence $\{B_m \subset C_m\}_m$ is
a periodic sequence of commuting squares of period 2
in the sense of Wenzl \cite{We}.  For a sufficiently large $m$,
we can make a basic construction $B_m \subset C_m \subset D_m$
so that $B\subset C\subset D=\bigvee_m D_m$ is also a basic
construction.  We can extend the definition of $D_m$ to small
$m$ so that the sequences $\{B_m \subset C_m \subset D_m\}_m$ is
a periodic sequence of commuting squares of period 2.
For a sufficiently large $m$, the graph of the Bratteli diagram
for $B_{2m}\subset C_{2m}$ stays the same and the graph for
$C_{2m}\subset D_{2m}$ is its reflection.  This graph can be
computed as in Fig.\ \ref{Brat-E7} in an elementary way
(see e.g.\ \cite{Ok}, \cite[Examples 11.25, 11.71]{EK3}),
so we also have the graph for $B_{2m}\subset D_{2m}$,
and we see that $D_0$ is ${\bf C}\oplus {\bf C}\oplus {\bf C}$
and the three minimal projections in $D_0$ correspond to
the 0th, 8th, and 16th vertices of A$_{17}$.
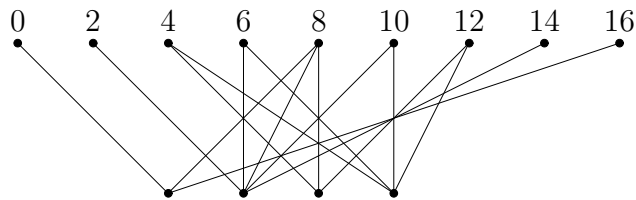
\begin{figure}[htb]
\begin{center}
\unitlength 1.0mm
\begin{picture}(100,40)
\thinlines 
\multiput(10,30)(10,0){9}{\circle*{1}}
\multiput(30,10)(10,0){4}{\circle*{1}}
\path(10,30)(30,10)
\path(20,30)(40,10)
\path(30,30)(50,10)
\path(30,30)(60,10)
\path(40,30)(40,10)
\path(40,30)(60,10)
\path(50,30)(30,10)
\path(50,30)(40,10)
\path(50,30)(50,10)
\path(60,30)(40,10)
\path(60,30)(60,10)
\path(70,30)(50,10)
\path(70,30)(60,10)
\path(80,30)(40,10)
\path(90,30)(30,10)
\put(10,33){\makebox(0,0){$0$}}
\put(20,33){\makebox(0,0){$2$}}
\put(30,33){\makebox(0,0){$4$}}
\put(40,33){\makebox(0,0){$6$}}
\put(50,33){\makebox(0,0){$8$}}
\put(60,33){\makebox(0,0){$10$}}
\put(70,33){\makebox(0,0){$12$}}
\put(80,33){\makebox(0,0){$14$}}
\put(90,33){\makebox(0,0){$16$}}
\end{picture}
\caption{The Bratteli diagram}
\label{Brat-E7}
\end{center}
\end{figure}
(The graph in Fig.\ \ref{Brat-E7} is actually
the principal graph of $B\subset C$ by \cite{Ok},
but this is not important here.)
Then we see that the Bratteli diagram for the sequence
$\{D_m\}_m$ starts with these three vertices and we have
the graph A$_{17}$ or a part of it as the Bratteli diagram at
each step, as in Fig.\ \ref{Brat}.
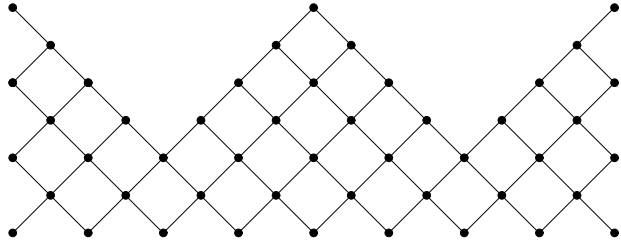
\begin{figure}[htb]
\begin{center}
\unitlength 1.0mm
\begin{picture}(100,50)
\thinlines 
\multiput(10,40)(40,0){3}{\circle*{1}}
\multiput(10,30)(10,0){2}{\circle*{1}}
\multiput(10,30)(10,0){2}{\circle*{1}}
\multiput(40,30)(10,0){3}{\circle*{1}}
\multiput(80,30)(10,0){2}{\circle*{1}}
\multiput(10,20)(10,0){9}{\circle*{1}}
\multiput(10,10)(10,0){9}{\circle*{1}}
\multiput(15,35)(30,0){2}{\circle*{1}}
\multiput(55,35)(30,0){2}{\circle*{1}}
\multiput(15,25)(10,0){2}{\circle*{1}}
\multiput(35,25)(10,0){4}{\circle*{1}}
\multiput(75,25)(10,0){2}{\circle*{1}}
\multiput(15,15)(10,0){8}{\circle*{1}}
\put(10,20){\line(1,-1){10}}
\put(10,30){\line(1,-1){20}}
\put(10,40){\line(1,-1){30}}
\put(35,25){\line(1,-1){15}}
\put(40,30){\line(1,-1){20}}
\put(45,35){\line(1,-1){25}}
\put(50,40){\line(1,-1){30}}
\put(75,25){\line(1,-1){15}}
\put(80,30){\line(1,-1){10}}
\put(85,35){\line(1,-1){5}}
\put(10,30){\line(1,1){5}}
\put(10,20){\line(1,1){10}}
\put(10,10){\line(1,1){15}}
\put(20,10){\line(1,1){30}}
\put(30,10){\line(1,1){25}}
\put(40,10){\line(1,1){20}}
\put(50,10){\line(1,1){15}}
\put(60,10){\line(1,1){30}}
\put(70,10){\line(1,1){20}}
\put(80,10){\line(1,1){10}}
\end{picture}
\caption{The Bratteli diagram for $\{D_m\}_m$}
\label{Brat}
\end{center}
\end{figure}
Each algebra $B_m$ is generated by the
Jones projections of the sequence $\{D_m\}_m$.

Similarly, if we choose A$_{17}$ as $\Gamma$ and $q$ be the
projection corresponding to the first vertex of A$_{17}$,
we get a periodic sequence $\{E_m\subset F_m\}_m$ of
commuting squares.  (Note that we start the numbering of
the vertices of A$_{17}$ with 0.)  It is well-known
that the resulting subfactor
$E=\bigvee_m E_m \subset \bigvee_m F_m =F$ is the Jones
subfactor \cite{J} with principal graph A$_{17}$.  We make
basic constructions of $E_m\subset F_m$ for 15 times in the
same way as above and
get a periodic sequence $\{E_m\subset G_m\}_m$ of
commuting squares. Let $\tilde q$ be a sum of three minimal
projections corresponding to the 0th, 8th, and 16th vertices
of A$_{17}$ in $G_0$.  Setting $\tilde E_m=\tilde q E_m$
and $\tilde G_m=q G_m q$, we get a periodic sequence of
commuting squares $\{\tilde E_m\subset \tilde G_m\}_m$
such that the resulting subfactor
$\bigvee_m \tilde E_m \subset \bigvee_m \tilde G_m$ is
isomorphic to $pP \subset p(Q_{15})p$ defined
in the first paragraph.

Now we see that the Bratteli diagram of the sequence
$\{\tilde G_m\}_m$ is the same as the one for
$\{D_m\}_m$ as in Fig.\ \ref{Brat}
and each algebra $\tilde E_m$ is generated by the Jones
projections for the sequence $\{\tilde G_m\}_m$.
This shows that the two periodic sequences of commuting
squares $\{B_m\subset D_m\}_m$ and
$\{\tilde E_m\subset \tilde G_m\}_m$ are isomorphic.
Thus the resulting subfactors $B\subset D$ and
$pP \subset p(Q_{15})p$ are also isomorphic. Since
the subfactor $B\subset D$ is a basic construction
of $B\subset C$, we conclude that the
subfactor $pP \subset p(Q_{15})p$ is also a
basic construction of some subfactor, as desired.
\end{proof}

\begin{remark}{\rm
With a different choice of $q$ corresponding to
another end vertex of E$_7$, we can also prove that
$[\la_0]\oplus[\la_6]\oplus[\la_{10}]\oplus[\la_{16}]$
for $SU(2)_{16}$ gives a dual canonical endomorphism
in a similar way.  This also produces
the E$_7$ modular invariant.

We can also choose D$_5$ as $\Gamma$ and $q$ to be a minimal
central projection corresponding to one of the two tail vertices
of D$_5$, and then the same method as in the above proof
shows that $[\la_0]\oplus[\la_4]$ for $SU(2)_6$ gives a
dual canonical endomorphism.  One can check that this
produces the  D$_5$ modular invariant.
}\end{remark}

In Lemma \ref{E7can} above, we have used the construction of
the GHJ-subfactor for E$_7$.  We can also apply the same
construction to E$_6$, E$_8$ as in \cite{GHJ}.  Note that
the principal graph \cite{Ok} of the GHJ-subfactor with
$\Gamma={\rm E}_6$ [resp.\ E$_8$] for the choice of $q$
corresponding to the vertex with the lowest Perron-Frobenius
eigenvector entry is the same as the principal graph,
Fig.\ 3 [resp.\ Fig.\ 6] in \cite{BE3},
of the subfactor arising from the
conformal inclusion $\SUz_{10}\subset\SOf_1$ 
[resp.\ $\SUz_{28}\subset(\Gtwo)_1$] studied in
\cite[Sect.\ 6.1]{BE3}.  It is then natural to expect that these
subfactor are indeed isomorphic (after tensoring a common 
injective factor of type \nolinebreak III$_1$).
For the E$_6$ case, a combinatorial unpublished argument of Rehren
shows that we have only two paragroups for the principal graph
in \cite[Fig.\ 3]{BE3} and these produce two mutually dual
subfactors.  This implies the desired isomorphism of our two
subfactors by \cite[Cor.\ 6.4]{P}, but it seems very hard to
obtain a similar argument for the E$_8$ case.  Here we prove the
desired isomorphism for both cases of E$_6$ and E$_8$. 

\begin{proposition}
\label{CI=GHJ}
The subfactor arising from the
conformal inclusion $\SUz_{10}\subset\SOf_1$
$[$resp.\ $\SUz_{28}\subset(\Gtwo)_1]$ is isomorphic to the
GHJ subfactor constructed as above for E$_6$ $[$resp.\ E$_8]$
tensored with a common injective factor of
type \nolinebreak III$_1$.
\end{proposition}

\begin{proof}
By \cite[Cor.\ 6.4]{P}, it is enough to prove that the two
subfactors have the same higher relative commutants.

Let $N\subset M$ be the subfactor arising from the
conformal inclusion and $\iota$ the inclusion map
$N\hookrightarrow M$.  We label $N$-$N$ morphisms as
$\la_0=\id, \la_1,\dots,\la_k$, where $k=10$ or $k=28$.
We set the finite dimensional $C^*$-algebras $A_{m,l}$,
$m\ge 0, l\ge -1$, to be as follows.  
(For $l=-1$, $m$ starts at 1.)
$$\left\{
\begin{array}{ll}
\Hom(\theta^{m/2} (\la_1 \bar\la_1)^{l/2},
\theta^{m/2} (\la_1 \bar\la_1)^{l/2}),&
(m:{\rm even}, l:{\rm even}),\\
\Hom(\theta^{m/2} (\la_1 \bar\la_1)^{(l-1)/2}{\la_1},
\theta^{m/2} (\la_1 \bar\la_1)^{(l-1)/2}{\la_1}),&
(m:{\rm even}, l:{\rm odd}),\\
\Hom(\iota \theta^{(m-1)/2} (\la_1 \bar\la_1)^{l/2},
\iota \theta^{(m-1)/2} (\la_1 \bar\la_1)^{l/2}),&
(m:{\rm odd}, l:{\rm even}),\\
\Hom(\iota \theta^{(m-1)/2} (\la_1 \bar\la_1)^{(l-1)/2}{\la_1},
\iota \theta^{(m-1)/2} (\la_1 \bar\la_1)^{(l-1)/2}{\la_1}),&
(m:{\rm odd}, l:{\rm odd}),\\
\Hom(\bar\iota \gamma^{(m-2)/2}, \bar\iota \gamma^{(m-2)/2}),&
(m:{\rm even}, l=-1),\\
\Hom(\gamma^{(m-1)/2}, \gamma^{(m-1)/2}),&
(m:{\rm odd}, l=-1).
\end{array}
\right.
$$
We then naturally have inclusions
$A_{m,l}\subset A_{m,l+1}$, and similarly embeddings
$\iota: A_{2m,l}\hookrightarrow A_{2m+1,l}$ as well as
$\bar \iota: A_{2m-1,l}\hookrightarrow A_{2m,l}$. 
With these, we have a double sequence of commuting squares.
Note that the sequence $\{A_{m,l}\}_{m,l\ge0}$ is a usual
double sequence of string algebras as in \cite[Chapter II]{O2}
(cf. \cite[Sect.\ 11.3]{EK3}) and we now have an extra sequence
$\{A_{m,-1}\}_{m\ge1}$ here.

Set $A_{m,\infty}$  to be the GNS-completions
of $\bigcup_{l=0}^\infty A_{m,l}$ with respect to the trace.
Then we have the Jones tower as
$$A_{0,\infty}\subset A_{1,\infty}\subset A_{2,\infty}\subset 
\cdots.$$
The Bratteli diagram of $\{A_{0,l}\}_l$ is given by
reflections of the Dynkin diagram of type A$_{11}$ or A$_{29}$,
so the algebra $A_{0,\infty}$ is generated by
the Jones projections.
The Bratteli diagram of $\{A_{1,l}\}_l$ is given by reflections of
the Dynkin diagram of type E$_6$ or E$_8$
since we know the fusion graph of $\la_1$ on the $M$-$N$
sectors, so the subfactor 
$A_{0,\infty}\subset A_{1,\infty}$ is isomorphic to the
GHJ-subfactor.  Then we next show
that the higher relative commutants of this subfactor
are given as
\begin{eqnarray*}
A'_{0,\infty}\cap A_{m,\infty}&=& A_{m,0},\\
A'_{1,\infty}\cap A_{m,\infty}&=& A_{m,-1},
\end{eqnarray*}
which are also the higher relative commutants of $N\subset M$
from the above definition, so the proof will be complete.

The definition of $\{A_{m,l}\}_{m,l}$ shows that
$A_{2m,0}$ and $A_{0,l}$ commute. Then Ocneanu's
compactness argument \cite[Sect.\ II.6]{O2}
(cf. \cite[Thm.\ 11.15]{EK3}) or Wenzl's
dimension estimate \cite[Thm.\ 1.6]{We} gives
$A_{m,0}= A'_{0,\infty}\cap A_{m,\infty}$.
We similarly have 
$A_{m,-1}\subset A'_{1,\infty}\cap A_{m,\infty}$.
In general, we have
$$\dim (A'_{0,\infty}\cap A_{2m+1,\infty})=
\dim (A'_{1,\infty}\cap A_{2m+2,\infty}),$$
so that we can compute
\[ \bearl
\dim\;\Hom(\iota \theta^m, \iota \theta^m)
= \dim\; A_{2m+1,0}
= \dim (A'_{0,\infty}\cap A_{2m+1,\infty}) \\[.4em]
\qquad = \dim (A'_{1,\infty}\cap A_{2m+2,\infty})
 \ge \dim\; A_{2m+2,-1}
 = \dim\; \Hom (\bar\iota \gamma^m,\bar\iota \gamma^m) \\[.4em]
\qquad = \dim\; \Hom(\iota \theta^m, \iota \theta^m) \,,
\eear \]
which shows equality
$A_{2m+2,-1} = A'_{1,\infty}\cap A_{2m+2,\infty}$.
We then have
$$ A'_{1,\infty}\cap A_{2m+1,\infty} \subset
(A'_{1,\infty}\cap A_{2m+2,\infty})\cap A_{2m+1,\infty}
 = A_{2m+2,-1} \cap A_{2m+1,\infty} = A_{2m+1,-1},$$
which completes the proof.
\end{proof}

Proposition \ref{CI=GHJ} implies in particular that the
graph in \cite[Fig.\ 7]{BE3} is also the dual principal
graph of the GHJ-subfactor arising from E$_8$.

\end{appendix}


\vspace{0.5cm}
\begin{footnotesize}
\noindent{\it Acknowledgment.}
Part of this work was done during visits of the third
author to the University of Wales Swansea and the
University of Wales Cardiff, visits of all the three
to Universit\`a di Roma ``Tor Vergata''
and visits of the first two authors to the
Australian National University, Canberra,
the University of Melbourne, the University of
Newcastle, the University of Tokyo and the
Research Institute for Mathematical Sciences, Kyoto.
We are indebted to R.\ Longo, L.\ Zsido, J.E.\ Roberts,
D.W.\ Robinson, P.A.\ Pearce, I.\ Raeburn, T.\ Miwa,
H.\ Araki and these institutions for their hospitality.
J.B.\ thanks R.\ Longo for an inspiring correspondence
on the subject of Subsect.\ \ref{Ising}.
We gratefully acknowledge the financial support of
the Australian National University, CNR (Italy), EPSRC
(U.K.), the EU TMR Network in Non-Commutative Geometry,
Grant-in-Aid for Scientific Research, Ministry of
Education (Japan), the Kanagawa Academy of Science
and Technology Research Grants,  the Universit\`a di
Roma ``Tor Vergata'', University of Tokyo,
and the University of Wales.
\end{footnotesize}

\begin{footnotesize}

\end{footnotesize}

\begin{thebibliography}{99}

\bibitem{A1}
H. Araki, 
{\sl On quasi-free states of CAR and Bogoliubov
automorphisms}, 
Publ. RIMS Kyoto Univ. {\bf 6} (1970/71) 385--442

\bibitem{Bo1}
J. B\"ockenhauer, 
{\sl Localized endomorphisms of the chiral Ising model}, 
Commun. Math. Phys. {\bf 177} (1996) 265--304

\bibitem{BE1}
J. B\"ockenhauer, D.E. Evans, 
{\sl Modular invariants, graphs and $\a$-induction for
nets of subfactors. I}, 
Commun. Math. Phys. {\bf 197} (1998) 361--386

\bibitem{BE2}
J. B\"ockenhauer, D.E. Evans, 
{\sl Modular invariants, graphs and $\a$-induction for
nets of subfactors. II},
Commun. Math. Phys. {\bf 200} (1999) 57--103

\bibitem{BE3}
J. B\"ockenhauer, D.E. Evans, 
{\sl Modular invariants, graphs and $\alpha$-induction for
nets of subfactors. III}, preprint, hep-th/9812110

\bibitem{BEK1}
J. B\"ockenhauer, D.E. Evans, Y. Kawahigashi,
{\sl On $\a$-induction, chiral generators and modular
invariants for subfactors}, preprint, math.OA/9904109

\bibitem{CIZ1}
A. Cappelli, C. Itzykson, J.-B. Zuber,
{\sl Modular invariant partition functions in
two dimensions}, Nucl. Phys. {\bf B280} (1987) 445--465

\bibitem{CIZ2}
A. Cappelli, C. Itzykson, J.-B. Zuber,
{\sl The $A$-$D$-$E$ classification of minimal and
$A^{(1)}_1$ conformal invariant theories},
Commun. Math. Phys. {\bf 113} (1987) 1--26

\bibitem{DZ1}
P. Di Francesco, J.-B. Zuber,
{\sl $SU(N)$ lattice integrable models associated with 
graphs}, Nucl. Phys. {\bf B338} (1990) 602--646

\bibitem{DZ2}
P. Di Francesco, J.-B. Zuber,
{\sl $SU(N)$ lattice integrable models and modular 
invariance}, in: S. Randjbar et al (eds.),
{\it Recent Developments in Conformal Field Theories},
Singapore: World Scientific 1990, pp.\ 179--215

\bibitem{DHR1}
S. Doplicher, R. Haag, J.E. Roberts,
{\sl Fields, observables and gauge transformations. I},
Commun. Math. Phys. {\bf 13} (1969), 1--23.

\bibitem{EK3}
D.E. Evans, Y. Kawahigashi,
{\it Quantum symmetries on operator algebras},
Oxford: Oxford University Press, 1998

\bibitem{FRS2}
K. Fredenhagen, K.-H. Rehren, B. Schroer,
{\sl Superselection sectors with braid group statistics
and exchange algebras. II}, 
Rev. Math. Phys. {\bf Special issue} (1992) 113--157

\bibitem{FG}
J. Fr\"ohlich, F. Gabbiani,
{\sl Braid statistics in local quantum theory}, 
Rev. Math. Phys. {\bf 2} (1990) 251--353

\bibitem{G1}
T. Gannon,
{\sl WZW commutants, lattices and level-one partition functions},
Nucl. Phys. {\bf B396} (1993) 708--736

\bibitem{G2}
T. Gannon,
{\sl The classification of affine $\SUd$ modular invariants},
Commun. Math. Phys. {\bf 161} (1994) 233--264

\bibitem{GS}
B. Gato-Rivera, A.N. Schellekens,
{\sl Complete classification of simple current automorphisms},
Nucl. Phys. {\bf B353} (1991) 519--537

\bibitem{GHJ}
F. Goodman, P. de la Harpe, V.F.R. Jones,
{\it Coxeter graphs and towers of algebras}, MSRI publications 14,
Berlin: Springer, 1989

\bibitem{GL}
D. Guido, R. Longo,
{\sl The conformal spin and statistics theorem},
Commun. Math. Phys. {\bf 181} (1996) 11--35

\bibitem{Itz}
C. Itzykson,
{\sl From the harmonic oscillator to the A-D-E
classification of conformal models},
Adv. Stud. in Pure Math. {\bf 19} (1989) 287--346

\bibitem{I0}
M. Izumi, 
{\sl Application of fusion rules to classification
of subfactors},
Publ. RIMS, Kyoto Univ. {\bf 27} (1991) 953--994

\bibitem{J}
V.F.R. Jones,
{\sl Index for subfactors},
Invent. Math. {\bf 72} (1983) 1--25

\bibitem{Kc}
V.G. Kac,
{\it Infinite dimensional Lie algebras},
3rd edition, Cambridge: Cambridge University Press, 1990

\bibitem{Kt}
A. Kato, 
{\sl Classification of modular invariant partition
functions in two dimensions},
Modern Phys. Lett {\bf A 2} (1987) 585--600

\bibitem {Kw}
T. Kawai,
{\sl On the structure of fusion rule algebras},
Phys. Lett. {\bf B217} (1989) 247-251

\bibitem {Ko}
H. Kosaki,
{\sl Extension of Jones theory on index to arbitrary factors},
J. Funct. Anal. {\bf 66} (1986) 123-140

\bibitem{L3}
R. Longo,
{\sl A Duality for Hopf algebras and for
subfactors I}, Commun. Math. Phys. {\bf 159} (1994) {133--150}

\bibitem{LR}
R. Longo, K.-H. Rehren,
{\sl Nets of subfactors},
Rev. Math. Phys.  {\bf 7} (1995) 567--597

\bibitem{MaS}
G. Mack, V. Schomerus, 
{\sl Conformal field algebras with quantum symmetry from
the theory of superselection sectors}, 
Commun. Math. Phys. {\bf 134} (1990) 139--196

\bibitem{MS}
G. Moore, N. Seiberg, 
{\sl Naturality in conformal field theory}, 
Nucl. Phys. {\bf B313} (1989) 16--40

\bibitem{N} W. Nahm,
{\sl Lie group exponents and $\SUz$ current algebras},
Commun. Math. Phys. {\bf 118} (1988) {171--176}

\bibitem{O2}
A. Ocneanu,
{\sl Quantum symmetry, differential geometry of finite
graphs and classification of subfactors},
Univ. of Tokyo Seminary Notes 45, 1991
(Notes recorded by Y. Kawahigashi)

\bibitem{O7}
A. Ocneanu,
{\sl Paths on Coxeter diagrams: From Platonic solids and
singularities to minimal models and subfactors}
(Notes recorded by S. Goto), 
in preparation.

\bibitem{Ok}
S. Okamoto, 
{\sl Invariants for subfactors arising from Coxeter
graphs}, in: H. Araki et al. (eds.),
{\it Current Topics in Operator Algebras},
Singapore: World Scientific 1991, pp.\ 84--103.

\bibitem{Pa}
V. Pasquier,
{\sl Etiology of IRF models},
Commun. Math. Phys. {\bf 118} (1988) 355--364

\bibitem{PZ}
V.B. Petkova, J.-B. Zuber,
{\sl From CFT to graphs},
Nucl. Phys. {\bf B463} (1996) 161--193

\bibitem{P}
S. Popa,
{\sl Classification of subfactors and of their endomorphisms},
CBMS Regional Conference Series, Am. Math. Soc. {\bf 86} (1995)

\bibitem{R0}
K.-H. Rehren,
{\sl Braid group statistics and their superselection rules},
in: D. Kastler (ed.), {\it The algebraic theory of
superselection sectors}, Palermo 1989,
Singapore: World Scientific 1990, pp. 333--355

\bibitem{R2}
K.-H. Rehren, 
{\sl Space-time fields and exchange fields},
Commun. Math. Phys. {\bf 132} (1990) 461--483

\bibitem{R4}
K.-H. Rehren, 
{\sl Chiral observables and modular invariants},
preprint, hep-th/9903262

\bibitem{SY}
A.N. Schellekens, S. Yankielowicz,
{\sl Extended chiral algebras and modular invariant
partition functions},
Nucl. Phys. {\bf B327} (1989) 673--703

\bibitem{W2}
A. Wassermann,
{\sl Operator algebras and conformal field theory III: Fusion
of positive energy representations of $SU(N)$ using bounded 
operators}, Invent. Math. {\bf 133} (1998) 467--538

\bibitem{We}
H. Wenzl, 
{\sl Hecke algebras of type $A_n$ and subfactors},
Invent. Math. {\bf 92} (1988) 345--383

\bibitem{X1}
F. Xu, 
{\sl New braided endomorphisms from conformal inclusions},
Commun. Math. Phys. {\bf 192} (1998) 347--403

\bibitem{X2}
F. Xu,
{\sl Applications of braided endomorphisms from conformal
inclusions},
Internat. Math. Research Notices,
(1998) 5--23, see also the erratum to Theorem 3.4 (1)
on page 437 of the same volume

\end{thebibliography}
\end{document}